\documentclass[reqno,11pt]{book}
\usepackage{psfig, amsmath, amstexnb, amsthm}
\usepackage{amssymb}  
\input dynsys.sty

\begin{document}

\thispagestyle{empty}

\begin{center}
{\Huge Perturbation Theory\\
 of Dynamical Systems\\}
\vspace{2cm}
{\Large
Nils Berglund\\
Department of Mathematics\\
ETH Z\"urich\\
8092 {\bf Z\"urich}\\
Switzerland\\
}
\vspace{2cm}
{\large Lecture Notes\\
Summer Semester 2001\\}
\vspace{2cm}
{\Large Version:\\
November 14, 2001}
\end{center}

\chapter*{Preface}

This text is a slightly edited version of lecture notes for a course I gave
at ETH, during the Summer term 2001, to undergraduate Mathematics and
Physics students. It covers a few selected topics from perturbation theory
at an introductory level. Only certain results are proved, and for some of
the most important theorems, sketches of the proofs are provided. 

Chapter~2 presents a topological approach to perturbations of planar vector
fields. It is based on the idea that the qualitative dynamics of most vector
fields does not change under small perturbations, and indeed, the set of all
these structurally stable systems can be identified. The most common
exceptional vector fields can be classified as well. 

In Chapter~3, we use the problem of stability of elliptic periodic orbits to
develop perturbation theory for a class of dynamical systems of dimension 3
and larger, including (but not limited to) integrable Hamiltonian systems.
This will bring us, via averaging and Lie-Deprit series, all the way to
KAM-theory.

Finally, Chapter~4 contains an introduction to singular perturbation theory,
which is concerned with systems that do not admit a well-defined limit when
the perturbation parameter goes to zero. After stating a fundamental result
on existence of invariant manifolds, we discuss a few examples of dynamic
bifurcations.  

An embarrassing number of typos from the first version has been corrected,
thanks to my students' attentiveness.  

\bigskip
\noindent
Files available at 
{\tt http://www.math.ethz.ch/$\sim$berglund}

\noindent
Please send any comments to 
{\tt berglund@math.ethz.ch}

\bigskip
\noindent
\begin{flushright}
Z\"urich, November 2001
\end{flushright}

\addtolength\textheight{2mm}
\tableofcontents
\goodbreak
\addtolength\textheight{-2mm}
\setcounter{page}{0}


\chapter{Introduction and Examples}
\label{ch_in}

The main purpose of perturbation theory is easy to describe. Consider for
instance an ordinary differential equation of the form
\begin{equation}
\label{i1}
\dot x = f(x,\eps),
\end{equation}
where $x\in\R^n$ defines the state of the system, $\eps\in\R$ is a
parameter, $f:\R^n\times\R\to\R^n$ is a sufficiently smooth vector field,
and the dot denotes derivation with respect to time. Assume that the
dynamics of the unperturbed equation
\begin{equation}
\label{i2}
\dot x = f(x,0)
\end{equation}
is well understood, meaning either that we can solve this system exactly,
or, at least, that we have a good knowledge of its long-time behaviour. What
can be said about solutions of \eqref{i1}, especially for $t\to\infty$, if
$\eps$ is sufficiently small? 

The answer to this question obviously depends on the dynamics of
\eqref{i2}, as well as on the kind of $\eps$-dependence of $f$. There are
several reasons to be interested in this problem:
\begin{itemiz}
\item	Only very few differential equations can be solved exactly, and if
we are able to describe the dynamics of a certain class of perturbations, we
will significantly increase the number of systems that we understand.

\item	In applications, many systems are described by differential
equations, but often with an imprecise knowledge of the function $f$; it is
thus important to understand the stability of the dynamics under small
perturbations of $f$. 

\item	If the dynamics of \eqref{i1} is qualitatively the same for some
set of values of $\eps$, we can construct equivalence classes of dynamical
systems, and it is sufficient to study one member of each class.
\end{itemiz}
Many different methods have been developed to study perturbation problems of
the form \eqref{i1}, and we will present a selection of them. Before doing
that for general systems, let us discuss a few examples. 


\section{One-Dimensional Examples}
\label{sec_in1}

One-dimensional ordinary differential equations are the easiest to study,
and in fact, their solutions never display a very interesting behaviour.
They provide, however, good examples to illustrate some basic facts of
perturbation theory.

\begin{example}
\label{ex_1D1}
Consider the equation
\begin{equation}
\label{1D1}
\dot x = f(x,\eps) = -x + \eps g(x).
\end{equation}
When $\eps=0$, the solution with initial condition $x(0)=x_0$ is 
\begin{equation}
\label{1D2}
x(t) = x_0 \e^{-t},
\end{equation}
and thus $x(t)$ converges to $0$ as $t\to\infty$ for any initial condition
$x_0$. Our aim is to find a class of functions $g$ such that the solutions
of \eqref{1D1} behave qualitatively the same for sufficiently small $\eps$.
We will compare different methods to achieve this.

\begin{enumerate}
\item	{\bf Method of exact solution:}
The solution of \eqref{1D1} with initial condition $x(0)=x_0$ satisfies
\begin{equation}
\label{1D3}
\int_{x_0}^{x(t)} \frac{\6x}{-x+\eps g(x)} = t.
\end{equation}
Thus the equation can be solved in principle by computing the integral and
solving with respect to $x(t)$. Consider the particular case $g(x)=x^2$.
Then the integral can be computed by elementary methods, and the result is 
\begin{equation}
\label{1D4}
x(t) = \frac{x_0\e^{-t}}{1-\eps x_0(1-\e^{-t})}.
\end{equation}
Analysing the limit $t\to\infty$, we find that
\begin{itemiz}
\item	if $x_0<1/\eps$, then $x(t)\to0$ as $t\to\infty$;
\item	if $x_0=1/\eps$, then $x(t)=x_0$ for all $t$;
\item	if $x_0>1/\eps$, then $x(t)$ diverges for a time $t^\star$ given by
$\e^{-t^\star}=1-1/(\eps x_0)$. 
\end{itemiz}
This particular case shows that the perturbation term $\eps x^2$ indeed has
a small influence as long as the initial condition $x_0$ is not too large
(and the smaller $\eps$, the larger the initial condition may be). The
obvious advantage of an exact solution is that the long-term behaviour can
be analysed exactly. There are, however, important disadvantages:
\begin{itemiz}
\item	for most functions $g(x)$, the integral in \eqref{1D3} cannot be
computed exactly, and it is thus difficult to make statements about general
classes of functions $g$;
\item	even if we manage to compute the integral, quite a bit of analysis
may be required to extract information on the asymptotic behaviour.
\end{itemiz}
Of course, even if we cannot compute the integral, we may still be able to
extract some information on the limit $t\to\infty$ directly from
\eqref{1D3}. In more than one dimension, however, formulae such as
\eqref{1D3} are not available in general, and thus alternative methods are
needed.

\item	{\bf Method of series expansion:}
A second, very popular method, is to look for solutions of the form
\begin{equation}
\label{1D5}
x(t) = x^0(t) + \eps x^1(t) + \eps^2 x^2(t) + \dots
\end{equation} 
with initial conditions $x^0(0)=x_0$ and $x^j(0)=0$ for $j\geqs 1$. By
plugging this expression into \eqref{1D1}, and assuming that $g$ is
sufficiently differentiable that we can write 
\begin{equation}
\label{1D6}
\eps g(x) = \eps g(x^0) + \eps^2 g'(x_0) x^1 + \dots
\end{equation}
we arrive at a hierarchy of equations that can be solved recursively: 
\begin{align}
\nonumber
\dot x^0 &= -x^0 
&& \Rightarrow &
x^0(t) &= x_0\e^{-t} \\
\label{1D7}
\dot x^1 &= -x^1 + g(x^0)
&& \Rightarrow &
x^1(t) &= \int_0^t \e^{-(t-s)} g(x^0(s))\6s \\
\nonumber
\dot x^2 &= -x^2 + g'(x^0)x^1
&& \Rightarrow &
x^2(t) &= \int_0^t \e^{-(t-s)} g'(x^0(s))x^1(s)\6s \\
\nonumber
&&& \dots &&
\end{align}
The advantage of this method is that the computations are relatively easy
and systematic. Indeed, at each order we only need to solve an inhomogeneous
linear equation, where the inhomogeneous term depends only on previously
computed quantities.

The disadvantages of this method are that it is difficult to control the
convergence of the series \eqref{1D5}, and that it depends on the smoothness
of $g$ (if $g$ has a finite number of continuous derivatives, the expansion
\eqref{1D5} has to be finite, with the last term depending on $\eps$). If
$g$ is analytic, it is possible to investigate the convergence of the
series, but this requires quite heavy combinatorial calculations.

\item	{\bf Iterative method:}
This variant of the previous method solves some of its disadvantages.
Equation \eqref{1D1} can be treated by the method of variation of the
constant: putting $x(t)=c(t)\e^{-t}$, we get an equation for $c(t)$ leading
to the relation
\begin{equation}
\label{1D8}
x(t) = x_0\e^{-t} + \eps \int_0^t \e^{-(t-s)}g(x(s))\6s.
\end{equation}
One can try to solve this equation by iterations, setting
$x^{(0)}(t)=x_0\e^{-t}$ and defining a sequence of functions recursively by
$x^{(n+1)}=\cT x^{(n)}$, where
\begin{equation}
\label{1D9}
(\cT x)(t) = x_0\e^{-t} + \eps\int_0^t \e^{-(t-s)} g(x(s)) \6s. 
\end{equation}
If $\cT$ is a contraction in an appropriate Banach space, then the sequence
$(x^{(n)})_{n\geqs0}$ will converge to a solution of \eqref{1D8}. It turns
out that if $g$ satisfies a uniform global Lipschitz condition 
\begin{equation}
\label{1D10}
\abs{g(x)-g(y)} \leqs K \abs{x-y} 
\qquad \forall x, y\in\R,
\end{equation}
and for the norm 
\begin{equation}
\label{1D11}
\norm{x}_t = \sup_{0\leqs s\leqs t} \abs{x(s)},
\end{equation}
we get from the definition of $\cT$ the estimate
\begin{equation}
\label{1D12}
\begin{split}
\abs{(\cT x)(t) - (\cT y)(t)} 
&\leqs \eps \int_0^t \e^{-(t-s)} K\abs{x(s)-y(s)}\6s \\
&\leqs \eps K (1-\e^{-t}) \norm{x-y}_t.
\end{split}
\end{equation}
Thus if $\eps<1/K$, $\cT$ is a contraction with contraction constant
$\lambda=\eps K (1-\e^{-t})$. The iterates $x^{(n)}(t)$ converge to the
solution $x(t)$, and we can even estimate the rate of convergence:
\begin{equation}
\label{1D13}
\begin{split}
\norm{x^{(n)}-x}_t &\leqs \sum_{m=n}^{\infty} \norm{x^{(m+1)}-x^{(m)}}_t \\
&\leqs \sum_{m=n}^{\infty} \lambda^m \norm{x^{(1)}-x^{(0)}}_t 
= \frac{\lambda^n}{1-\lambda} \norm{x^{(1)}-x^{(0)}}_t.
\end{split}
\end{equation}
Thus after $n$ iterations, we will have approximated $x(t)$ up to order
$\eps^n$. 

The advantage of the iterative method is that it allows for a precise
control of the convergence, under rather weak assumptions on $g$ (if $g$
satisfies only a local Lipschitz condition, as in the case $g(x)=x^2$, the
method can be adapted by restricting the domain of $\cT$, which amounts to
restricting the domain of initial conditions $x_0$). 

However, this method usually requires more calculations than the previous
one, which is why one often uses a combination of both methods: iterations
to prove convergence, and expansions to compute the terms of the series. It
should also be noted that the fact that $\cT$ is contracting or not will
depend on the unperturbed vector field $f(x,0)$. 

\item	{\bf Method of changes of variables:}
There exists still another method to solve \eqref{1D1} for small $\eps$,
which is usually more difficult to implement, but often easier to extend to
other unperturbed vector fields and higher dimension. 

Let us consider the effect of the change of variables $x=y+\eps h(y)$, where
the function $h$ will be chosen in such a way as to make the equation for
$y$ as simple as possible. Replacing $x$ by $y+\eps h(y)$ in \eqref{1D1}, we
get 
\begin{equation}
\label{1D14}
\dot y + \eps h'(y) \dot y = -y - \eps h(y) + \eps g(y+\eps h(y)).
\end{equation}
Ideally, the transformation should remove the perturbation term $\eps g(x)$.
Let us thus impose that 
\begin{equation}
\label{1D15}
\dot y = - y.
\end{equation}
This amounts to requiring that
\begin{equation}
\label{1D16}
h'(y) = \frac 1y h(y) - \frac 1y g(y+\eps h(y)).
\end{equation}
If such a function $h$ can be found, then the solution of \eqref{1D1} is
determined by the relations
\begin{equation}
\label{1D17}
\begin{split}
x(t) &= y(t) + \eps h(y(t)) = y_0 \e^{-t} + \eps h(y_0\e^{-t}) \\
x_0 &= y_0 + \eps h(y_0).
\end{split}
\end{equation}
In particular, we would have that $x(t)\to \eps h(0)$ as $t\to\infty$. 
We can try to construct a solution of \eqref{1D16} by the method of
variation of the constant. Without $g$, $h(y)=y$ is a solution. Writing
$h(y)=c(y)y$ and solving for $c(y)$, we arrive at the relation
\begin{equation}
\label{1D18}
h(y) = -y \int_a^y \frac 1{z^2} g(z+\eps h(z)) \6z 
\bydef (\cT h)(y),
\end{equation}
where $a$ is an arbitrary constant. It turns out that if $g$ satisfies the
Lipschitz condition \eqref{1D10} and $\abs{a}\to\infty$, then $\cT$ is a
contraction if $\eps<1/K$, and thus $h$ exists and can be computed by
iterations. 

For more general $f(x,0)$ it is usually more difficult to find $h$, and one
has to construct it by successive approximations. 

\item	{\bf Geometrical method:}
The geometrical approach gives less detailed informations than the previous
methods, but can be extended to a larger class of equations and requires
very little computations. 

In dimension one, the basic observation is that $x(t)$ is (strictly)
increasing if and only if $f(x(t),\eps)$ is (strictly) positive. Thus
$x(t)$ either converges to a zero of $f$, or it is unbounded. Hence the
asymptotic behaviour of $x(t)$ is essentially governed by the zeroes and
the sign of $f$. Now we observe that if $g$ satisfies the Lipschitz
condition \eqref{1D10}, then for all $y>x$ we have 
\begin{equation}
\label{1D19}
\begin{split}
f(y,\eps) - f(x,\eps) 
&= -y + x + \eps \brak{g(y)-g(x)} \\
&\leqs -(1-\eps K) (y-x), 
\end{split}
\end{equation}
and thus $f$ is monotonously decreasing if $\eps<1/K$. Since
$f(x,\eps)\to\mp\infty$ as $x\to\pm\infty$, it follows that $f$ has only one
singular point, which attracts all solutions. 

If $g$ satisfies only a local Lipschitz condition, as in the case
$g(x)=x^2$, one can show that for small $\eps$ there exists a point near or
at $0$ attracting solutions starting in some \nbh. 
\end{enumerate}
\end{example}

We have just seen various methods showing that the equation $\dot x=-x$ is
robust under quite a general class of small perturbations, in the sense
that solutions of the perturbed system are also attracted by a single point
near the origin, at least when they start not to far from $0$. The system
$\dot x=-x$ is called \defwd{structurally stable}. 

We will now give an example where the situation is not so simple. 

\begin{example}
\label{ex_1D2}
Consider the equation 
\begin{equation}
\label{1D20}
\dot x = f(x,\eps) = -x^3 + \eps g(x).
\end{equation}
If $\eps=0$, the solution with initial condition $x(0)=x_0$ is given by 
\begin{equation}
\label{1D21}
x(t) = \frac{x_0}{\sqrt{1+2x_0^2 t}},
\end{equation}
and thus $x(t)$ converges to $0$ as $t\to\infty$, though much slower than in
Example \ref{ex_1D1}. Let us start by considering the particular case
$g(x)=x$. The perturbation term being linear, it is a good idea to look for
a solution of the form $x(t)=c(t)\e^{\eps t}$, and we find in this way
\begin{equation}
\label{1D22}
x(t) = \frac{x_0\e^{\eps t}}{\sqrt{1+x_0^2 \frac{\e^{2\eps t}-1}\eps}}.
\end{equation}
If $\eps\leqs 0$, then $x(t)\to 0$ as $t\to \infty$ as before. 
If $\eps>0$, however, it turns out that 
\begin{equation}
\label{1D23}
\lim_{t\to\infty} x(t) = 
\begin{cases}
-\sqrt\eps &  \text{if $x_0<0$} \\
0 &  \text{if $x_0=0$} \\
\sqrt\eps &  \text{if $x_0>0$.} 
\end{cases}
\end{equation}
This behaviour differs drastically from the behaviour for $\eps\leqs 0$.
Moreover, it is not obvious how the term $\sqrt\eps$ could be obtained from
an expansion of $x$ in powers of $\eps$. In fact, since $x(t)$ is a function
of $\eps t$, any such expansion will only converge for finite values of $t$,
and most perturbative methods introduced in Example \ref{ex_1D1} will fail
in the limit $t\to\infty$. 

\begin{figure}
 \centerline{\psfig{figure=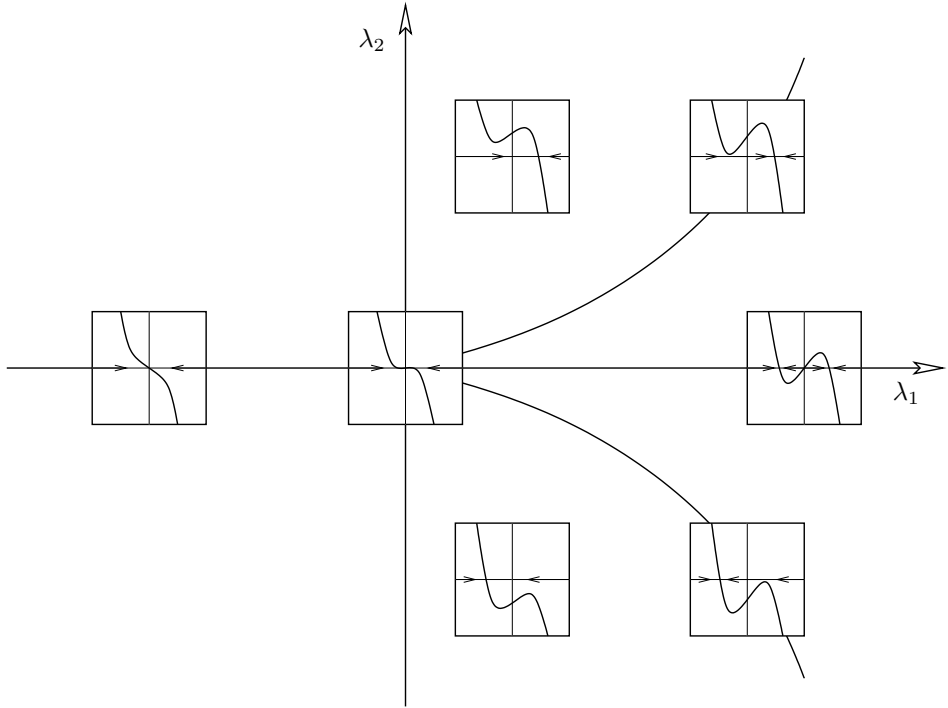,width=130mm}}
 \figtext{
 	\writefig	12.9	4.8	$\lambda_1$
 	\writefig	5.8	9.5	$\lambda_2$
 }
 \captionspace
 \caption[]{Graphs of the function $F_{\lambda_1,\lambda_2}(x) = -x^3 +
 \lambda_1 x + \lambda_2$ as a function of $\lambda_1$ and
 $\lambda_2$. The lines $4\lambda_1^3 = 27 \lambda_2^2$, $\lambda_1>0$ are
 lines of a codimension $1$ bifurcation, and the origin is a codimension $2$
 bifurcation point.}
 \label{fig_unfold1D}
\end{figure}

On the other hand, the geometrical approach still works. Let us first
consider the particular two-parameter family 
\begin{equation}
\label{1D24}
\dot x = F_{\lambda_1,\lambda_2}(x) = -x^3 + \lambda_1 x + \lambda_2.
\end{equation}
If $\lambda_2=0$, $F_{\lambda_1,0}$ vanishes only at $x=0$ if
$\lambda_1\leqs 0$, while $F_{\lambda_1,0}$ has three zeroes located at
$-\sqrt\lambda_1$, $0$ and $\sqrt\lambda_1$ if $\lambda_1>0$. This explains
in particular the asymptotic behaviour \eqref{1D23}. Adding the constant
term $\lambda_2$ will move the graph of $F$ vertically. If $\lambda_1\leqs
0$, nothing special happens, but if $\lambda_1>0$, two zeroes of $F$ will
disappear for sufficiently large $\abs{\lambda_2}$ (\figref{fig_unfold1D}). 
In the limiting case, the graph of $F$ is tangent to the $x$-axis at some
point $x_\star$. We thus have the conditions
\begin{equation}
\label{1D25}
\begin{split}
F_{\lambda_1,\lambda_2}(x_\star) &= -x_\star^3 + \lambda_1 x_\star +
\lambda_2 = 0 \\
F'_{\lambda_1,\lambda_2}(x_\star) &= -3 x_\star^2 + \lambda_1 = 0. 
\end{split}
\end{equation} 
Eliminating $x_\star$ from these equations, we get the relation 
\begin{equation}
\label{1D26}
4\lambda_1^3 = 27 \lambda_2^2
\end{equation}
for the critical line. In summary, 

\begin{itemiz}
\item	if $27\lambda_2^2<4\lambda_1^3$, the system has two stable
equilibrium points $x^\pm$ and one unstable equilibrium point
$x^0\in(x^-,x^+)$, with 
\begin{equation}
\label{1D27}
\lim_{t\to\infty}x(t) = 
\begin{cases}
x^- & \text{if $x(0)<x^0$} \\
x^0 & \text{if $x(0)=x^0$} \\
x^+ & \text{if $x(0)>x^0$;}
\end{cases}
\end{equation}
\item	if $\lambda_1>0$ and $3\sqrt3 \lambda_2 = 2\lambda_1^{3/2}$, there
is a stable equilibrium $x^+$ and an equilibrium $x^-$ such that 
\begin{equation}
\label{1D28}
\lim_{t\to\infty}x(t) = 
\begin{cases}
x^- & \text{if $x(0)\leqs x^-$} \\
x^+ & \text{if $x(0)>x^-$;}
\end{cases}
\end{equation}
\item	if $\lambda_1>0$ and $3\sqrt3 \lambda_2 = -2\lambda_1^{3/2}$, a
similar situation occurs with the roles of $x^-$ and $x^+$ interchanged;
\item	in all other cases, there is a unique equilibrium point attracting
all solutions.
\end{itemiz}

What happens for general functions $g$? The remarkable fact is the
following: for any sufficiently smooth $g$ and any sufficiently small
$\eps$, there exist $\lambda_1(g,\eps)$ and $\lambda_2(g,\eps)$ (going
continuously to zero as $\eps\to0$) such that the dynamics of $\dot
x=-x^3+\eps g(x)$ is qualitatively the same, near $x=0$, as the dynamics of
\eqref{1D24} with $\lambda_1=\lambda_1(g,\eps)$ and
$\lambda_2=\lambda_2(g,\eps)$.

One could have expected the dynamics to depend, at least, on quadratic terms
in the Taylor expansion of $g$. It turns out, however, that these terms can
be removed by a transformation $x=y+\gamma y^2$. 

The system $\dot x=-x^3$ is called \defwd{structurally unstable} because
different perturbations, even arbitrarily small, can lead to different
dynamics. It is called \defwd{singular of codimension 2} because every
(smooth and small) perturbation is equivalent to a member of the
2-parameter family \eqref{1D24}, which is called an \defwd{unfolding} (or a
\defwd{versal deformation}) of the singularity. 
\end{example}

This geometrical approach to the perturbation of vector fields is rather
well established for two-dimensional systems, where all structurally stable
vector fields, and all singularities of codimension 1 have been classified.
We will discuss this classification in Chapter 2. Little is known on
higher-dimensional systems, for which alternative methods have to be used. 


\section{Forced Oscillators}
\label{sec_inf}

In mechanical systems, solutions are not attracted by stable equilibria, but
oscillate around them. This makes the question of long-term stability an
extremely difficult one, which was only solved about 50 years ago by the
so-called Kolmogorov--Arnol'd--Moser (KAM) theory.

We introduce here the simplest example of this kind. 

\begin{example}
\label{ex_oscillator}
Consider an oscillator which is being forced by a small-amplitude periodic
term:
\begin{equation}
\label{os1}
\ddot x = f(x) + \eps \sin(\w t).
\end{equation}
We assume that $x=0$ is a stable equilibrium of the unperturbed equation,
that is, $f(0)=0$ and $f'(0)=-\w_0^2<0$. 

Let us first consider the linear case, i.e.\ the forced harmonic oscillator:
\begin{equation}
\label{os2}
\ddot x = -\w_0^2 x + \eps\sin(\w t).
\end{equation} 
The theory of linear ordinary differential equations states that the general
solution is obtained by superposition of a particular solution with the
general solution of the homogeneous system $\ddot x=-\w_0^2 x$, given by 
\begin{equation}
\label{os3}
x_{\math{h}}(t) = x_{\math{h}}(0) \cos(\w_0 t) + 
\frac{\dot x_{\math{h}}(0)}{\w_0} \sin(\w_0 t).
\end{equation}
Looking for a particular solution which is periodic with the same period as
the forcing term, we find 
\begin{equation}
\label{os4}
x_{\math{p}}(t) = \frac{\eps}{\w_0^2-\w^2}\sin(\w t),
\end{equation}
and thus the general solution of \eqref{os2} is given by
\begin{equation}
\label{os5}
x(t) = x(0)\cos(\w_0 t) + \frac 1{\w_0} \Bigbrak{\dot x(0) -
\frac{\eps\w}{\w_0^2-\w^2}} \sin(\w_0 t) + 
\frac{\eps}{\w_0^2-\w^2}\sin(\w t). 
\end{equation}
This expression is only valid for $\w\neq\w_0$, but it admits a limit as
$\w\to\w_0$, given by
\begin{equation}
\label{os6}
x(t) = \Bigbrak{x(0)-\frac{\eps t}{2\w_0}} \cos(\w_0 t) + 
\frac{\dot x(0)}{\w_0} \sin(\w_0 t).
\end{equation}
The solution \eqref{os5} is periodic if $\w/\w_0$ is rational, and
\defwd{quasiperiodic} otherwise. Its maximal amplitude is of the order
$\eps/(\w_0^2-\w^2)$ when $\w$ is close to $\w_0$. In the case $\w=\w_0$,
$x(t)$ is no longer quasiperiodic, but its amplitude grows proportionally to
$\eps t$. This effect, called \defwd{resonance}, is a first example showing
that the long-term behaviour of an oscillator can be fundamentally changed
by a small perturbation. 

Before turning to the nonlinear case, we remark that equation \eqref{os2}
takes a particularly simple form when written in the complex variable
$z=\w_0x+\icx y$. Indeed, $z$ satisfies the equation 
\begin{equation}
\label{os7}
\dot z = -\icx\w_0 z + \icx\eps\sin(\w t),
\end{equation}
which is easily solved by the method of variation of the constant. 

Let us now consider equation \eqref{os1} for more general functions $f$. If
$\eps=0$, the system has a constant of the motion. Indeed, let $V(x)$ be a
\defwd{potential function}, i.e.\ such that $V'(x)=-f(x)$ (for
definiteness, we take $V(0)=0$). Then 
\begin{equation}
\label{os8}
\dtot{}{t} \Bigpar{\frac12 \dot x^2 + V(x)} = \dot x \ddot x + V'(x) \dot x
= \dot x(\ddot x-f(x)) = 0. 
\end{equation}
By the assumptions on $f$, $V(x)$ behaves like $\frac12\w_0^2x^2$ near
$x=0$, and thus the level curves of $H(x,\dot x)=\frac12 \dot x^2 + V(x)$
are closed in the $(x,\dot x)$-plane near $x=\dot x=0$, which means that the
origin is stable. 

We can exploit this structure of phase space to introduce variables in which
the equations of motion take a simpler form. For instance if $x=x(H,\psi)$,
$y=y(H,\psi)$, $\psi\in[0,2\pi)$, is a parametrization of the level curves
of $H$, then $\dot H=0$ and $\dot\psi$ can be expressed as a function of $H$
and $\psi$. One can do even better, and introduce a parametrization by an
angle $\ph$ such that $\dot\ph$ is constant on each level curve. In other
words, there exists a coordinate transformation $(x,y)\mapsto(I,\ph)$ such
that 
\begin{equation}
\label{os9}
\begin{split}
\dot I &= 0 \\
\dot \ph &= \Omega(I). 
\end{split}
\end{equation}
$I$ and $\ph$ are called \defwd{action-angle variables}. $I=0$
corresponds to $x=y=0$ and $\Omega(0)=-\w_0$. The variable $z=I\e^{\icx\ph}$
satisfies the equation
\begin{equation}
\label{os10}
\dot z = \icx \Omega(\abs{z}) z
\qquad\Rightarrow\qquad
z(t) = \e^{\icx\Omega(\abs{z(0)})} z(0).
\end{equation}
Now if $\eps\neq0$, there will be an additional term in the equation,
depending on $z$ and its complex conjugate $\cc z$:
\begin{equation}
\label{os11}
\dot z = \icx \Omega(\abs z)z + \eps\sin(\w t) g(z,\cc z). 
\end{equation}
If $\Omega$ were a constant, one could try to solve this equation by
iterations as we did in Example \ref{ex_1D1}. However, the operator $\cT$
constructed in this way is not necessarily contracting because of the
imaginary coefficient in $\dot z$. In fact, it turns out that the method of
changes of variables works for certain special values of
$\Omega(\abs{z})$. 

\begin{figure}
 \centerline{\psfig{figure=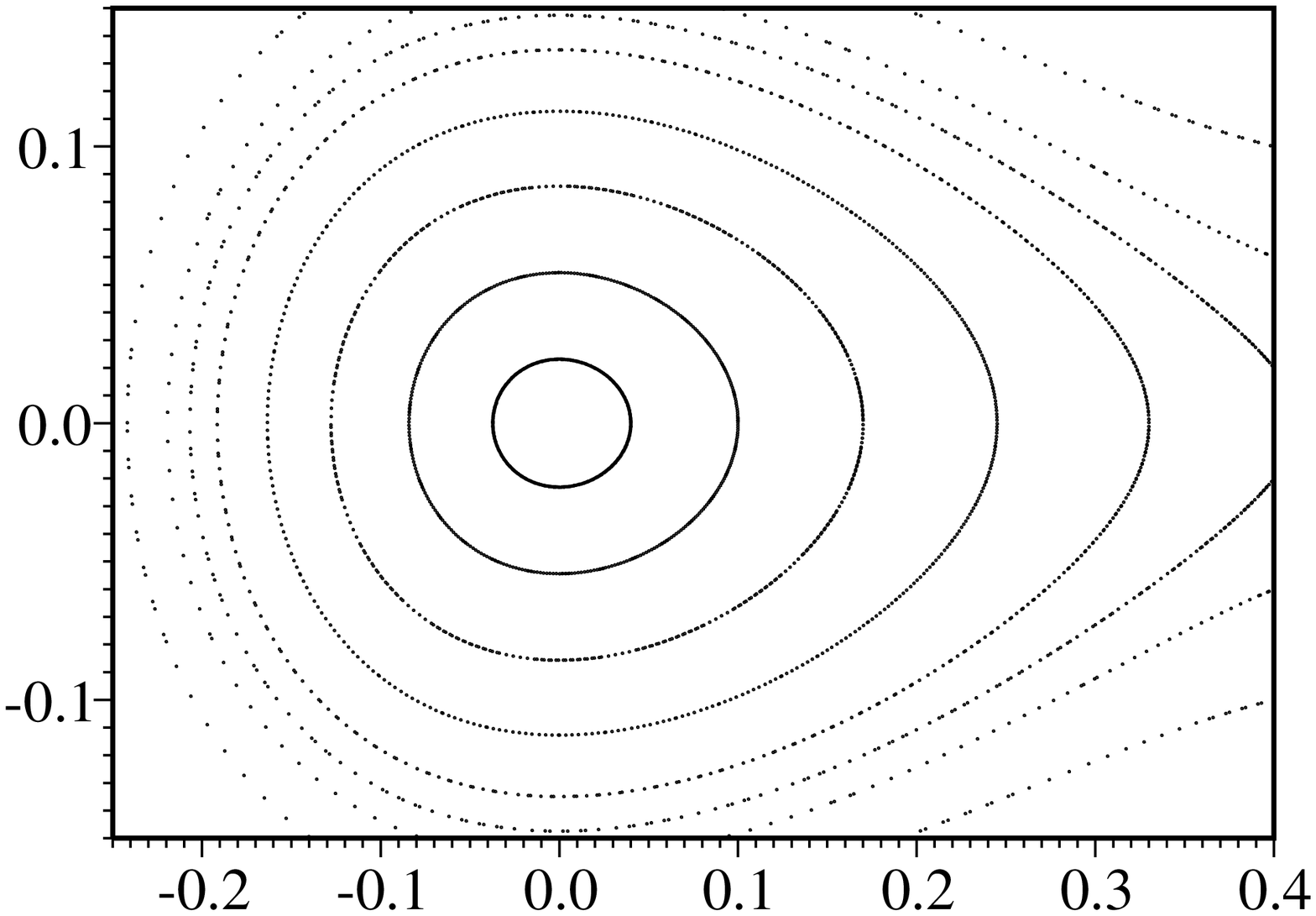,height=50mm}
 \psfig{figure=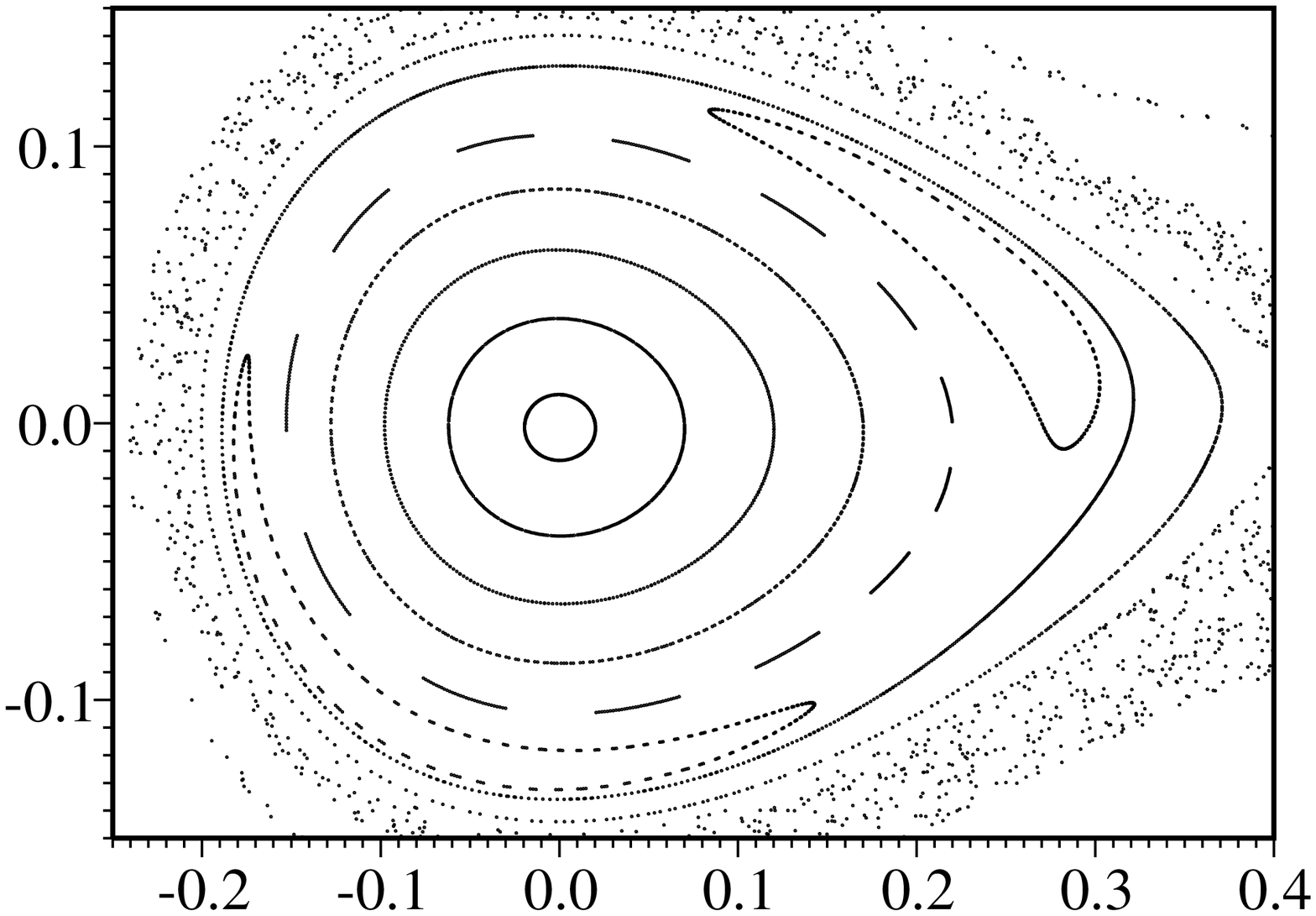,height=50mm}}
 \vspace{3mm}
 \centerline{\psfig{figure=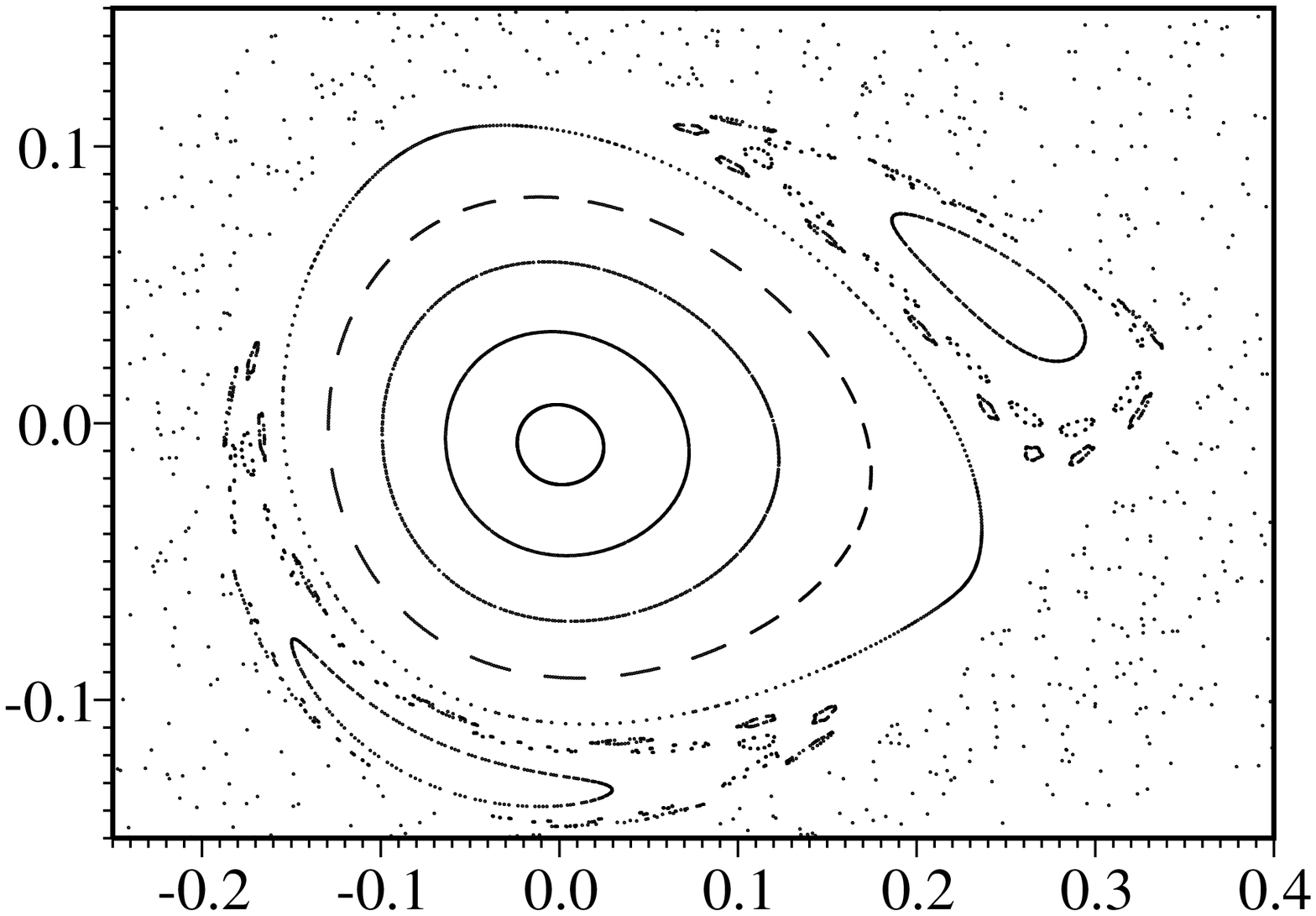,height=50mm}
 \psfig{figure=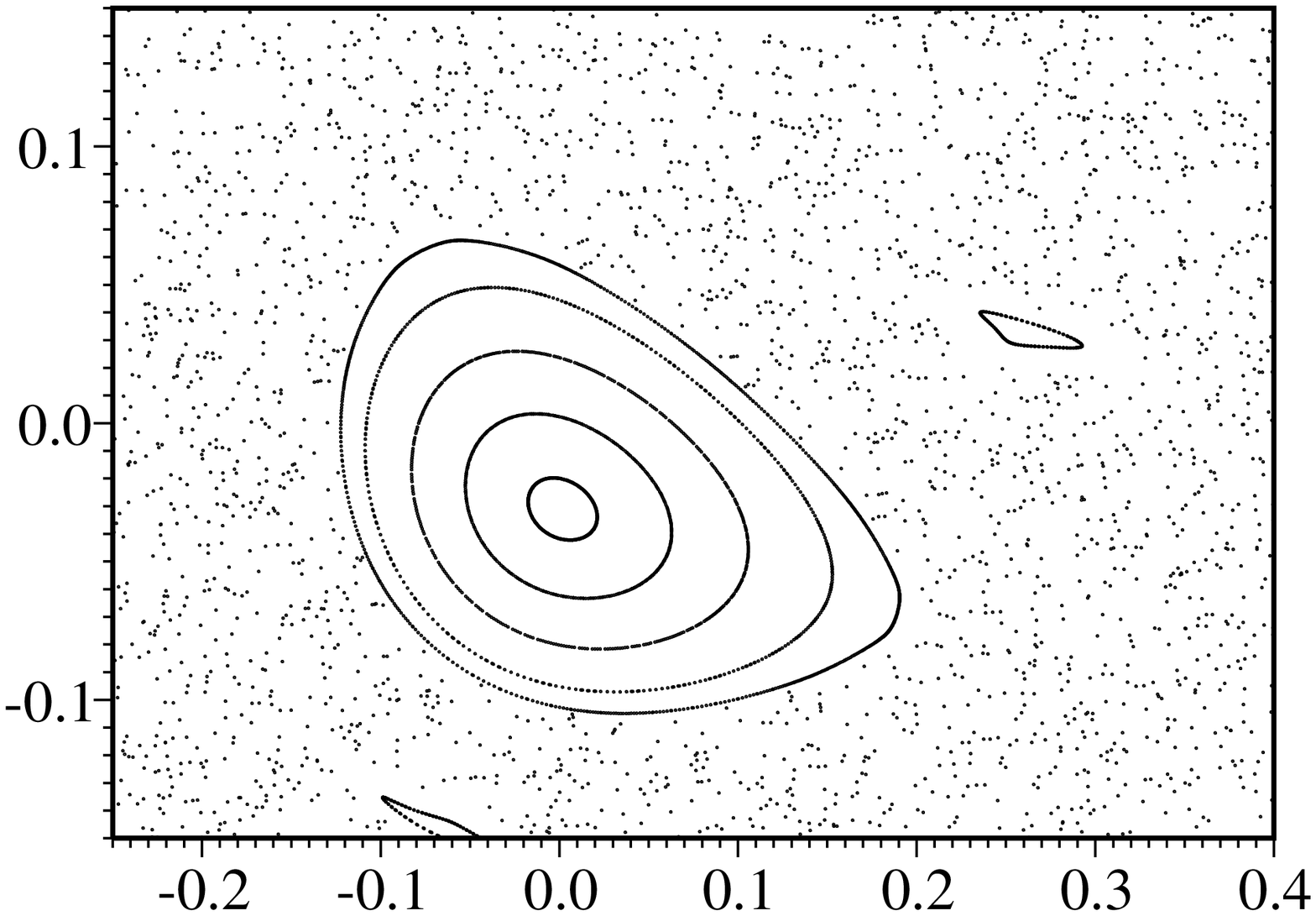,height=50mm}}
 \captionspace
 \caption[]{Poincar\'e maps of equation \eqref{os11} in the case
 $f(x)=-\w_0^2 x + x^2 - \frac12 x^3$, with $\w_0=0.6$ and for increasing
 values of $\eps$: from left to right and top to bottom, $\eps=0$,
 $\eps=0.001$, $\eps=0.005$ and $\eps=0.02$. }
 \label{fig_forced}
\end{figure}

Let us now find out what solutions of equation \eqref{os11} look like
qualitatively. Since the right-hand side is $2\pi/\w$-periodic in $t$, it
will be sufficient to describe the dynamics during one period. Indeed, if
we define the \defwd{Poincar\'e map}
\begin{equation}
\label{os12}
P_\eps: z(0) \mapsto z(2\pi/\w), 
\end{equation}
then $z$ at time $n(2\pi/\w)$ is simply obtained by applying $n$ times the
Poincar\'e map to $z(0)$. Let us now discuss a few properties of this map. 
If $\eps=0$, we know from \eqref{os10} that $\abs{z}$ is constant, and
thus 
\begin{equation}
\label{os13a}
P_0(z)=\e^{2\pi\icx\Omega(\abs z)/\w}z.
\end{equation}
In particular, 
\begin{equation}
\label{os13}
P_0(0)=0, \qquad P_0'(0) = \e^{-2\pi\icx \w_0/\w}. 
\end{equation}
For $\eps>0$, the implicit function theorem can be used to show that
$P_\eps$ has a fixed point near $0$, provided $\Omega/\w$ is not an integer
(which would correspond to a strong resonance). One obtains that for
$\Omega/\w\not\in\Z$, there exists a $z^\star(\eps)=\Order{\eps}$ such that 
\begin{equation}
\label{os14}
P_\eps(z^\star) = z^\star, \qquad
P_\eps'(z^\star) = \e^{2\pi\icx\th}, \qquad
\th = -\frac{\w_0}\w + \Order{\eps}.
\end{equation}
Thus if $\z=z-z^\star$, the dynamics will be described by a map of the form 
\begin{equation}
\label{os15}
\z\mapsto \e^{2\pi\icx\th}\z + G(\z,\cc \z,\eps), 
\end{equation}
where $G(\z,\cc \z,\eps)=\Order{\abs{\z}^2}$. In other words, the dynamics
is a small perturbation of a rotation in the complex plane. Although the map
\eqref{os15} looks quite simple, its dynamics can be surprisingly
complicated. 

\figref{fig_forced} shows phase portraits of the Poincar\'e map for a
particular function $f$ and increasing values of $\eps$. They are obtained
by plotting a large number of iterates of a few different initial
conditions. For $\eps=0$, these orbits live on level curves of the constant
of the motion $H$. As $\eps$ increases, some of these invariant curves seem
to survive, so that the point in their center (which is $z^\star$) is
surrounded by an \lq\lq elliptic island\rq\rq. This means that stable
motions still exist. However, more and more of the invariant curves are
destroyed, so that the island is gradually invaded by chaotic orbits. A
closer look at the island shows that it is not foliated by invariant
curves, but contains more complicated structures, also called resonances. 
\end{example}

One of the aims of Chapter 3 will be to develop methods to study perturbed
systems similar to \eqref{os1}, and to find conditions under which stable
equilibria of the unperturbed system remain stable when a small perturbation
is added. 


\section{Singular Perturbations: The Van der Pol Oscillator}
\label{sec_ins}

Singular perturbation theory considers systems of the form $\dot x =
f(x,\eps)$ in which $f$ behaves singularly in the limit $\eps\to0$. A simple
example of such a system is the Van der Pol oscillator in the large damping
limit. 

\begin{example}
\label{ex_vdPol}
The Van der Pol oscillator is a second order system with nonlinear damping,
of the form 
\begin{equation}
\label{vdP1}
\ddot x + \alpha(x^2-1)\dot x + x = 0.
\end{equation}
The special form of the damping (which can be realized by an electric
circuit) has the effect of decreasing the amplitude of large oscillations,
while increasing the amplitude of small oscillations. 

\begin{figure}
 \centerline{\psfig{figure=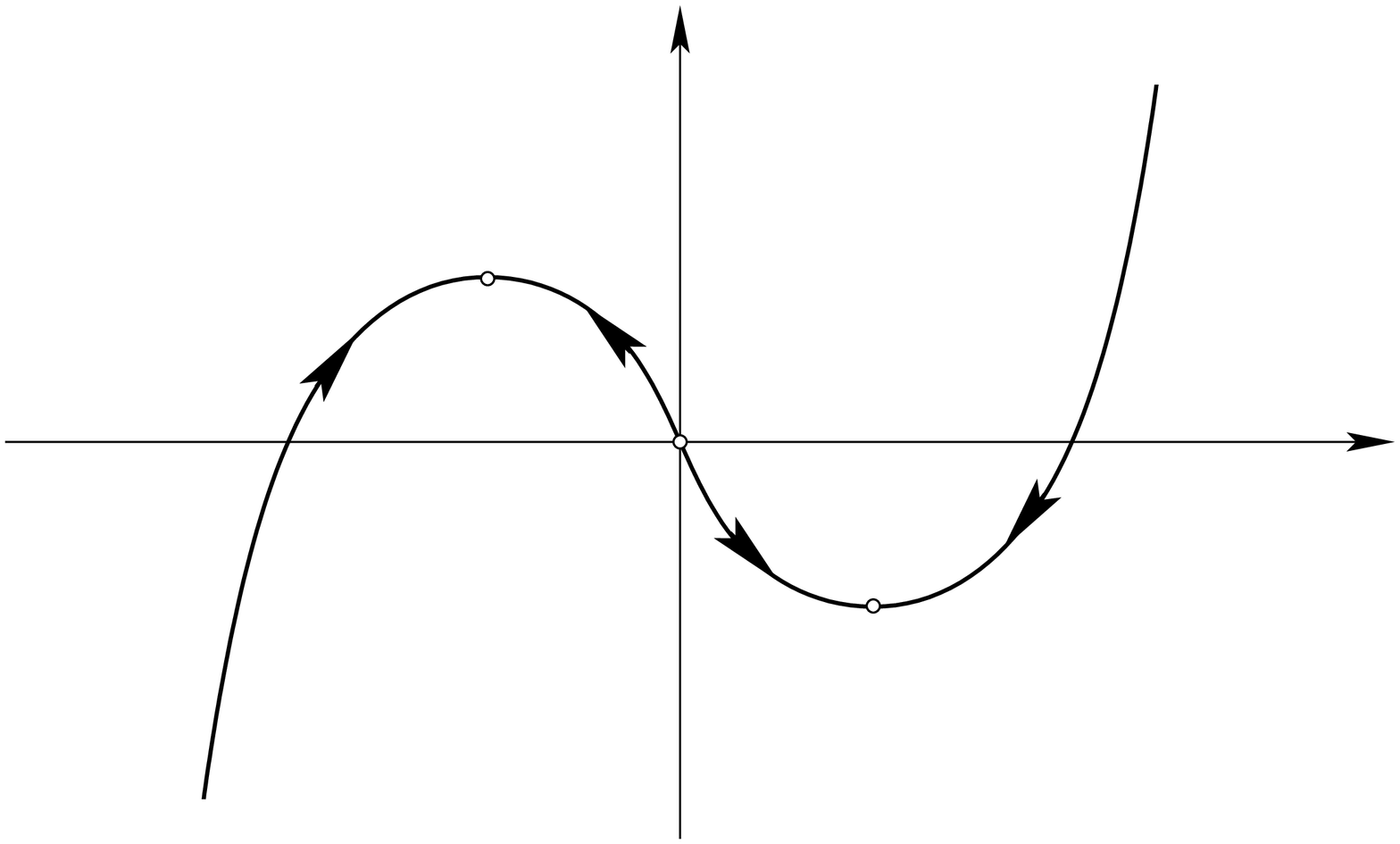,height=40mm}
 \hspace{5mm}
 \psfig{figure=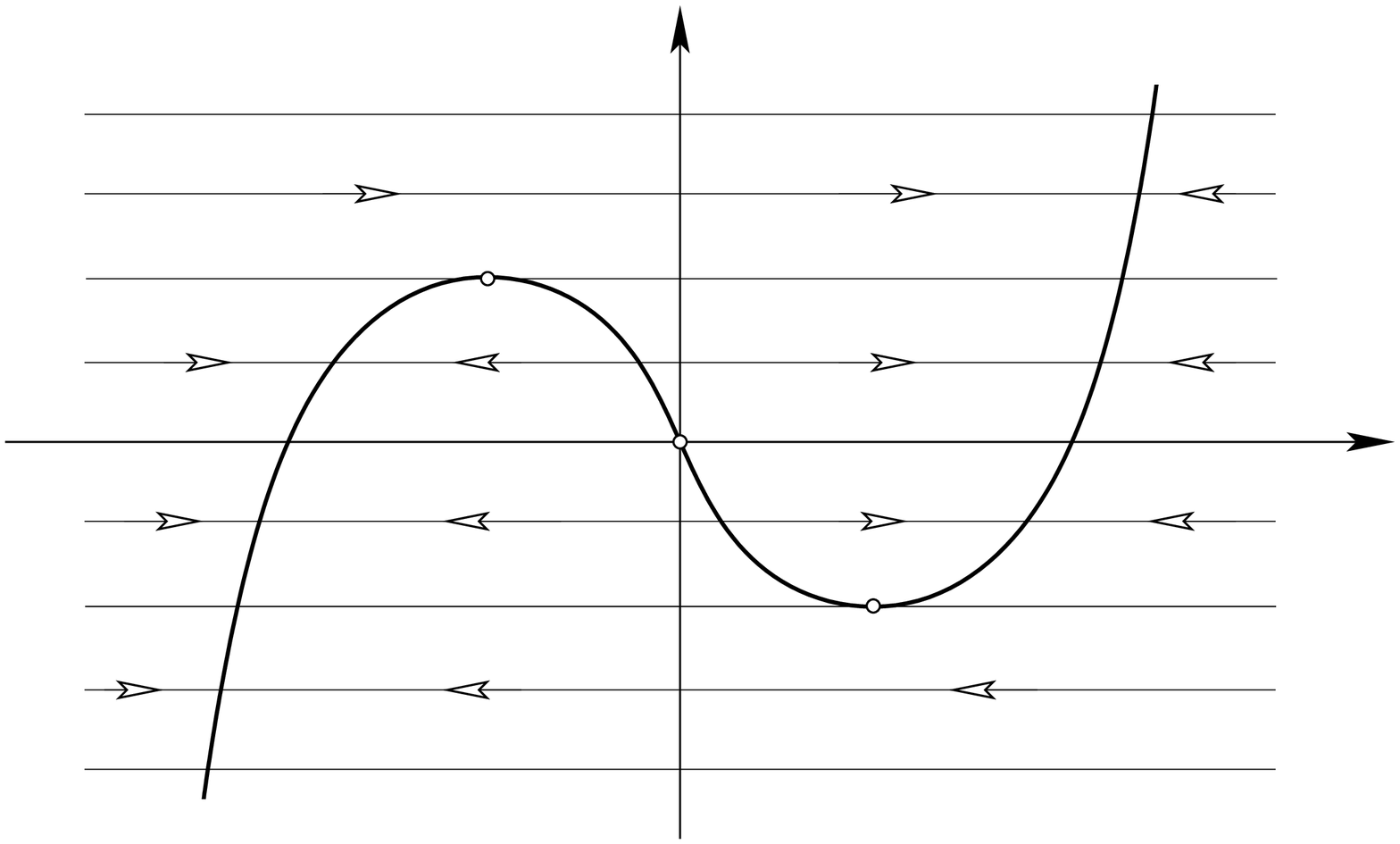,height=40mm}}
 \figtext{
 	\writefig	0.0	4.2	(a)
 	\writefig	3.9	4.1	$y$
 	\writefig	6.5	2.1	$x$
 	\writefig	7.2	4.2	(b)
 	\writefig	11.1	4.1	$y$
 	\writefig	13.7	2.1	$x$
 }
 \captionspace
 \caption[]{Behaviour of the Van der Pol equation in the singular limit
 $\eps\to0$, (a) on the slow time scale $t'=\sqrt\eps \mskip1.5mu t$, given by
 \eqref{vdP4}, and (b) on the fast time scale $t''=t/\sqrt\eps$, see
 \eqref{vdP6}.}
 \label{fig_vdP1}
\end{figure}

We are interested in the behaviour for large $\alpha$. There are several
ways to write \eqref{vdP1} as a first order system. For our purpose, a
convenient representation is 
\begin{equation}
\label{vdP2}
\begin{split}
\dot x &= \alpha \Bigpar{y+x-\frac{x^3}3} \\
\dot y &= -\frac x\alpha. 
\end{split}
\end{equation}
One easily checks that this system is equivalent to \eqref{vdP1} by
computing $\ddot x$. If $\alpha$ is very large, $x$ will move quickly,
while $y$ changes slowly. To analyse the limit $\alpha\to\infty$, we
introduce a small parameter $\eps=1/\alpha^2$ and a \lq\lq slow time\rq\rq\
$t'=t/\alpha=\sqrt\eps \mskip1.5mu t$. Then \eqref{vdP2} can be rewritten
as
\begin{equation}
\label{vdP3}
\begin{split}
\eps\dtot{x}{t'} &= y + x - \frac{x^3}3 \\
\dtot{y}{t'} &= -x. 
\end{split}
\end{equation}
In the limit $\eps\to0$, we obtain the system
\begin{equation}
\label{vdP4}
\begin{split}
0 &= y + x - \frac{x^3}3 \\
\dtot{y}{t'} &= -x, 
\end{split}
\end{equation}
which is no longer a system of differential equations. In fact, the
solutions are constrained to move on the curve $\cC: y=\frac13 x^3 - x$,
and eliminating $y$ from the system we have 
\begin{equation}
\label{vdP5}
-x = \dtot y{t'} = (x^2-1) \dtot x{t'} 
\qquad\Rightarrow\qquad 
\dtot x{t'} = -\frac x{x^2-1}. 
\end{equation}
The dynamics is shown in \figref{fig_vdP1}a. 
Another possibility is to introduce the \lq\lq fast time\rq\rq\
$t''=\alpha t=t/\sqrt\eps$. Then \eqref{vdP2} becomes
\begin{equation}
\label{vdP6}
\begin{split}
\dtot x{t''} &= y + x - \frac{x^3}3 \\
\dtot y{t''} &= -\eps x.
\end{split}
\end{equation}
In the limit $\eps\to 0$, we get the system
\begin{equation}
\label{vdP7}
\begin{split}
\dtot x{t''} &= y + x - \frac{x^3}3 \\
\dtot y{t''} &= 0.
\end{split}
\end{equation}
In this case, $y$ is a constant and acts as a parameter in the equation for
$x$. Some orbits are shown in \figref{fig_vdP1}b. 

\begin{figure}
 \centerline{\psfig{figure=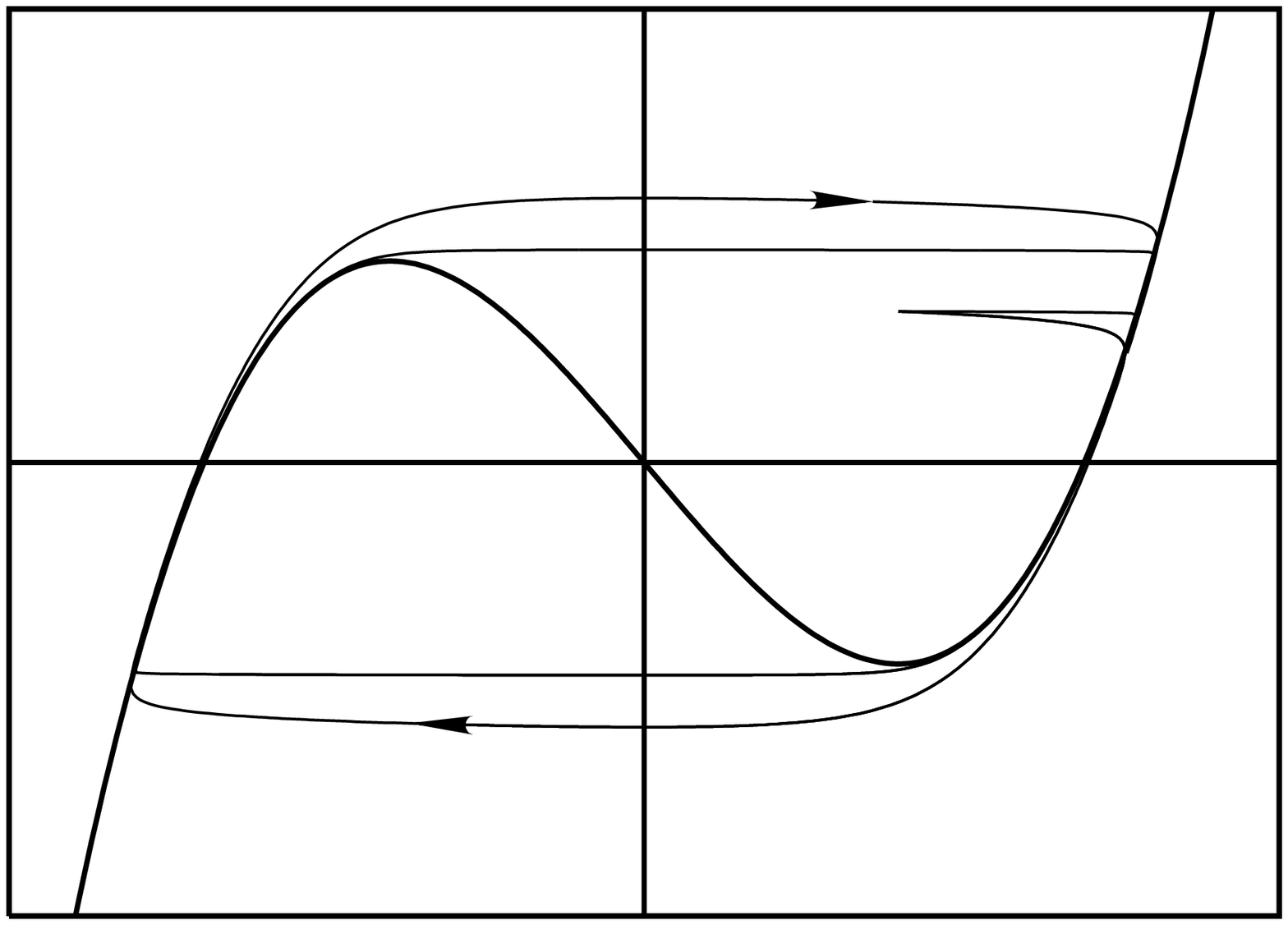,height=40mm}
 \hspace{8mm}
 \psfig{figure=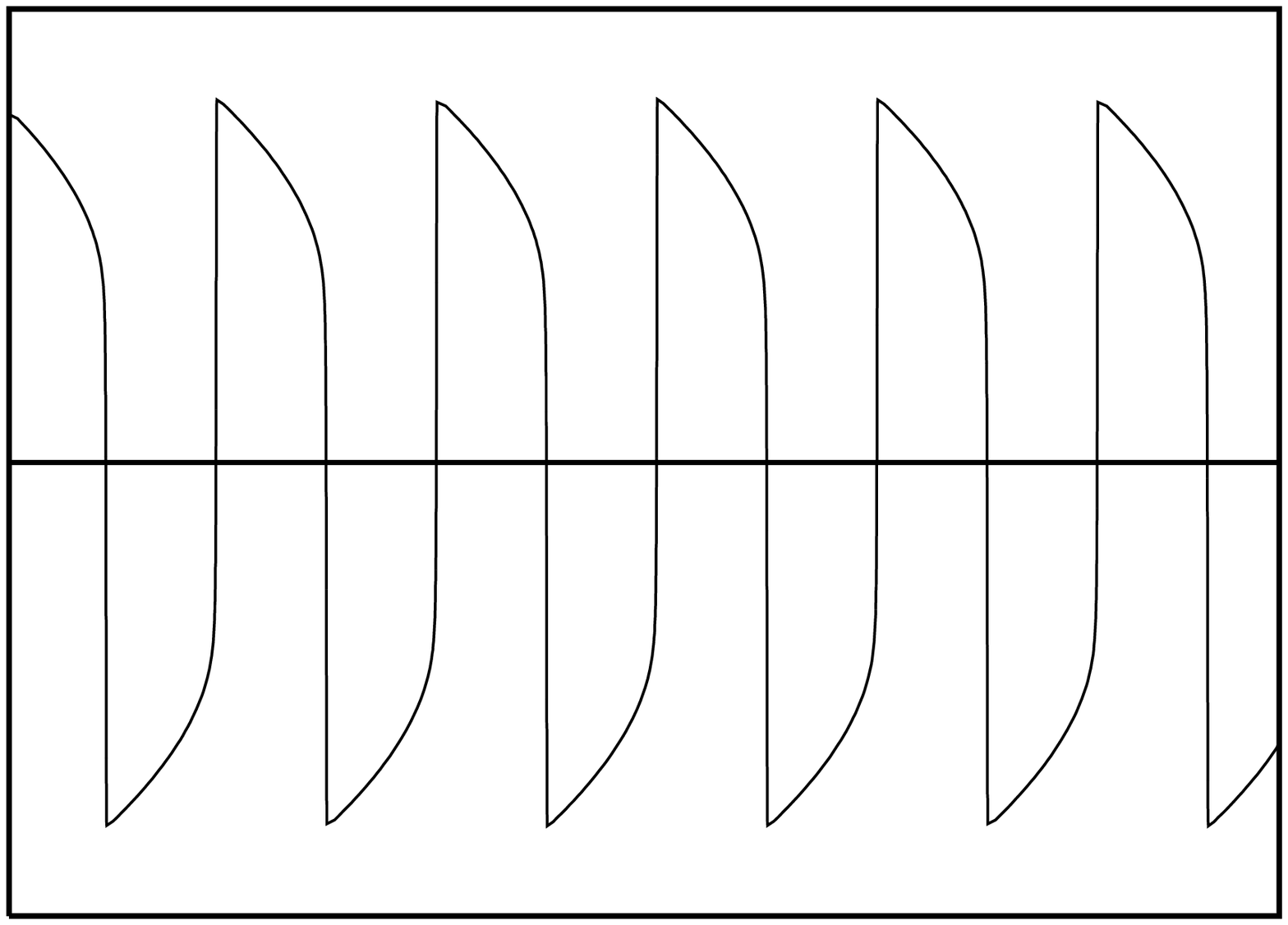,height=40mm}}
 \figtext{
 	\writefig	0.6	4.2	(a)
 	\writefig	3.8	4.2	$y$
 	\writefig	6.5	2.2	$x$
 	\writefig	7.2	4.2	(b)
 	\writefig	8.0	4.2	$x$
 	\writefig	13.7	2.4	$t$
 }
 \captionspace
 \caption[]{(a) Two solutions of the Van der Pol equations \eqref{vdP2}
 (light curves) for the same initial condition $(1,0.5)$, for $\alpha=5$
 and $\alpha=20$. The heavy curve is the curve $\cC: y=\frac13 x^3 - x$.
 (b) The graph of $x(t)$ ($\alpha=20$) displays relaxation oscillations.}
 \label{fig_vdP2}
\end{figure}

Of course, the systems \eqref{vdP2}, \eqref{vdP3} and \eqref{vdP6} are
strictly equivalent for $\eps>0$. They only differ in the singular limit
$\eps\to0$.  The dynamics for small but positive $\eps$ can be understood
by sketching the vector field. Let us note that
\begin{itemiz}
\item	$\dot x$ is positive if $(x,y)$ lies above the curve $\cC$ and
negative when it lies below; this curve separates the plane into regions
where $x$ moves to the right or to the left, and the orbit must cross the
curve vertically;
\item	$\tdtot yx$ is very small unless the orbit is close to the curve
$\cC$, so that the orbits will be almost horizontal except near this curve;
\item	orbits move upward if $x<0$ and downward if $x>0$. 
\end{itemiz}
The resulting orbits are shown in \figref{fig_vdP2}a. An orbit starting
somewhere in the plane will first approach the curve $\cC$ on a nearly
horizontal path, in a time $t$ of order $1/\alpha$. Then it will track the
curve at a small distance until it reaches a turning point, after a time
$t$ of order $\alpha$. Since the equations forbid it to follow $\cC$ beyond
this point, the orbit will jump to another branch of $\cC$, where the
behaviour repeats. The graph of $x(t)$ thus contains some parts with a
small slope and others with a large slope (\figref{fig_vdP2}b). This
phenomenon is called a \defwd{relaxation oscillation}.
\end{example}

The system \eqref{vdP3} is called a \defwd{slow-fast} system, and is a
particular case of a singularly perturbed equation. We will discuss some
properties of these systems in Chapter 4. 


\chapter{Bifurcations and Unfolding}
\label{ch_bu}

The qualitative theory of two-dimensional vector fields is relatively well
developed, since all structurally stable systems and all singularities of
codimension $1$ have been identified. Thus we know which systems are
well-behaved under small perturbations, and most of the other systems are
qualitatively well understood. 

In this chapter we will consider ordinary differential equations of the form
\begin{equation}
\label{bu1}
\dot x = f(x),
\qquad f\in\cC^k(\cD,\R^2),
\end{equation}
where $\cD$ is an open domain in $\R^2$ and the differentiability $k$ is at
least $1$ (we will often need a higher differentiability). $f$ is called a
\defwd{vector field}. Equivalently, equation \eqref{bu1} can be written in
components as 
\begin{equation}
\label{bu2}
\begin{split}
\dot x_1 &= f_1(x_1,x_2) \\
\dot x_2 &= f_2(x_1,x_2). 
\end{split}
\end{equation}
In fact, one may also consider the case where $\cD$ is a compact,
two-dimensional manifold, such as the $2$-sphere $\fS^2$ or the $2$-torus
$\T^2$. This actually simplifies some of the problems, but the main results
are essentially the same. 

We will start by establishing some properties of invariant sets of the ODE
\eqref{bu1}, before proceeding to the classification of structurally stable
and unstable vector fields. 


\section{Invariant Sets of Planar Flows}
\label{sec_bui}

Results from basic analysis show that the condition $f\in\cC^1(\cD,\R^2)$
is sufficient to guarantee the existence of a unique solution $x(t)$ of
\eqref{bu1} for every initial condition $x(0)\in\cD$. This solution can be
extended to a maximal open $t$-interval $I\ni0$; $I$ may not be equal to
$\R$, but in this case $x(t)$ must tend to the boundary of $\cD$ as $t$
approaches the boundary of $I$.

\begin{definition}
\label{def_bui1}
Let $x_0\in\cD$ and let $x:I\to\cD$ be the unique solution of \eqref{bu1}
with initial condition $x(0)=x_0$ and maximal interval of existence $I$. 
\begin{itemiz}
\item	the \defwd{orbit} through $x_0$ is the set
$\setsuch{x\in\cD}{x=x(t),t\in I}$; 
\item	the \defwd{flow} of \eqref{bu1} is the map 
\begin{equation}
\label{bui1}
\ph : (x_0,t) \mapsto \ph_t(x_0) = x(t).
\end{equation}
\end{itemiz}
\end{definition} 

The uniqueness property implies that 
\begin{equation}
\label{bui2}
\ph_0(x_0) = x_0 
\qquad\qquad\text{and}\qquad\qquad
\ph_t(\ph_s(x_0)) = \ph_{t+s}(x_0)
\end{equation}
for all $x_0\in\cD$ and all $t,s\in\R$ for which these quantities are
defined. 

\begin{definition}
\label{def_bui2}
A set $\cS\in\cD$ is called 
\begin{itemiz}
\item	\defwd{positively invariant} if $\ph_t(\cS)\subset\cS$ for all
$t\geqs0$;
\item	\defwd{negatively invariant} if $\ph_t(\cS)\subset\cS$ for all
$t\leqs0$;
\item	\defwd{invariant} if $\ph_t(\cS)=\cS$ for all $t\in\R$.
\end{itemiz}
\end{definition}

Note that if $\cS$ is positively invariant, then $\ph_t(x_0)$ exists for
all $t\in[0,\infty)$ whenever $x_0\in\cS$, and similar properties hold in
the other cases. A sufficient condition for $\cS$ to be positively
invariant is that the vector field $f$ be directed inward $\cS$ on the
boundary $\partial\cS$. 

We shall now introduce some particularly important invariant sets. 

\begin{definition}
\label{def_bui3}
Assume that $\ph_t(x)$ exists for all $t\geqs0$. The \defwd{$\w$-limit set}
of $x$ for $\ph_t$, denoted $\w(x)$, is the set of points $y\in\cD$ such
that there exists a sequence $\set{t_n}_{n\geqs0}$ such that $t_n\to\infty$
and  $\ph_{t_n}(x)\to y$ as $n\to\infty$. 
The \defwd{$\alpha$-limit set} $\alpha(x)$ is defined in a similar way,
with $t_n\to-\infty$. 
\end{definition}

The terminology is due to the fact that $\alpha$ and $\w$ are the first and
last letters of the Greek alphabet. 

\begin{prop}
\label{prop_bui1}
Let $\cS$ be positively invariant and compact. Then for any $x\in\cS$, 
\begin{enum}
\item	$\w(x)\neq\emptyset$;
\item	$\w(x)$ is closed;
\item	$\w(x)$ is invariant under $\ph_t$;
\item	$\w(x)$ is connected. 
\end{enum}
\end{prop}
\begin{proof}\hfill
\begin{enum}
\item	Choose a sequence $t_n\to\infty$ and set $x_n=\ph_{t_n}(x)$. Since
$\cS$ is compact, $\set{x_n}$ has a convergent subsequence, and its limit
is a point in $\w(x)$. 
\item	Pick any $y\notin\w(x)$. There must exist a \nbh\ $\cU$ of $y$ and
$T\geqs0$ such that $\ph_t(x)\notin\cU$ for all $t\geqs T$. This shows that
the complement of $\w(x)$ is open, and thus $\w(x)$ is closed. 
\item	Let $y\in\w(x)$ and let $t_n\to\infty$ be an increasing sequence
such that $\ph_{t_n}(x)\to y$ as $n\to\infty$. For any $s\in\R$,
\eqref{bui2} implies that $\ph_s(\ph_{t_n}(x))=\ph_{s+t_n}(x)$ exists for
all $s\geqs -t_n$. Taking the limit $n\to\infty$, we obtain by continuity
that $\ph_s(y)$ exists for all $s\in\R$. Moreover, since $\ph_{s+t_n}(x)$
converges to $\ph_s(y)$ as $t_n\to\infty$, we conclude that
$\ph_s(y)\in\w(x)$, and therefore $\w(x)$ is invariant. 
\item	Suppose by contradiction that $\w(x)$ is not connected. Then there
exist two open sets $\cU_1$ and $\cU_2$ such that $\overline\cU_1 \cap
\overline\cU_2 = \emptyset$, $\w(x) \subset \cU_1\cup\cU_2$, and
$\w(x)\cap\cU_i\neq\emptyset$ for $i=1,2$. By continuity of $\ph_t(x)$,
given $T>0$, there exists $t>T$ such that
$\ph_t(x)\in\cK\defby\cS\setminus(\cU_1\cup\cU_2)$. We can thus find a
sequence $t_n\to\infty$ such that $\ph_{t_n}(x)\in\cK$, which must admit a
subsequence converging to some $y\in\cK$. But this would imply $y\in\w(x)$,
a contradiction. 
\qed
\end{enum}
\renewcommand{\qed}{}
\end{proof}

Typical examples of $\w$-limit sets are attracting equilibrium points and
periodic orbits. We will start by examining these in more detail before
giving the general classification of $\w$-limit sets. 


\subsection{Equilibrium Points}
\label{ssec_bueq}

\begin{definition}
\label{def_bueq1}
$x^\star\in\cD$ is called an \defwd{equilibrium point} (or a \defwd{singular
point}) of the system $\dot x=f(x)$ if $f(x^\star)=0$. Equivalently, we have
$\ph_t(x^\star)=x^\star$ for all $t\in\R$, and hence $\set{x^\star}$ is
invariant. 
\end{definition}

We would like to characterize the dynamics near $x^\star$. Assuming
$f\in\cC^2$, we can Taylor-expand $f$ to second order around $x^\star$: if
$y=x-x^\star$, then 
\begin{equation}
\label{bueq1}
\dot y = f(x^\star+y) = Ay + g(y),
\end{equation}
where $A$ is the \defwd{Jacobian matrix} of $f$ at $x^\star$, defined by 
\begin{equation}
\label{bueq2}
A = \dpar fx(x^\star) \defby
\begin{pmatrix}
\vrule height 10pt depth 14pt width 0pt
\displaystyle \dpar{f_1}{x_1}(x^\star) & 
\displaystyle \dpar{f_1}{x_2}(x^\star) \\
\vrule height 12pt depth 10pt width 0pt
\displaystyle \dpar{f_2}{x_1}(x^\star) & 
\displaystyle \dpar{f_2}{x_2}(x^\star)
\end{pmatrix}
\end{equation}
and $g(y)=\Order{\norm{y}^2}$ as $y\to 0$. The matrix $A$ has eigenvalues
$a_1$, $a_2$, which are either both real, or complex conjugate. 

\begin{definition}
\label{def_bueq2}
The equilibrium point $x^\star$ is called
\begin{itemiz}
\item	\defwd{hyperbolic} if $\re a_1\neq 0$ and $\re a_2\neq 0$;
\item	\defwd{non-hyperbolic} if $\re a_1 = 0$ or $\re a_2 = 0$;
\item	\defwd{elliptic} if $\re a_1 =\re a_2 = 0$ (but one usually requires
$\im a_1\neq 0$ and $\im a_2\neq 0$).
\end{itemiz}
\end{definition}

\begin{figure}
 \centerline{\psfig{figure=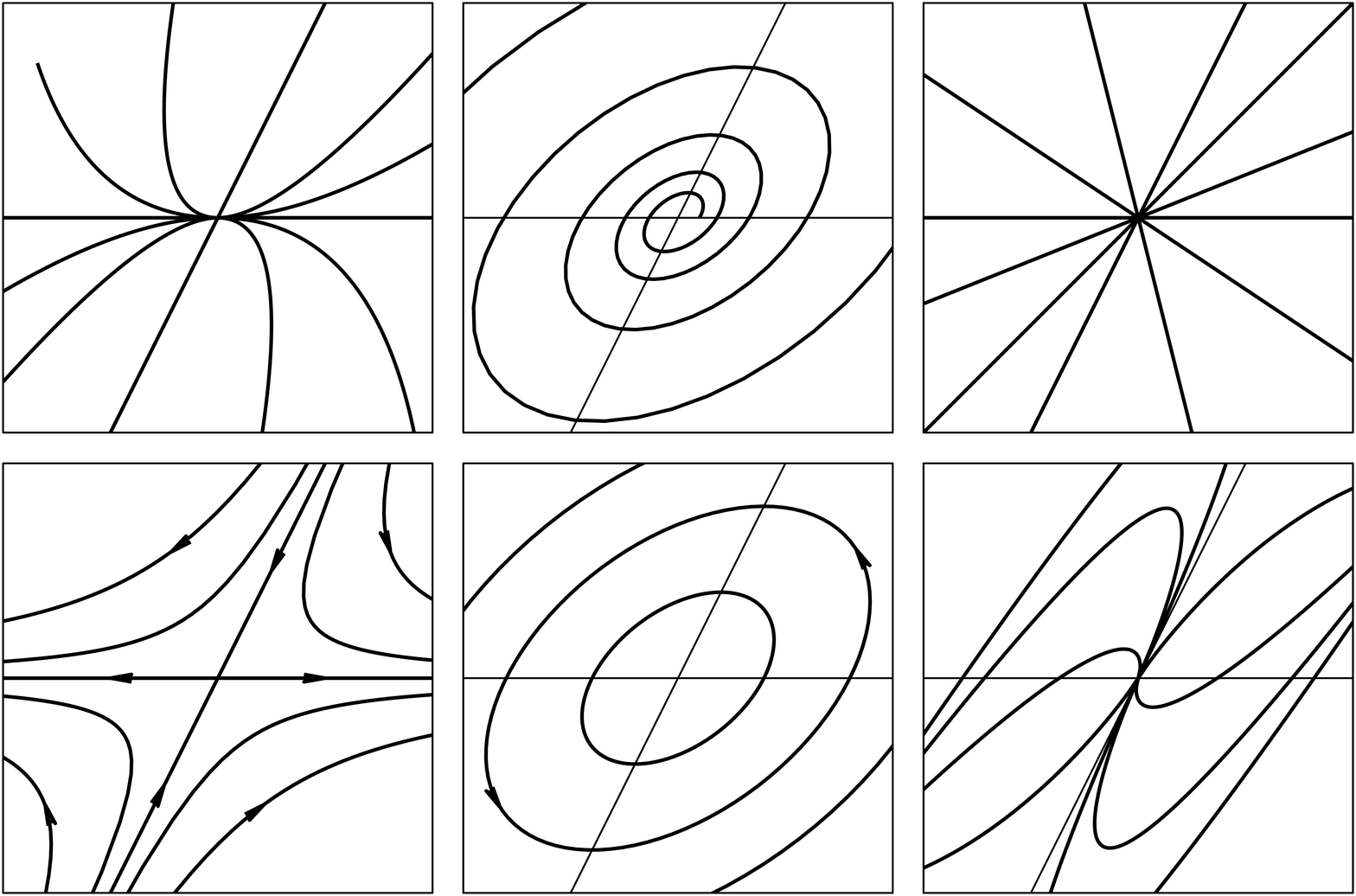,height=60mm,clip=t}}
 \figtext{
 	\writefig	3.0	6.2	a
 	\writefig	3.0	3.0	b
 	\writefig	6.1	6.2	c
 	\writefig	6.1	3.0	d
 	\writefig	9.2	6.2	e
 	\writefig	9.2	3.0	f
 }
 \captionspace
 \caption[Phase portraits of a linear two--dimensional system]
 {Phase portraits of a linear two--dimensional system: (a) node,
 (b) saddle, (c) focus, (d) center, (e) degenerate node, (f) improper node.}
\label{fig_slin}
\end{figure}

In order to understand the meaning of these definitions, let us start by
considering the \defwd{linearization} of \eqref{bueq1}: 
\begin{equation}
\label{bueq3}
\dot y = Ay
\qquad\qquad\Rightarrow\qquad\qquad
y(t) = \e^{At} y(0),
\end{equation}
where the \defwd{exponential of $At$} is defined by the absolutely
convergent series
\begin{equation}
\label{bueq4}
\e^{At} \defby \sum_{k=0}^\infty \frac{t^k}{k!} A^k. 
\end{equation}
The simplest way to compute $\e^{At}$ is to choose a basis in which $A$ is
in Jordan canonical form. There are four qualitatively different (real)
canonical forms:
\begin{equation}
\label{bueq5}
\begin{pmatrix}
a_1 & 0 \\ 0 & a_2
\end{pmatrix}
\qquad
\begin{pmatrix}
a & -\w \\ \w & a
\end{pmatrix}
\qquad
\begin{pmatrix}
a & 1 \\ 0 & a
\end{pmatrix}
\qquad
\begin{pmatrix}
a & 0 \\ 0 & a
\end{pmatrix},
\end{equation}
where we assume that $a_1\neq a_2$ and $\w\neq0$. They correspond to the
following situations:
\begin{enum}
\item	If $a_1$ and $a_2$ are both real and different, then 
\begin{equation}
\label{bueq6a}
\e^{At} = 
\begin{pmatrix}
\e^{a_1t} & 0 \\ 0 & \e^{a_2t}
\end{pmatrix}
\qquad
\Rightarrow
\qquad
\begin{matrix}
y_1(t) = \e^{a_1t} y_1(0) \phantom{.}\\
y_2(t) = \e^{a_2t} y_2(0).
\end{matrix}
\end{equation}
The orbits are curves of the form $y_2=cy_1^{a_2/a_1}$, and $x^\star$ is
called a \defwd{node} if $a_1$ and $a_2$ have the same sign, and a
\defwd{saddle} if they have opposite signs (\figref{fig_slin}a and b). 

\item	If $a_1 = \cc{a_2} = a+\icx\w$ with $\w\neq0$, then 
\begin{equation}
\label{bueq6b}
\e^{At} = \e^{at}
\begin{pmatrix}
\cos\w t & -\sin\w t \\ \sin\w t & \cos\w t
\end{pmatrix}
\quad
\Rightarrow
\quad
\begin{matrix}
y_1(t) = \e^{at} (y_1(0) \cos\w t - y_2(0) \sin\w t) \phantom{.}\\
y_2(t) = \e^{at} (y_1(0) \sin\w t + y_2(0) \cos\w t).
\end{matrix}
\end{equation}
The orbits are spirals or ellipses, and $x^\star$ is called a \defwd{focus}
is $a\neq0$ and a \defwd{center} if $a=0$ (\figref{fig_slin}c and d). 

\item	If $a_1=a_2=a$ and $A$ admits two independent eigenvectors, then 
\begin{equation}
\label{bueq6c}
\e^{At} = \e^{at}
\begin{pmatrix}
1 & 0 \\ 0 & 1
\end{pmatrix}
\qquad
\Rightarrow
\qquad
\begin{matrix}
y_1(t) = \e^{at} y_1(0) \phantom{.}\\
y_2(t) = \e^{at} y_2(0).
\end{matrix}
\end{equation}
The orbits are straight lines, and $x^\star$ is called a \defwd{degenerate
node} (\figref{fig_slin}e). 

\item	If $a_1=a_2=a$ but $A$ admits only one independent eigenvector, then 
\begin{equation}
\label{bueq6d}
\e^{At} = \e^{at}
\begin{pmatrix}
1 & t \\ 0 & 1
\end{pmatrix}
\qquad
\Rightarrow
\qquad
\begin{matrix}
y_1(t) = \e^{at}(y_1(0) + y_2(0)t) \phantom{.}\\
y_2(t) = \e^{at} y_2(0).\phantom{({}+ y_2(0))}
\end{matrix}
\end{equation}
$x^\star$ is called an \defwd{improper node} (\figref{fig_slin}f). 
\end{enum}

The above discussion shows that the asymptotic behaviour depends only on
the sign of the real parts of the eigenvalues of $A$. This property carries
over to the nonlinear equation \eqref{bueq1} if  $x^\star$ is hyperbolic,
while nonhyperbolic equilibria are more sensitive to nonlinear
perturbations. We will demonstrate this by considering first the case of
nodes/foci, and then the case of saddles. 

\begin{prop}
\label{prop_bueq1}
Assume that $\re a_1<0$ and $\re a_2<0$. Then there exists a \nbh\ $\cU$ of
$x^\star$ such that $\w(x)=x^\star$ for all $x\in\cU$. 
\end{prop}
\begin{proof}
The proof is a particular case of a theorem due to Liapunov. Consider the
case where $a_1$ and $a_2$ are real and different. Then system \eqref{bueq1}
can be written in appropriate coordinates as
\[
\begin{split}
\dot y_1 &= a_1 y_1 + g_1(y) \\
\dot y_2 &= a_2 y_2 + g_2(y),
\end{split}
\]
where $\abs{g_j(y)}\leqs M\norm{y}^2$, $j=1,2$, for sufficiently
small $\norm{y}$ and for some constant $M>0$. Given a solution $y(t)$ with
initial condition $y(0)=y_0$, we define the (Liapunov) function
\[
V(t) = \norm{y(t)}^2 = y_1(t)^2 + y_2(t)^2.
\]
Differentiating with respect to time, we obtain 
\[
\begin{split}
\dot V(t) &= 2 \bigbrak{y_1(t)\dot y_1(t) + y_2(t)\dot y_2(t)} \\
&= 2 \bigbrak{a_1 y_1(t)^2 + a_2 y_2(t)^2 + y_1(t) g_1(y(t)) +
y_2(t) g_2(y(t))}.
\end{split}
\]
If, for instance, $a_1>a_2$, we conclude from the properties of the $g_j$
that 
\[
\dot V(t) \leqs \bigbrak{a_1 + M(\abs{y_1(t)} + \abs{y_2(t)})} V(t).
\]
Since $a_1<0$, there exists a constant $r>0$ such that term in brackets is
strictly negative for $0\leqs\norm{y(t)}\leqs r$. This shows that the disc
of radius $r$ is positively invariant. Any solution starting in this disc
must converge to $y=0$, since otherwise we would contradict the fact that
$\dot V<0$ if $y\neq0$. 
The other cases can be proved in a similar way. 
\end{proof}

\begin{cor}
\label{cor_bueq}
Assume that $\re a_1>0$ and $\re a_2>0$. Then there exists a \nbh\ $\cU$ of
$x^\star$ such that $\alpha(x)=x^\star$ for all $x\in\cU$. 
\end{cor}
\begin{proof}
By changing $t$ into $-t$, we transform the situation into the previous one.
\end{proof}

\begin{figure}
 \centerline{\psfig{figure=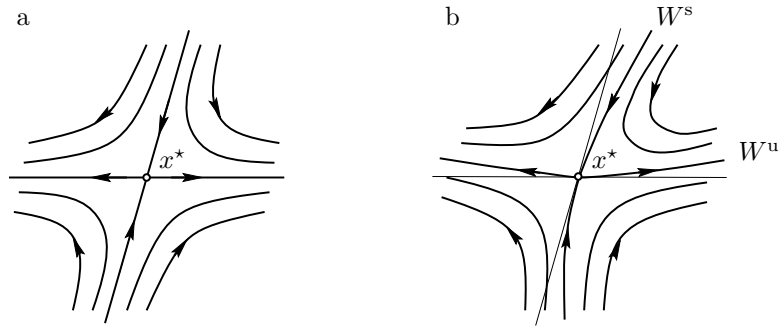,height=40mm,clip=t}}
 \figtext{
 	\writefig	2.7	4.5	a
 	\writefig	8.4	4.5	b
 	\writefig	4.6	2.6	$x^\star$
 	\writefig	10.35	2.6	$x^\star$
 	\writefig	12.3	2.7	$\Wglo{u}$
 	\writefig	11.2	4.5	$\Wglo{s}$
 }
 \captionspace
 \caption[]
 {Orbits near a hyperbolic fixed point: (a) orbits of the
 linearized system, (b) orbits of the nonlinear system with local stable
 and unstable manifolds.}
\label{fig_hyperbol}
\end{figure}

The remaining hyperbolic case is the saddle point. In the linear case, we
saw that the eigenvectors define two invariant lines, on which the motion is
either contracting or expanding. This property can be generalized to the
nonlinear case.

\begin{theorem}[Stable Manifold Theorem]
\label{thm_bueq}
Assume that $a_1<0<a_2$. Then
\begin{itemiz}
\item	there exists a curve $\Wglo{s}(x^\star)$, tangent in $x^\star$ to
the eigenspace of $a_1$, such that $\w(x)=x^\star$ for all
$x\in\Wglo{s}(x^\star)$; $\Wglo{s}$ is called the \defwd{stable manifold} of
$x^\star$;
\item	there exists a curve $\Wglo{u}(x^\star)$, tangent in $x^\star$ to
the eigenspace of $a_2$, such that $\alpha(x)=x^\star$ for all
$x\in\Wglo{u}(x^\star)$; $\Wglo{u}$ is called the \defwd{unstable manifold} of
$x^\star$.
\end{itemiz}
\end{theorem}

The situation is sketched in \figref{fig_hyperbol}. We conclude that there
are three qualitatively different types of hyperbolic points: \defwd{sinks},
which attract all orbits from a \nbh, \defwd{sources}, which repel all
orbits from a \nbh, and \defwd{saddles}, which attract orbits from one
direction, and repel them into another one. These three situations are
robust to small perturbations. 


\subsection{Periodic Orbits}
\label{ssec_bupo}

\begin{definition}
\label{def_bupo1}
A \defwd{periodic orbit} is an orbit forming a closed curve $\Gamma$ in
$\cD$. Equivalently, if $x_0\in\cD$ is not an equilibrium point, and
$\ph_T(x_0)=x_0$ for some $T>0$, then the orbit of $x_0$ is a periodic orbit
with period $T$. $T$ is called the \defwd{least period} of $\Gamma$ is
$\ph_t(x_0)\neq x_0$ for $0<t<T$. 
\end{definition}

We denote by $\gamma(t)=\ph_t(x_0)$ a periodic solution living on $\Gamma$.
In order to analyse the dynamics near $\Gamma$, we introduce the variable
$y=x-\gamma(t)$ which satisfies the equation 
\begin{equation}
\label{bupo1}
\dot y = f(\gamma(t)+y) - f(\gamma(t)). 
\end{equation}
We start again by considering the linearization of this equation around
$y=0$, given by 
\begin{equation}
\label{bupo2}
\dot y = A(t) y
\qquad\qquad
\text{where $\displaystyle A(t) = \dpar fx(\gamma(t))$.}
\end{equation}
The solution can be written in the form $y(t)=U(t)y(0)$, where the
\defwd{principal solution} $U(t)$ is a matrix-valued function solving the
linear equation 
\begin{equation}
\label{bupo3}
\dot U = A(t) U, 
\qquad\qquad
U(0) = \one.
\end{equation}
Such an equation is difficult to solve in general. In the present case,
however, there are several simplifications. First of all, since
$A(t)=A(t+T)$ for all $t$, Floquet's theorem allows us to write 
\begin{equation}
\label{bupo4}
U(t) = P(t) \e^{Bt},
\end{equation}
where $B$ is a constant matrix, and $P(t)$ is a $T$-periodic function of
time, with $P(0)=\one$. The asymptotic behaviour of $y(t)$ thus depends only
on the eigenvalues of $B$, which are called \defwd{characteristic
exponents} of $\Gamma$. $U(T)=\e^{BT}$ is called the \defwd{monodromy
matrix}, and its eigenvalues are called the \defwd{characteristic
multipliers}. 

\begin{prop}
\label{prop_bupo1}
The characteristic exponents of $\Gamma$ are given by 
\begin{equation}
\label{bupo5}
0 
\qquad\qquad
\text{and}
\qquad\qquad
\frac1T \int_0^T \Bigpar{\dpar{f_1}{x_1}+\dpar{f_2}{x_2}}(\gamma(t)) \6t.
\end{equation}
\end{prop}
\begin{proof}
We first observe that 
\[
\dtot{}{t} \dot\gamma(t) = \dtot{}{t} f(\gamma(t)) = \dpar
fx(\gamma(t))\dot\gamma(t) = A(t)\dot\gamma(t).
\]
Thus $\dot\gamma(t)$ is a particular solution of \eqref{bupo2} and can be
written, by \eqref{bupo4}, as $\dot\gamma(t)=P(t)\e^{Bt}\dot\gamma(0)$. 
It follows by periodicity of $\gamma$ that 
\[
\dot\gamma(0) = \dot\gamma(T) = P(T)\e^{BT}\dot\gamma(0) =
\e^{BT}\dot\gamma(0). 
\]
In other words, $\dot\gamma(0)$ is an eigenvector of $\e^{BT}$ with
eigenvalue $1$, which implies that $B$ has an eigenvalue equal to $0$. 

Next we use the fact that the principal solution satisfies the relation
\[
\dtot{}{t} \det U(t) = \Tr A(t) \det U(t),
\]
as a consequence of the fact that $\det(\one+\eps B)=1+\eps\Tr B
+\Order{\eps^2}$ as $\eps\to0$. Since $\det U(0)=1$, we obtain 
\[
\det\e^{BT} = \det U(T) = \exp\Bigset{\int_0^T \Tr A(t)\6t}. 
\]
But $\det\e^{BT}$ is the product of the characteristic multipliers, and thus
$\log\det\e^{BT}$ is $T$ times the sum of the characteristic exponents, one
of which we already know is equal to $0$. The result follows from the
definition of $A$. 
\end{proof}

In the linear case, the second characteristic exponent in \eqref{bupo5}
determines the stability: if it is negative, then $\Gamma$ attracts other
orbits, and if it is positive, then $\Gamma$ repels other orbits. 
The easiest way to extend these properties to the nonlinear case is
geometrical. It relies heavily on a \lq\lq geometrically obvious\rq\rq\
property of planar curves, which is notoriously hard to prove. 

\begin{theorem}[Jordan Curve Theorem]
\label{thm_bupo1}
A closed curve in $\R^2$ which does not intersect itself separates $\R^2$
into two connected components, one bounded, called the \defwd{interior} of
the curve, the other unbounded, called the \defwd{exterior} of the curve. 
\end{theorem}

Let $\Sigma$ be an arc (a piece of a differentiable curve) in the plane. We
say that $\Sigma$ is \defwd{transverse} to the vector field $f$ if at any
point of $\Sigma$, the tangent vector to $\Sigma$ and $f$ are not
collinear. 

\begin{figure}
 \centerline{\psfig{figure=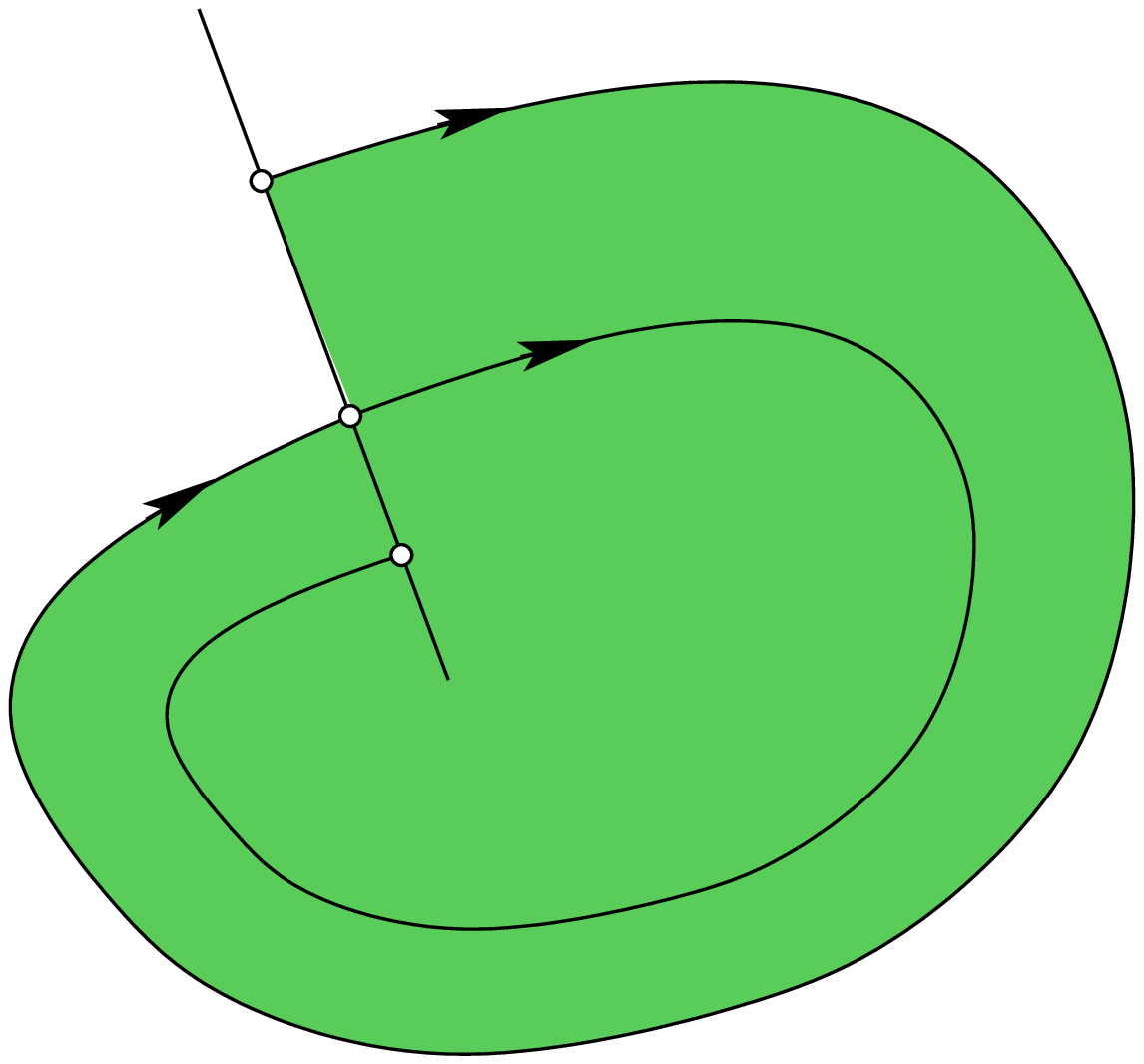,height=40mm,clip=t}
 \hspace{15mm}
 \psfig{figure=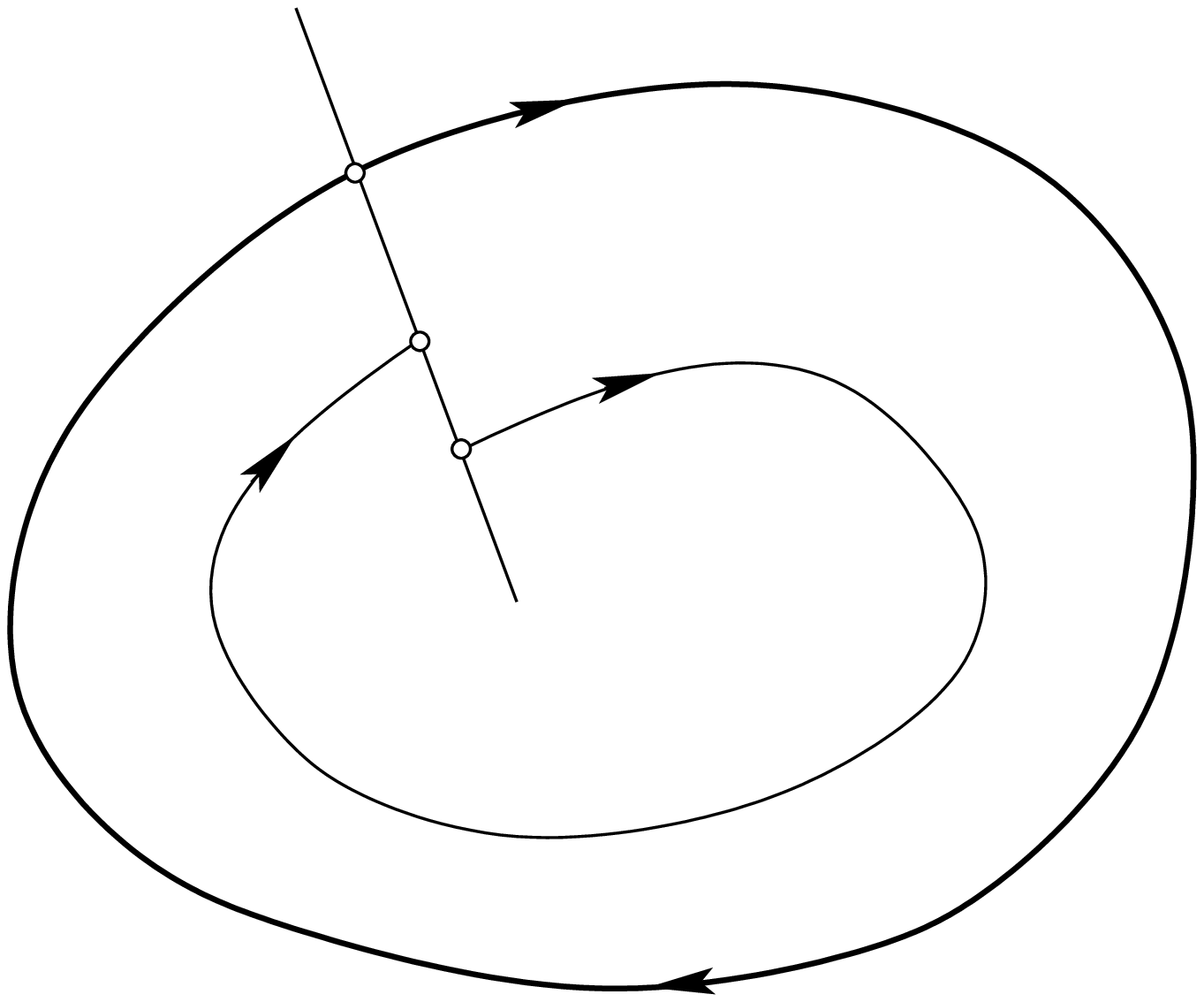,height=40mm,clip=t}}
 \figtext{
 	\writefig	1.0	4.2	(a)
 	\writefig	7.5	4.2	(b)
 	\writefig	2.9	4.2	$\Sigma$
  	\writefig	4.8	4.2	$\cC$
	\writefig	4.4	2.1	$\cN$
 	\writefig	2.6	3.6	$x_0$
  	\writefig	2.9	2.95	$x_1$
 	\writefig	3.6	2.3	$x_2$
 	\writefig	9.4	4.2	$\Sigma$
	\writefig	12.0	4.0	$\Gamma$
 	\writefig	8.95	3.75	$x^\star$
  	\writefig	9.8	3.15	$\Pi(x)$
 	\writefig	9.55	2.5	$x$
 }
 \captionspace
 \caption[]
 {(a) Proof of Lemma~\ref{lem_bupo1}: The orbit of $x_0$ defines a
 positively invariant set and crosses the transverse arc $\Sigma$ in a
 monotone sequence of points. (b) Definition of a Poincar\'e map $\Pi$
 associated with the periodic orbit $\Gamma$.}
\label{fig_poincare}
\end{figure}

\begin{lemma}
\label{lem_bupo1}
Assume the vector field is transverse to $\Sigma$. Let $x_0\in\Sigma$ and
denote by $x_1, x_2, \dots$ the successive intersections, as time increases,
of the orbit of $x_0$ with $\Sigma$ (as long as they exist). Then for every
$i$, $x_i$ lies between $x_{i-1}$ and $x_{i+1}$ on $\Sigma$. 
\end{lemma}
\begin{proof}
Assume the orbit of $x_0$ intersects $\Sigma$ for the first (positive) time
in $x_1$ (otherwise there is nothing to prove). Consider the curve $\cC$
consisting of the orbit between $x_0$ and $x_1$, and the piece of $\Sigma$
between the same points (\figref{fig_poincare}a). Then $\cC$ admits an
interior $\cN$, which is either positively invariant or negatively
invariant. In the first case, we conclude that $x_2$, if it exists, lies
inside $\cN$. In the second case, $x_2$ must lie outside $\cN$ because
otherwise the negative orbit of $x_2$ would have to leave $\cN$. 
\end{proof}

Assume now that $\Sigma$ is a small segment, intersecting the periodic orbit
$\Gamma$ transversally at exactly one point $x^\star$. For each
$x\in\Sigma$, we define a \defwd{first return time}
\begin{equation}
\label{bupo6}
\tau(x) = \inf\bigsetsuch{t>0}{\ph_t(x)\in\Sigma} \in (0,\infty].
\end{equation}
By continuity of the flow, $\tau(x)<\infty$ for $x$ sufficiently close to
$x^\star$ (\figref{fig_poincare}b). The \defwd{Poincar\'e map} of $\Gamma$
associated with $\Sigma$ is the map 
\begin{equation}
\label{bupo7}
\fctndef{\Pi}{\setsuch{x\in\Sigma}{\tau(x)<\infty}}{\Sigma}
{x}{\ph_{\tau(x)}(x).}
\end{equation}

\begin{theorem}
\label{thm_bupo2}
The Poincar\'e map has the following properties:
\begin{itemiz}
\item	$\Pi$ is a monotone increasing $\cC^1$ map in a \nbh\ of $x^\star$;
\item	the multipliers of $\Gamma$ are $1$ and $\Pi'(x^\star)$;
\item	if $\Pi'(x^\star)<1$, then there exists a \nbh\ $\cU$ of $\Gamma$
such that $\w(x)=\Gamma$ for every $x\in\cU$;
\item	if $\Pi'(x^\star)>1$, then there exists a \nbh\ $\cU$ of $\Gamma$
such that $\alpha(x)=\Gamma$ for every $x\in\cU$.
\end{itemiz}
\end{theorem}
\begin{proof}\hfill
\begin{itemiz}
\item	The differentiability of $\Pi$ is a consequence of the implicit
function theorem. Define some ordering on $\Sigma$, and consider two points
$x<y$ on $\Sigma$. We want to show that $\Pi(x)<\Pi(y)$. As in the proof of
Lemma~\ref{lem_bupo1}, the orbit of each point defines a positively or
negatively invariant set $\cN$, resp.\ $\cM$. Assume for instance that
$\cN$ is positively invariant and $y\in\cN$. Then $\Pi(y)\in\cN$ and the
conclusion follows. The other cases are similar. 

\item	We already saw in Proposition~\ref{prop_bupo1} that
$\e^{BT}f(x^\star)=f(x^\star)$. Let $v$ be a tangent vector to $\Sigma$ of
unit length, $a$ a sufficiently small scalar, and consider the relation
\[
\Pi(x^\star+av) = \ph_{\tau(x^\star+av)}(x^\star+av).
\]
Differentiating this with respect to $a$ and evaluating at $a=0$, we get 
\[
\Pi'(x^\star)v = \e^{BT} v + \tau'(x^\star) f(x^\star).
\]
This shows that in the basis $(f(x^\star),v)$, the monodromy matrix has the
representation 
\[
\e^{BT} = 
\begin{pmatrix}
1 & -\tau'(x^\star) \\
0 & \Pi'(x^\star)
\end{pmatrix}.
\]

\item	Assume that $\Pi'(x^\star)<1$. Then $\Pi$ is a contraction in some
\nbh\ of $x^\star$. This implies that there exists an open interval
$(x_1,x_2)$ containing $x^\star$ which is mapped strictly into itself. Thus
the orbits of $x_1$ and $x_2$, between two successive intersections with
$\Sigma$, define a positively invariant set $\cN$. We claim that
$\w(x)=\Gamma$ for any $x\in\cN$. 
Let $\Sigma_0$ denote the piece of $\Sigma$ between $x_1$ and $x_2$.  If
$x\in\Sigma_0$, then the iterates $\Pi^n(x)$ converge to $x^\star$ as
$n\to\infty$, which shows that $x^\star\in\w(x)$. It also shows that
$x^\star$ is the only point of $\Sigma$ which belongs to $\w(x)$. If
$x\notin\Sigma_0$, there exists by continuity of the flow a time
$t\in(0,T)$ such that $x\in\ph_t(\Sigma_0)$. Hence the sequence
$\set{\Pi^n\circ\ph_{T-t}(x)}$ converges to $x^\star$. Using a similar
combination of $\ph_t$ and $\Pi^n$, one can construct a sequence converging
to any $y\in\Gamma$.  

\item	The case $\Pi'(x^\star)>1$ can be treated in the same way by
considering $\Pi^{-1}$. 
\qed
\end{itemiz}
\renewcommand{\qed}{}
\end{proof}

\begin{definition}
\label{def_bupo2}
The periodic orbit $\Gamma$ is \defwd{hyperbolic} if $\Pi'(x^\star)\neq
1$.  In other words, $\Gamma$ is hyperbolic if 
\begin{equation}
\label{bupo8}
\int_0^T \Bigpar{\dpar{f_1}{x_1}+\dpar{f_2}{x_2}}(\gamma(t)) \6t \neq 0.
\end{equation}
\end{definition}


\subsection{The Poincar\'e--Bendixson Theorem}
\label{ssec_bupb}

\begin{figure}
 \centerline{\psfig{figure=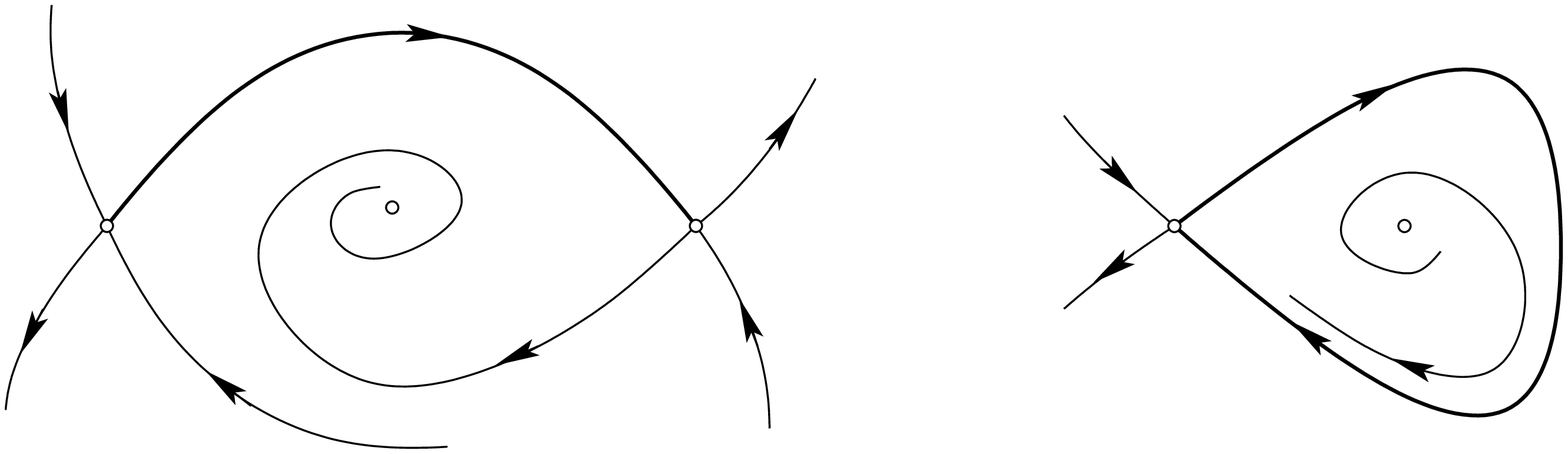,height=40mm,clip=t}}
 \figtext{
 	\writefig	0.5	4.5	a
 	\writefig	2.0	2.45	$x^\star_1$
  	\writefig	6.9	2.45	$x^\star_2$
 	\writefig	2.3	4.1	$\Wglo{u}(x^\star_1)$
 	\writefig	5.7	3.7	$\Wglo{s}(x^\star_2)$
	\writefig	9.0	4.5	b
 	\writefig	10.45	2.7	$x^\star$
 	\writefig	11.2	3.8	$\Wglo{u}(x^\star)$
 	\writefig	11.2	0.9	$\Wglo{s}(x^\star)$
 }
 \captionspace
 \caption[]
 {Examples of saddle connections: (a) heteroclinic connection, (b)
 homoclinic connection.}
\label{fig_connection}
\end{figure}

Up to now, we have encountered equilibrium points and periodic orbits as
possible limit sets. There exists a third kind of orbit that may be part of
a limit set. 

\begin{definition}
\label{def_bupb1}
Let $x^\star_1$ and $x^\star_2$ be two saddles. Assume that their stable and
unstable manifolds intersect, i.e.\ there exists a point
$x_0\in\Wglo{u}(x^\star_1)\cap\Wglo{s}(x^\star_2)$. Then the orbit $\Gamma$
of $x_0$ must be contained in $\Wglo{u}(x^\star_1)\cap\Wglo{s}(x^\star_2)$
and is called a \defwd{saddle connection}. This implies that
$\alpha(x)=x^\star_1$ and $\w(x)=x^\star_2$ for any $x\in\Gamma$. The
connection is called \defwd{heteroclinic} if $x^\star_1\neq x^\star_2$ and
\defwd{homoclinic} if $x^\star_1=x^\star_2$.
\end{definition}

Examples of saddle connections are shown in \figref{fig_connection}. 
The remarkable fact is that this is the exhaustive list of all
possible limit sets in two dimensions, as shown by the following theorem.

\begin{figure}
 \centerline{\psfig{figure=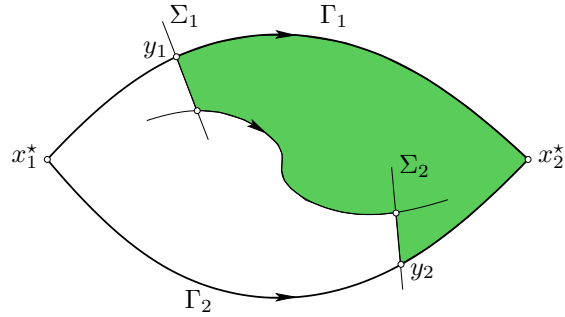,height=40mm,clip=t}}
 \figtext{
 	\writefig	3.7	2.45	$x^\star_1$
 	\writefig	10.7	2.45	$x^\star_2$
 	\writefig	5.45	3.9	$y_1$
 	\writefig	9.0	0.95	$y_2$
 	\writefig	5.8	4.3	$\Sigma_1$
 	\writefig	8.85	2.3	$\Sigma_2$
 	\writefig	7.8	4.3	$\Gamma_1$
 	\writefig	6.0	0.5	$\Gamma_2$
 }
 \captionspace
 \caption[]
 {If the $\w$-limit set contained two saddle connections $\Gamma_1,
 \Gamma_2$, admitting the same equilibria as $\alpha$- and $\w$-limit sets,
 the shaded set would be invariant.}
\label{fig_Bendixson}
\end{figure}

\begin{theorem}[Poincar\'e--Bendixson]
\label{thm_bupb}
Let $\cS$ be a compact positively invariant region containing a finite
number of equilibrium points. For any $x\in\cS$, one of the following
possibilities holds:
\begin{enum}
\item	$\w(x)$ is an equilibrium point;
\item	$\w(x)$ is a periodic orbit;
\item	$\w(x)$ consists of a finite number of equilibrium points
$x^\star_1,\dots,x^\star_n$ and orbits $\Gamma_k$ with
$\alpha(\Gamma_k)=x^\star_i$ and $\w(\Gamma_k)=x^\star_j$. 
\end{enum}
\end{theorem}
\begin{proof}
The proof is based on the following observations:
\begin{enum}
\item	{\bf Let $\Sigma\subset\cS$ be an arc transverse to the vector
field. Then $\w(x)$ intersects $\Sigma$ in one point at most.}

Assume by contradiction that $\w(x)$ intersects $\Sigma$ in two points
$y,z$. Then there exist sequences of points $\set{y_n}\in\Sigma$ and
$\set{z_n}\in\Sigma$ in the orbit of $x$ such that $y_n\to y$ and $z_n\to z$
as $n\to\infty$. However, this would contradict Lemma~\ref{lem_bupo1}. 

\item	{\bf If $\w(x)$ does not contain equilibrium points, then
it is a periodic orbit.}

Choose $y\in\w(x)$ and $z\in\w(y)$. The definition of $\w$-limit sets
implies that $z\in\w(x)$.  Since $\w(x)$ is closed, invariant, and contains
no equilibria, $z\in\w(x)$ is not an equilibrium. We can thus construct a
transverse arc $\Sigma$ through $z$. The orbit of $y$ intersects $\Sigma$
in a monotone sequence $\set{y_n}$ with $y_n\to z$ as $n\to\infty$. Since
$y_n\in\w(x)$, we must have $y_n=z$ for all $n$ by point 1. Hence the orbit
of $y$ must be a closed curve.

Now take a transverse arc $\Sigma'$ through $y$. By point 1., $\w(x)$
intersects $\Sigma'$ only at $y$. Since $\w(x)$ is connected, invariant, and
without equilibria, it must be identical with the periodic orbit through
$y$. 

\item	{\bf Let $x^\star_1$ and $x^\star_2$ be distinct equilibrium points
contained in $\w(x)$. There exists at most one orbit $\Gamma\in\w(x)$ such
that $\alpha(y)=x^\star_1$ and $\w(y)=x^\star_2$ for all $y\in\Gamma$.}

Assume there exist two orbits $\Gamma_1, \Gamma_2$ with this property.
Choose points $y_1\in\Gamma_1$ and $y_2\in\Gamma_2$ and construct
transverse arcs $\Sigma_1, \Sigma_2$ through these points. Since $\Gamma_1,
\Gamma_2 \in \w(x)$, there exist times $t_2>t_1>0$ such that
$\ph_{t_j}(x)\in\Sigma_j$ for $j=1,2$. But this defines an invariant region
(\figref{fig_Bendixson}) that does not contain $\Gamma_1, \Gamma_2$, a
contradiction. 

\item	{\bf If $\w(x)$ contains only equilibrium points, then it must
consist of a unique equilibrium.}

This follows from the assumption that there are only finitely many
equilibria and $\w(x)$ is connected.

\item	{\bf Assume $\w(x)$ contains both equilibrium points and points that
are not equilibria. Then any point in $\w(x)$ admits equilibria as $\alpha$-
and $\w$-limit sets.}

Let $y$ be a point in $\w(x)$ which is not an equilibrium. If $\alpha(y)$
and $\w(y)$ were not equilibria, they must be periodic orbits by point 2.
But this is impossible since $\w(x)$ is connected and contains equilibria. 
\qed 
\end{enum}
\renewcommand{\qed}{}
\end{proof}

One of the useful properties of limit sets is that many points share the
same limit set, and thus the number of these sets is much smaller than the
number of orbits.

\begin{definition}
\label{def_bupb2}
Let $\w_0=\w(x)$ be the $\w$-limit set of a point $x\in\cD$. The
\defwd{basin of attraction} of $\w_0$ is the set 
$\cA(\w_0) = \setsuch{y\in\cD}{\w(y)=\w_0}.$ 
\end{definition}

\begin{figure}[hb]
 \centerline{\psfig{figure=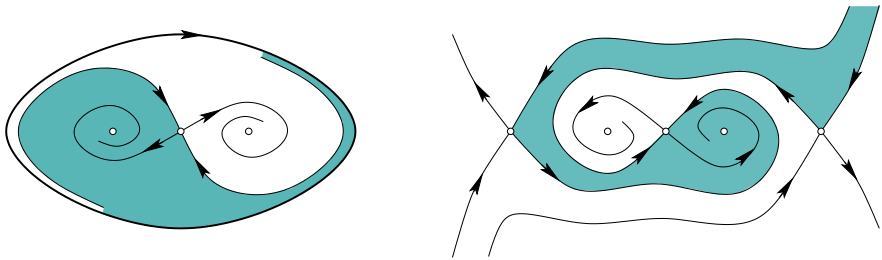,height=37mm,clip=t}}
 \captionspace
 \caption[]
 {Examples of basins of attraction.}
\label{fig_basins}
\end{figure}

We leave it as an exercise to show that if $\cS$ is a positively invariant
set containing $\w_0$, then the sets $\cA(\w_0)\cap\cS$,
$\cS\setminus\cA(\w_0)$ and $\partial\cA(\w_0)\cap\cS$ are positively
invariant. As a consequence, $\cS$ can be decomposed into a union of
disjoint basins of attraction, whose boundaries consist of periodic orbits
and stable manifolds of equilibria (see \figref{fig_basins}). 

\section{Structurally Stable Vector Fields}
\label{sec_ssvf}

Loosely speaking, a vector field $f$ is \defwd{structurally stable} if any
sufficiently small perturbation does not change the qualitative dynamics of
$\dot x =f(x)$. 

Obviously, to make this a precise definition, we need to specify what we
mean by \lq\lq small perturbation\rq\rq\ and \lq\lq same qualitative
dynamics\rq\rq. After having done this, we will state a general result,
mainly due to Peixoto, which gives necessary and sufficient conditions for
a vector field $f$ to be structurally stable. 


\subsection{Definition of Structural Stability}
\label{ssec_ssdef}

We start by defining what we mean by a small perturbation. To avoid certain
technical difficulties, we will work on a {\em compact} phase space
$\cD\subset\R^2$. 

We denote by $\cX^k(\cD)$ the space of $\cC^k$ vector fields on $\cD$,
pointing inwards on the boundary of $\cD$, so that $\cD$ is positively
invariant. On $\cX^k(\cD)$, we introduce the norm
\begin{equation}
\label{ssd1}
\norm{f}_{\cC^k} \defby 
\sup_{x\in\cD} \;\max_{j=1,2} \;
\max_{\substack{0\leqs p_1+p_2\leqs k\\p_1,p_2\geqs0}}\;
\biggabs{\dpar{^{p_1+p_2}f_j}{x_1^{p_1}\partial x_2^{p_2}}},
\end{equation}
called the \defwd{$\cC^k$-norm}. Note that $\norm{\cdot}_{\cC^0}$ is simply
the sup norm. $\norm{f}_{\cC^k}$ is small if the sup norms of $f_1, f_2$
and all their derivatives up to order $k$ are small. The resulting topology
in $\cX^k(\cD)$ is the \defwd{$\cC^k$-topology}. Note that
$(\cX^k(\cD),\norm{\cdot}_{\cC^k})$ is a Banach space. 

We will say that the vector field $g\in\cX^k(\cD)$ is a small perturbation
of $f\in\cX^k(\cD)$ if $\norm{g-f}_{\cC^k}$ is small, for a certain $k$.
Taking $k=0$ yields a notion of proximity which is too weak (i.e.\ too large
\nbh s in the function space $\cX^k(\cD)$). Indeed, consider for instance
the one-dimensional example $f(x)=-x$, $g(x,\eps)=-x+\eps\sqrt{\abs x}$,
$\cD=[-1,1]$. Then $\norm{g-f}_{\cC^0} = \eps$, but for every $\eps>0$, $g$
has two equilibria at $0$ and $\eps^2$, while $f$ has only one equilibrium. 
To obtain interesting results, we will need to work with $k=1$ at least,
and sometimes with larger $k$.

We now introduce a way to compare qualitative dynamics. 

\begin{definition}
\label{def_ssd1}
Two vector fields $f$ and $g:\cD\to\R^2$ are said to be \defwd{topologically
equivalent} if there is a homeomorphism $h:\cD\to\cD$ (i.e. $h$ is
continuous, bijective, and has a continuous inverse) taking orbits of $f$
onto orbits of $g$, and preserving direction of time. 
\end{definition}

In other words, if we denote by $\ph_t$ and $\psi_t$ the flows of $f$ and
$g$ respectively, the vector fields are topologically equivalent if there
exists a homeomorphism $h$ and a continuous, monotonously increasing
bijection $\tau:\R\to\R$ such that 
\begin{equation}
\label{ssd2}
\psi_{\tau(t)}(x) = h \circ \ph_t \circ h^{-1}(x)
\end{equation}
for all $x\in\cD$ and all $t\geqs 0$. 

Since $h$ transforms orbits of $f$ into orbits of $g$, both vector fields
also have homeomorphic invariant sets, and in particular the same type of
$\alpha$- and $\w$-limit sets. Thus, if we understand the qualitative
dynamics of $f$ (in particular, asymptotically for $t\to\pm\infty$), we
also understand the dynamics of $g$, even if we don't know the
homeomorphism $h$ explicitly. 

We can also define \defwd{differentiable equivalence} by requiring that $h$
be $\cC^1$ (or $\cC^k$ for some $k\geqs1$). This, however, often turns out
to be too strong a requirement to be helpful in classifying vector fields.
Consider for instance the linear system
\begin{equation}
\label{ssd3}
\dot x = 
\begin{pmatrix}
1 & 0 \\ 0 & 1+\eps
\end{pmatrix}
x.
\end{equation}
The orbits are of the form $x_2=x_1^{1+\eps}$. For $\eps=0$, they belong to 
straight lines through the origin, while for $\eps>0$, they are tangent to
the $x_1$-axis. Thus, the systems for $\eps=0$ and $\eps>0$ are
topologically equivalent, but not differentiably equivalent. On the other
hand, topological equivalence will only distinguish between sources, sinks
and saddles, but not between nodes and foci. 

We will choose the following definition of structural stability:

\begin{definition}
\label{def_ssd2}
The vector field $f\in\cX^k(\cD)$, $k\geqs 1$, is \defwd{structurally
stable} if there exists $\eps>0$ such that every $g\in\cX^k(\cD)$ with
$\norm{g-f}_{\cC^k}<\eps$ is topologically equivalent to $f$. 
\end{definition}

This definition may seem a bit arbitrary. If we only allow perturbations
which are small in the $\cC^1$-topology, why do we not require
differentiable  equivalence? The reason has to do with the size of
equivalence classes. If we were to allow perturbations which are only small
in the $\cC^0$-topology, virtually no vector field would be structurally
stable, because one can find continuous perturbations which change the
orbit structure. On the other hand, differentiable equivalence is such a
strong requirement that it would lead to a very large number of relatively
small equivalence classes. Definition~\ref{def_ssd2} turns out to have
exactly the right balance to yield a useful result, which we now state. 


\subsection{Peixoto's Theorem}
\label{ssec_ssp}

A characterization of structurally stable vector fields was first given by
Andronov and Pontrjagin in 1937, and generalized by Peixoto in 1959. We
state here a particular version of Peixoto's result. 

\begin{theorem}[Peixoto]
\label{thm_peixoto}
A vector field $f\in\cX^1(\cD)$ is structurally stable if and only if it
satisfies the following properties:
\begin{enum}
\item	$f$ has finitely many equilibrium points, all being hyperbolic;
\item	$f$ has finitely many periodic orbits, all being hyperbolic;
\item	$f$ has no saddle connections.
\end{enum}
Moreover, the set of all structurally stable vector fields is dense in
$\cX^1(\cD)$. 
\end{theorem}

Note that $f$ is allowed to have no equilibrium or no periodic orbit at all
(although, if $\cD$ is positively invariant, we know that there exists a
non-empty $\w$-limit set). We will not give a full proof of this result,
which is quite involved, but we shall discuss the main ideas of the proof
below. First of all, a few remarks are in order:
\begin{itemiz}
\item	The conditions for $f$ to be structurally stable are remarkably
simple, even though they are not always easy to check, especially the
absence of saddle connections.

\item	The fact that the structurally stable vector fields are dense is a
very strong result. It means that given any structurally {\em unstable}
vector field, one can find an arbitrarily small perturbation which makes it
structurally stable (this is no longer true in three dimensions). It also
means that a \lq\lq typical\rq\rq\ vector field in $\cX^k(\cD)$ will
contain only hyperbolic equilibria and hyperbolic periodic orbits, and no
saddle connections. 

\item	This result allows to classify vector fields with respect to
topological equivalence. This can be done, roughly speaking, by associating
a graph with each vector field: the vertices of the graph correspond to
$\alpha$- and $\w$-limit sets, and the edges to orbits connecting these
limit sets (this is not always sufficient to distinguish between different
classes of vector fields, and some extra information has to be added to the
graph).

\item	Peixoto has given various generalizations of this result, to domains
$\cD$ which are not positively invariant (in that case, one has to impose
certain conditions on the behaviour near the boundary), and to general
compact two-dimensional manifolds (where an additional condition on the
limit sets is required). 
\end{itemiz}
Before proceeding to the sketch of the proof, let us recall the implicit
function theorem:

\begin{theorem}
\label{thm_ift}
Let $\cN$ be a \nbh\ of $(x^\star,y^\star)$ in $\R^n\times\R^m$. Let
$\Phi:\cN\to\R^n$ be of class $\cC^k$, $k\geqs 1$, and satisfy
\begin{align}
\label{ssp1}
\Phi(x^\star,y^\star) &= 0, \\
\det \dpar \Phi x (x^\star,y^\star) &\neq 0.
\label{ssp2}
\end{align}
Then there exists a \nbh\ $\cU$ of $y^\star$ in $\R^m$ and a unique
function $\ph\in\cC^k(\cU,\R^n)$ such that 
\begin{align}
\label{ssp3}
\ph(y^\star) &= x^\star, \\
\Phi(\ph(y),y) &= 0 \qquad\text{for all $y\in\cU$}.
\label{ssp4}
\end{align}
\end{theorem}

\begin{figure}
 \centerline{\psfig{figure=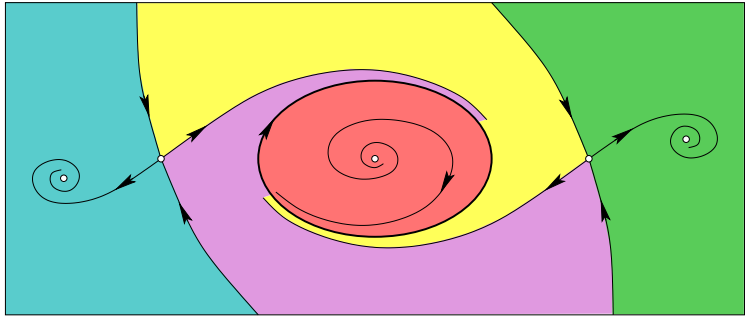,height=45mm,clip=t}}
 \captionspace
 \caption[]
 {Example of a vector field with five canonical regions. Any sufficiently
 small perturbation will have similar canonical regions, and the
 homeomorphism $h$ is constructed separately in a \nbh\ of each sink, in a
 \nbh\ of the periodic orbit, and in each remaining part of a canonical
 region.}
\label{fig_canonical}
\end{figure}

\begin{proof}[{\sc Some ideas of the proof of Peixoto's Theorem}]
\hfill

\noindent
To prove that the three given conditions are sufficient for structural
stability, we have to prove the following: Let $f$ satisfy Conditions 1.--3.,
and let $g\in\cX^1(\cD)$ with $\norm{g-f}_{\cC^1}$ small enough. Then there
exists a homeomorphism $h$ taking orbits of $f$ onto orbits of $g$. It will
be convenient to introduce the one-parameter family
\[
F(x,\eps) = f(x) + \eps \Bigpar{\frac{g(x)-f(x)}{\eps_0}},
\qquad 
\eps_0 = \norm{g-f}_{\cC^1},
\]
which satisfies $F(x,0)=f(x)$, $F(x,\eps_0)=g(x)$ and
$\norm{F(\cdot,\eps)-f}_{\cC^1}=\eps$. The main steps of the proof are the
following:
\begin{enum}
\item	Let $x^\star$ be an equilibrium point of $f$. Then we have
$F(x^\star,0)=0$ by definition, and $\dpar Fx(x^\star,0)= \dpar
fx(x^\star)$ is the Jacobian matrix of $f$ at $x^\star$. Since $x^\star$
is hyperbolic by assumption, this matrix has a nonvanishing determinant.
Hence the implicit function theorem yields the existence, for sufficiently
small $\eps$, of a unique equilibrium point of $F(x,\eps)$ in a \nbh\ of
$x^\star$. This equilibrium is also hyperbolic by continuous dependence of
the eigenvalues of $\dpar Fx$ on $\eps$. 

\item	Let $\Gamma$ be a periodic orbit of $f$. Let $\Sigma$ be a
transverse arc, intersecting $\Gamma$ at $x_0$, and let $\Pi_0$ be the
associated Poincar\'e map. For sufficiently small $\eps$, $\Sigma$ is also
transverse to $F(x,\eps)$, and a Poincar\'e map $\Pi_\eps$ can be defined in
some subset of $\Sigma$. Consider the function 
\[
\Phi(x,\eps) = \Pi_\eps(x) - x.
\] 
Then $\Phi(x_0,0)=0$ and $\dpar\Phi x(x_0,0) = \Pi_0'(x_0)-1 \neq 0$, since
$\Gamma$ is assumed to be hyperbolic. Thus the implicit function theorem can
be applied again to show that $\Pi_\eps$ admits a fixed point for
sufficiently small $\eps$, which corresponds to a periodic orbit. 

\item	Since $f$ has no saddle connections, the Poincar\'e--Bendixson
Theorem shows that $f$ and $g$ have the same  $\alpha$- and $\w$-limit
sets, up to a small deformation. We call \defwd{canonical region} a maximal
open subset of points in $\cD$ sharing the same limit sets (where the
boundary $\partial\cD$ is considered as an $\alpha$-limit set). Each
canonical region is the intersection of the basin of attraction of an
$\w$-limit set, and the \lq\lq basin of repulsion\rq\rq\ of an
$\alpha$-limit set (\figref{fig_canonical}). Its boundary, called a
\defwd{separatrix}, consists of equilibria, periodic orbits, and
stable/unstable manifolds of saddles. One can show that these manifolds
connect the same limit sets even after a small perturbation of $f$, and
thus $f$ and $g$ have homeomorphic canonical regions.

\item	The last step is to construct the homeomorphism $h$. This is done
separately in \nbh s of all limit sets, and in the remaining part of each
canonical region, with appropriate matching on the boundaries. In this way,
the sufficiency of Conditions 1.--3.\ has been proved.
\end{enum}
To prove that the given conditions are also necessary, it suffices to prove
that if any of them is violated, then one can find $g\in\cX^1(\cD)$ with
$\norm{g-f}_{\cC^1}$ arbitrarily small, such that $f$ and $g$ are not
topologically equivalent. We will do this in more detail in the next
section. The basic idea is that nonhyperbolic equilibria lead to local
bifurcations: the number of limit sets near such an equilibrium can change
under a small perturbation. Also, nonhyperbolic periodic orbits may
disappear or duplicate when they are perturbed. And saddle connections can
be broken by a small transversal perturbation.

Finally, to show that structurally stable vector fields are dense, one
proceeds as follows: Let $f$ be structurally unstable and $\eps>0$ an
arbitrarily small number. Then one constructs a structurally stable vector
field $g$ such that $\norm{f-g}<\eps$. The construction proceeds by
successive small deformations of $f$, which remove the violated conditions
one by one. 
\end{proof}

\section{Singularities of Codimension 1}
\label{sec_sc1}

Having obtained a rather precise characterization of structurally stable
vector fields, we would now like to understand the structure of the set
$\cS_0$ of structurally {\em unstable} vector fields in $\cX^k(\cD)$ (with
$k$ sufficiently large). Why should we want to do this, if \lq\lq
typical\rq\rq\ vector fields are structurally stable? There are at least two
reasons:
\begin{itemiz}
\item	The set $\cS_0$ constitutes boundaries between different equivalence
classes of structurally stable vector fields, so that elements in $\cS_0$
describe transitions between these classes.
\item	While typical vector fields will not belong to $\cS_0$, it is quite
possible that one-parameter families $f_\lambda$ (i.e., curves in
$\cX^k(\cD)$) has members in $\cS_0$. 
\end{itemiz}
By Peixoto's theorem, a structurally unstable vector field will display at
least one of the following properties:
\begin{itemiz}
\item	it has a non-hyperbolic equilibrium;
\item	it has a non-hyperbolic periodic orbit;
\item	it has a saddle connection.
\end{itemiz}
At least the first two properties can be expressed as conditions of the
form $H(f)=0$, where $H:\cX^k(\cD)\to\R$ is, for instance, an eigenvalue of
a stability matrix. This suggests that $\cS_0$ is composed of \lq\lq
hypersurfaces\rq\rq\ in the infinite-dimensional space $\cX^k(\cD)$. To
make this idea a bit more precise, let us first consider some
finite-dimensional examples:
\begin{enum}
\item	If $H:\R^3\to\R$ is given by $H(x)=x_1x_2$, the set $\cS=H^{-1}(0)$
is composed of the two planes $x_1=0$ and $x_2=0$. Most points in $\cS$
admit a \nbh\ in which $\cS$ is a nice \lq\lq interface\rq\rq\ separating
two regions. In the vicinity of the line $x_1=x_2=0$, however, $\cS$
separates four regions of $\R^3$. Note that the gradient of $H$ vanishes on
this line, and only there.
\item	If $H:\R^3\to\R$ is given by $H(x)=x_1^2 + x_2^2 - x_3^2$, the set
$\cS=H^{-1}(0)$ is a cone. Near any point of $\cS$, $\cS$ separates two
three-dimensional regions, except at the origin, where three such regions
meet. Again, the gradient of $H$ vanishes only at this point. 
\end{enum} 
In the infinite-dimensional case of $\cX^k(\cD)$, we will make the following
definition:
\begin{definition}
\label{def_sc01}
A set $\cS\subset\cX^k(\cD)$ is called a \defwd{$\cC^r$ submanifold of
codimension $1$} if there exists an open set $\cU\subset\cX^k(\cD)$ and a
function $H\in\cC^r(\cU,\R)$ such that $\Fd H(f)\neq 0$ in $\cU$ and
$\cS=\setsuch{f\in\cU}{H(f)=0}$. 
\end{definition}
Here $\Fd H(f)$ should be understood as the Fr\'echet derivative, that is,
a linear operator such that for all sufficiently small $g$, 
\begin{equation}
\label{sc01}
\begin{split}
&H(f+g) = H(f) + \Fd H(f) g + R(g), \\
&\text{with}\quad
\lim_{\norm{g}_{\cC^k}\to 0} \frac{\norm{R(g)}_{\cC^k}}{\norm{g}_{\cC^k}}
= 0.
\end{split}
\end{equation}
One can show that the Fr\'echet derivative is given by the following limit,
provided it exists and is continuous in $f$ and $g$:
\begin{equation}
\label{sc02}
\Fd H(f) g = \lim_{\eps\to0} \frac{H(f+\eps g)-H(f)}\eps.
\end{equation}
The set $\cS_0$ of all structurally unstable vector fields is not a
submanifold of codimension $1$, because, loosely speaking, it contains
points where $\Fd H$ vanishes. We can, however, single out those points of
$\cS_0$ in a \nbh\ of which $\cS_0$ behaves like such a manifold. A way to
do this is to start with the following definition.

\begin{definition}
\label{def_sc02}
A structurally unstable vector field $f\in\cS_0$ is called \defwd{singular of
codimension $1$} if there exists a \nbh\ $\cU_0$ of $f$ in $\cS_0$ (in the
induced topology) such that every $g\in\cU_0$ is topologically equivalent to
$f$. We will denote by $\cS_1$ the set of all singular vector fields of
codimension $1$.
\end{definition}

Then the following result holds:

\begin{theorem}
\label{thm_sc01}
Let $k\geqs4$. Then $\cS_1$ is a $\cC^{k-1}$ submanifold of codimension $1$
of $\cX^k(\cD)$ and is open in $\cS_0$ (in the induced topology). 
\end{theorem}

An important consequence of this result is that if $f\in\cS_1$, then there
exists a one-parameter family $F_\lambda$ of vector fields (depending
smoothly on $\lambda$) such that for all $g\in\cX^k(\cD)$ with
$\norm{g-f}_{\cC^k}$ small enough, $g$ is topologically equivalent to
$F_\lambda$ for some $\lambda\in\R$. The family $F_\lambda$ is called an
\defwd{unfolding} of the singularity $f$. 

\begin{figure}
 \centerline{\psfig{figure=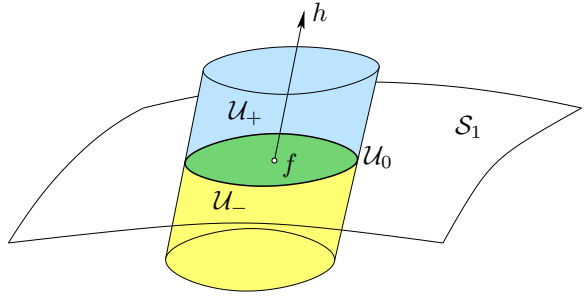,height=40mm,clip=t}}
 \figtext{
 	\writefig	8.3	2.3	$\cU_0$
 	\writefig	7.2	2.2	$f$
 	\writefig	9.5	2.7	$\cS_1$
 	\writefig	7.6	4.2	$h$
 	\writefig	6.3	1.7	$\cU_-$
 	\writefig	6.5	2.9	$\cU_+$
 }
 \captionspace
 \caption[]
 {Schematic representation of $\cX^k(\cD)$ near a singular vector field $f$ of
 codimension $1$. In a small \nbh\ of $f$, the manifold $\cS_1$ of singular
 vector fields divides $\cX^k(\cD)$ into two components $\cU_+$ and $\cU_-$,
 belonging to two different equivalence classes of structurally stable
 vector fields.}
\label{fig_codim1}
\end{figure}

To see that this is true, we choose a vector field $h$ satisfying the
\defwd{transversality condition} $\Fd H(f) h>0$ (such an $h$ must exist
since $\Fd H(f)$ is not identically zero). Let $\eps$ be a sufficiently
small number, and define the sets 
\begin{equation}
\label{sc03}
\begin{split}
\cU_0 &= \setsuch{g\in\cS_1}{\norm{g-f}_{\cC^k} < \eps} \\ 
\cU_+ &= \setsuch{g_\lambda = g_0+\lambda h}
{g_0\in\cU_0, 0< \lambda < \eps } \\ 
\cU_- &= \setsuch{g_\lambda = g_0+\lambda h}
{g_0\in\cU_0, -\eps< \lambda < 0 }. 
\end{split}
\end{equation}
Note that by definition, $H=0$ in $\cU_0$. 
Since, for sufficiently small $\eps$, 
\begin{equation}
\label{sc04}
\dtot{}{\lambda} H(g_\lambda) = \Fd H(g_\lambda) h>0
\end{equation}
by continuity of $\Fd H(\cdot)$, we conclude that $H>0$ in $\cU_+$ and $H<0$
in $\cU_-$. Hence all $g\in\cU_+$ are structurally stable, and since $\cU_+$
is clearly open and connected, they must be topologically equivalent.
Similarly, all $g\in\cU_-$ are topologically equivalent. This shows that any
family $g_\lambda$ is an unfolding of $f$, since it has members in $\cU_-$,
$\cU_0$ and $\cU_+$. 

Unfoldings are mainly used to describe the dynamics in \nbh s of
nonhyperbolic equilibria, and their perturbations. They are closely related
to the so-called \defwd{singularity theory} of functions, also known as
\defwd{catastrophe theory}. 

In the sequel, we shall give some ideas of the proof of
Theorem~\ref{thm_sc01} by enumerating all possible singularities of
codimension $1$, which can be classified into five types. In each case, we
will construct a function $H:\cU\to\R$, defined in a \nbh\ $\cU$ of a given
structurally unstable vector field $f$. This function has the following
properties: all $g\in\cU_0=H^{-1}(0)$ are topologically equivalent
to $f$, and the sets $\cU_\pm=H^{-1}(\R_\pm)$ belong to two different
equivalence classes of structurally stable vector fields. 


\subsection{Saddle--Node Bifurcation of an Equilibrium}
\label{ssec_snb}

A first type of codimension $1$ singular vector field occurs when a
nonhyperbolic equilibrium point is present, such that the linearization
around this point admits $0$ as a simple eigenvalue. In appropriate
coordinates, the vector field can be written as 
\begin{equation}
\label{snb1}
\begin{split}
\dot x_1 &= g_1(x_1,x_2) \\
\dot x_2 &= a x_2 + g_2(x_1,x_2), 
\qquad\qquad a\neq 0,
\end{split}
\end{equation}
where $g_1, g_2$ and their derivatives all vanish at the origin. We shall
further assume that $g_1$ and $g_2$ are of class $\cC^2$. The stable
manifold theorem (Theorem~\ref{thm_bueq}) admits the following
generalization:

\begin{theorem}[Center Manifold Theorem]\hfill
\label{thm_centermanif}
\begin{itemiz}
\item	If $a<0$, \eqref{snb1} admits a unique invariant curve $\Wglo{s}$,
tangent to the $x_2$-axis at the origin, such that $\w(x)=0$ for all
$x\in\Wglo{s}$; $\Wglo{s}$ is called the \defwd{stable manifold} of $0$.
\item	If $a>0$, \eqref{snb1} admits a unique invariant curve $\Wglo{u}$,
tangent to the $x_2$-axis at the origin, such that $\alpha(x)=0$ for all
$x\in\Wglo{u}$; $\Wglo{u}$ is called the \defwd{unstable manifold} of $0$.
\item	There exists an invariant curve $\Wglo{c}$, tangent to the
$x_1$-axis at the origin, called a \defwd{center manifold}.
\end{itemiz}
\end{theorem}

The center manifold is not necessarily unique. However, any center manifold
can be described, for sufficiently small $x_1$, by an equation of the form
$x_2=h(x_1)$. Inserting this relation into \eqref{snb1}, we obtain the
equation 
\begin{equation}
\label{snb2}
a h(x_1) + g_2(x_1,h(x_1)) = h'(x_1) g_1(x_1,h(x_1)),
\end{equation}
which must be satisfied by $h$. One can show that all solutions of this
equation admit the same Taylor expansion around $x_1=0$. For our purposes,
it will be sufficient to know that $h(x_1)=\Order{x_1^2}$ as $x_1\to0$. 

If $y=x_2-h(x_1)$ describes the distance of a general solution to the
center manifold, one deduces from $\eqref{snb1}$ and $\eqref{snb2}$ that $y$
satisfies an equation of the form 
\begin{equation}
\label{snb3}
\dot y = G(y, x_1)y, \qquad
G(y,x_1) = a + \Order{x_1}
\end{equation}
for small $x_1$. Thus $\Wglo{c}$ attracts nearby orbits if $a<0$
and repels them if $a>0$, so that there can be no equilibrium points
outside $\Wglo{c}$. 

\begin{figure}
 \centerline{\psfig{figure=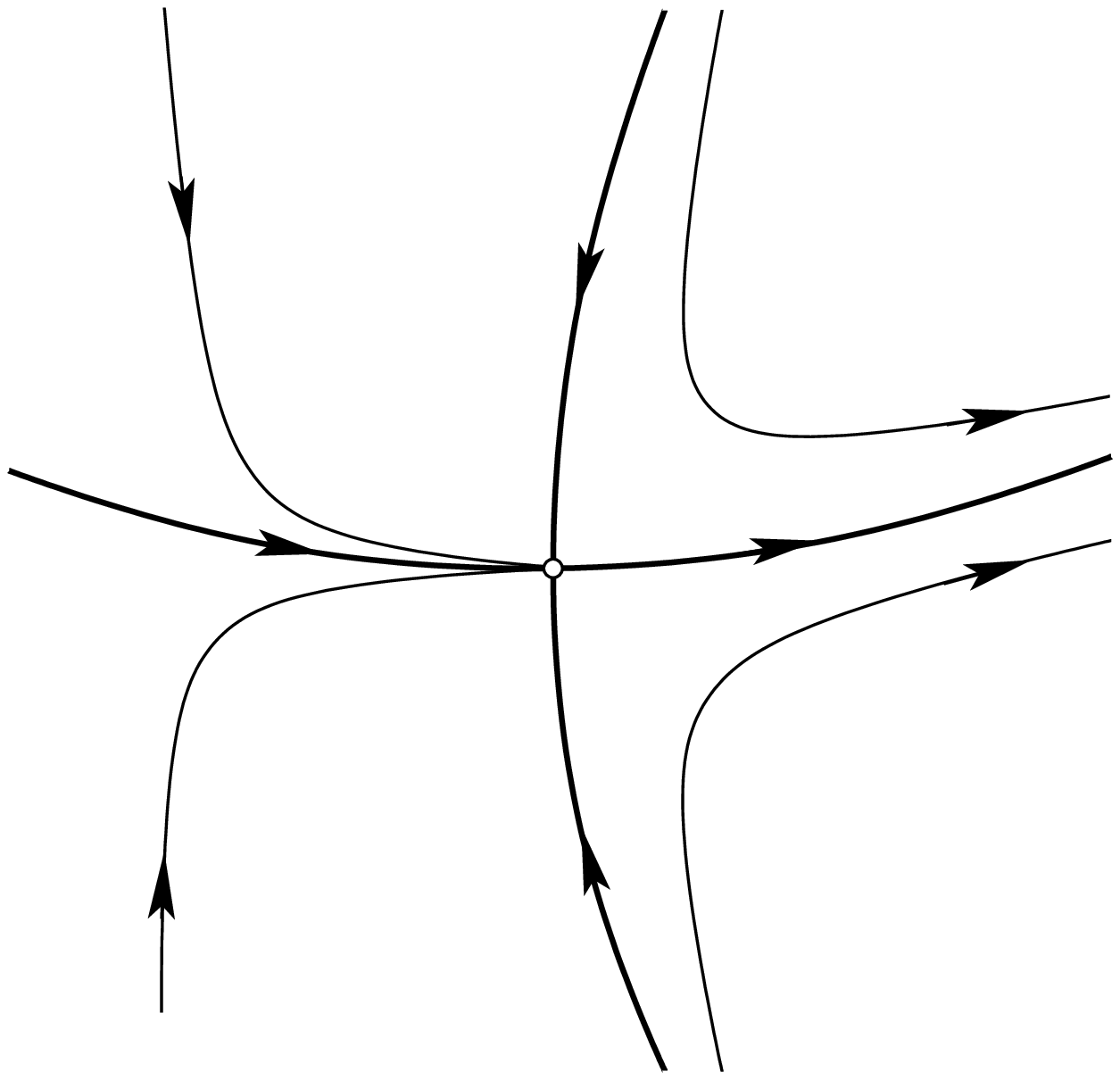,height=40mm,clip=t}
 \hspace{15mm}
 \psfig{figure=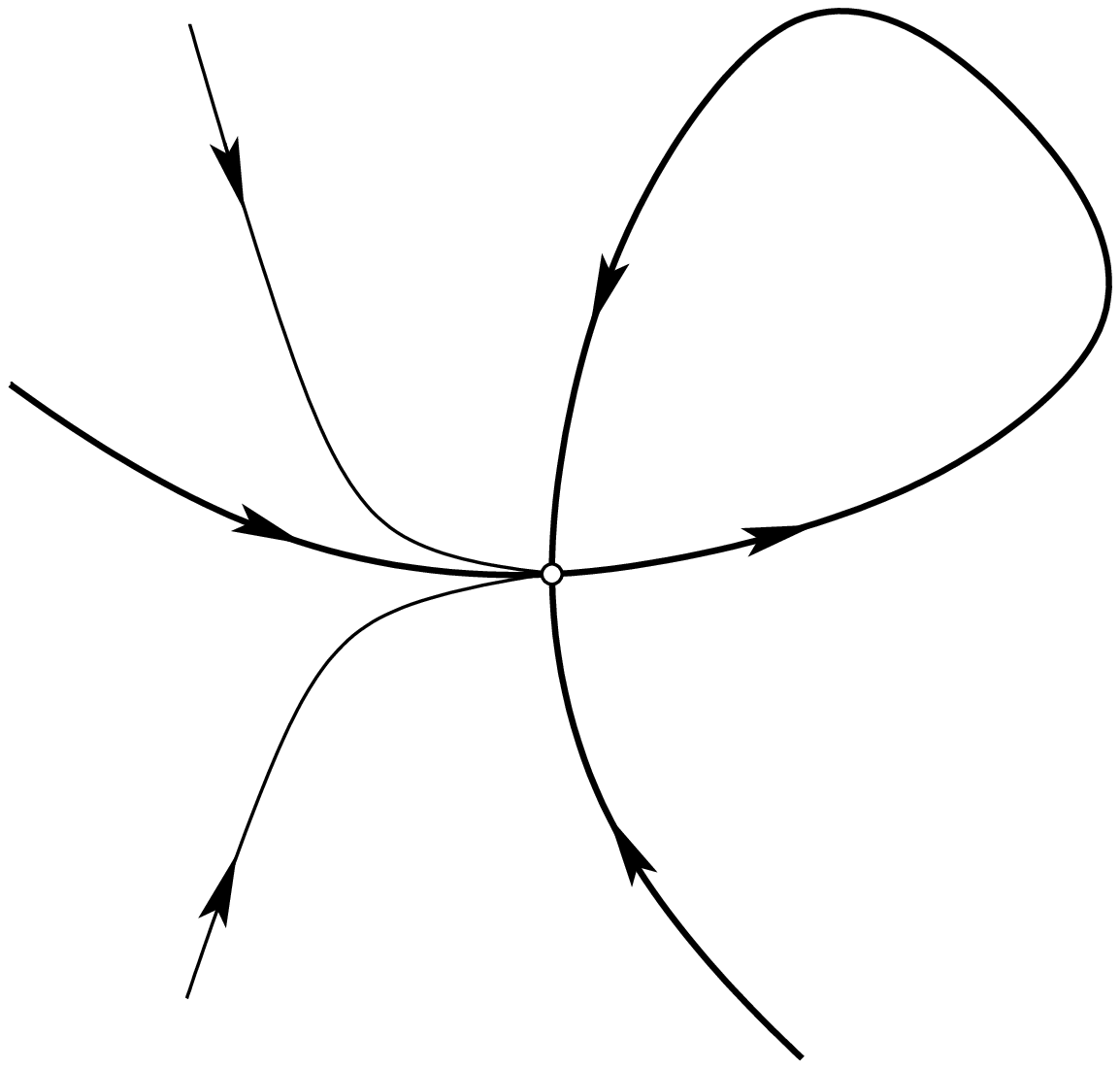,height=40mm,clip=t}}
 \figtext{
 	\writefig	1.5	4.2	(a)
 	\writefig	7.5	4.2	(b)
  	\writefig	4.1	3.9	$\Wglo{s}$
	\writefig	6.5	2.6	$\Wglo{c}$
	\writefig	4.5	2.1	$x^\star$
	\writefig	12.0	4.0	$\Wglo{c}\cap\Wglo{s}$
 	\writefig	10.4	2.05	$x^\star$
 }
 \captionspace
 \caption[]
 {(a) Local phase portrait near an elementary saddle--node $x^\star$, in a
 case where $a<0$. (b) If the center and stable manifolds intersect, then
 the vector field is not singular of codimension $1$, because one can find
 small perturbations admitting a saddle connection.}
\label{fig_saddlenode1}
\end{figure}

In order to understand the qualitative dynamics near the origin, it is thus
sufficient to understand the dynamics on the center manifold. It is governed
by the equation
\begin{equation}
\label{snb4}
\dot x_1 = g_1(x_1,h(x_1)) = c x_1^2 + \order{x_1^2}, 
\qquad\qquad
c = \frac12 \dpar{^2g_1}{x_1^2}(0,0).
\end{equation}
If $c\neq 0$, which is the typical case, then the orbits on $\Wglo{c}$ will
be attracted by the origin from one side, and repelled from the other side
(\figref{fig_saddlenode1}a). 

\begin{definition}
\label{def_snb1}
Let $x^\star$ be a nonhyperbolic equilibrium point of $f$. Assume that the
linearization of $f$ at $x^\star$ has the eigenvalues $0$ and $a\neq0$, and
that the equation \eqref{snb4} on its center manifold satisfies $c\neq0$.
Then $x^\star$ is called an \defwd{elementary saddle--node}. 
\end{definition}

Now we would like to examine the effect of a small perturbation of $f$.
Consider to that effect a one-parameter family $f(x,\lambda)$ such that
$f(\cdot,0)=f$, and write the system in the form
\begin{align}
\nonumber
\dot x_1 &= f_1(x_1,x_2,\lambda), &
f_1(x_1,x_2,0) &= g_1(x_1,x_2), \\
\label{snb5}
\dot x_2 &= f_2(x_1,x_2,\lambda), &
f_2(x_1,x_2,0) &= a x_2 + g_2(x_1,x_2), \\
\nonumber
\dot\lambda &= 0. &&
\end{align}
The point $(0,0,0)$ is a non-hyperbolic equilibrium point of this extended
system, with eigenvalues $(0,a,0)$. The center manifold theorem shows the
existence of an invariant manifold, locally described by the equation
$x_2=h(x_1,\lambda)$, which attracts nearby orbits if $a<0$ and repels
them if $a>0$. The dynamics on the center manifold is described by 
\begin{equation}
\label{snb6}
\dot x_1 = f_1(x_1,h(x_1,\lambda),\lambda) \bydef F(x_1,\lambda),
\end{equation}
where $F\in\cC^2$ satisfies the relations 
\begin{equation}
\label{snb7}
F(0,0)=0, \qquad
\dpar F{x_1}(0,0) = 0, \qquad
\dpar{^2F}{x_1^2}(0,0) = 2c \neq 0.
\end{equation}
The graph of $F(x_1,0)$ has a quadratic tangency with the $x_1$-axis at the
origin. For small $\lambda$, the graph of $F(x_1,\lambda)$ can thus have
zero, one or two intersections with the $x_1$ axis. 

\begin{prop}
\label{prop_snb}
For small $\lambda$, there exists a differentiable function
$H(\lambda)$ such that
\begin{itemiz}
\item	if $H(\lambda)>0$, then $F(x_1,\lambda)$ has no equilibria near
$x_1=0$;
\item	if $H(\lambda)=0$, then $F(x_1,\lambda)$ has an isolated
non-hyperbolic equilibrium point near $x_1=0$;
\item	if $H(\lambda)<0$, then $F(x_1,\lambda)$ has two hyperbolic
equilibrium points of opposite stability near $x_1=0$.
\end{itemiz}
\end{prop} 
\begin{proof}
Consider the function
\[
G(x_1,\lambda) = \dpar F{x_1}(x_1,\lambda).
\]
Then $G(0,0)=0$ and $\dpar G{x_1}(0,0)=2c\neq0$ by \eqref{snb7}. Thus the
implicit function theorem shows the existence of a unique differentiable
function $\ph$ such that $\ph(0)=0$ and $G(\ph(\lambda),\lambda)=0$ for all
sufficiently small $\lambda$. 

Consider now the function $K(y,\lambda)=F(\ph(\lambda)+y,\lambda)$. Taylor's
formula shows that
\begin{align*}
K(y,\lambda) &= F(\ph(\lambda),\lambda) + y^2\bigbrak{c+R_1(y,\lambda)} \\
\dpar Ky(y,\lambda) &= y \bigbrak{2c+R_2(y,\lambda)},
\end{align*}
for some continuous functions $R_1, R_2$ which vanish at the origin.  Let
$H(\lambda)=F(\ph(\lambda),\lambda)/c$. Then $H(0)=0$, and if $\lambda$
and $y$ are sufficiently small that $\abs{R_1(y,\lambda)}<c/2$, then 
\[
H(\lambda) + \frac12 y^2 \leqs \frac{K(y,\lambda)}c \leqs H(\lambda) +
\frac32 y^2. 
\]
Hence $K(y,\lambda)$ vanishes twice, once, or never, depending on the sign
of $H(\lambda)$. The expression for $\dpar Ky$ shows that $K$ is
monotonous in $y$ for small positive or small negative $y$, which means that
there are no other equilibria near the origin. 
\end{proof}

\begin{example}\hfill
\label{ex_snb}
\begin{itemiz}
\item	If $F(x_1,\lambda)=\lambda+x_1^2$, then $H(\lambda)=\lambda$.
There are two equilibria for $\lambda<0$ and no equilibria for $\lambda>0$.
This is the generic case, called the \defwd{saddle--node bifurcation}. The
family $F(x_1,\lambda)$ constitutes an unfolding of the singularity
$f(x_1)=x_1^2$. 
\item	If $F(x_1,\lambda)=\lambda^2-x_1^2$, then $H(\lambda)=-\lambda^2$.
In this case, there are two equilibrium points for all $\lambda\neq0$,
because the family $F(x_1,\lambda)$ is tangent to the manifold $\cS_1$.
This case is called the \defwd{transcritical bifurcation}. 
\item	If $F(x_1,\lambda)=\lambda x_1-x_1^2$, then
$H(\lambda)=-\lambda^2/4$. This is again a transcritical bifurcation. 
\end{itemiz}
\end{example}

\begin{figure}
 \centerline{\psfig{figure=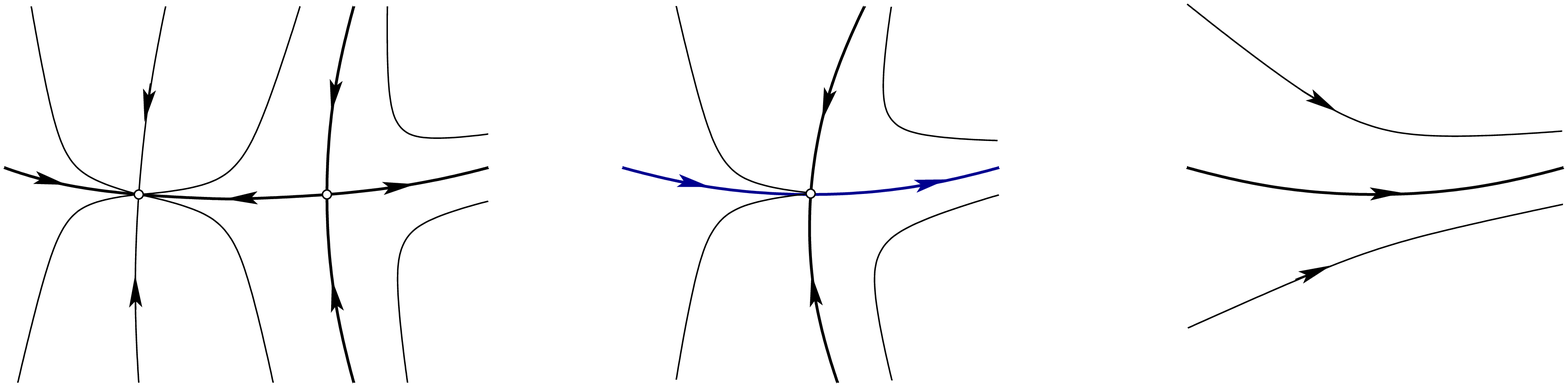,height=35mm,clip=t}}
 \centerline{\psfig{figure=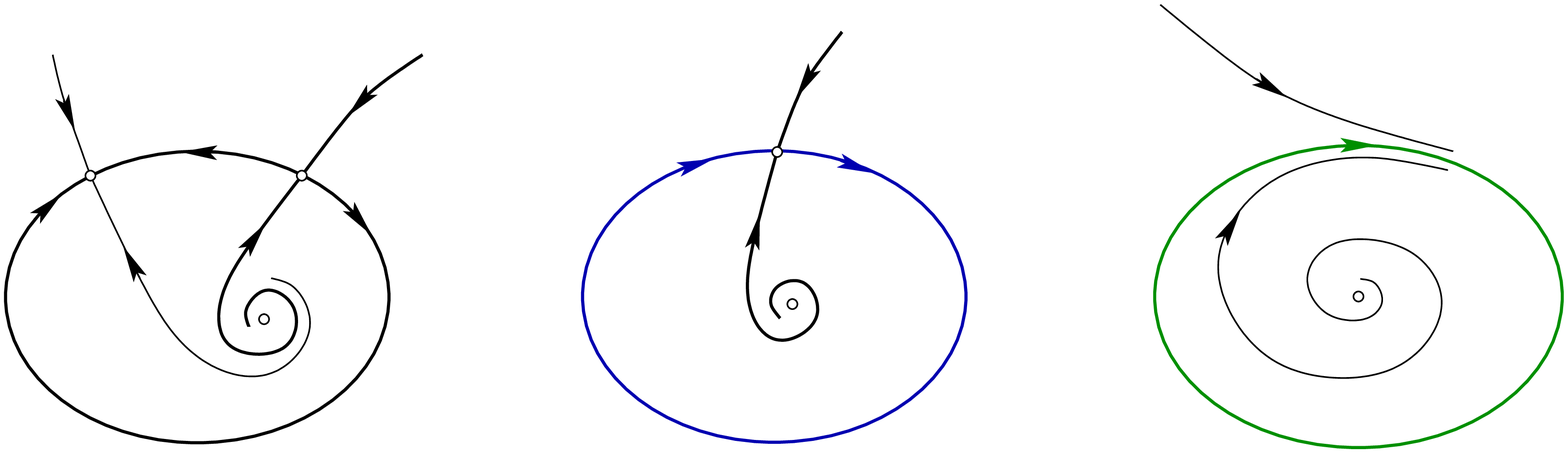,height=40mm,clip=t}}
 \figtext{
 	\writefig	0.2	7.7	(a)
 	\writefig	0.2	3.7	(b)
 }
 \captionspace
 \caption[]
 {(a) Unfolding of a saddle--node in a case where the center manifold
 $\Wglo{c}$ does not form a loop: the system has a saddle and a node if
 $H(\lambda)<0$, an elementary saddle--node if $H(\lambda)=0$ and there
 is no equilibrium point if $H(\lambda)>0$. (b) If the center manifold
 does form a loop, then a periodic orbit appears for $H(\lambda)>0$.}
\label{fig_saddlenode2}
\end{figure}

\goodbreak

The existence of the scalar function $H(\lambda)$, which determines
entirely the topology of the orbits near the equilibrium point, shows that
$f$ is indeed a singularity of codimension $1$, in a sufficiently small
\nbh\ of the equilibrium. It admits the {\em local} unfolding 
\begin{equation}
\label{snb8}
\begin{split}
\dot x_1 &= \lambda + x_1^2 \\
\dot x_2 &= a x_2.
\end{split}
\end{equation}
For positive $H$, there are no equilibrium points, while for negative
$H$, the system admits a saddle and a node (\figref{fig_saddlenode2}a). 

The {\em global} behaviour of the perturbed vector field $f(x,\lambda)$
depends on the global behaviour of the center manifold of the saddle--node
at $\lambda=0$. The result is the following.

\begin{theorem}
\label{thm_snb}
Assume that $f\in\cX^k(\cD)$, $k\geqs 2$, admits an elementary saddle--node
$x^\star$ as an equilibrium point. Then $f$ is singular of codimension $1$
if and only if 
\begin{itemiz}
\item	all other equilibria and periodic orbits of $f$ are hyperbolic;
\item	$f$ has no saddle connections;
\item	the manifolds $\Wglo{s,u}(x^\star)$ and $\Wglo{c}(x^\star)$ have no
common point and do not go to a saddle point (\figref{fig_saddlenode1}a).
\end{itemiz}

If $\Wglo{c}$ does not form a loop, then any small perturbation of $f$
admits, in a small \nbh\ of $x^\star$, either a saddle and a node, or an
elementary saddle--node, or no equilibria at all
(\figref{fig_saddlenode2}a), all other limit sets being topologically the
same.  

If $\Wglo{c}$ does form a loop, then any small perturbation of $f$ admits
either a saddle and a node near $x^\star$, or an elementary saddle--node
near $x^\star$, or a hyperbolic periodic orbit (\figref{fig_saddlenode2}b),
all other limit sets being topologically the same.  
\end{theorem}


\subsection{Hopf Bifurcation}
\label{ssec_hopf}

\begin{figure}
 \centerline{\psfig{figure=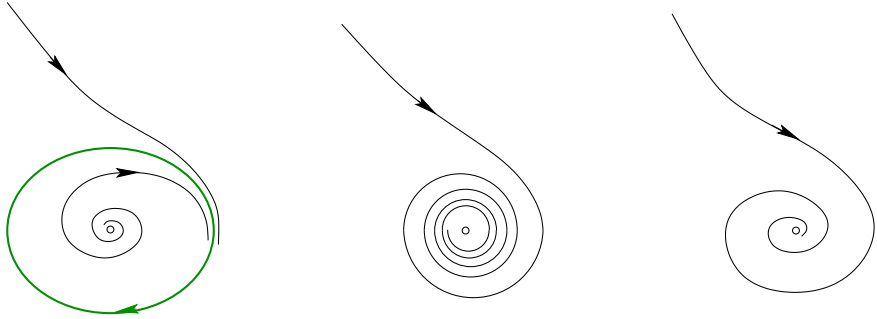,height=45mm,clip=t}}
 \figtext{
 }
 \captionspace
 \caption[]
 {The unfolding of an elementary composed focus is a
 Poincar\'e--Andronov--Hopf bifurcation. Depending on the sign of $a'(0)$,
 the direction of arrows may be reversed.}
\label{fig_hopf}
\end{figure}

We consider now another type of nonhyperbolic equilibrium point, namely we
assume that $f\in\cC^3$ vanishes at $x^\star$ and that the linearization of
$f$ at $x^\star$ has imaginary eigenvalues $\pm\icx\w_0$ with $\w_0\neq0$.
In appropriate coordinates, we can thus write
\begin{equation}
\label{hopf1}
\begin{split}
\dot x_1 &= -\w_0 x_2 + g_1(x_1,x_2) \\
\dot x_2 &= \w_0 x_1 + g_2(x_1,x_2),
\end{split}
\end{equation}
where $g_1, g_2$ are again nonlinear terms, that is, $g_1, g_2$ and their
derivatives vanish for $x_1=x_2=0$. The easiest way to describe the
dynamics is to introduce the complex variable $z=x_1+\icx x_2$. This amounts
to carrying out the change of variables
\begin{equation}
\label{hopf2}
x_1 = \frac{z+\cc z}2 = \re z, 
\qquad\qquad
x_2 = \frac{z-\cc z}{2\icx} = \im z.
\end{equation}
The system \eqref{hopf1} becomes
\begin{equation}
\label{hopf3}
\dot z = \icx\w_0 z + G(z,\cc z),
\end{equation}
where $G$ admits a Taylor expansion of the form 
\begin{equation}
\label{hopf4}
G(z,\cc z) = \sum_{2\leqs n+m \leqs 3} G_{nm}z^n\cc z^m + \order{\abs{z}^3}. 
\end{equation}
The next step is to simplify the nonlinear term $G$, using the theory of
\defwd{normal forms}. Observe to this end that if $h(z,\cc z)$ is a
homogeneous polynomial of degree $2$, then the variable $w=z+h(z,\cc z)$
satisfies the equation
\begin{equation}
\label{hopf5}
\dot w = \icx \w_0 z + \icx \w_0 \Bigbrak{z\dpar hz-\cc z\dpar h{\cc z}} +
G(z,\cc z) + \dpar hz G(z,\cc z) + \dpar h{\cc z} \cc G(z,\cc z).
\end{equation}
Now assume that $h$ satisfies the equation
\begin{equation}
\label{hopf6}
\icx \w_0 \Bigbrak{h - z\dpar hz + \cc z\dpar h{\cc z}} 
= G_2(z,\cc z) \defby \sum_{n+m=2} G_{nm}z^n\cc z^m.
\end{equation}
Then we see that 
\begin{equation}
\label{hopf7}
\dot w = \icx \w_0 w + G(z,\cc z) - 
G_2(z,\cc z) + \dpar hz G(z,\cc z) + \dpar h{\cc z} \cc G(z,\cc z).
\end{equation}
All terms of order $2$ have disappeared from this equation, while new terms
of order $3$ have appeared. The right-hand side can be expressed as a
function of $w$ and $\cc w$ without generating any new terms of order $2$.
The same procedure can be used to eliminate terms of order $3$ as well, by
solving an equation analogous to \eqref{hopf6}. Let us now consider this
equation in more detail. If $h$ is a sum of terms of the form $h_{nm}z^n\cc
z^m$, we get for each $(n,m)$ the equation 
\begin{equation}
\label{hopf8}
\icx\w_0 (1-n+m) h_{nm} = G_{nm}. 
\end{equation}
This equation can be solved if $m\neq n-1$. Hence the only term which cannot
be eliminated is $(n,m)=(2,1)$, that is, the term proportional to
$\abs{z}^2z$. We can thus reduce \eqref{hopf3} to the form
\begin{equation}
\label{hopf9}
\dot w = \icx \w_0 w + c_{21}\abs{w}^2 w + \order{\abs{w}^3},
\end{equation}
where the coefficient $c_{21}$ can be expressed as a function of the
$G_{nm}$. If we now introduce polar coordinates $w=r\e^{\icx\ph}$,
substitute in \eqref{hopf9} and take the real and imaginary parts, we end up
with the system
\begin{equation}
\label{hopf10}
\begin{split}
\dot r &= \re c_{21} r^3 + R_1(r,\ph)\\
\dot\ph &= \w_0 + \im c_{21}r^2 + R_2(r,\ph),
\end{split}
\end{equation}
where $R_1=\order{r^3}$ and $R_2=\order{r^2}$. We conclude that if $\re
c_{21}<0$, solutions starting sufficiently close to the origin will spiral
slowly towards $(0,0)$ ($r$ decreases like $1/\sqrt t$)
(\figref{fig_hopf}). If $\re c_{21}>0$, solutions slowly spiral outwards. 

\begin{definition}
\label{def_hopf}
Assume $x^\star$ is an equilibrium point such that the linearization of $f$
at $x^\star$ admits imaginary eigenvalues $\pm\icx\w_0$, $\w_0\neq0$, and
such that the normal form \eqref{hopf10} satisfies $\re c_{21}\neq0$. Then
$x^\star$ is called an \defwd{elementary composed focus}. 
\end{definition}

Let us now consider perturbations of the form $\dot x=f(x,\lambda)$ of
\eqref{hopf1}. Observe that since $f(0,0)=0$ and $\dpar fx(0,0)$ has no
vanishing eigenvalue, the implicit function theorem can be applied to show
the existence of a unique equilibrium $x^\star(\lambda)$ near $0$ for small
$\lambda$. The linearization at $x^\star(\lambda)$ has eigenvalues
$a(\lambda)\pm\icx\w(\lambda)$, where $a(0)=0$ and $\w(0)=\w_0$. In
appropriate coordinates, we thus have 
\begin{equation}
\label{hopf11}
\begin{split}
\dot x_1 &= a(\lambda)x_1 -\w(\lambda) x_2 + g_1(x_1,x_2,\lambda) \\
\dot x_2 &= \w(\lambda) x_1 + a(\lambda)x_2 + g_2(x_1,x_2,\lambda).
\end{split}
\end{equation}
The normal form can be computed in the same way, and becomes in polar
coordinates
\begin{equation}
\label{hopf12}
\begin{split}
\dot r &= a(\lambda)r + \re c_{21}(\lambda) r^3 + R_1(r,\ph,\lambda)\\
\dot\ph &= \w(\lambda) + \im c_{21}(\lambda)r^2 + R_2(r,\ph,\lambda).
\end{split}
\end{equation}
Then the following result holds:

\begin{theorem}
\label{thm_hopf}
Assume that $\re c_{21}(0)\neq0$ and let $H(\lambda)=a(\lambda)/\re
c_{21}(0)$. Then for sufficiently small $\lambda$,
\begin{itemiz}
\item	if $H(\lambda)>0$, $f$ has an isolated hyperbolic equilibrium
near $0$;
\item	if $H(\lambda)=0$, $f$ has an elementary composed focus near $0$;
\item	if $H(\lambda)<0$, $f$ has a hyperbolic equilibrium and a
hyperbolic periodic orbit near $0$, with opposite stability. 
\end{itemiz}
\end{theorem}

The proof should be clear if $R_1$ is identically zero. If not, this result
can be proved by computing $\tdtot r\ph$ and examining the properties of the
Poincar\'e map associated with the section $\ph=0$. 

We conclude that if $f$ admits an elementary composed focus, then $f$ is
singular of codimension $1$ (provided all other limit sets are hyperbolic
and there are no saddle connections). The family of equations $\dot z =
(\lambda+\icx\w_0)z + c_{21}\abs{z}^2z$ is a local unfolding of the
singularity. 


\subsection{Saddle--Node Bifurcation of a Periodic Orbit}
\label{ssec_snp}

A third class of structurally unstable vector fields consists of the vector
fields which admit a non-hyperbolic periodic orbit $\Gamma$. Recall that if
$\Pi$ is the Poincar\'e map associated with a transverse section $\Sigma$
through $x^\star\in\Gamma$, then $\Gamma$ is non-hyperbolic if
$\Pi'(x^\star)=1$, which means that the graph of $\Pi(x)$ is tangent to the
diagonal at $x^\star$. 

\begin{figure}
 \centerline{\psfig{figure=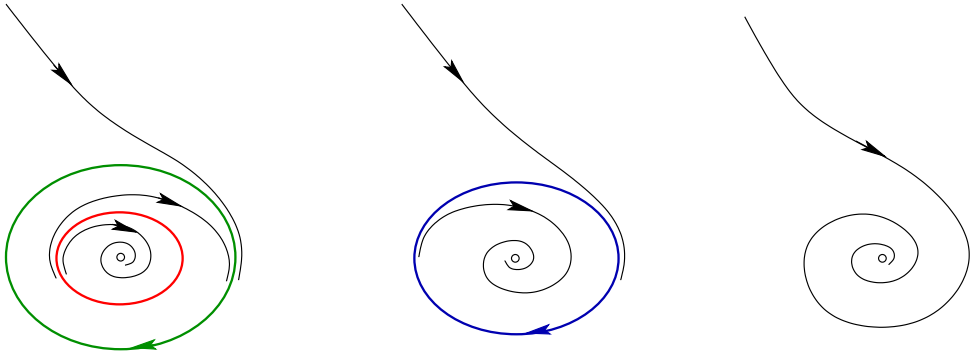,height=50mm,clip=t}}
 \figtext{
 }
 \captionspace
 \caption[]
 {The unfolding of a quasi-hyperbolic periodic orbit $\Gamma$ is a
 saddle--node bifurcation of periodic orbits: A small perturbation will
 have either two hyperbolic periodic orbits, or one quasi-hyperbolic
 periodic orbit, or no periodic orbit near $\Gamma$. The focus in the center
 of the pictures does not take part in the bifurcation, and can be replaced
 by a more complicated limit set.}
\label{fig_snorbit}
\end{figure}

\begin{definition}
\label{def_snp}
A periodic orbit $\Gamma$ is called \defwd{quasi-hyperbolic} if its
Poincar\'e map satisfies $\Pi'(x^\star)=1$ and $\Pi''(x^\star)\neq 0$. 
\end{definition}

The condition $\Pi''(x^\star)\neq 0$ means that the graph of $\Pi(x)$ has a
quadratic tangency with the diagonal. Hence, the orbit $\Gamma$ attracts
other orbits from one side, and repels them from the other side
(\figref{fig_snorbit}). 

Let us now consider a one-parameter family of perturbations $f(x,\lambda)$
of the vector field. For sufficiently small $\lambda$, $f(x,\lambda)$ will
still be transverse to $\Sigma$, and it admits a Poincar\'e map
$\Pi(x,\lambda)$ defined in some \nbh\ of $(x^\star,0)$. Periodic orbits
correspond to fixed points of $\Pi$, and are thus solutions of the equation 
$\Pi(x,\lambda)-x=0$. 

Note that $\Pi(x,0)-x$ is tangent to the $x$-axis at $x=x^\star$. This
situation is very reminiscent of the saddle--node bifurcation of equilibrium
points, see \eqref{snb6} and \eqref{snb7}. Indeed,
Proposition~\ref{prop_snb} applied to $\Pi(x,\lambda)-x=0$ yields the
existence of a function $H(\lambda)$ which controls the number of fixed
points of $\Pi$. 

\begin{theorem}
\label{thm_sno}
Let $f(x,0)$ admit a quasi-hyperbolic periodic orbit $\Gamma$. Then there is
a function $H(\lambda)$ such that, for sufficiently small $\lambda$,
\begin{itemiz}
\item	if $H(\lambda)>0$, then $f(x,\lambda)$ admits no periodic orbit near
$\Gamma$;
\item	if $H(\lambda)=0$, then $f(x,\lambda)$ admits a quasi-hyperbolic
periodic orbit near $\Gamma$ ;
\item	if $H(\lambda)<0$, then $f(x,\lambda)$ admits two hyperbolic
periodic orbits of opposite stability near $\Gamma$.
\end{itemiz}
Moreover, $f(x,0)$ is singular of codimension $1$ if and only if there do
not exist two invariant manifolds of saddles which approach $\Gamma$, one for
$t\to\infty$, the other for $t\to-\infty$.
\end{theorem}

\begin{figure}[b]
 \centerline{\psfig{figure=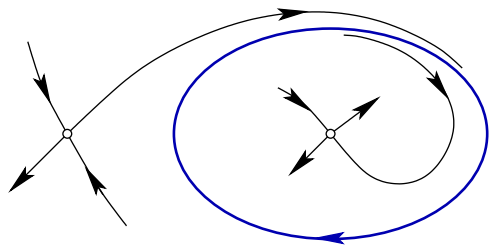,height=35mm,clip=t}}
 \figtext{
 }
 \captionspace
 \caption[]
 {Partial phase portrait near a quasi-hyperbolic periodic orbit which is not
 a codimension $1$ singularity: a small perturbation may admit a saddle
 connection.}
\label{fig_snbad}
\end{figure}

The reason for this last condition is that if such manifolds exist, then
one can construct perturbations which admit a saddle connection, and are
thus neither structurally stable, nor topologically equivalent to $f$  
(\figref{fig_snbad}).

\goodbreak


\subsection{Global Bifurcations}
\label{ssec_gb}

The remaining singular vector fields of codimension $1$ are those which
admit a saddle connection. Such saddle connections are indeed very sensitive
to small perturbations.

Consider first the case of a heteroclinic saddle connection between two
saddle points $x^\star_1$ and $x^\star_2$. To simplify the notations, let
us choose coordinates in such a way that $x^\star_1=(0,0)$ and
$x^\star_2=(1,0)$, and that the saddle connection belongs to the
$x_1$-axis. We write the perturbed system as 
\begin{equation}
\label{gb1}
\begin{split}
\dot x_1 &= f_1(x_1,x_2) + \lambda g_1(x_1,x_2,\lambda) \\
\dot x_2 &= f_2(x_1,x_2) + \lambda g_2(x_1,x_2,\lambda),
\end{split}
\end{equation}
where, in particular, $f_2(x_1,0)=0$ for all $x_1\in[0,1]$. We want to find
a condition for this system to admit a saddle connection for $\lambda\neq0$.
Note first that since $x^\star_1$ and $x^\star_2$ are hyperbolic, the
perturbed system \eqref{gb1} still admits two saddles near $x^\star_1$ and
$x^\star_2$ for small $\lambda$. Carrying out an affine coordinate
transformation, we may assume that they are still located at $(0,0)$ and
$(1,0)$. 

\begin{figure}
 \centerline{\psfig{figure=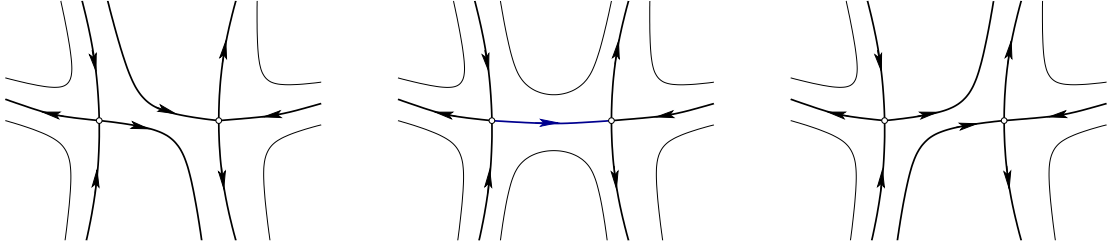,height=35mm,clip=t}}
 \figtext{
 }
 \captionspace
 \caption[]
 {Breaking a heteroclinic saddle connection.}
\label{fig_heteroclinic}
\end{figure}

Assume that the saddle connection has an equation of the form $x_2 =
\lambda h_1(x_1) + \Order{\lambda^2}$, with $h_1(0)=h_1(1)=0$ (such an
expansion can be proved to hold for sufficiently small $\lambda$). By
Taylor's formula, we find that on this curve,
\begin{equation}
\label{gb2}
\begin{split}
\dot x_1 &= f_1(x_1,0) + \Order{\lambda} \\
\dot x_2 &= \lambda \dpar{f_2}{x_2}(x_1,0)h_1(x_1) + \lambda g_2(x_1,0,0) 
+ \Order{\lambda^2}. 
\end{split}
\end{equation}
Thus $h_1$ must satisfy the linear equation
\begin{equation}
\label{gb3}
h'_1(x_1) = \lim_{\lambda\to0} \frac{\dot x_2}{\lambda\dot x_1} = 
\frac1{f_1(x_1,0)} \biggbrak{\dpar{f_2}{x_2}(x_1,0)h_1(x_1) + g_2(x_1,0,0)},
\end{equation}
which admits the solution (using $h_1(0)=0$ as initial value)
\begin{equation}
\label{gb4}
h_1(x_1) = \int_0^{x_1} \exp\biggset{\int_y^{x_1} \frac1{f_1(z,0)}
\dpar{f_2}{x_2}(z,0)\6z} \frac{g_2(y,0,0)}{f_1(y,0)}\6y.
\end{equation}
This relation can also be written in the more symmetric form
\begin{equation}
\label{gb5}
h_1(x_1) = \frac1{f_1(x_1,0)} 
\int_0^{x_1} \exp\biggset{\int_y^{x_1} \frac1{f_1}
\Bigbrak{\dpar{f_1}{x_1}+\dpar{f_2}{x_2}}(z,0)\6z} g_2(y,0,0)\6y.
\end{equation}
Note that only the last factor in the integrand depends on the perturbation
term $g$. Here we have only computed the first order term of
$h(x_1,\lambda) = \lambda h_1(x_1)+\Order{\lambda^2}$. Including higher
order terms, the requirement $h(1,\lambda)=0$ provides us with a
codimension $1$ condition of the form $H(\lambda)=0$. This kind of integral
condition is closely related to the method of \defwd{Melnikov functions},
which allows to estimate the splitting between manifolds. 

If the condition $h(1,\lambda)=0$ is not satisfied, then the unstable
manifold of $x^\star_1$ and the stable manifold of $x^\star_2$ will avoid
each other, and the sign of $h$ will determine which manifold lies above the
other (\figref{fig_heteroclinic}). (Strictly speaking, $\Wglo{u}(x^\star_1)$
may not be the graph of a function when $x_1$ approaches $1$, and one would
rather compute estimations of $\Wglo{u}(x^\star_1)$ for $0<x_1<1/2$ and of
$\Wglo{s}(x^\star_2)$ for $1/2<x<1$, and determine their splitting at
$1/2$.) 

\medskip

\begin{figure}
 \centerline{\psfig{figure=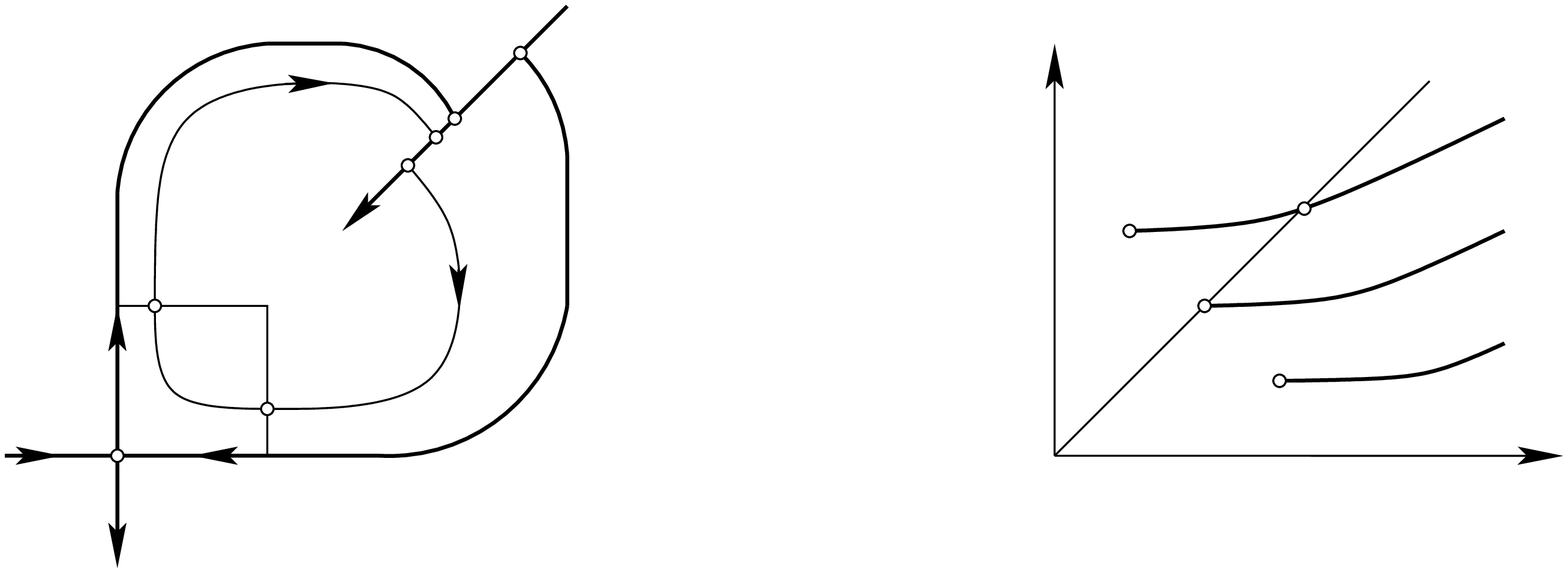,height=44mm,clip=t}}
 \figtext{
 	\writefig	2.1	1.6	$x^\star$
 	\writefig	3.7	1.1	$\Wglo{s}(x^\star)$
 	\writefig	1.2	2.8	$\Wglo{u}(x^\star)$
 	\writefig	5.9	4.6	$\Sigma$
 	\writefig	5.0	4.6	$x_-$
 	\writefig	5.0	3.8	$x_+$
 	\writefig	3.7	2.1	$(\eps,x^0_2)$
 	\writefig	2.9	2.8	$(x^0_1,\eps)$
 	\writefig	13.2	1.5	$x$
 	\writefig	13.0	2.3	$\Delta<0$
 	\writefig	13.0	3.2	$\Delta=0$
 	\writefig	13.0	4.1	$\Delta>0$
 	\writefig	9.6	4.4	$\Pi(x)$
 }
 \vspace{1mm}
 \caption[]
 {(a) Definition of the first return map $\Pi$ associated with a section
 $\Sigma$ transverse to the stable and unstable manifolds of $x^\star$. 
 (b) Behaviour of the first return map in a case where $\abs{a_1}>a_2$ for
 different signs of the splitting $\Delta = x_+- x_-$. $\Pi$ admits a fixed
 point when the splitting is positive, which corresponds to a periodic
 orbit. If $\abs{a_1}<a_2$, then $\Pi$ admits a fixed point when the
 splitting is negative.} 
\label{fig_homoclinic1}
\end{figure}

Let us finally consider the other kind of saddle connection, also called a
\defwd{homoclinic loop} $\Gamma$. In order to describe the flow near
$\Gamma$, let us introduce a transverse arc $\Sigma$ through a point
$x_0\in\Gamma$. We order the points on $\Sigma$ in the inward direction
(\figref{fig_homoclinic1}a). Orbits starting slightly \lq\lq above\rq\rq\
$x_0$ with respect to this ordering will return to $\Sigma$, defining a
first return map $\Pi$ for small $x>x_0$. 

For small $\lambda\neq0$, the stable and unstable manifold of $x^\star$
will intersect $\Sigma$ at points $x_-(\lambda)$ and $x_+(\lambda)$, and
their distance $\Delta(\lambda) = x_+(\lambda)- x_-(\lambda)$ may be
estimated as in the case of a heteroclinic saddle connection. A first
return map $\Pi(x,\lambda)$ can be defined for sufficiently small $x>x_-$
(\figref{fig_homoclinic1}a). 

Our aim will now be to determine the behaviour of $\Pi$ in the limit $x\to
x_-$. One can show that in this limit, the first return map depends
essentially on the behaviour of the orbits as they pass close to $x^\star$.
Assume that the linearization of $f$ at $x^\star$ has eigenvalues
$a_1<0<a_2$, and choose a coordinate system where this linearization is
diagonal. We consider an orbit connecting the points $(\eps,x^0_2)$ and
$(x^0_1,\eps)$ with $0<x^0_2<\eps\ll1$, and would like to determine $x^0_1$ as
a function of $x^0_2$. In a linear approximation, we have
\begin{align}
\nonumber
\dot x_1 &= -\abs{a_1}x_1 
&& \Rightarrow & x_1(t) &= \e^{-\abs{a_1}t}\eps \\
\label{gb6}
\dot x_2 &= a_2x_2 
&& \Rightarrow & x_2(t) &= \e^{a_2t}x^0_2.
\end{align} 
Requiring that $x_2(t)=\eps$ and eliminating $t$, we obtain 
\begin{equation}
\label{gb7}
\frac{x^0_1}{x^0_2} = \frac{\e^{-\abs{a_1}\log(\eps/x^0_2)/a_2}\eps}{x^0_2} =
\Bigpar{\frac\eps{x^0_2}}^{1-\abs{a_1}/a_2}. 
\end{equation}
\goodbreak
\noindent
We conclude that 
\begin{equation}
\label{gb8}
\lim_{x^0_2\to0} \frac{x^0_1}{x^0_2} = 
\begin{cases}
0 & \text{if $\abs{a_1}>a_2$} \\
\infty & \text{if $\abs{a_1}<a_2$.}
\end{cases}
\end{equation}
By controlling the error terms a bit more carefully, one can show that
$\lim_{x\to x_-}\Pi'(x)$ has the same behaviour, as long as $\abs{a_1}/a_2$
is bounded away from $1$. 

\begin{definition}
\label{def_gb}
Let $x^\star$ be a saddle point with eigenvalues $a_1<0<a_2$ satisfying
$a_1+a_2\neq0$. An orbit $\Gamma$ such that $\w(y)=\alpha(y)=x^\star$ for
all $y\in\Gamma$ is called an \defwd{elementary homoclinic loop}. 
\end{definition}

\begin{figure}
 \centerline{\psfig{figure=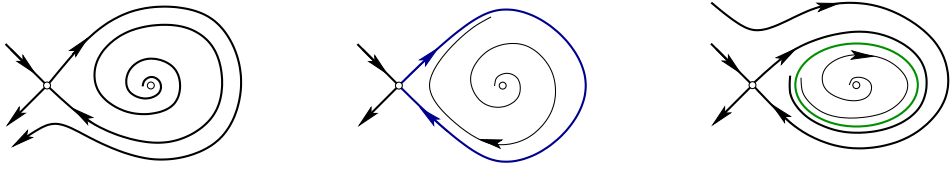,height=25mm,clip=t}}
 \figtext{
 }
 \captionspace
 \caption[]
 {Unfolding an elementary homoclinic loop $\Gamma$: any small perturbation
 will either admit a hyperbolic periodic orbit close to $\Gamma$, or an
 elementary homoclinic loop, or no invariant set near $\Gamma$.}
\label{fig_homoclinic2}
\end{figure}

In summary, we know that as $x\to x_-$, the first return map of an
elementary homoclinic loop satisfies
\begin{align}
\nonumber
\lim_{x\to x_-} \Pi(x) &= x_- + \Delta(\lambda) \\
\label{gb9}
\lim_{x\to x_-} \Pi'(x) &=
\begin{cases}
0 & \text{if $\abs{a_1}>a_2$} \\
\infty & \text{if $\abs{a_1}<a_2$.}
\end{cases}
\end{align}
The shape of $\Pi$ for different values of $\Delta$ is sketched in
\figref{fig_homoclinic1}b, in a case where $\abs{a_1}>a_2$. We see that if
$\Delta>0$, $\Pi$ must admit a fixed point, which corresponds to a periodic
orbit (\figref{fig_homoclinic2}). If $\abs{a_1}<a_2$, there is a periodic
orbit for $\Delta<0$.

\begin{theorem}
\label{thm_gb}
Let $f$ admit an elementary homoclinic loop $\Gamma$. Given a
one-parameter family of perturbations $f(x,\lambda)$, there exists a
function $H(\lambda) = \Delta(\lambda) \log(\abs{a_1}/a_2)$ (where $\Delta$
describes the splitting of stable and unstable manifolds) such that
\begin{itemiz}
\item	if $H(\lambda)<0$, then there are no invariant sets near $\Gamma$;
\item	if $H(\lambda)=0$, then $f(x,\lambda)$ admits an elementary
homoclinic loop near $\Gamma$;
\item	if $H(\lambda)>0$, then $f(x,\lambda)$ admits a hyperbolic 
periodic orbit near $\Gamma$.
\end{itemiz}
If, moreover, all other equilibria and periodic orbits of $f$ are
hyperbolic, there are no saddle connections, and no stable or unstable
manifold of another saddle point tends to $\Gamma$, then $f$ is a
singularity of codimension $1$. 
\end{theorem}

This completes the list of all structurally unstable vector fields of
codimension $1$: $f\in\cS_1$ if and only if it satisfies the hypotheses of
Peixoto's theorem, except for one limit set, chosen from the following
list:
\begin{enum}
\item	an elementary saddle-node (satisfying the hypotheses of
Theorem~\ref{thm_snb});
\item	an elementary composed focus;
\item	a quasi-hyperbolic periodic orbit (satisfying the hypotheses of
Theorem~\ref{thm_sno});
\item	a heteroclinic saddle connection;
\item	an elementary homoclinic loop (satisfying the hypotheses of
Theorem~\ref{thm_gb}).
\end{enum}

\goodbreak


\section{Local Bifurcations of Codimension 2}
\label{sec_sc2}

Once we have understood the set $\cS_1$ of singular vector fields of
codimension $1$, the next step would be to characterize the remaining
structurally unstable vector fields in $\cS_0\setminus\cS_1$. In analogy
with the codimension $1$ case, we can define a set $\cS_2\subset
\cS_0\setminus\cS_1$ of singular vector fields of codimension $2$ in two
equivalent ways:
\begin{enum}
\item	a vector field $f\in\cS_0\setminus\cS_1$ is singular of codimension
$2$ if any sufficiently small perturbation $g$ of $f$ which still belongs to
$\cS_0\setminus\cS_1$ is topologically equivalent to $f$;
\item	a vector field $f\in\cS_0\setminus\cS_1$ is singular of codimension
$2$ if there exists a two-parameter family $F_{\lambda_1\lambda_2}$ of
vector fields such that {\em any} sufficiently small perturbation of $f$ is
topologically equivalent to a member of the family $F_{\lambda_1\lambda_2}$.
\end{enum}
More generally, we could define singular vector fields of higher
codimension, and thus decompose $\cS_0$ into a disjoint union of components
\begin{equation}
\label{sc2_1}
\cS_0 = \cS_1 \cup \cS_2 \cup \dots \cup \cS_k \cup \dots
\end{equation}
The classification of vector fields in $\cS_k$ becomes, however,
increasingly difficult as $k$ grows, and is not even complete in the case
$k=2$. The discussion in the previous section shows that the following
vector fields are candidates for belonging to $\cS_2$: 
\begin{enum}
\item	vector fields containing {\em two} non-hyperbolic equilibrium
points, periodic orbits or saddle connections (or any combination of them); 
\item	vector fields with a non-hyperbolic equilibrium or periodic orbit
having more than one vanishing term in their Taylor series;
\item	vector fields with a non-hyperbolic equilibrium point having a
double zero eigenvalue;
\item	structurally unstable vector fields which may admit a saddle
connection after a small perturbation (as in \figref{fig_saddlenode1}b and
\figref{fig_snbad}).
\end{enum}
The first type is not very interesting, and the fourth type is very hard to
analyse. We will discuss below the simplest examples of types 2.\
and 3. 


\subsection{Pitchfork Bifurcation}
\label{ssec_pb}

Let $f\in\cC^3$ have a non-hyperbolic equilibrium point $x^\star$, such that
the Jacobian matrix $\dpar fx(x^\star)$ admits the eigenvalues $0$ and
$a\neq0$. Unlike in the case of the saddle--node bifurcation, let us assume
that the dynamics on the center manifold is governed by the equation
$\dot{x_1}=F(x_1)$ with 
\begin{equation}
\label{pb1}
F(0) = \dpar Fx(0) = \dpar{^2F}{x^2}(0) = 0, \qquad
\dpar{^3F}{x^3}(0) = 6c \neq 0.
\end{equation} 
Then $\dot{x_1}= c x_1^3 + \order{x_1^3}$, and thus $x^\star$ will attract
nearby solutions on the center manifold from both sides if $c<0$, and repel
them if $c>0$, but at a much slower rate than in the hyperbolic case. 

\begin{figure}
 \centerline{\psfig{figure=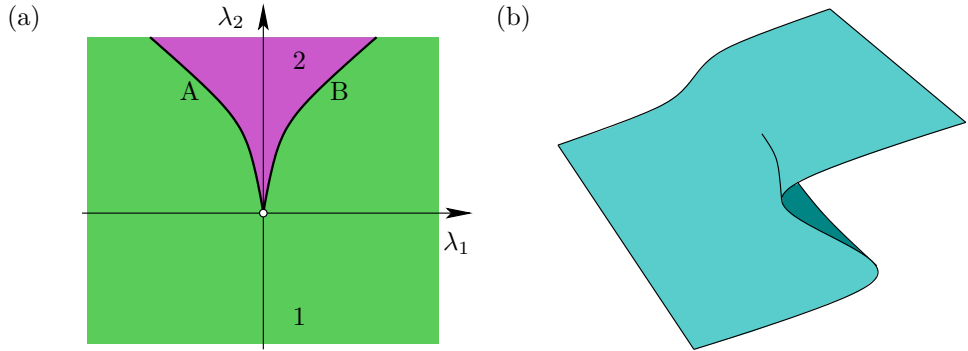,height=50mm,clip=t}}
 \figtext{
 	\writefig	0.5	5.0	(a)
 	\writefig	7.0	5.0	(b)
	\writefig	4.3	1.0	$1$
	\writefig	4.3	4.4	$2$
	\writefig	4.8	4.0	B
	\writefig	2.8	4.0	A
	\writefig	6.3	2.0	$\lambda_1$
	\writefig	3.3	5.0	$\lambda_2$
 }
 \captionspace
 \caption[]
 {(a) Bifurcation diagram of the family $\dot x_1=\lambda_1 + \lambda_2 x_1
 - x_1^3$. The lines $A$ and $B$ correspond to saddle--node bifurcations.
 (b) The surface $\lambda_1 + \lambda_2 x_1 - x_1^3=0$ has a fold in the
 space $(\lambda_1,\lambda_2,x_1)$.}
\label{fig_pitchfork0}
\end{figure}

A perturbation $f(x,\lambda)$ of $f$ will admit a center manifold on which
the dynamics is governed by an equation $\dot{x_1}=F(x_1,\lambda)$. We
claim that for sufficiently small $\lambda$ and $x_1$, this system is
topologically equivalent to a member of the $2$-parameter family 
\begin{equation}
\label{pb2}
\dot x_1 = F_{\lambda_1\lambda_2}(x_1) = \lambda_1 + \lambda_2 x_1 \pm
x_1^3,
\end{equation}
where the signs $\pm$ correspond to the cases $c>0$ or $c<0$. Let us focus
on the $-$ case, since the other case can be obtained by inverting the
direction of time. As we saw in Example~\ref{ex_1D2}, the dynamics for
$(\lambda_1,\lambda_2)\neq(0,0)$ can be of three types
(\figref{fig_unfold1D}):
\begin{enum}
\item	if $4\lambda_2^3<27\lambda_1^2$, $F_{\lambda_1\lambda_2}$ has
exactly one equilibrium point, which is hyperbolic and stable; 
\item	if $4\lambda_2^3>27\lambda_1^2$, $F_{\lambda_1\lambda_2}$ has three
equilibrium points, all being hyperbolic, and two of them stable while the
third one is unstable;
\item	if $4\lambda_2^3=27\lambda_1^2$, $F_{\lambda_1\lambda_2}$ has two
equilibrium points; one of them is hyperbolic and stable, and the other one
is nonhyperbolic with a quadratic tangency to the $x_1$-axis, 
corresponding to an elementary saddle--node of the associated
two-dimensional system.
\end{enum}
The singular lines $4\lambda_2^3=27\lambda_1^2$ form a cusp in the
$(\lambda_1,\lambda_2)$-plane (\figref{fig_pitchfork0}). A one-parameter
family of vector fields crossing one of these curves away from the cusp
undergoes a saddle--node bifurcation. Crossing the origin from one side of
the cusp to the other one changes the number of equilibria from $1$ to $3$
and is called a \defwd{pitchfork bifurcation}. 

\begin{figure}
 \centerline{\psfig{figure=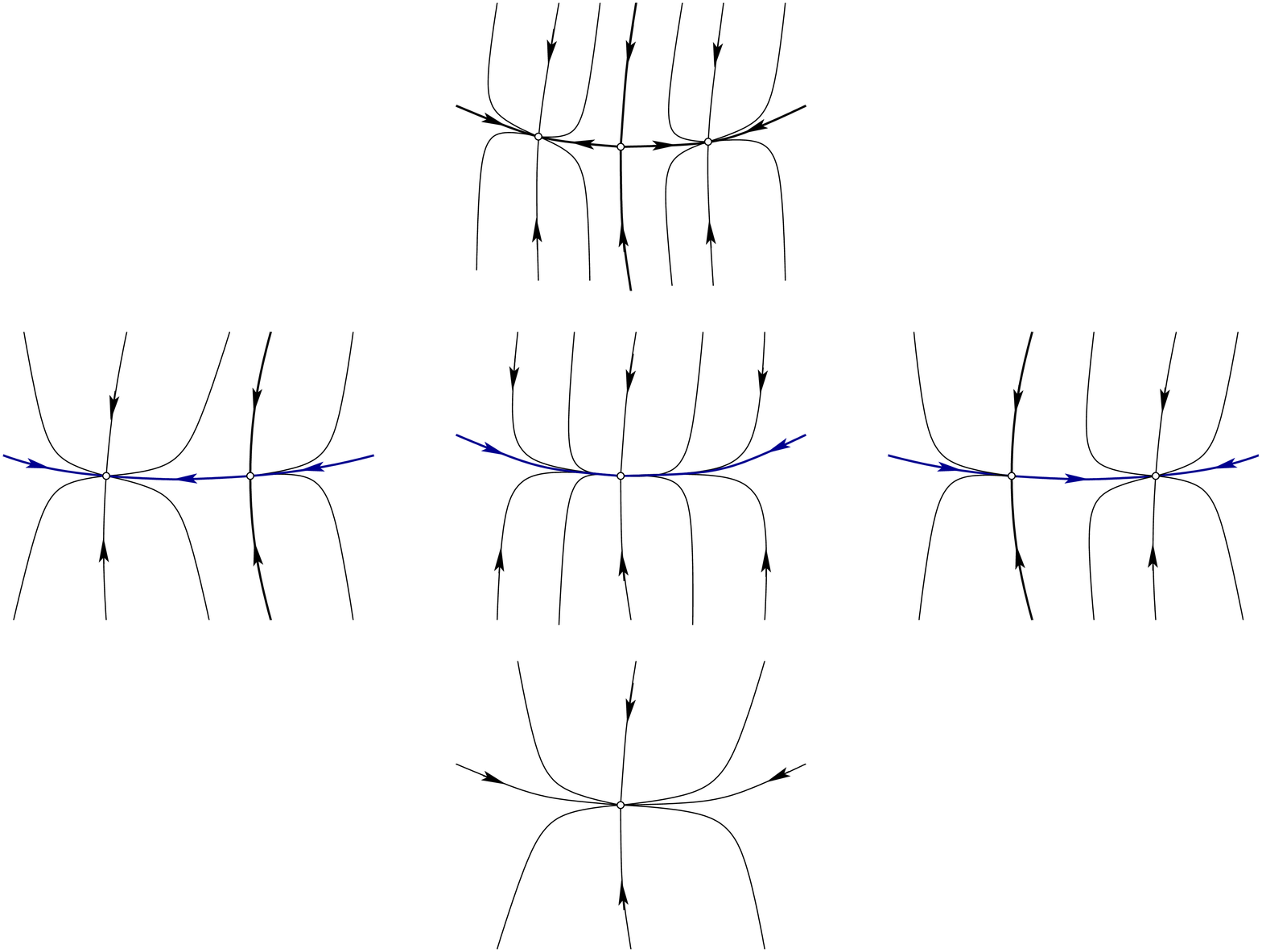,width=130mm,clip=t}}
 \figtext{
 	\writefig	0.5	6.7	A
 	\writefig	5.3	6.7	O
	\writefig	9.5	6.7	B
	\writefig	5.3	10.0	2
	\writefig	5.3	3.4	1
 }
 \captionspace
 \caption[]
 {Unfolding a singular vector field (O) equivalent to $(-x_1^3,-x_2)$.
 Structurally stable perturbations have either one hyperbolic equilibrium
 (1) or three hyperbolic equilibria (2). They are separated by singular
 vector fields containing an elementary saddle--node and a hyperbolic
 equilibrium (A,B).}
\label{fig_pitchfork}
\end{figure}

\begin{prop}
\label{prop_pb}
The two-parameter family \eqref{pb2} is a local unfolding of the singular
vector field whose restriction to the center manifold satisfies \eqref{pb1}.
\end{prop}
\begin{proof}
The function 
\[
G(x_1,\lambda) = \dpar F{x_1}(x_1,\lambda)
\]
satisfies the hypotheses of $F$ in Proposition~\ref{prop_snb}, that we used
in the case of the saddle--node bifurcation. There we showed the existence
of a function $H(\lambda)$ such that, in a \nbh\ of $x_1=0$, $G$ vanishes
twice if $H<0$, once if $H=0$ and never if $H>0$. It follows that $F$ has
two extrema if $H<0$, is monotonous with an inflection point if $H=0$, and
is strictly monotonous if $H>0$. This shows in particular that $H$ plays
the role of $\lambda_2$ in the unfolding \eqref{pb2}. 

Now let $\ph(\lambda)$ be the extremum of $G$ (see
Proposition~\ref{prop_snb}), i.e.\ the inflection point of $F$. Then the
Taylor expansion of $F$ around $\ph(\lambda)$ can be written as 
\[
\frac1cF(\ph(\lambda)+y,\lambda) = \frac1cF(\ph(\lambda),\lambda) 
+ \frac1c\dpar F{x_1}(\ph(\lambda),\lambda)y + y^3\bigbrak{1+R(y,\lambda)},
\]
where $R(y,\lambda)$ is a continuous function with $R(0,0)=0$. The
coefficient of $y$ is exactly $H(\lambda)= \lambda_2$. Since
$F(\ph(\lambda),\lambda)/c$ controls the vertical position of the graph of
$F$, it plays the same r\^ole as $\lambda_1$ (one can show that $F$ has the
same behaviour as $F_{\lambda_1\lambda_2}$ for a $\lambda_1$ of the form  
$c^{-1}F(\ph(\lambda),\lambda)[1+\rho(\lambda)]$, with $\rho(\lambda)\to0$
as $\lambda\to0$).  
\end{proof}


\subsection{Takens--Bogdanov Bifurcation}
\label{ssec_tb}

Consider now the case where $f$ has a non-hyperbolic equilibrium point
$x^\star$, such that the linearization of $f$ at $x^\star$ admits $0$ as a
double eigenvalue. Then $\dpar fx(x^\star)$ can admit one of the two Jordan
canonical forms
\begin{equation}
\label{tb1}
\begin{pmatrix}
0 & 1 \\ 0 & 0
\end{pmatrix}
\qquad\text{or}\qquad
\begin{pmatrix}
0 & 0 \\ 0 & 0
\end{pmatrix}.
\end{equation}
The second case does not correspond to a codimension $2$ bifurcation. This
can be roughly understood as follows: Any $2\times2$ matrix is equivalent to
a triangular matrix 
\begin{equation}
\label{tb2}
\begin{pmatrix}
\lambda_1 & \lambda_3 \\ 0 & \lambda_2
\end{pmatrix}.
\end{equation}
None of the systems obtained by setting two among the three $\lambda_j$
equal to zero is topologically equivalent to the system with all three
$\lambda_j$ equal to zero, and thus two parameters do not suffice to
describe all possible perturbations of $f$. 

\begin{figure}
 \centerline{\psfig{figure=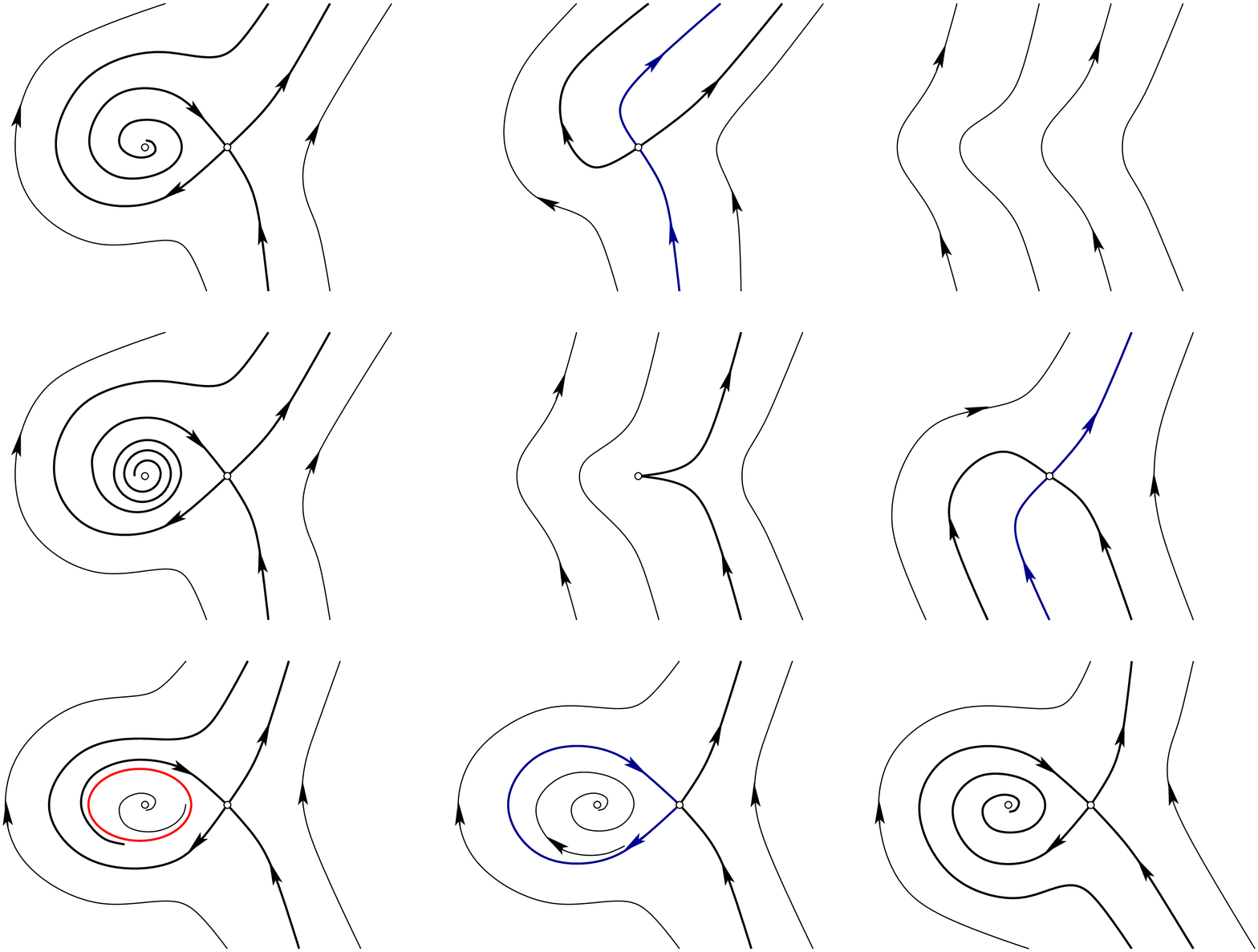,width=130mm,clip=t}}
 \figtext{
 	\writefig	0.7	10.0	$2$
 	\writefig	0.7	6.7	B
 	\writefig	0.7	3.4	$3$
 	\writefig	5.5	10.0	A
 	\writefig	5.5	6.7	O
 	\writefig	5.5	3.4	C
 	\writefig	9.7	10.0	$1$
 	\writefig	9.7	6.7	D
 	\writefig	9.7	3.4	$4$
 }
 \captionspace
 \caption[]
 {Unfolding the singular vector field \eqref{tb8} give rise to four
 structurally stable vector fields ($1$--$4$), and to four singular vector
 fields of codimension $1$ (A--D).}
\label{fig_takens}
\end{figure}

Consider now the first case. In appropriate coordinates, we can write 
\begin{equation}
\label{tb3}
\begin{split}
\dot x_1 &= x_2 + g_1(x_1,x_2) \\
\dot x_2 &= g_2(x_1,x_2),
\end{split}
\end{equation}
where $g_1, g_2$ and all their first order derivatives vanish at the
origin. These nonlinear terms can be simplified by a change of variables,
and the general theory describing how to do this is the theory of
\defwd{normal forms}. We will only illustrate it in this particular case.
First, note that the Jacobian of the change of variables
$(x_1,x_2)\mapsto(x_1,x_2+g_1(x_1,x_2))$ is close to unity near
$(x_1,x_2)=0$, and transforms the system into
\begin{equation}
\label{tb4}
\begin{split}
\dot x_1 &= x_2 \\
\dot x_2 &= g(x_1,x_2),
\end{split}
\end{equation}
where an easy calculation shows that $g$ has the same properties as $g_1$
and $g_2$. Note that this system is equivalent to the second order equation 
\begin{equation}
\label{tb5}
\ddot x_1 = g(x_1,\dot x_1),
\end{equation} 
called a \defwd{Li\'enard equation}. As we did in the case of the Hopf
bifurcation, we can try to simplify the term $g$ further by introducing the
new variable $y_1=x_1+h(x_1,\dot x_1)$, where $h$ is a homogeneous
polynomial of degree $2$. Differentiating $y_1$ and using \eqref{tb5}, we
obtain 
\begin{equation}
\label{tb6}
\dot y_1 = \dot x_1 \Bigpar{1+\dpar h{x_1}} + \dpar h{\dot x_1}g(x_1,\dot x_1).
\end{equation}
If we differentiate this relation again, we obtain several terms, but most
of them are of order $3$ (or larger):
\begin{equation}
\label{tb7}
\ddot y_1 = g(x_1,\dot x_1) + \dot x_1^2 \dpar {^2h}{x_1^2} + \text{terms of
order $3$}.
\end{equation}
Now if $h(x_1,\dot x_1) = h_1 x_1^2 + h_2 x_1\dot x_1 + h_3\dot x_1^2$, we
see that we can choose $h_1$ in such a way as to eliminate the term
proportional to $\dot x_1^2$ of $g$, while $h_2$ and $h_3$ have no
influence on quadratic terms. The other second order terms of $g$ cannot be
eliminated, and are called \defwd{resonant}. Expressing the remaining terms
as functions of $y_1, y_2$, we obtain the normal form of \eqref{tb4} to
second order
\begin{equation}
\label{tb8}
\begin{split}
\dot y_1 &= y_2 \\
\dot y_2 &= c_1 y_1^2 + c_2 y_1 y_2 + R(y_1,y_2),
\end{split}
\end{equation}
where $c_1$ and $c_2$ are Taylor coefficients of $g$ and $R$ is a remainder
of order $3$. Takens has proved the nontrivial fact that if $c_1\neq 0$,
then the remainder $R$ does not change the topology of the orbits near the
origin. If $c_1\neq0$ and $c_2\neq0$, one can reduce the problem, by
scaling $y_1$, $y_2$ and $t$, to one of the cases $c_1=1$ and $c_2=\pm1$.
We shall consider here the case $c_2=1$, that is, 
\begin{equation}
\label{tb9}
\begin{split}
\dot y_1 &= y_2 \\
\dot y_2 &= y_1^2 + y_1 y_2.
\end{split}
\end{equation}
The case $c_2=-1$ is obtained by inverting the direction of time (and
changing the sign of $y_2$). The system \eqref{tb9} is not straightforward
to analyse by itself, but one can check that there exist two orbits of the
form 
\begin{equation}
\label{tb10}
y_2 = \pm \Bigpar{\frac23}^{1/2} y_1^{3/2} \brak{1+\rho(y_1)}, 
\qquad
y_1>0,
\qquad
\rho(y_1)=\Order{\sqrt{y_1}\mskip1.5mu},
\end{equation}
which form a cusp (\figref{fig_takens}--O). Takens and Bogdanov have shown
that the following two-parameter family of vector fields constitutes a
local unfolding of the singular vector field \eqref{tb9}:
\begin{equation}
\label{tb11}
\begin{split}
\dot y_1 &= y_2 \\
\dot y_2 &= \lambda_1 + \lambda_2 y_2 + y_1^2 + y_1 y_2.
\end{split}
\end{equation}
(Various equivalent unfoldings are found in the literature).  Proving this
lies beyond the scope of the present course. Let us, however, examine the
dynamics of \eqref{tb11} as a function of $\lambda_1$ and $\lambda_2$. If
$\lambda_1>0$, there are no equilibrium points. If $\lambda_1<0$, there are
two equilibria 
\begin{equation}
\label{tb12}
x^\star_\pm = (\pm\sqrt{-\lambda_1},0).
\end{equation}
The linearization of $f$ at $x^\star_+$ has the eigenvalues 
\begin{equation}
\label{tb13}
\frac{\lambda_2+\sqrt{-\lambda_1}}2 \pm 
\sqrt{\frac{(\lambda_2+\sqrt{-\lambda_1})^2}4 + 2\sqrt{-\lambda_1}},
\end{equation}
and thus $x^\star_+$ is always a saddle. As for $x^\star_-$, the
eigenvalues are 
\begin{equation}
\label{tb14}
\frac{\lambda_2-\sqrt{-\lambda_1}}2 \pm 
\sqrt{\frac{(\lambda_2-\sqrt{-\lambda_1})^2}4 - 2\sqrt{-\lambda_1}},
\end{equation}
which shows that $x^\star_-$ is a node or a focus, unless
$\lambda_2=\sqrt{-\lambda_1}$. At $\lambda_2=\sqrt{-\lambda_1}$, there is a
Hopf bifurcation since the eigenvalues \eqref{tb14} are imaginary. One can
check that the normal form is such that an unstable periodic orbit appears
for $\lambda_2<\sqrt{-\lambda_1}$. 

With these informations, one can already draw some phase portraits. We know
that there are no equilibrium points for $\lambda_1>0$
(\figref{fig_takens}--1). Crossing the $\lambda_2$-axis corresponds to a
saddle--node bifurcation, producing a saddle and a source if $\lambda_2>0$
(\figref{fig_takens}--2), and a saddle and a sink if $\lambda_2<0$
(\figref{fig_takens}--4). The source expels an unstable periodic orbit
when the curve $\lambda_2=\sqrt{-\lambda_1}$ is crossed
(\figref{fig_takens}--3). 

\begin{figure}
 \centerline{\psfig{figure=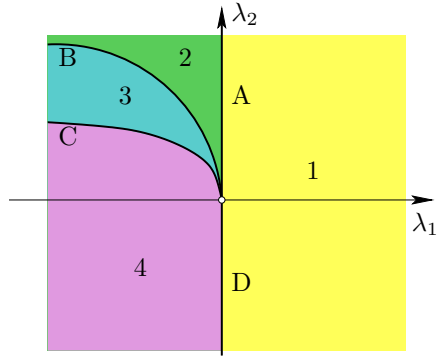,height=50mm,clip=t}}
 \figtext{
 	\writefig	9.9	2.3	$\lambda_1$
 	\writefig	7.5	5.1	$\lambda_2$
	\writefig	8.5	3.0	$1$
	\writefig	6.8	4.5	$2$
	\writefig	6.0	4.0	$3$
	\writefig	6.2	1.7	$4$
	\writefig	7.5	4.0	A
	\writefig	5.2	4.5	B
	\writefig	5.2	3.45	C
	\writefig	7.5	1.5	D
 }
 \captionspace
 \caption[]
 {Bifurcation diagram of the Takens--Bogdanov singularity. The lines A and D
 correspond to saddle--node bifurcations, B to a Hopf bifurcation, and C to
 a homoclinic loop bifurcation.}
\label{fig_takens0}
\end{figure}

There must be another bifurcation allowing to connect the vector fields
slightly below $\lambda_2=\sqrt{-\lambda_1}$, which admit a periodic orbit,
and slightly to the left of $\lambda_1=0$, $\lambda_2<0$, which have no
periodic orbit. The periodic orbit turns out to disappear in a homoclinic
bifurcation, occurring on the curve
\begin{equation}
\label{tb15}
\lambda_2 = \frac57 \sqrt{-\lambda_1} \bigbrak{1+\Order{\abs{\lambda_1}^{1/4}}}, 
\qquad \lambda_1<0,
\end{equation}
see \figref{fig_takens0}. 
We omit the proof of this fact, which is rather elaborate and relies on a
rescaling, transforming the system into a close to Hamiltonian form, and on
the computation of a Melnikov function, allowing to determine the splitting
of the stable and unstable manifolds of the saddle. 

It is remarkable that the double-zero eigenvalue, which is of purely local
nature, gives rise to a global bifurcation. Bogdanov has proved that there
are no other bifurcations in a \nbh\ of $(\lambda_1,\lambda_2)=0$, and that
the periodic orbit is unique. Thus any small perturbation $g$ of the singular
vector field \eqref{tb8} can be of one of the following structurally stable
types (\figref{fig_takens}):
\begin{enum}
\item	$g$ has no equilibria or periodic orbits in a \nbh\ of the origin;
\item	$g$ admits a saddle and a source;
\item	$g$ admits a saddle, a sink and an unstable hyperbolic periodic
orbit;
\item	$g$ admits a saddle and a sink;
\end{enum}
or of one of the following codimension $1$ singular types:
\begin{enum}
\item[A.]	$g$ has a (\lq\lq semi-unstable\rq\rq ) elementary saddle--node;
\item[B.]	$g$ has a saddle and an elementary composed focus;
\item[C.]	$g$ has a saddle admitting an elementary homoclinic loop;
\item[D.]	$g$ has a (\lq\lq semi-stable\rq\rq ) elementary saddle--node.
\end{enum}
Observe how the codimension $2$ singularity organizes several manifolds of
codimension $1$ singularities in its vicinity. 

Other singular vector fields have been studied. For instance, the system 
\begin{equation}
\label{tb16}
\begin{split}
\dot x_1 &= x_2 \\
\dot x_2 &= -x_1^2x_2 \pm x_1^3
\end{split}
\end{equation}
is not a singularity of codimension $2$ (it is conjectured to have
codimension $4$). However, Takens has shown that this vector field has
codimension $2$ {\em within} the class of vector fields which are invariant
under the transformation $(x_1,x_2)\mapsto(-x_1,-x_2)$. It admits the
unfolding
\begin{equation}
\label{tb17}
\begin{split}
\dot x_1 &= x_2 \\
\dot x_2 &= \lambda_1 + \lambda_2 x_1 -x_1^2x_2 \pm x_1^3
\end{split}
\end{equation}

\begin{exercise}
\label{exo_takens}
Find all equilibria of \eqref{tb17}, determine their stability and the
bifurcation curves. Try to sketch the phase portraits and to locate global
bifurcations (the cases $+$ and $-$ are different). 
\end{exercise}


\section{Remarks on Dimension $3$}
\label{sec_dim3}

We have just seen that the classification of two-dimensional vector fields
is quite well developed. What is the situation in dimension $3$? 

In a sufficiently small \nbh\ of equilibrium points, the situation is very
similar. Hyperbolic equilibrium points are {\em locally} insensitive to
small perturbations. Codimension $1$ bifurcations of equilibria are the
same as in dimension $2$, i.e., saddle--node and Hopf bifurcations. There
is an additional candidate for a bifurcation of codimension $2$,
corresponding to a zero eigenvalue and a pair of imaginary eigenvalues,
which is more difficult to analyze. 

The dynamics near periodic orbits can be described by a two-dimensional
Poincar\'e map. There are two new types of codimension $1$ bifurcations:
The period doubling bifurcation, which produces a new periodic orbit with
twice the period of the old one, and a Hopf bifurcation, which can produce
an invariant torus. Invariant tori constitute a new type of $\w$-limit set. 

Up to this point, the complications due to the additional dimension seem
manageable, and one could hope for a similar classification of
three-dimensional vector fields as of two-dimensional ones. In particular,
one has introduced a class of so-called \defwd{Morse--Smale systems}, which
generalize to higher dimensions the vector fields satisfying Peixoto's
theorem (in particular they have finitely many equilibria and periodic
orbits, all hyperbolic). It was then conjectured that a vector field is
structurally stable if and only if it is Morse--Smale, and that
structurally stable vector fields are dense in the space of all vector
fields.  

These conjectures turned out to be false. In particular, Smale and others
have constructed examples of vector fields that are structurally unstable,
and such that any sufficiently small perturbation of them is still
structurally unstable. Thus structurally {\em stable} vector fields are no
longer dense in the space of all vector fields. 

\begin{figure}
 \centerline{\psfig{figure=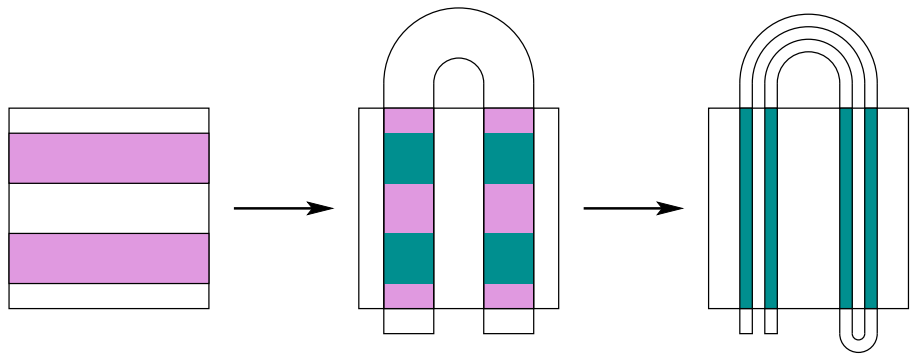,height=50mm,clip=t}}
 \figtext{
 }
 \captionspace
 \caption[]
 {Smale's horseshoe map. Shaded horizontal rectangles are mapped to vertical
 ones. There is a countable infinity of periodic orbits, and an uncountable
 infinity of non-periodic orbits. Arbitrarily close points will become
 separated after a sufficient number of iterations.}
\label{fig_horseshoe}
\end{figure}

On the other hand, Morse--Smale systems are structurally stable but the
converse is false. If $\Sigma$ is a surface transverse to the flow, the
dynamics of orbits returning to $\Sigma$ is described by a two-dimensional
map. Smale has shown that if this map behaves qualitatively like the
\defwd{horseshoe map} of \figref{fig_horseshoe}, then it admits an
invariant Cantor set. This set contains a countable infinity of hyperbolic
periodic orbits and an uncountable infinity of non-periodic orbits, and thus
it is not a Morse--Smale system. Nevertheless, the vector field restricted
to this invariant set is structurally stable. 

The invariant set of the horseshoe map exhibits very complicated dynamics
(which is sensitive to small changes in the initial conditions), but the
set is not attracting. Nevertheless, it has a strong influence on the
dynamics in a \nbh. Later it was shown that there also exist $\w$-limit
sets with such a complicated structure, called \defwd{strange attractors}.
These constitute a new type of $\w$-limit set, and are generally {\em not}
structurally stable. 

\begin{figure}[b]
 \centerline{\psfig{figure=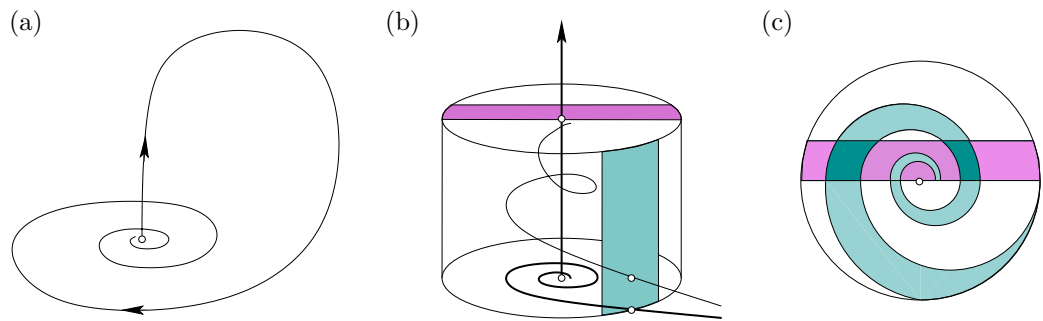,width=140mm,clip=t}}
 \figtext{
 	\writefig	0.5	4.5	(a)
 	\writefig	5.5	4.5	(b)
 	\writefig	10.5	4.5	(c)
 }
 \captionspace
 \caption[]
 {\v Silnikov's example: (a) A saddle--focus with a homoclinic loop. (b)
 Consider a small cylinder whose bottom lies on the stable manifold of the
 saddle--focus. A rectangle on the side of the cylinder is the image of
 some strip on the top of the cylinder (the dimensions of the strip depend
 on the choice of radius and height of the cylinder). (c) The image of the
 rectangle is a spiral, which can be shown, under suitable hypotheses on
 the eigenvalues, to intersect the strip on the top of the cylinder in a
 horseshoe-like way. Thus the Poincar\'e map with respect to the top of the
 cylinder contains a horseshoe.}
\label{fig_Silnikov}
\end{figure}

Codimension $1$ bifurcations of equilibria and periodic orbits do not
produce complicated invariant sets. One could thus wonder where strange
attractors come from, since it should be possible to produce these
complicated systems by continuous deformation of a simple system. It turns
out that horseshoes can appear in bifurcations involving a saddle
connection. Such an example was first given by \v Silnikov: Consider a
saddle--focus, i.e., an equilibrium whose linearization admits one real
eigenvalue $\alpha$ and two complex conjugate eigenvalues
$\beta\pm\icx\gamma$ (where $\alpha$ and $\beta$ have opposite sign). If
this equilibrium admits a saddle connection, \v Silnikov has shown that
under the condition $\abs{\beta}<\abs{\alpha}$, a perturbation of the flow
admits a countable set of horseshoes (\figref{fig_Silnikov}).

The theory developed in the two-dimensional case is still useful in the
vicinity of equilibria and periodic orbits, but the global behaviour of
three-dimensional flows may be substantially more difficult to analyze. We
will encounter some specific examples in the next chapter.  


\section*{Bibliographical comments}

Most of the material discussed here, including the properties of
equilibrium points and periodic orbits, limit sets, local and global
bifurcations are treated in standard monographs such as \cite{GH} and
\cite{Wiggins}. The classification of codimension $1$ bifurcations is
discussed more specifically in \cite{HK}, though with only few proofs.

The proof of Peixoto's theorem is quite nicely exposed in \cite{PP}. 
Further developments can be found in \cite{Pe62} (two-dimensional compact
manifolds) and \cite{Pe73} (classification of flows by graphs). A more
advanced description of results related to structural stability and
genericity is found in \cite{Palis}. 

The Bogdanov--Takens bifurcation is originally described in \cite{Takens}
and \cite{Bogdanov}. The theory of unfoldings (or versal deformations) and
relations to singularity theory are discussed in \cite{Arnold}. There exist
many partial results on local bifurcations of higher codimension, see for
instance \cite{KKR}.  

The fact that structurally stable vector fields are not dense was first
shown in \cite{Smale} and \cite{PP68}.  The example of \v Silnikov
originally appeared in \cite{Silnikov}, and is also discussed in \cite{GH}
and \cite{Wiggins}.


\chapter{Regular Perturbation Theory}
\label{ch_rpt}

In this chapter, we will consider non--autonomous ordinary differential
equations of the form 
\begin{equation}
\label{rpt1}
\dot x = f(x) + \eps g(x,t,\eps),
\qquad\qquad x\in\cD\subset\R^n,
\end{equation}
where $f$ and $g$ are of class $\cC^r$, $r\geqs1$, $g$ is $T$-periodic in
$t$, and $\eps$ is a small parameter. We will be mainly interested in the
following question: If $x^\star$ is a stable equilibrium point of the
unperturbed autonomous system $\dot x=f(x)$, does the perturbed system
\eqref{rpt1} admit a stable $T$-periodic orbit in a \nbh\ of $x^\star$? We
will see below that this question frequently arises in applications. 

Let us first recall the definition of stability:

\begin{definition}
\label{def_rpt}
A periodic solution $\gamma(t)$ of \eqref{rpt1} is called \defwd{stable} if
for any $\delta>0$, there exists $\eta>0$ such that, if
$\norm{x(0)-\gamma(0)}\leqs\eta$, then $\norm{x(t)-\gamma(t)}\leqs\delta$
for all $t\in\R$. 
\end{definition} 

It is easy to show that \eqref{rpt1} admits a stable periodic solution if
$x^\star$ is a sink:

\begin{prop}
\label{prop_rpt}
Assume that all eigenvalues of $A=\dpar fx(x^\star)$ have a strictly
negative real part. Then there exist constants $\eps_0,c_0>0$ such that,
whenever $\abs\eps\leqs\eps_0$, \eqref{rpt1} admits a unique stable
$T$-periodic solution $\gamma(t)$ such that $\norm{\gamma(t)-x^\star}\leqs
c_0\eps$ for all $t\in\R$. 
\end{prop}
\begin{proof}
Let $\ph^\eps_t$ denote the flow of \eqref{rpt1}, i.e., the solutions are
written as $x(t)=\ph^\eps_t(x(0))$. We are interested in the application
$\ph^\eps_T$, which can be interpreted as the Poincar\'e map of
\eqref{rpt1} in the extended phase space $\set{(x,t)\in\cD\times\R/T\Z}$,
taking the plane $t=0$ as a Poincar\'e section. Indeed, since the equation
is $T$-periodic in $t$, we can identify the planes $t=0$ and $t=T$.

If $\eps=0$, we simply have $\dpar{}x\ph^0_t(x^\star)=\e^{At}$. Thus,
if we introduce the function $\Phi(x,\eps)=\ph^\eps_T(x)-x$, we have the
relations
\[
\Phi(x^\star,0)=0,
\qquad\qquad
\dpar\Phi x(x^\star,0)=\e^{AT}-\one.
\]
Since all eigenvalues of $A$ have a negative real part, $\e^{AT}$ has all
its eigenvalues inside the unit circle, which means in particular that
$\dpar{}x\Phi (x^\star,0)$ has no eigenvalue equal to zero. Thus the implicit
function theorem implies the existence, for $\eps$ sufficiently small, of a
unique $\cC^r$ function $x^\star(\eps)$ such that
$\ph^\eps_T(x^\star)=x^\star$. This fixed point of the Poincar\'e map
corresponds to a periodic orbit of \eqref{rpt1}. 

For sufficiently small $\eps$, the eigenvalues of
$\dpar{}x\ph^\eps_T(x^\star(\eps))$ still lie inside the unit circle, which
means that $\ph^\eps_T$ is a contraction in some \nbh\ of $x^\star(\eps)$.
This implies the uniqueness of the periodic orbit. It also shows that
$\norm{x(t)-\gamma(t)}$ is smaller than $\norm{x(0)-\gamma(0)}$ for
$t/T\in\N$. One can show that $\norm{x(t)-\gamma(t)}$ is smaller than a
constant times $\norm{x(0)-\gamma(0)}$ for {\em all} $t$ by a standard
argument on differential inequalities, see Corollary~\ref{cor_rpb} below.
This implies the stability of $\gamma$.
\end{proof}

A similar argument shows that if $x^\star$ is a hyperbolic equilibrium point
of $f$, then \eqref{rpt1} admits a periodic orbit for sufficiently small
$\eps$, which is stable if and only if $x^\star$ is stable. 

As we saw in Example~\ref{ex_oscillator}, one often has to do with cases
where $x^\star$ is elliptic. If the linearization $A$ of $f$ at $x^\star$
has imaginary eigenvalues, but none of them is a multiple of $2\pi\icx/T$,
the implicit function theorem can still be applied to prove the existence
of a unique periodic orbit near $x^\star$. Analysing the stability of this
orbit, however, turns out to be considerably more difficult, and will be
one of the main goals of this chapter. Along the way, we will develop some
methods which are useful in other problems as well. 

Before discussing these methods, we will recall some properties of
Hamiltonian systems, which form an important source of equations of the form
\eqref{rpt1}, and provided the original motivation for their study. 


\section{Preliminaries}
\label{sec_rpp}


\subsection{Hamiltonian Systems}
\label{ssec_rph}

Let $H\in\cC^{r+1}(\cD,\R)$, $r\geqs 1$, be defined in an open subset $\cD$
of $\R^{2n}$ (or on a $2n$-dimensional differentiable manifold). The integer
$n$ is called the \defwd{number of degrees of freedom}, and $H$ is the
\defwd{Hamilton function} or \defwd{Hamiltonian}. We denote elements of
$\cD$ by $(q,p)=(q_1,\dots,q_n,p_1,\dots,p_n)$. The \defwd{canonical
equations} associated with $H$ are defined as 
\begin{equation}
\label{rph1}
\dot q_i = \dpar H{p_i}, \qquad\qquad
\dot p_i = -\dpar H{q_i}, \qquad\qquad
i=1,\dots,n,
\end{equation}
and define a \defwd{Hamiltonian flow} on $\cD$. 
The simplest example of a Hamiltonian is 
\begin{equation}
\label{rph2}
H = \frac12 \norm{p}^2 + V(q),
\end{equation}
whose canonical equations
\begin{equation}
\label{rph3}
\dot q_i = p_i, \qquad\qquad
\dot p_i = -\dpar V{q_i}(q)
\end{equation}
are those of a particle in a potential $V$. The advantage of the canonical
equations \eqref{rph1} is that they are invariant under a large class of
coordinate transformations, and thus adapted to mechanical systems with
constraints. 

The \lq\lq skew-symmetric\rq\rq\ structure of the canonical equations
\eqref{rph1}, called \defwd{symplectic structure}, has some important
consequences.

\begin{prop}\hfill
\label{prop_rph1}
\begin{itemiz}
\item	$H$ is a constant of the motion.
\item	The Hamiltonian flow preserves the volume element
$\6q\6p=\6q_1\dots\6q_n\6p_1\dots\6p_n$. 
\item	A curve $\gamma=\setsuch{q(t),p(t)}{t_1\leqs t\leqs t_2}$
corresponds to a solution of \eqref{rph1} if and only if the \defwd{action
integral}
\begin{equation}
\label{rph4}
\cS(\gamma) = 
\int_\gamma p\cdot\6q - H\6t \defby \int_{t_1}^{t_2} \bigbrak{p(t)\cdot
\dot q(t) - H(q(t),p(t))} \6t 
\end{equation}
is stationary with respect to variations of $\gamma$ (with fixed end
points). 
\end{itemiz}
\end{prop}
\goodbreak
\begin{proof}\hfill
\begin{enum}
\item	By a direct calculation,
\[
\dtot{}t H(q(t),p(t)) = \sum_{i=1}^n \Bigpar{\dpar H{q_i}\dot q_i + \dpar
H{p_i}\dot p_i} = 0.
\]

\item	In the extended phase space $\cD\times\R$, the canonical equations
define the vector field $f(q,p,t)=(\tdpar Hp,-\tdpar Hq,1)$. Let $M(0)$ be
a compact set in $\cD$, and $M(t)$ the image of $M(0)$ under the flow
$\ph_t$. Define a cylinder $\cC=\setsuch{(M(s),s)}{0\leqs s\leqs t}$. Then
we have, by Gauss' theorem,
\[
\int_{M(t)}\6q\6p - \int_{M(0)}\6q\6p = 
\int_{\partial\cC} f\cdot \6n = \int_\cC (\nabla\cdot f)\6p\6q\6t = 0,
\] 
where $\6n$ is a normal vector to $\partial\cC$, and we have used the fact
that $f\cdot \6n=0$ on the sides of the cylinder. 

\item	Let $\delta=\setsuch{\x(t),\y(t)}{t_1\leqs t\leqs t_2}$ be a curve
vanishing for $t=t_1$ and $t=t_2$. The first variation of the action in the
direction $\delta$ is defined as 
\[
\begin{split}
\dtot{}\eps \cS(\gamma+\eps\delta) \Bigevalat{\eps=0} 
&= \int_{t_1}^{t_2} \Bigpar{p\dot\x + \dot q\y - \dpar Hq\x - \dpar Hp\y}\6t
\\
&= \int_{t_1}^{t_2} \x\Bigpar{-\dot p - \dpar Hq} + \y\Bigpar{\dot q-\dpar
Hp} \6t,
\end{split}
\]
where we have integrated the term $p\dot\x$ by parts. Now it is clear that
this expression vanishes for all perturbations $\delta$ if and only if
$(q(t),p(t))$ satisfies the canonical equations \eqref{rph1}. 
\qed
\end{enum}
\renewcommand{\qed}{}
\end{proof}

A useful notation is the \defwd{Poisson bracket} of two differentiable
functions $f,g:\cD\to\R$, defined as 
\begin{equation}
\label{rph5}
\poisson fg = \sum_{i=1}^n \dpar f{q_i}\dpar g{p_i} - \dpar f{p_i}\dpar
g{q_i}.
\end{equation}
Then we have, in particular, 
\begin{equation}
\label{rph6}
\dtot ft = \sum_{i=1}^n \dpar f{q_i}\dot q_i + \dpar f{p_i} \dot p_i =
\poisson fH.
\end{equation}
Since $\poisson{}{}$ is an antisymmetric bilinear form, we immediately
recover that $\dot H = \poisson HH=0$. A function $J:\cD\to\R$ is a
\defwd{constant of motion} if and only if $\poisson JH=0$. One can show that
if $J_1$ and $J_2$ are two constants of motion, then $\poisson{J_1}{J_2}$ is
also a constant of motion. $J_1$ and $J_2$ are \defwd{in involution} if
$\poisson{J_1}{J_2}=0$. 

For one degree of freedom ($n=1$), orbits of the Hamiltonian flow are simply
level curves of $H$. With more degrees of freedom, the manifolds
$H=\text{\it constant}$ are invariant, but the motion within each manifold
may be complicated. If there are additional constants of motion, however,
there will be invariant manifolds of smaller dimension, in which the motion
may be easier to analyse. To do this, we would like to choose coordinates
which exploit these constants of motion, but keeping the symplectic
structure of the differential equation. Transformations which achieve this
are called \defwd{canonical transformations}. 

Let $S(q,Q)$ be twice continuously differentiable in an open domain of
$\R^{2n}$, such that
\begin{equation}
\label{rph7}
\det \dpar{^2S}{q\partial Q} \neq 0
\qquad\qquad\forall q,Q. 
\end{equation}
$S$ is called a \defwd{generating function}. We introduce
\begin{equation}
\label{rph8}
\begin{split}
p(q,Q) &= \dpar Sq(q,Q), \\
P(q,Q) &= -\dpar SQ(q,Q). 
\end{split}
\end{equation}
By the assumption \eqref{rph7}, we can invert the relation $p=p(q,Q)$ with
respect to $Q$, and thus the relations \eqref{rph8} define a transformation
$(q,p)\mapsto(Q,P)$ in some open domain in $\R^{2n}$. 

\begin{prop}
\label{prop_rph2}
The transformation $(q,p)\mapsto(Q,P)$ is volume-preserving and canonical.
The new equations of motion are given by 
\begin{equation}
\label{rph9}
\dot Q_i = \dpar K{P_i}, \qquad\qquad
\dot P_i = -\dpar K{Q_i}, \qquad\qquad
i=1,\dots n,
\end{equation}
where $K(Q,P)=H(q,p)$ is the old Hamiltonian expressed in the new variables.
\end{prop}
\begin{proof}
Note that $\6S=\sum_i p_i\6q_i - P_i\6Q_i$. For every $i$ and a given
compact $M\subset\R^2$, we have by Green's theorem
\[
\int_M \6q_i\6p_i - \6Q_i\6P_i = \int_{\partial M} p_i\6q_i - P_i\6Q_i =
\int_{\partial M} \dpar S{q_i}\6q_i + \dpar S{Q_i}\6Q_i = \int_{\partial
M}\6S, 
\]
but this last integral vanishes since $\partial M$ consists of one or
several closed curves. This proves that $\6q_i\6p_i=\6Q_i\6P_i$ for every
$i$. Furthermore, for any solution of the canonical equations, we have 
\[
\int_{t_1}^{t_2} P\cdot\6Q-H\6t = \int_{t_1}^{t_2} (p\cdot\6q - H\6t) -
\int_{t_1}^{t_2}\6S.
\]
The first integral on the right-hand side is stationary, and the second
depends only on the end points. Thus the integral of $P\6Q-H\6t$ is
stationary with respect to variations of the curve, and by
Proposition~\ref{prop_rph1}, $(Q,P)$ satisfies the canonical equations
\eqref{rph9}. 
\end{proof}

Consider for the moment the case $n=1$. As we have seen, the orbits are
level curves of $H$. It would thus be useful to construct a canonical
transformation such that one of the new variables is constant on each
invariant curve. Such variables are called \defwd{action--angle
variables}. 

Assume that the level curve $H^{-1}(h)=\setsuch{(q,p)}{H(q,p)=h}$ is
bounded. Then its \defwd{action} is defined by 
\begin{equation}
\label{rph10}
I(h) = \frac1{2\pi} \int_{H^{-1}(h)} p\6q.
\end{equation}
Assume further that $I'(h)\neq 0$ on some interval, so that the map
$I=I(h)$ can be inverted to give $h=h(I)$. If $p(h,q)$ is the solution of
$H(q,p)=h$ with respect to $p$, we can define a generating function 
\begin{equation}
\label{rph11}
S(I,q) = \int_{q_0}^q p(h(I),q')\6q'.
\end{equation}
(The function $p(h,q)$ is not uniquely defined in general, but the
definition \eqref{rph11} makes sense for any parametrization of the curve
$H(q,p)=h$). The transformation from $(q,p)$ to action--angle variables
$(I,\ph)$ is defined implicitly by 
\begin{equation}
\label{rph12}
p(I,q) = \dpar Sq, \qquad\qquad
\ph(I,q) = \dpar SI.
\end{equation} 
Note that the normalization in \eqref{rph10} has been chosen in such a way
that $\ph\to2\pi-$ as $q\to q_0-$, and thus $\ph$ is indeed an angle. The
Hamiltonian in action--angle variables takes the form
\begin{equation}
\label{rph13}
K(\ph,I) = h(I),
\end{equation} 
and the associated canonical equations are simply 
\begin{equation}
\label{rph14}
\dot\ph = \dpar KI = h'(I), \qquad\qquad
\dot I = -\dpar K\ph = 0. 
\end{equation}
The first relation means that the angle $\ph$ is an \lq\lq equal-time\rq\rq\
parametrization of the level curve. Indeed, we have 
\begin{equation}
\label{rph15}
\ph(t) = \int_{q_0}^{q(t)} \dpar pH(h(I),q')\6q' \dtot hI = 
\int_{q_0}^{q(t)} \frac{\6q'}{\dot q}\dtot hI = (t-t_0)\dtot hI. 
\end{equation}

In the case $n>1$, there is no reason for there to be other constants of
motion than $H$ (and functions of $H$). There are, however, a certain number
of cases where additional constants of motion exist, for instance when $H$
is invariant under a (continuous) symmetry group (Noether's theorem). A
Hamiltonian with $n$ degrees of freedom is called \defwd{integrable} if it
has $n$ constants of motion $J_1,\dots,J_n$ such that 
\begin{itemiz}
\item	the $J_i$ are in involution, i.e., $\poisson{J_i}{J_k}=0$ for all
$(i,k)$;
\item	the gradients of the $J_i$ are everywhere linearly independent in
$\R^{2n}$.
\end{itemiz}
Arnol'd has shown that if a given level set of $J_1,\dots,J_n$ is compact,
then it typically has the topology of an $n$-dimensional torus. One can then
introduce action--angle variables in a similar way as for $n=1$, and the
Hamiltonian takes the form 
\begin{equation}
\label{rph16}
K(\ph,I) = H_0(I). 
\end{equation}
Consider now a perturbation of the integrable Hamiltonian, of the form
\begin{equation}
\label{rph17}
H(\ph,I,\eps) = H_0(I) + \eps H_1(I,\ph,\eps). 
\end{equation}
The associated canonical equations are
\begin{equation}
\label{rph18}
\dot \ph_j = \dpar{H_0}{I_j} + \eps \dpar{H_1}{I_j}, \qquad\qquad
\dot I_j = - \eps \dpar{H_1}{\ph_j}.
\end{equation}
We can now make the link with the system $\dot x=f(x)+\eps g(x,t,\eps)$
introduced at the beginning of this chapter. A first possibility is to take
$x=(\ph,I)$, $f=(\tdpar{H_0}I,0)$ and $g=(\tdpar{H_1}I,-\tdpar{H_1}\ph)$.
Then $g$ is not time-dependent. There are, however, other possibilities. 

\begin{example}
\label{ex_rph1}
Assume for instance that $\tdpar{H_0}{I_1}\neq 0$ for all $I$. If
$\tdpar{H_1}{I_1}$ is bounded, then $\ph_1$ will be a monotonous function of
$t$ for sufficiently small $\eps$. We can thus replace $t$ by $\ph_1$ and
consider the $2n-1$ equations 
\begin{align}
\nonumber
\dtot{\ph_k}{\ph_1} &= \frac{\tdpar{H_0}{I_k} + \eps \tdpar{H_1}{I_k}}
{\tdpar{H_0}{I_1} + \eps \tdpar{H_1}{I_1}}, 
&& k=2,\dots,n, \\
\dtot{I_j}{\ph_1} &= -\eps \frac{\tdpar{H_1}{\ph_j}}
{\tdpar{H_0}{I_1} + \eps \tdpar{H_1}{I_1}}, 
&& j=1,\dots,n.
\label{rph19}
\end{align}
We thus obtain a system of the form $\dot x=f(x)+\eps g(x,t,\eps)$ for
$x=(\ph_2,\dots,\ph_n,I_1,\dots,I_n)$, and with
$f=((\tdpar{H_0}I)/(\tdpar{H_0}{I_1}),0)$. The remainder $g$
depends periodically on $\ph_1$ because $\ph_1$ is an angle. 
\end{example}

\begin{example}
\label{ex_rph2}
Assume that $\gamma(t)=(q(t),p(t))$ is a $T$-periodic solution of a
Hamiltonian system (which need not be close to integrable). The variable
$y=x-\gamma(t)$ satisfies an equation of the form
\begin{equation}
\label{rph20}
\dot y = A(t)y + g(y,t),
\end{equation}
where 
\begin{equation}
\label{rph21}
A(t) = 
\begin{pmatrix}
\tdpar{^2H}{q\partial p} & \tdpar{^2H}{p^2} \\
-\tdpar{^2H}{q^2} & -\tdpar{^2H}{q\partial p}
\end{pmatrix}\biggevalat{\gamma(t)}
\end{equation}
and there are constants $M,\delta>0$ such that $\norm{g(y,t)}\leqs
M\norm{y}^2$ for $\norm{y}\leqs\delta$. The equation \eqref{rph20} is not
yet in the form we are looking for, but we are going to transform it now. 
Let $U(t)$ be the solution of 
\begin{equation}
\label{rph22}
\dot U = A(t) U, \qquad\qquad U(0)=\one.
\end{equation}
Since $A$ is $T$-periodic,  we know by Floquet's theorem that
$U(t)=P(t)\e^{tB}$ for some $T$-periodic matrix $P(t)$. By substitution
into \eqref{rph22}, we see that 
\begin{equation}
\label{rph23}
\dot P = A(t)P - PB,
\end{equation}
and thus the change of variables $y=P(t)z$ yields the equation 
\begin{equation}
\label{rph24}
\dot z = Bz + P^{-1} g(Pz,t). 
\end{equation}
Finally we carry out a rescaling $z=\eps w$, which means that we zoom on a
small \nbh\ of the periodic orbit. Here $\eps$ is a small positive
parameter, and the change of variables is defined for all $\eps>0$. $w$
obeys the equation 
\begin{equation}
\label{rph25}
\dot w = B w + \eps G(w,t,\eps),
\end{equation}
which is exactly of the form \eqref{rpt1}, with
\begin{equation}
\label{rph26}
\begin{split}
G(w,t,\eps) &= \frac1{\eps^2} P^{-1} g(\eps Pw,t), \\
\norm{G(w,t,\eps)} &\leqs \frac1{\eps^2} \norm{P^{-1}} M \norm{\eps Pw}^2 
\leqs M \norm{P^{-1}} \norm{P} \norm{w}^2.
\end{split}
\end{equation}
In fact, \eqref{rph20} and \eqref{rph24} can both be considered as small
perturbations of a linear system, deriving from a (possibly time-dependent)
quadratic Hamiltonian. Since a quadratic form can always be diagonalized,
these linear equations are integrable, and thus the motion in a small \nbh\
of a periodic orbit can be viewed as a perturbed integrable system. 
\end{example}


\subsection{Basic Estimates}
\label{ssec_rpb}

Let us return to the equation 
\begin{equation}
\label{rpb1}
\dot x = f(x) + \eps g(x,t,\eps),
\qquad\qquad x\in\cD\subset\R^n,
\end{equation}
where $f$ and $g$ are of class $\cC^r$ for some $r\geqs1$.  Our aim is to
describe the difference between the dynamics for $\eps=0$ and $\eps>0$. Let
us recall a fundamental theorem on ordinary differential equations, and the
dependence of their solutions on parameters. 

\begin{theorem}
\label{thm_rpb1}
Let $F(x,t,\lambda)\in\cC^r(\cD_0,\R^n)$, $r\geqs 1$, for an open set
$\cD_0\subset\R^n\times\R\times\R^p$. Then the solution of $\dot{x} =
F(x,t,\lambda)$ with initial condition $x(t_0) = x_0$ is a $\cC^r$
function of $x_0, t_0, t$ and $\lambda$ on its domain of existence.
\end{theorem}	

It follows that $x(t)$ admits an expansion of the form 
\begin{equation}
\label{rpb2}
x(t) = x^0(t) + \eps x^1(t) + \eps^2 x^2(t) + \dots + \eps^r x^r(t,\eps)
\end{equation}
in some intervals of $t$ and $\eps$. In particular, $x^0(t)$ is a solution
of the unperturbed equation $\dot x = f(x)$, and thus the solution of
\eqref{rpb1} stays close to $x^0(t)$ on compact time intervals. 

A useful tool to estimate the time-dependence of the error $x(t)-x^0(t)$ is
the following lemma from the theory of differential inequalities.

\begin{lemma}[Gronwall's inequality]
\label{lem_Gronwall}
Let $\ph$, $\alpha$ and $\beta$ be continuous real-valued functions on
$[a,b]$, with $\beta$ non-negative, and assume that
\begin{equation}
\label{rpb3}
\ph(t) \leqs \alpha(t) + \int_a^t \beta(s)\ph(s) \6s 
\qquad \forall t\in[a,b].
\end{equation}
Then
\begin{equation}
\label{rpb4}
\ph(t) \leqs \alpha(t) + \int_a^t \beta(s)\alpha(s) 
\e^{\int_s^t\beta(u)\6u}\6s 
\qquad \forall t\in[a,b].
\end{equation}
\end{lemma} 
\begin{proof}
Let 
\[
R(t) = \int_a^t \beta(s)\ph(s) \6s.
\]
Then $\ph(t)\leqs\alpha(t)+R(t)$ for all $t\in[a,b]$ and, since
$\beta(t)\geqs0$,  
\[
\dtot Rt(t) = \beta(t)\ph(t) \leqs \beta(t)\alpha(t) + \beta(t)R(t).
\]
Let $B(s)=\int_a^s\beta(u)\6u$. Then 
\[
\dtot {}{s} \e^{-B(s)}R(s) 
= \bigbrak{R'(s)-\beta(s)R(s)}\e^{-B(s)}
\leqs \beta(s)\alpha(s)\e^{-B(s)}
\]
and thus, integrating from $a$ to $t$,
\[
\e^{-B(t)}R(t) \leqs \int_a^t \beta(s)\alpha(s)\e^{-B(s)}\6s.
\]
We obtain the conclusion by multiplying this expression by $\e^{B(t)}$ and
inserting the result into \eqref{rpb3}.
\end{proof}

\begin{cor}
\label{cor_rpb}
Assume $f$ admits a uniform Lipschitz constant $K$ in $\cD$, and $\norm{g}$
is uniformly bounded by $M$ in $\cD\times\R\times[0,\eps_0]$. If
$x(0)=x^0(0)$, then 
\begin{equation}
\label{rpb5}
\norm{x(t)-x^0(t)} \leqs \frac{\eps M}K \bigbrak{\e^{Kt}-1} 
\qquad\qquad \forall t\geqs 0.
\end{equation}
\end{cor}
\begin{proof}
Let $y=x - x^0(t)$. Then 
\[
\dot y = f(x^0(t) + y) - f(x^0(t)) + \eps g(x(t),t,\eps).
\]
Integrating between $0$ and $t$, we obtain   
\[
y(t) = \int_0^t \bigbrak{f(x^0(s) + y(s)) - f(x^0(s))} \6s 
+ \eps \int_0^t g(x(s),s,\eps) \6s. 
\]
Taking the norm, we arrive at the estimate
\[
\begin{split}
\norm{y(t)} &\leqs 
\int_0^t \norm{f(x^0(s) + y(s)) - f(x^0(s))} \6s
+ \eps \int_0^t \norm{g(x(s),s,\eps)} \6s \\
&\leqs \int_0^t K \norm{y(s)} \6s + \eps Mt. 
\end{split}
\]
Applying Gronwall's inequality, we obtain \eqref{rpb5}, as desired. 
\end{proof}

The estimate \eqref{rpb5} shows that for a given $t$, we can make the
deviation $\norm{x(t)-x^0(t)}$ small by taking $\eps$ small. However, since
$\e^{Kt}$ grows without bound as $t\to\infty$, we cannot make any statements
on stability. Of course, \eqref{rpb5} is only an upper bound. The following
example shows, however, that it cannot be improved without further
assumptions on $f$. 

\begin{example}
\label{ex_rpb}
Consider, as a special case of \eqref{rph25}, the equation 
\begin{equation}
\label{rpb6}
\dot x = B x + \eps g(t),
\end{equation}
where $g$ is $T$-periodic in $t$. Let us assume that $B$ is in diagonal
form, with (possibly complex) eigenvalues $b_1,\dots,b_n$. Then the
solution can be written, for each component, as 
\begin{equation}
\label{rpb7}
x_j(t) = \e^{b_jt} x_j(0) + \eps \int_0^t \e^{b_j(t-s)} g_j(s) \6s. 
\end{equation}
If the $g_j$ are written in Fourier series as 
\begin{equation}
\label{rpb8}
g_j(s) = \sum_{k\in\Z} c_{jk}\e^{\icx\w ks}, 
\qquad\qquad \w = \frac{2\pi}T, 
\end{equation}
we have to compute integrals of the form 
\begin{equation}
\label{rpb9}
\int_0^t \e^{b_j(t-s)} \e^{\icx\w ks} \6s = 
\begin{cases}
\displaystyle\frac{\e^{\icx\w kt} - \e^{b_jt}}{\icx\w k-b_j}
& \text{if $b_j \neq \icx\w k$,} \\
\vrule height 16pt depth 7pt width 0pt
t \e^{b_jt}
& \text{if $b_j = \icx\w k$.} 
\end{cases}
\end{equation}
If $\re b_j\neq 0$, solutions grow or decrease exponentially fast, just like
for $\eps=0$. If $\re b_j=0$, however, the origin is a stable equilibrium
point of the unperturbed system, while the perturbed system displays the
phenomenon of \defwd{resonance}: If $b_j$ is a multiple of $\icx\w$, the
deviation grows linearly with $t$; if not, it remains bounded, but the
amplitude may become large if the denominator $\icx\w k-b_j$ is small for
some $k$. 
\end{example}


\section{Averaging and Iterative Methods}
\label{sec_avi}

Corollary \ref{cor_rpb} gives a rather rough estimate on the difference
between solutions of the perturbed and unperturbed equations. There exist a
number of iterative methods that allow to improve these estimates, by
replacing the unperturbed equation by a better approximation of the
perturbed one, which still remains easier to solve. 


\subsection{Averaging}
\label{ssec_av}

The method of averaging allows to deal with differential equations
depending periodically on time, by comparing their solutions with a
simpler, autonomous equation. The method applies to systems of the form 
\begin{equation}
\label{av1}
\dot x = \eps g(x,t,\eps),
\end{equation}
where $g$ is $T$-periodic in $t$. 

We are, in fact, interested in more general equations of the form 
\begin{equation}
\label{av2}
\dot x = f(x) + \eps g(x,t,\eps),
\end{equation}
but it is often possible to transform this equation into \eqref{av1}.
Indeed, let $h(x,t,\eps)\in\cC^r(\cD\times\R/T\Z\times[0,\eps_0],\R^n)$ be
$T$-periodic in $t$, and invertible in $x$ near $x=0$. Then the variable
$y=h(x,t,\eps)$ satisfies 
\begin{equation}
\label{av3}
\dot y = \dpar ht + \dpar hx f(x) + \eps \dpar hx g(x,t,\eps). 
\end{equation}
Thus, if we manage to find an $h$ such that 
\begin{equation}
\label{av4}
\dpar ht + \dpar hx f(x) = \Order{\eps},
\end{equation}
then $y$ satisfies an equation of the form \eqref{av1}. This is, in fact,
possible if the unperturbed equation has periodic solutions with a period
close to $T$, i.e., near resonance. 

\begin{example}
\label{ex_av1}
Consider the equation 
\begin{equation}
\label{av5}
\dot x = B(\w_0) x + \eps g(x,t,\eps), 
\qquad\qquad
B(\w_0) = 
\begin{pmatrix}
0 & 1 \\ -\w_0^2 & 0
\end{pmatrix}.
\end{equation}
The solutions of the unperturbed equation are of the form 
\begin{equation}
\label{av6}
x^0(t) = \e^{B(\w_0)t} x^0(0), 
\qquad\qquad
\e^{B(\w_0)t} = 
\begin{pmatrix}
\cos(\w_0t) & \frac1{\w_0}{\sin(\w_0t)} \\ -\w_0\sin(\w_0t) & \cos(\w_0t)
\end{pmatrix},
\end{equation}
and are thus periodic with period $2\pi/\w_0$. Now let $\w=2\pi/T$ and take
$h(x,t)=\e^{-B(\w)t}x$. We obtain that 
\begin{equation}
\label{av7}
\dpar ht + \dpar hx f(x) =\e^{-B(\w)t} \bigbrak{B(\w_0)-B(\w)} x.
\end{equation}
Thus if $\w_0$ is close to $\w$, that is, if $\w^2=\w_0^2 + \Order{\eps}$,
then \eqref{av7} is indeed of order $\eps$. More generally, taking
$h(x,t)=\e^{-B(\w/k)t}x$ for some $k\in\Z$ allows to treat cases where
$\w^2=k^2\w_0^2+\Order{\eps}$. This trick is known as the \defwd{van der
Pol transformation}. 
\end{example}

Let us now return to the equation 
\begin{equation}
\label{av8}
\dot x = \eps g(x,t,\eps).
\end{equation}
We assume that $g\in\cC^r(\cD\times\R/T\Z\times[0,\eps_0],\R^n)$, where
$\cD$ is a compact subset of $\R^n$ and $r\geqs2$. Note that if $M$ is an
upper bound for $\norm{g}$ in $\cD$, then $\norm{x(t)-x(0)}\leqs \eps Mt$.
Thus the trivial approximation $\dot x=0$ of \eqref{av8} is good to order
$\eps$ for times of order $1$ (this is, in fact, a particular case of
Corollary \ref{cor_rpb} in the limit $K\to0$). A better approximation is
given by the averaged equation:

\begin{definition}
\label{def_av}
The \defwd{averaged system} associated with \eqref{av8} is the autonomous
equation 
\begin{equation}
\label{av9}
\dot y^0 = \eps \avrg{g}(y^0), 
\qquad\qquad
\avrg{g}(x) = \frac1T \int_0^T g(x,t,0)\6t. 
\end{equation}
\end{definition}

The averaging theorem shows that \eqref{av9} is a good approximation of
\eqref{av8} up to times of order $1/\eps$:

\begin{theorem}
\label{thm_avrg}
There exists a $\cC^r$ change of coordinates $x=y+\eps w(y,t)$, with $w$
a $T$-periodic function of $t$, transforming \eqref{av8} into
\begin{equation}
\label{av10}
\dot y = \eps \avrg{g}(y) + \eps^2 g_1(y,t,\eps),
\end{equation}
with $g_1$ also $T$-periodic in $t$. Moreover, 
\begin{itemiz}
\item	if $x(t)$ and $y^0(t)$ are solutions of the original and averaged
systems \eqref{av8} and \eqref{av9} respectively, with initial conditions
$x(0)=y^0(0)+\Order{\eps}$, then $x(t)=y^0(t)+\Order{\eps}$ for times $t$ of
order $1/\eps$;
\item	if $y^\star$ is an equilibrium point of \eqref{av9} such that
$\dpar{}y \avrg{g}(y^\star)$ has no eigenvalue equal to zero, then
\eqref{av8} admits a $T$-periodic solution $\gamma_\eps(t) =
y^\star+\Order{\eps}$ for sufficiently small $\eps$; if, in addition,
$\dpar{}y \avrg{g}(y^\star)$ has no imaginary eigenvalue, then
$\gamma_\eps(t)$ has the same type of stability as $y^\star$ for
sufficiently small $\eps$.  
\end{itemiz}
\end{theorem}

\goodbreak
\begin{proof}\hfill
\begin{itemiz}
\item	We begin by constructing the change of variables explicitly. We
first split $g$ into its averaged and oscillating part:
\[
g(x,t,\eps) = \avrg{g}(x) + g^0(x,t,\eps), 
\qquad\qquad
\avrg{g^0}(x)=0.
\]
Inserting $x=y+\eps w(y,t)$ into \eqref{av8}, we get 
\[
\dot x = \Bigbrak{\one+\eps\dpar wy} \dot y + \eps \dpar wt(y,t) = 
\eps\avrg{g}(y+\eps w) + \eps g^0(y+\eps w,t,\eps),
\]
and expanding $\dot y$ into Taylor series, we obtain 
\[
\dot y = \eps\bigbrak{\avrg{g}(y) + g^0(y,t,0) - \dpar wt(y,t)} +
\Order{\eps^2},
\]
where the remainder of order $\eps^2$ is a $T$-periodic function of $t$, the
Taylor expansion of which can be computed if necessary. Now choosing 
\[
w(y,t) = \int_0^t g^0(y,s,0) \6s
\]
yields the desired form for $\dot y$. Note that $w$ is indeed $T$-periodic
in $t$ since $g^0$ has average zero. 

\item	The difference $z(t)=y(t)-y^0(t)$ satisfies the equation 
\[
\dot z = \eps \bigbrak{\avrg{g}(y) - \avrg{g}(y^0)} + \eps^2 g_1(y,t,\eps), 
\]
and thus we can proceed as in the proof of Corollary~\ref{cor_rpb}. Since 
\[
z(t) = z(0) + \int_0^t \Bigset{\eps \bigbrak{\avrg{g}(y(s)) -
\avrg{g}(y^0(s))}  + \eps^2 g_1(y(s),s,\eps)} \6s, 
\]
taking $K$ as a Lipschitz constant for $\avrg{g}$ and $M$ as a bound for
$\norm{g_1}$ (which must exist since $\cD$ is compact), we arrive at the
estimate
\[
\norm{z(t)} \leqs \norm{z(0)} + \int_0^t \eps K \norm{z(s)}\6s + \eps^2Mt, 
\]
and Gronwall's inequality provides us with the bound 
\[
\norm{z(t)} \leqs \norm{z(0)} \e^{\eps Kt} +  
\frac{\eps M}K \bigpar{\e^{\eps Kt}-1}. 
\]
Since we have assumed that $\norm{z(0)}=\Order{\eps}$, the error made by the
averaging approximation remains of order $\eps$ for $t\leqs 1/(K\eps)$. 

\item	The assertion on the periodic orbit is proved in a similar way as
we did in Proposition~\ref{prop_rpt}. Let $\ph_t$ denote the flow of
the transformed equation \eqref{av10}, and let $\ph^0_t$ denote the flow of
the averaged equation \eqref{av9} ($\ph^0$ is obtained from $\ph$ by
letting $g_1=0$, {\em not} by letting $\eps=0$). By the fundamental
Theorem~\ref{thm_rpb1}, these flows may be Taylor-expanded in $\eps$ up to
order $\eps^2$. Setting $t=T$, we obtain the Poincar\'e maps of both systems
\[
\begin{split}
\ph^0_T(y) &\bydef P^0(y,\eps) 
= y + \eps P^0_1(y) + \eps^2 P^0_2(y,\eps) \\
\ph_T(y) &\bydef P(y,\eps) 
= y + \eps P_1(y) + \eps^2 P_2(y,\eps).
\end{split}
\]
By our estimate on $\norm{z(t)}$, we must have $P_1(y)=P^0_1(y)$.  If
$y^\star$ is an equilibrium of \eqref{av10} (that is, if
$\avrg{g}(y^\star)=0$), it is a fixed point of $P^0$, and in particular 
\[
P^0_1(y^\star) = 0.
\]
If we denote $\dpar{}{y}\avrg{g}(y^\star)$ by $A$, we also have 
\[
\dpar{P^0}y (y^\star,\eps) = \e^{\eps AT} = \one + \eps AT + \Order{\eps^2} 
\qquad\Rightarrow\qquad
\dpar{P^0_1}y (y^\star) = AT. 
\]
Now if we define the function 
\[
\Phi(y,\eps) = P_1(y) + \eps P_2(y,\eps),
\]
the above relations imply that 
\[
\Phi(y^\star,0) = 0,
\qquad\qquad
\dpar\Phi{y}(y^\star,0) = AT.
\]
Thus if $A$ has no vanishing eigenvalue, the implicit function theorem
implies the existence, for small $\eps$, of a fixed point $y^\star(\eps)$ of
$P$. This fixed point corresponds to a periodic orbit $\gamma_\eps$, which
remains in an $\eps$-\nbh\ of $y^\star$ because of our estimate on
$\norm{z}$. Since $x=y+\eps w(y,t)$, the corresponding $x$-coordinates of
the periodic orbit are also $\eps$-close to $y^\star$. 

Finally, if $A$ has no imaginary eigenvalue, then the linearization of $P^0$
at $y^\star$ has no eigenvalue on the unit circle, and this remains true,
for sufficiently small $\eps$, for the linearization of $P$ at
$y^\star(\eps)$. Thus the fixed points of $P$ and $P^0$ have the same type
of stability, and this property carries over to the periodic orbit.
\qed
\end{itemiz}
\renewcommand{\qed}{}
\end{proof}

The averaged approximation of a periodic system is considerably easier to
study, since it does not depend on time. In order to find periodic orbits of
the initial system, it is sufficient to find equilibrium points of the
averaged system. The approximation is reliable on a larger time scale than
the unperturbed system. However, it does not give any information on
stability if the averaged system has an elliptic equilibrium. 

\begin{example}
\label{ex_av2}
Consider a small perturbation of a two-degree-of-freedom integrable
system, whose Hamiltonian depends only on one action:
\begin{equation}
\label{av11}
H(\ph_1,\ph_2,I_1,I_2,\eps) = H_0(I_1) + \eps H_1(\ph_1,\ph_2,I_1,I_2,\eps).
\end{equation}
This situation arises, for instance, in the Kepler problem of a planet
orbiting the Sun. If $\eps=0$, $\ph_1$ rotates with constant speed
$\Omega(I_1) = H_0'(I_1)$ (as does the area swept by the planet in its
orbit), while all other variables are fixed. If $\Omega(I_1)\neq 0$, we may
write the equations of motion for sufficiently small $\eps$ as 
\begin{equation}
\label{av12}
\begin{split}
\dtot{\ph_2}{\ph_1} &= \eps
\frac{\tdpar{H_1}{I_2}}{\Omega(I_1)+\eps\tdpar{H_1}{I_1}} \\
\dtot{I_j}{\ph_1} &= -\eps
\frac{\tdpar{H_1}{\ph_j}}{\Omega(I_1)+\eps\tdpar{H_1}{I_1}}, 
\qquad\qquad j=1,2.
\end{split}
\end{equation}
The averaged version of this system is 
\begin{equation}
\label{av13}
\begin{split}
\dtot{\ph_2}{\ph_1} &= \frac\eps{\Omega(I_1)} 
\frac1{2\pi} \int_0^{2\pi} 
\dpar{H_1}{I_2}(\ph_1,\ph_2,I_1,I_2,0) \6\ph_1 \\
\dtot{I_2}{\ph_1} &= -\frac\eps{\Omega(I_1)} 
\frac1{2\pi} \int_0^{2\pi} 
\dpar{H_1}{\ph_2}(\ph_1,\ph_2,I_1,I_2,0) \6\ph_1 \\
\dtot{I_1}{\ph_1} &= -\frac\eps{\Omega(I_1)} 
\frac1{2\pi} \int_0^{2\pi} 
\dpar{H_1}{\ph_1}(\ph_1,\ph_2,I_1,I_2,0) \6\ph_1 = 0.
\end{split}
\end{equation}
Now we can make an important observation. Consider the \defwd{averaged
Hamiltonian}
\begin{equation}
\label{av14}
\begin{split}
\avrg{H}(\ph_2,I_2;I_1,\eps) 
&= \frac1{2\pi} \int_0^{2\pi} H(\ph_1,\ph_2,I_1,I_2,\eps) \6\ph_1 \\
&= H_0(I_1) + \eps \avrg{H_1}(\ph_2,I_2;I_1,\eps).
\end{split}
\end{equation}
The ordering of the variables in \eqref{av14} is meant to be suggestive,
since we will obtain a closed system of equations for $\ph_2,I_2$. Indeed, 
the associated canonical equations are
\begin{align}
\nonumber
\dot\ph_1 &= \Omega(I_1) + \eps\dpar{\avrg{H_1}}{I_1}(\ph_2,I_2;I_1,\eps) & 
\dot I_1 &= 0 \\
\label{av15}
\dot\ph_2 &= \eps\dpar{\avrg{H_1}}{I_2}(\ph_2,I_2;I_1,\eps) & 
\dot I_2 &= -\eps\dpar{\avrg{H_1}}{\ph_2}(\ph_2,I_2;I_1,\eps). 
\end{align}
At first order in $\eps$, these equations are equivalent to \eqref{av13}. We
conclude that for Hamiltonian systems of the form \eqref{av11}, averaging
can be performed directly on the Hamiltonian, with respect to the \lq\lq
fast\rq\rq\ variable $\ph_1$. 

Note that $I_1$ and $\avrg{H}$ are constants of the motion for the averaged
system. By the averaging theorem, they vary at most by an amount of order
$\eps$ during time intervals of order $1/\eps$. Such quantities are called
\defwd{adiabatic invariants}. 

Remark also that since the averaged Hamiltonian no longer depends on
$\ph_1$, the only dynamical variables are $\ph_2$ and $I_2$, while $I_1$
plays the r\^ole of a parameter. The averaged Hamiltonian thus has one
degree of freedom and is integrable. One can generalize the method to more
degrees of freedom. In that case, the averaged system has one degree of
freedom less than the original one. There also exist generalizations to
systems with several fast variables. 
\end{example}

The method of averaging of Hamiltonian systems with respect to fast
variables is of major importance in applications. For instance, it allowed
Laskar to integrate the motion of the Solar System over several hundred
million years, while numerical integrations of the original system are only
reliable on time spans of a few hundred to a few thousand years. 


\subsection{Iterative Methods}
\label{ssec_itm}

We have seen that the averaging procedure allows us to obtain a better
approximation of the perturbed system, valid for longer times. A natural
question is whether this method can be pushed further, in order to yield
still better approximations. This is indeed possible to some extent.
Consider for instance the equation 
\begin{equation}
\label{it1}
\dot y = \eps \avrg{g}(y) + \eps^2 g_1(y,t,\eps)
\end{equation}
of the transformed system in Theorem~\ref{thm_avrg}. Assuming sufficient
differentiability, and proceeding in the same way as when we constructed
the change of variables $x=y+\eps w(y,t)$, we can construct a new change of
variables $y=y_2+\eps^2 w_2(y_2,t)$ such that 
\begin{equation}
\label{it2}
\dot y_2 = \eps \avrg{g}(y_2) + \eps^2 \avrg{g_1}(y_2) + \eps^3 g_2(y_2,t,\eps).
\end{equation}
The associated averaged system is 
\begin{equation}
\label{it2a}
\dot y^0_2 = \eps \avrg{g}(y^0_2) + \eps^2 \avrg{g_1}(y^0_2),
\end{equation}
and if $K_2$ and $M_2$ denote, respectively, a Lipschitz constant for
$\avrg{g}+\eps\avrg{g_1}$ and an upper bound on $\norm{g_2}$, one easily
obtains from Gronwall's inequality that 
\begin{equation}
\label{it2b}
\norm{y_2(t)-y^0_2(t)} \leqs \Bigpar{\norm{y_2(0)-y^0_2(0)} +
\frac{\eps^2M_2}{K_2}} \e^{\eps K_2t}. 
\end{equation}
Thus for an appropriate initial condition, the error remains of order
$\eps^2$ for times smaller than $1/(K_2\eps)$. This procedure can be
repeated as long as the system is differentiable enough, and after $r-1$
steps we get exact and approximated systems of the form 
\begin{equation}
\label{it2c}
\begin{split}
\dot y_r &= \eps \avrg{g}(y_r) + \dots + \eps^r \avrg{g_{r-1}}(y_r) +
\eps^{r+1} g_r(y_r,t,\eps)\\
\dot y^0_r &= \eps \avrg{g}(y^0_r) + \dots + \eps^r \avrg{g_{r-1}}(y^0_r),
\end{split}
\end{equation}
satisfying 
\begin{equation}
\label{it2d}
\norm{y_r(t)-y^0_r(t)} \leqs \Bigpar{\norm{y_r(0)-y^0_r(0)} +
\frac{\eps^{r}M_r}{K_r}} \e^{\eps K_rt}. 
\end{equation}
For sufficiently small $\norm{y_r(0)-y^0_r(0)}$, this error remains of order
$\eps^{r}$ for times of order $1/(K_r\eps)$. One could thus wonder, if the
initial system is analytic, whether this procedure can be repeated
infinitely often, in such a way that the time-dependence is eliminated
completely. It turns out that this is {\em not} possible in general. The
reason is that the bounds $M_r$ and $K_r$ can grow rapidly with $r$
(typically, like $r!$), so that the amplitude of the error \eqref{it2d}
becomes large, while the time interval on which the approximation is valid
shrinks to zero. 

Once again, we are forced to conclude that the method of averaging provides
very useful approximations on a given bounded time scale, up to some power
in $\eps$, but fails in general to describe the long-term dynamics,
especially in the \nbh\ of elliptic periodic orbits. We will return to this
point in the next section. 

Similar properties hold for Hamiltonian systems, for which these iterative
methods were first developed. The main difference is that certain components
of the averaged equations may vanish (see for instance \eqref{av15} in
Example~\ref{ex_av2}), which implies that the averaged approximation may be
valid for longer time intervals. Consider a Hamiltonian system of the form 
\begin{equation}
\label{it3}
H(I,\ph,\eps) = H_0(I) + \eps H_1(I,\ph) + \eps^2 H_2(I,\ph) + \dotsb
\end{equation}
appearing for instance when describing the dynamics of the Solar System.  A
method to simplify this system, which is due to von Zeipel, proceeds as
follows:
\begin{itemiz}
\item	Let $W(\ph,J,\eps) = \ph\cdot J+\eps W_1(\ph,J)$ be a generating
function. It defines a near-identity canonical transformation
$(\ph,I)\mapsto(\psi,J)$ implicitly by 
\begin{equation}
\label{it4}
I = \dpar W\ph = J + \eps\dpar{W_1}\ph, \qquad\qquad
\psi = \dpar WJ = \ph + \eps\dpar{W_1}J.
\end{equation}

\item	Invert the relations \eqref{it4} to obtain $I$ and $\ph$ as
functions of $J$ and $\psi$, expressed as series in $\eps$. Express $H$ as a
function of the new variables. Choose $W_1$ in such a way that the
coefficient $K_1(J,\psi)$ of order $\eps$ is as simple as possible, e.g.\
depends only on $J$. 

\item	Repeat this procedure in order to simplify higher order terms. 
\end{itemiz}

The aim is to produce, after a certain number of iterations, new variables
$(J,\psi)$ in which the Hamiltonian takes a simple form such as  
\begin{equation}
\label{it5}
K(J,\psi,\eps) = K_0(J) + \eps K_1(J) + \dots + \eps^r K_r(J) +
\eps^{r+1}K_{r+1}(J,\psi,\eps),
\end{equation}
so that $\dot J=\Order{\eps^{r+1}}$, which means that $J$ is an adiabatic
invariant on the time scale $1/\eps^{r+1}$. Ideally, if we were able to
eliminate all $\psi$-dependent terms, we would have produced an integrable
system by a change of variables. 

Such a transformation does {\em not} exist in general. But even computing
an approximation of the form \eqref{it5} to an integrable system by von
Zeipel's method, which would be useful for finite but long time scales,
turns out to be extremely difficult in practice. The main difficulty is to
invert the relations \eqref{it4}, which requires a lot of computations.
Nevertheless, the method was used in the nineteenth century to compute the
orbits of the Moon and planets to high order in $\eps$. 


\subsection{Lie--Deprit Series}
\label{ssec_lds}

Andr\'e Deprit has introduced a different method, allowing to compute
simplified Hamiltonians of the form \eqref{it5} in a very systematic way,
and avoiding the time-consuming inversions. The method has been developed
for analytic Hamiltonians of the form 
\begin{equation}
\label{ld1}
H(q,p,\eps) = \sum_{n\geqs0} \frac{\eps^n}{n!} H_n(q,p).
\end{equation}
Instead of simplifying $H$ by successive canonical transformations of the
form \eqref{it4}, which mix old and new variables, only one canonical
transformation is carried out, which does not mix the variables. It is
defined by a function 
\begin{equation}
\label{ld2}
W(q,p,\eps) = \sum_{n\geqs0} \frac{\eps^n}{n!} W_n(q,p),
\end{equation}
where the $W_n$ have to be determined. The idea is to construct a
near-identity transformation by letting $\eps$ play the r\^ole of time: $q$
and $p$ are considered as functions of $Q$ and $P$ via the initial value
problem 
\begin{align}
\nonumber
\dtot q\eps &= \dpar Wp(q,p,\eps) & 
q(0) &= Q \\
\dtot p\eps &= -\dpar Wq(q,p,\eps) & 
p(0) &= P. 
\label{ld3}
\end{align}
For $\eps=0$, this is simply the identity transformation, while for
$0<\eps\ll1$, it is a near-identity canonical transformation because it has
the structure of a Hamiltonian flow (if $p\cdot\6q-H\6t$ and
$p\cdot\6q-W\6\eps$ are stationary, then $P\cdot\6Q-H\6t$ is stationary). 

The major advantage is that no functions need to be inverted. The main
problem is now, given the $W_n$, to express $H(q,p,\eps)$ as a function of
$Q$, $P$ and $\eps$. This can be done in a straightforward way, using the
notion of a Lie derivative. 

Consider first a system with Hamiltonian $H(q,p)$. The \defwd{Lie
derivative} generated by $H$ is the map
\begin{equation}
\label{ld4}
L_H: f\mapsto \poisson fH.
\end{equation}
It is convenient to introduce the notations 
\begin{equation}
\label{ld5}
L^0_H(f) = f, \qquad
L^1_H(f) = L_H(f), \qquad
L^k_H(f) = L_H(L^{k-1}_H(f)), \quad k\geqs 2.
\end{equation}
Then for any analytic function $f(q,p)$ we have 
\begin{align}
\nonumber
\dtot{}t f(q(t),p(t)) &= 
\poisson{f}H (q(t),p(t)) = L_H(f)(q(t),p(t)) \\
\label{ld6}
\dtot{^2}{t^2} f(q(t),p(t)) &= 
\poisson{\dtot{}t f}H = L_H(L_H(f)) = L^2_H(f),
\end{align}
and so on. By Taylor's formula, we thus obtain the relation 
\begin{equation}
\label{ld7}
f(q(t),p(t)) = \sum_{k\geqs 0} \frac{t^k}{k!} L^k_H(f)(q(0),p(0)),
\end{equation}
valid for all $t$ sufficiently small for this series to converge.
Symbolically, we can write this relation as 
\begin{equation}
\label{ld8}
f(q(t),p(t)) = \exp\set{tL_H}(f)(q(0),p(0)).
\end{equation}
Consider now a case where $f(q,p,t)$ depends explicitly on time. The
situation is similar, with the difference that 
\begin{equation}
\label{ld9}
\dtot{}{t} f(q(t),p(t),t) = \poisson fH + \dpar ft.
\end{equation}
Thus is we introduce the notation
\begin{equation}
\label{ld10}
\Delta_H(f) = \poisson fH + \dpar ft,
\end{equation}
and define $\Delta^k_H$ in a similar way as $L^k_H$, we can write 
\begin{equation}
\label{ld11}
f(q(t),p(t),t) = \sum_{k\geqs 0} \frac{t^k}{k!} \Delta^k_H(f) (q(0),p(0),0)
= \exp\set{t\Delta_H}(f) (q(0),p(0),0). 
\end{equation}
With these notations, the solution of \eqref{ld3} can be represented as 
\begin{equation}
\label{ld12}
\begin{split}
q(Q,P,\eps) &= \exp\set{\eps\Delta_W}(Q) \\
p(Q,P,\eps) &= \exp\set{\eps\Delta_W}(P).
\end{split}
\end{equation}
For a general analytic function $f(q,p,\eps)$, we have 
\begin{equation}
\label{ld13}
F(Q,P,\eps) \defby f(q(Q,P,\eps),p(Q,P,\eps),\eps) 
= \exp\set{\eps\Delta_W}(f)(Q,P,0).
\end{equation}
The inverse transformation of \eqref{ld12} is given by 
\begin{equation}
\label{ld14}
\begin{split}
Q(q,p,\eps) &= \exp\set{-\eps\Delta_W}(q) \\
P(q,p,\eps) &= \exp\set{-\eps\Delta_W}(p).
\end{split}
\end{equation}
We are now prepared to state and prove a general result, allowing to
determine the effect of a change of variables defined by \eqref{ld3} on a
function of $p$ and $q$, such as the Hamiltonian. 

\begin{prop}
\label{prop_ld}
Consider an analytic function 
\begin{equation}
\label{ld15}
f(q,p,\eps) = \sum_{n\geqs 0}\frac{\eps^n}{n!} f_n(q,p). 
\end{equation}
Let $(Q,P)$ be new variables defined by \eqref{ld3} with 
\begin{equation}
\label{ld16}
W(q,p,\eps) = \sum_{n\geqs 0}\frac{\eps^n}{n!} W_n(q,p). 
\end{equation}
Then the expression of $f$ in the new variables is given by 
\begin{equation}
\label{ld17}
F(Q,P,\eps) \defby f(q(Q,P,\eps),p(Q,P,\eps),\eps) 
= \sum_{n\geqs 0}\frac{\eps^n}{n!} F_n(Q,P),
\end{equation}
where the $F_n$ are determined iteratively as follows. Let
$f^0_n(Q,P)=f_n(Q,P)$. Define functions $f^k_n(Q,P)$ recursively by the
relation 
\begin{equation}
\label{ld18}
f^k_n(Q,P) = f^{k-1}_{n+1}(Q,P) + \sum_{m=0}^n \binom nm
\poisson{f^{k-1}_{n-m}}{W_{m}}(Q,P).
\end{equation}
Then $F_n(Q,P) = f^n_0(Q,P)$. 
\end{prop}
\begin{proof}
The idea is to compute $\Delta^k_W(f)$ by induction, and then to use
\eqref{ld13}. We first observe that 
\[
\begin{split}
\dpar{}{\eps} f(Q,P,\eps) &= \sum_{n\geqs0} \frac{\eps^n}{n!} f_{n+1}(Q,P)
\\ 
L_W f(Q,P,\eps) &= \sum_{n\geqs0} \sum_{m\geqs0} \frac{\eps^n}{n!} 
\frac{\eps^m}{m!} \poisson{f_n}{W_m}(Q,P) \\
&= \sum_{k\geqs0} \frac{\eps^k}{k!} \sum_{m=0}^k \frac{k!}{m!(k-m)!}
\poisson{f_{n-m}}{W_{m}}(Q,P). 
\end{split}
\]
This shows that 
\[
\Delta_W f(Q,P,\eps) = \sum_{n\geqs0} \frac{\eps^n}{n!} f^1_n(Q,P),
\]
where 
\[
f^1_n(Q,P) = f_{n+1}(Q,P) + \sum_{m=0}^n \binom nm
\poisson{f_{n-m}}{W_{m}}(Q,P). 
\]
In particular, for $\eps=0$ we find
$F_1(Q,P)=\Delta_W(f)(Q,P,0)=f^1_0(Q,P)$.  Now by induction, one obtains in
a similar way that 
\[
\Delta^k_W f(Q,P,\eps) = \sum_{n\geqs0} \frac{\eps^n}{n!} f^k_n(Q,P),
\]
where 
\[
f^k_n(Q,P) = f^{k-1}_{n+1}(Q,P) + \sum_{m=0}^n \binom nm
\poisson{f^{k-1}_{n-m}}{W_{m}}(Q,P). 
\]
Since, for $\eps=0$, $F_k(Q,P)=\Delta_W^k(f)(Q,P,0)=f^k_0(Q,P)$, the
assertion is proved. 
\end{proof}

The iterative scheme may by represented by the \defwd{Lie triangle}: 

\begin{center}
\begin{tabular}{ccccccccccc}
 & & $F_0$ & & $F_1$ & & $F_2$ & & $F_3$ & & \dots \\
 & &  $\parallel$  & &  $\parallel$  & &  $\parallel$  & &  $\parallel$  \\
\vrule height 14pt depth 7pt width 0pt
$f_0$ & $=$ & $f^0_0$ & $\longrightarrow$ & $f^1_0$ & $\longrightarrow$ &
$f^2_0$ & $\longrightarrow$ & $f^3_0$ & $\longrightarrow$ & \dots \\
 & & & $\nearrow$ & & $\nearrow$ & & $\nearrow$ \\
$f_1$ & $=$ & $f^0_1$ & $\longrightarrow$ & $f^1_1$ & $\longrightarrow$ &
$f^2_1$ & $\longrightarrow$ & \dots \\
 & & & $\nearrow$ & & $\nearrow$ \\
$f_2$ & $=$ & $f^0_2$ & $\longrightarrow$ & $f^1_2$ & $\longrightarrow$ & 
\dots  \\
 & & & $\nearrow$ \\
$f_3$ & $=$ & $f^0_3$ & $\longrightarrow$ & \dots  \\
\\
\dots
\end{tabular}
\end{center}

In practice, one proceeds as follows. The Hamiltonian $H(q,p,\eps)$ is
expanded in $\eps$ to some finite order $r$, and provides the entries $H_n$
in the left column of the triangle. One can then compute the $H^k_n$ in one
diagonal after the other, and the upper row will contain the components
$K_n$ of the Hamiltonian in the new variables. Each diagonal depends only
on previously computed quantities and one additional term $W_n$. 

The first two relations are 
\begin{align}
\nonumber
K_0(Q,P) &= H^0_0(Q,P) = H_0(Q,P) \\
K_1(Q,P) &= H^1_0(Q,P) = H_1(Q,P) + \poisson{H_0}{W_0}(Q,P).
\label{ld19}
\end{align}
The second relation allows us to determine $W_0$, which we choose in such a
way that $K_1$ is as simple as possible. The second diagonal of the triangle
can then be computed by the relations
\begin{align}
\nonumber
H^1_1 &= H_2 + \poisson{H_1}{W_0} +
\poisson{H_0}{W_1} \\ 
H^2_0 &= H^1_1 +
\poisson{H^1_0}{W_0},
\label{ld20}
\end{align}
so that we obtain 
\begin{equation}
\label{ld21}
K_2 = H_2 + 2\poisson{H_1}{W_0} +
\poisson{H_0}{W_1} + \poisson{\poisson{H_0}{W_0}}{W_0}.
\end{equation}
Again, this allows us to choose $W_1$ in such a way as to simplify $K_2$ as
much as possible. This procedure is very systematic, and thus suitable for
computer algebra. 

\begin{example}
\label{ex_ld}
Consider again a Hamiltonian as in Example~\ref{ex_av2}:
\begin{equation}
\label{ld22}
H(\ph_1,\ph_2,I_1,I_2,\eps) = H_0(I_1) + \sum_{n\geqs1} \frac{\eps^n}{n!}
H_n(\ph_1,\ph_2,I_1,I_2), 
\end{equation}
where $\Omega(I_1)=H_0'(I_1)\neq 0$. We would like to determine $W$ in such
a way as to eliminate the dependence on $\ph_1$ of the $H_n$. From
\eqref{ld19} and \eqref{ld20}, we get 
\begin{equation}
\label{ld23}
\begin{split}
K_0 &= H_0 \\
K_1 &= H_1 + \poisson{H_0}{W_0} \\
K_2 &= H_2 + 2\poisson{H_1}{W_0} +
\poisson{H_0}{W_1} + \poisson{\poisson{H_0}{W_0}}{W_0}.
\end{split}
\end{equation}
Since $H_0$ depends only on $I_1$, the second relation can be written as 
\begin{equation}
\label{ld24}
K_1(\psi,J) = H_1(\psi,J) - \Omega(J_1) \dpar{W_0}{\psi_1}(\psi,J). 
\end{equation}
We want to find $W_0$ periodic in $\psi_1$ such that $K_1$ no longer depends
on $\psi_1$. This can be achieved by taking 
\begin{equation}
\label{ld25}
W_0(\psi_1,\psi_2,J_1,J_2) = \frac1{\Omega(J_1)} \int_0^{\psi_1}
\bigbrak{H_1(\psi_1',\psi_2,J_1,J_2) - \avrg{H_1}(\psi_2,J_1,J_2)}
\6\psi_1',
\end{equation}
where $\avrg{H_1}$ denotes the average of $H_1$ over $\psi_1$. As a result,
we have 
\begin{equation}
\label{ld26}
K_1(\psi_1,\psi_2,J_1,J_2) = \avrg{H_1}(\psi_2,J_1,J_2).
\end{equation}
Thus in this particular case, the Lie--Deprit method to first order is
identical with the averaging method. Continuing to higher orders, we can
produce a new Hamiltonian of the form 
\begin{equation}
\label{ld27}
\begin{split}
K(\psi_1,\psi_2,J_1,J_2,\eps) ={}& K_0(J_1) + \eps K_1(\psi_2,J_1,J_2) +
\dots + \eps^r K_r(\psi_2,J_1,J_2) \\
&{}+ \eps^{r+1} R(\psi_1,\psi_2,J_1,J_2,\eps). 
\end{split}
\end{equation}
Writing down the canonical equations, we find in particular that
$\dot{J_1}=\Order{\eps^{r+1}}$, and thus $J_1$ is an adiabatic invariant on
the time scale $1/\eps^r$. 
\end{example}



\section{Kolmogorov--Arnol'd--Moser Theory}
\label{sec_kam}

So far, we have discussed a few iterative methods allowing to solve
approximately differential equations depending on a small parameter. By $r$
successive changes of variables, the original equation is transformed into a
perturbation of order $\eps^{r+1}$ of a solvable equation. The solution of
the transformed equation remains, on some time scale, $\eps^{r+1}$-close to
the solution of its approximation. Going back to original variables, we thus
obtain an expansion of the form
\begin{equation}
\label{kam1}
x(t) = x^0(t) + \eps x^1(t) + \dots + \eps^r x^r(t) +
\eps^{r+1}R_{r+1}(t,\eps),
\end{equation}
where the functions $x^0,\dots,x^r$ can be computed, and $\norm{R}$ can be
bounded on some time interval. 

For a number of purposes, expansions of the form \eqref{kam1} are
sufficient. If one wants to understand the long-time behaviour of the
system, however, they are not of much use. So the question arises whether
the limit $r\to\infty$ can be taken, in such a way that the remainder
$\eps^{r+1}R_{r+1}$ disappears. In certain cases, it turns out that
although the $x^j(t)$ have more or less constant amplitude for a number of
$j$, they ultimately grow very quickly with $j$, preventing the convergence
of the series. 

Poincar\'e already pointed out the difference between a mathematician's
definition of convergence of a series, and a more pragmatic one, useful in
applications. Consider the two formal series
\begin{align}
\label{kam2a}
S_1 &= \sum_{j\geqs0} \frac{(1000)^j}{j!} 
= 1 + 10^3 + \frac12 \cdot 10^6 + \frac16 \cdot 10^9 + \dotsb 
= \e^{1000}\\
S_2 &= \sum_{j\geqs0} \frac{j!}{(1000)^j} 
= 1 + 10^{-3} + 2 \cdot 10^{-6} + 6 \cdot 10^{-9} + \dotsb 
\end{align}
The first series is convergent, while the second is only formal. However,
the fact that $S_1$ converges is pretty useless in practice since the first
thousand terms increase, while the partial sums of the divergent series
$S_2$ hardly change for the first thousand terms. 

In a number of cases, it has been shown that the series \eqref{kam1}
diverges for some initial conditions. There are thus two possible approaches
to the problem: 
\begin{enum}
\item	Admitting that \eqref{kam1} may diverge for some initial conditions,
find the optimal $r$ such that the remainder $\eps^{r+1}R_{r+1}(t,\eps)$ is
as small as possible. This approach was first developed by Nekhoroshev and
Neishtadt, and we will discuss an example in the last chapter. 

\item	Find initial conditions for which the series \eqref{kam1} converges.
This was first achieved by Kolmogorov in 1954 for a class of analytic
near-integrable Hamiltonians (he did only publish an outline of the proof),
and generalized by Arnol'd in 1961. In 1962, Moser obtained a similar result
for near-integrable mappings of {\em finite} differentiability. 
\end{enum}
Results of the second category are known as Kolmorogov--Arnol'd--Moser
(or KAM) theory. They have been considerably developed in the last 40 years,
and are still an active area of research. 

In the sequel, we shall illustrate the problems and results in a relatively
simple case, which can be described by Moser's theorem. After that, we will
mention some other results of KAM theory. 


\subsection{Normal Forms and Small Denominators}
\label{ssec_nfsd}

Let us consider a two-degree-of-freedom Hamiltonian $H(q_1,q_2,p_1,p_2)$,
admitting a periodic solution $\gamma(t)$. Our aim is to determine the
stability of $\gamma$. We will assume for simplicity that $q_1$ is an
angular variable such that 
\begin{equation}
\label{nfs1}
\dot q_1 = \dpar H{p_1} > 0
\end{equation}
in a \nbh\ of $\gamma$ (this can be achieved by a local canonical
transformation). Note that Condition \eqref{nfs1} implies that the function
$H$ can be inverted with respect to $p_1$ near $\gamma$, i.e., we have 
\begin{equation}
\label{nfs2}
p_1 = P(H,q_1,q_2,p_2). 
\end{equation}
Since $H$ is a constant of the motion, trajectories are contained in a
three-dimensional level set of $H$, which can be parametrized by $q_1$,
$q_2$ and $p_2$ (the last coordinate $p_1$ being determined by
\eqref{nfs2}). Condition \eqref{nfs1} also implies that the vector field is
transverse to hyperplanes $q_1=\text{\it constant}$. We can thus take a
Poincar\'e section, say, at $q_1=0\equiv2\pi$, and consider the Poincar\'e
map 
\begin{equation}
\label{nfs3}
\Pi: (q_2,p_2,q_1=0) \mapsto (q_2,p_2,q_1=2\pi),
\end{equation}
which is, in effect, two-dimensional, and depends on $H$ as on a parameter.
Note, moreover, that 
\begin{equation}
\label{nfs4}
\begin{split}
\dtot{q_2}{q_1} &= \frac{\tdpar H{p_2}}{\tdpar H{p_1}} = -\dpar P{p_2} \\
\dtot{p_2}{q_1} &= -\frac{\tdpar H{q_2}}{\tdpar H{p_1}} = \dpar P{q_2}. 
\end{split}
\end{equation}
$P$ plays the r\^ole of a (time-dependent) Hamiltonian, which shows in
particular that $\Pi$ is area-preserving. 

Let $x=(q_2,p_2)$ and denote by $x^\star$ the coordinates of the
intersection of $\gamma$ with the section $q_1=0$. Then
$\Pi(x^\star)=x^\star$, and since $\Pi$ is area-preserving, the Jacobian of
$\Pi$ is equal to $1$. There are three possibilities:
\begin{enum}
\item	$\dpar\Pi x(x^\star)$ has real eigenvalues $a_1\neq\pm1$ and
$a_2=1/a_1$. Then $x^\star$ is a hyperbolic saddle point of $\Pi$, and thus
$\gamma$ is unstable because nearby orbits are repelled along the unstable
manifold. 

\item	$\dpar\Pi x(x^\star)$ has complex eigenvalues $\e^{\pm2\pi\icx\th}$
of unit module, $\th\notin\Z$, and $x^\star$ is called an \defwd{elliptic}
equilibrium. This is the case we want to study. 

\item	$\dpar\Pi x(x^\star)$ has eigenvalues $+1$ or $-1$, and $x^\star$ is
called a \defwd{parabolic} equilibrium. This situation arises in
bifurcations, we will not discuss it here. 
\end{enum}
In an appropriate basis, the linearization of $\Pi$ at $x^\star$ is a
rotation of angle $2\pi\th$. We can introduce a complex variable $z$ for
which the Poincar\'e map becomes
\begin{equation}
\label{nfs5}
z \mapsto \hat z = \e^{2\pi\icx\th}z + g(z,\cc z),
\end{equation}
where $g$ admits a Taylor expansion of the form 
\begin{equation}
\label{nfs6}
g(z,\cc z) = \sum_{n+m=2}^{r} g_{nm}z^n\cc z^m + \order{\abs{z}^r}
\end{equation}
if the original Hamiltonian is of class $\cC^r$. 

If the term $g$ were absent, the map \eqref{nfs5} would be a rotation, and
thus the periodic solution $\gamma$ would be stable. A natural thing to do
is thus try to eliminate $g$ by a change of variables, i.e., compute a
\defwd{normal form} of \eqref{nfs5}. Consider a transformation of the form 
\begin{equation}
\label{nfs7}
z = \z + c \z^j \cc\z^k, 
\qquad\qquad
2\leqs j+k\leqs r.
\end{equation}
Since $g(z,\cc z)=g(\z,\cc\z) + \Order{\abs{\z}^{j+k+1}}$, the map
\eqref{nfs5} becomes 
\begin{equation}
\label{nfs8}
\hat\z + c {\hat\z}^j \cc{\hat\z}^k 
= \e^{2\pi\icx\th}\z + c \e^{2\pi\icx\th} \z^j \cc\z^k + g(\z,\cc\z) +
\Order{\abs{\z}^{j+k+1}}.  
\end{equation}
This shows that $\hat\z = \e^{2\pi\icx\th}\z + \Order{\abs{\z}^2}$.
Substituting this in the term $c {\hat\z}^j \cc{\hat\z}^k$ and grouping like
powers of $\z$, we get  
\begin{equation}
\label{nfs9}
\hat\z = \e^{2\pi\icx\th}\z + c \bigpar{\e^{2\pi\icx\th} -
\e^{2\pi\icx(j-k)\th}}
\z^j\cc\z^k + \sum_{n+m=2}^{r} g_{nm}\z^n\cc \z^m + \Order{\abs{\z}^{j+k+1}}.  
\end{equation}
We see that we can kill the term proportional to $\z^j \cc\z^k$ if we choose
\begin{equation}
\label{nfs10}
c = \frac{g_{jk}}{\e^{2\pi\icx(j-k)\th} - \e^{2\pi\icx\th}} = 
\frac{g_{jk}\e^{-2\pi\icx\th}}{\e^{2\pi\icx(j-k-1)\th}-1}. 
\end{equation}
The transformation creates new terms, but they are of higher order in
$\abs{\z}$. We can thus simplify the map \eqref{nfs5} by eliminating first
terms of order $2$, then terms of order $3$, and so on up to order $r$. 

However, the equation \eqref{nfs10} for $c$ can only be solved under the
condition
\begin{equation}
\label{nfs11}
\e^{2\pi\icx(j-k-1)\th}\neq1.
\end{equation}
Terms $g_{jk}z^j\cc z^k$ which violate this condition are called
\defwd{resonant}. There are different kinds of resonant terms:
\begin{itemiz}
\item	terms for which $k=j-1$ are always resonant, they are of the form
$g_{k+1\,k}\abs z^{2k}z$; 
\item	if $\th=p/q$ is rational, terms for which $j-k-1$ is a multiple of
$q$ are also resonant, they are of the form
$g_{jj+1+nq}\abs{z}^{2j}z^{1+nq}$ (the smaller $q$, the more such resonant
terms exist).
\end{itemiz}
One can show, in fact, that the origin is an unstable fixed point of the
map if $\th=p/q$ for some $q\leqs 4$. On the other hand, if $\th$ is
irrational, the normal form of \eqref{nfs5} is 
\begin{equation}
\label{nfs12}
\hat\z = \e^{2\pi\icx\th}\z + C_1 \abs{\z}^2\z + \dots + 
C_m \abs{\z}^{2m}\z + \order{\abs{\z}^r}, 
\qquad
m = \biggbrak{\frac r2 - 1},
\end{equation}
where $\brak{\cdot}$ denotes the integer part. Such a map is called a
\defwd{Birkhoff normal form}.

However, even if $g$ is analytic and $\th$ is irrational, it is not clear
whether all terms which are not of the form $\abs{\z}^{2j}\z$ can be
eliminated. The reason is that the remainder at some order $j+k+1$ will
involve terms depending on $c$ given by \eqref{nfs10}, the denominator of
which is of the form $\e^{2\pi\icx q\th}-1$ for some integer $q$. Since the
numbers $\e^{2\pi\icx q\th}$ are dense on the unit circle, the quantity
$\e^{2\pi\icx q\th}-1$ may become arbitrarily small for large $q$. This is
the problem of the \defwd{small denominators}. 

For now, let us assume that the differentiability $r$ is at least $4$ and
that 
\begin{equation}
\label{nfs13}
\bigpar{\e^{2\pi\icx\th}}^q \neq 1
\qquad\qquad
\text{for $q=1,2,3,4$.}
\end{equation}
Then the Birkhoff normal form is 
\begin{equation}
\label{nfs14}
\hat\z = \e^{2\pi\icx\th}\z + C_1\abs{\z}^2\z + \order{\abs{\z}^4}. 
\end{equation}
We now introduce polar coordinates $\z=I\e^{\icx\ph}$ (where $I$ will play
the r\^ole of an action variable). In order to find the expression of the
map in polar coordinates, we first compute 
\begin{equation}
\label{nfs15}
\hat I^2 = \abs{\hat\z}^2 = I^2 + 2\re(\e^{-2\pi\icx\th}C_1) I^4 +
\order{I^5}
\end{equation}
which determines $\hat I$. This allows us to compute 
\begin{equation}
\label{nfs16}
\e^{\icx\hat\ph} = \e^{\icx\ph+2\pi\icx\th}
\bigpar{1+\icx\im(\e^{-2\pi\icx\th}C_1) I^2} + \order{I^3}. 
\end{equation}
Taking the logarithm and expanding in $I$, we find that the map in polar
coordinates takes the form 
\begin{equation}
\label{nfs17}
\begin{split}
\hat\ph &= \ph + 2\pi\th + \im(\e^{-2\pi\icx\th}C_1) I^2 + \order{I^3} \\
\hat I &= I + \re(\e^{-2\pi\icx\th}C_1) I^3 + \order{I^4}.
\end{split}
\end{equation}
If $\re(\e^{-2\pi\icx\th}C_1)$ were negative, we could conclude that $\hat
I<I$ for sufficiently small $I$, which would imply the asymptotic stability
of the orbit. This, however, is incompatible with the fact that the map is
area-preserving, since it would imply that the image of a small disc is
strictly contained in the disc. For similar reasons,
$\re(\e^{-2\pi\icx\th}C_1)$ cannot be positive, and we conclude that
$\re(\e^{-2\pi\icx\th}C_1)=0$. We will write the map \eqref{nfs17} in the
form 
\begin{equation}
\label{nfs18}
\begin{split}
\hat\ph &= \ph + \Omega(I) + f(\ph,I) \\
\hat I &= I + g(\ph,I),
\end{split}
\end{equation}
where 
\begin{equation}
\label{nfs19}
\Omega(I) = 2\pi\th + \frac{\e^{-2\pi\icx\th}C_1}{\icx} I^2, 
\qquad
f=\order{I^3}, \qquad g=\order{I^4}. 
\end{equation}
If the original system satisfies higher differentiability and non-resonance
conditions, the map in polar coordinates can also be written in the form
\eqref{nfs18}, with accordingly more terms in $\Omega$ and smaller
remainders $f$ and $g$. 

The smallness of $f$ and $g$ near $I=0$ can also be expressed by
introducing a small parameter $\eps$ and a scaled variable $J=I/\eps$. We
will, however, stick to the representation \eqref{nfs18}. 


\subsection{Diophantine Numbers}
\label{ssec_diop}

We have seen that if we want to eliminate the remainders from maps of the
form \eqref{nfs5}, \eqref{nfs14} or \eqref{nfs18}, we will produce small
denominators of the form $\e^{2\pi\icx q\th}-1$ with $q\in\Z_*$. The
successive changes of variables can only be expected to converge if these
denominators are not too small, which amounts to requiring that $\th$ is
badly approximated by rationals. Such numbers are called
\defwd{Diophantine}. 

\begin{definition}
\label{def_diop1}
A number $\w\in\R\setminus\Q$ is called \defwd{Diophantine of type
$(C,\nu)$} for some real numbers $C>0$, $\nu\geqs1$ if 
\begin{equation}
\label{diop1}
\Bigabs{\w-\frac pq} \geqs \frac C{\abs{q}^{1+\nu}} 
\end{equation}
for all relatively prime $(p,q)\in\Z\times\Z_*$.  Diophantine numbers of
type $(C,1)$ are said to be of \defwd{constant type}. 
\end{definition}

Note that Condition \eqref{diop1} becomes harder to fulfill when $C$
increases or $\nu$ decreases. One can, of course, wonder whether Diophantine
numbers do exist at all. The answer is yes, for $C$ small enough. Algebraic
numbers are examples of Diophantine numbers. 

\begin{definition}
\label{def_diop2}
An irrational number $\w$ is called \defwd{algebraic of order $n\geqs2$} if
there exists a polynomial with integer coefficients 
\begin{equation}
\label{diop2}
P(z) = a_n z^n + a_{n-1}z^{n-1} + \dots + a_1z + a_0, 
\qquad\qquad
a_0,\dots,a_{n}\in\Z,\; a_n\neq0, 
\end{equation}
such that $P(\w)=0$. Algebraic numbers of order $2$ are called
\defwd{quadratic}. 
\end{definition}

For instance, $\sqrt2$ is a quadratic irrational. The following result, due
to Liouville, states that algebraic numbers are indeed Diophantine. 

\begin{theorem}[Liouville]
\label{thm_diop}
Let $\w$ be an algebraic irrational of order $n$, and let $k$, $0\leqs
k\leqs n-2$, be such that 
\begin{equation}
\label{diop3}
P(\w)=0, \quad
P'(\w)=0, \quad
\dots \quad
P^{(k)}(\w) = 0, \quad
P^{(k+1)}(\w) \neq 0.
\end{equation}
(Note that $P^{(n-1)}(\w)=n!\w\neq0$). Then $\w$ is Diophantine of type
$(C,\frac n{k+1}-1)$ for some $C>0$. In particular, quadratic irrationals
are of constant type. 
\end{theorem}
\begin{proof}
The zeroes of polynomials being isolated, there exists a $\delta>0$ such
that $P(x)\neq0$ for $0<\abs{x-\w}<\delta$. It is sufficient to check
Condition \eqref{diop1} for $\abs{\w-p/q}<\delta$, since for other rationals
$p/q$, it is trivially satisfied with $C=\delta$. 
Since 
\[
q^n P\Bigpar{\frac pq} = \sum_{j=0}^n a_j \Bigpar{\frac pq}^j q^n \in\Z,  
\]
the fact that $P(p/q)\neq0$ implies that 
\[
\Bigabs{q^n P\Bigpar{\frac pq}} \geqs 1. 
\]
On the other hand, Taylor's formula shows that there exists a $z$ between
$\w$ and $p/q$ such that 
\[
P\Bigpar{\frac pq} = \underbrace{P(\w) + \sum_{j=1}^k \frac{P^{(j)}(\w)}{j!}
\Bigpar{\frac pq-\w}^j}_{=0} + \frac{P^{(k+1)}(z)}{(k+1)!} \Bigpar{\frac
pq-\w}^{k+1}. 
\]
This shows that there exists a constant $M$, depending only on $P$ and
$\delta$, such that 
\[
\Bigabs{P\Bigpar{\frac pq}} \leqs M \Bigabs{\frac pq-\w}^{k+1}, 
\]
and thus
\[
\Bigabs{\w-\frac pq} \geqs \biggpar{\frac1M \Bigabs{P\Bigpar{\frac
pq}}}^{1/(k+1)} \geqs \Bigpar{\frac1{Mq^n}}^{1/(k+1)}.
\]
The result follows with $C=\min\set{\delta,M^{-1/(k+1)}}$. 
\end{proof}

Note that $k=0$ is always possible, so that an algebraic number of order
$n$ is always Diophantine of type $(C,n-1)$. This result is certainly not
optimal, except for $n=2$.  There are no Diophantine numbers of type
$(C,\nu)$ for $\nu<1$, but Roth showed that all algebraic irrationals are
Diophantine of type $(C(\eps),1+\eps)$ for any $\eps>0$. 

There are Diophantine numbers which are not algebraic, in fact Diophantine
numbers form a (Cantor) set of positive measure for $C$ small enough. This
can be seen by noting that Diophantine numbers in $[0,1]$ are obtained by
excluding intervals of length $C$ around $0$ and $1$, of length
$C/2^{1+\nu}$ around $1/2$, of length $C/3^{1+\nu}$ around $1/3$ and $2/3$,
and so on. Since there are less than $q$ rationals of the form $p/q$ in
$[0,1]$, the total measure of numbers in $[0,1]$ which are {\em not}
Diophantine of type $(C,\nu)$ is smaller than 
\begin{equation}
\label{diop4}
\sum_{q\geqs 1} q \frac C{q^{1+\nu}} = \sum_{q\geqs 1} \frac C{q^\nu} =
C\z(\nu),
\end{equation}
where the Riemann Zeta Function $\z(s)=\sum_{n\geqs1}n^{-s}$ is bounded for
$s>1$. Thus Diophantine numbers of type $(C,\nu)$ have positive measure for 
$\nu>1$ and $C<1/\z(\nu)$. In fact, the upper bound \eqref{diop4} can be
replaced by an exact value, using properties of the number $\phi(q)$ of
integers $0<p<q$ relatively prime with $q$ ($\phi(q)$ is called the 
\defwd{Euler function}). 

If $\w$ is a Diophantine irrational of type $(C,\nu)$, we can find a lower
bound for the small denominator in \eqref{nfs10}. Indeed, for all
$q\in\Z_*$, we have 
\begin{equation}
\label{diop5}
\abs{\e^{2\pi\icx q\w}-1}^2 = 2(1-\cos(2\pi q\w)) = 4\sin^2(\pi q\w). 
\end{equation}
Let $p$ be the integer closest to $q\w$. Since $\abs{\sin x} >
\frac2\pi\abs{x}$ for $\abs{x}\leqs\frac\pi2$ and $\abs{\sin x}$ is
$\pi$-periodic in $x$, we have 
\begin{equation}
\label{diop6}
\abs{\e^{2\pi\icx q\w}-1} = 2\abs{\sin(\pi q\w-\pi p)} 
\geqs 4 \abs{q\w - p} \geqs 4\frac C{\abs{q}^\nu}. 
\end{equation}


\subsection{Moser's Theorem}
\label{ssec_moser}

Let us now return to the map
\begin{equation}
\label{mos1}
\begin{split}
\hat\ph &= \ph + \Omega(I) + f(\ph,I) \\
\hat I &= I + g(\ph,I).
\end{split}
\end{equation}
If $f$ and $g$ are identically zero, the dynamics is a simple rotation: $I$
is constant while $\ph$ increases at each iteration by an amount depending
only on $I$. The curves $I=\text{{\it constant}}$ are invariant curves of
the map (which are closed since $\ph$ is periodic). If this map were
structurally stable, we should be able to find a change of variables
$(\ph,I)\mapsto(\psi,J)$ such that the dynamics of the perturbed system
\eqref{mos1} reduces to a rotation for $(\psi,J)$. This is not possible in
general, but Moser's theorem asserts that it is possible for some initial
conditions, living on a Cantor set. 

\begin{theorem}[Moser]
\label{thm_Moser}
Assume that the area-preserving map \eqref{mos1} is $\cC^r$, $r\geqs4$, in a
strip $a\leqs I\leqs b$. Assume that $\Omega(I)$ satisfies the \defwd{twist
condition}
\begin{equation}
\label{mos2}
\dtot\Omega I \geqs W > 0 
\qquad\qquad\text{for $a\leqs I\leqs b$.}
\end{equation}
We introduce the $\cC^r$-norm 
\begin{equation}
\label{mos3}
\norm{f}_{\cC^r} = \sup_{\ph,a\leqs I\leqs b} \;\max_{i+j\leqs r} 
\biggabs{\dpar{^{i+j}f}{\ph^i\partial I^j}}. 
\end{equation}
Then for every $\delta>0$, there exists a number $\eps>0$, depending on
$\delta$ and $r$, with the following property. If
$\w\in[\Omega(a)+C,\Omega(b)-C]$ is Diophantine of type $(C,\nu)$ for some
$\nu\in(1,\frac{r-1}2)$, $C>0$, and $\norm{f}_{\cC^r}+\norm{g}_{\cC^r} <
\eps WC^2$, then the map \eqref{mos1} admits an invariant curve of the
form 
\begin{equation}
\label{mos4}
\begin{array}{rcl}
\ph &=& \psi + u(\psi,\w) \\
I &=& \Omega^{-1}(2\pi\w) + v(\psi,\w)
\end{array}
\qquad\qquad
0\leqs\psi\leqs2\pi,
\end{equation}
where $u$ and $v$ are $2\pi$-periodic in $\psi$ and differentiable, with 
\begin{equation}
\label{mos5}
\norm{u}_{\cC^1} + \norm{v}_{\cC^1} < \delta.
\end{equation}
The dynamics on this invariant curve is given by the circle map 
\begin{equation}
\label{mos6}
\psi \mapsto \hat\psi = \psi + 2\pi\w.
\end{equation}
($\w$ is called the \defwd{rotation number} of the invariant curve). 
\end{theorem}

In fact, it is not necessary to require that the map be area-preserving. It
suffices to satisfy the weaker \defwd{intersection property}: every curve
sufficiently close to $I=\text{{\it constant}}$ should intersect its image. 
The initial result published by Moser required a differentiability
$r\geqs333$, and this condition was subsequently weakened by R\"ussman,
Douady and Herman. In fact, the condition $r>3$ is sufficient, meaning that
the map should be $\cC^3$ and satisfy a H\"older condition. 

Note that larger values of $r$ allow for larger values of $\nu$, and thus
there will be more invariant curves for a given perturbation size $\eps$.
Also, $C$ can be taken larger, so that invariant curves subsist for larger
perturbations. 

We will now give an outline of the proof, and postpone a more detailed
discussion of the consequences of Moser's theorem to the next subsection.
In this outline, we will not try to establish the optimal relations between
$\eps$, $\delta$, $r$, $C$ and $\tau$, but only give an idea how the
functions $u$ and $v$ can be constructed under suitable hypotheses. 

\begin{proof}[{\sc Some ideas of the proof of Moser's theorem}]\hfill
\begin{enum}
\item	{\bf Preliminary transformation:} 
First note that the average $\avrg{f}(I)$ can be assumed to vanish, since
otherwise we can include this average into $\Omega$. Area preservation
implies that if $f$ and $g$ are of order $\eps$, then the average of $g$ is
of order $\eps^2$ (to see this, compute the Jacobian of the map
\eqref{mos1} and average over $\ph$). For simplicity, we will assume that
$g$ also has average zero, although this assumption is not necessary. 

Since $\Omega'(I)\neq0$, we can simplify the map by using the variable
$\Omega(I)$ instead of $I$. Then $g$ will be replaced by
$\Omega(I+g)-\Omega(I)$. To keep notations simple, we do not introduce new
variables, and write the system in the form 
\[
\begin{split}
\hat\ph &= \ph + I + f(\ph,I) \\
\hat I &= I + g(\ph,I),
\end{split}
\]
where $f$ and $g$ satisfy the former bounds with $W=1$. 

\item	{\bf A difference equation:}
We will need to solve difference equations of the form 
\[
(\cL_\w u)(\psi) \defby u(\psi+2\pi\w) - u(\psi) = h(\psi),
\]
where $h$ is a given $2\pi$-periodic function with zero average. We will
denote the solution $u$ of this equation with zero average by
$u=\cL_\w^{-1}h$. 

Consider the Fourier series of $h$, 
\[
h(\psi) = \sum_{k\in\Z_*} h_k \e^{\icx k\psi}, 
\qquad\qquad
h_k = \frac1{2\pi} \int_0^{2\pi} \e^{-\icx k\psi}h(\psi)\6\psi. 
\]
Then the solution of the difference equation is given by 
\[
u(\psi)=(\cL_\w^{-1}h)(\psi) 
= \sum_{k\in\Z_*} \frac{h_k}{\e^{2\pi\icx k\w}-1} \e^{\icx k\psi}
\]
At this point, the small denominators have shown up once again. We will
need to bound the norm of $u$, as a function of the norm of $h$. 
Note that if $h\in\cC^r$, then integration by parts yields
\[
h_k = \frac1{2\pi\icx k} \int_0^{2\pi} \e^{-\icx k\psi}h'(\psi)\6\psi = \dots 
= \frac1{2\pi(\icx k)^r} \int_0^{2\pi} \e^{-\icx k\psi}h^{(r)}(\psi)\6\psi, 
\]
and thus $\abs{h_k}\leqs \abs{k}^{-r} \norm{h}_{\cC^r}$. 
Using the property \eqref{diop6} of Diophantine numbers, we obtain 
\[
\norm{\cL_\w^{-1}h}_{\cC^0} = \sup_\psi\abs{u(\psi)} 
\leqs \sum_{k\in\Z_*} \frac{\norm{h}_{\cC^r}}{\abs{k}^r}
\frac{\abs{k}^\nu}{4C} 
\leqs \frac{\z(r-\nu)}{2C} \norm{h}_{\cC^r},
\] 
provided $r>\nu$. Note that a larger differentiability $r$ allows for larger
values of $\nu$, and thus a larger set of $\w$'s. The effect of the small
denominators is now merely that the factor $\z(r-\nu)/2C$ may be large, but
it is a constant in the problem. Note that under stronger assumptions on
$r-\nu$, we can bound higher derivatives of $u$. 

\item	{\bf Change of variables:}
We would like to find a change of variables 
\[
\begin{split}
\ph &= \psi + u(\psi,\w) \\
I &= 2\pi\w + v(\psi,\w)
\end{split}
\]
transforming the map into 
\[
\begin{split}
\hat\psi &= \psi+2\pi\w \\
\hat\w &= \w.
\end{split}
\]
Inserting this into our map $(\ph,I)\mapsto(\hat\ph,\hat I)$, we find 
\[
\hat I = 2\pi\hat\w + v(\hat\psi,\hat\w) = 
2\pi\w + v(\psi,\w) + g(\psi+u(\psi,\w),\w+v(\psi,\w)).
\]
We now obtain a condition for $v$ by requiring that $\hat\psi=\psi+2\pi\w$
and $\hat\w=\w$. Proceeding in a similar way with the equation for $\ph$,
we end up with the system
\[
\begin{split}
v(\psi+2\pi\w,\w) - v(\psi,\w) &= g(\psi+u(\psi,\w),\w+v(\psi,\w)) \\
u(\psi+2\pi\w,\w) - u(\psi,\w) &= v(\psi,\w) + f(\psi+u(\psi,\w),\w+v(\psi,\w)).
\end{split}
\]
These functional equations are very hard to solve, but we may try to find a
solution by iterations. Recall that $u$ and $v$ are expected to be small. As
a first approximation, we may neglect the terms $u$ and $v$ appearing in the
argument of $f$ and $g$, so that we are reduced to the difference equation
discussed above. Thus we use as a first approximation 
\[
\begin{split}
v_0(\psi,\w) &= \cL_\w^{-1} g(\psi,\w) \\
u_0(\psi,\w) &= \cL_\w^{-1} \bigbrak{f(\psi,\w)+v_0(\psi,\w)}. 
\end{split}
\]
Note that since $\norm{g}_{\cC^r}$ and $\norm{f}_{\cC^r}$ are of order 
$\eps C^2$, our bound on $\cL_\w^{-1}$ implies that $\norm{v_0}_{\cC^0}$ is
of order $\eps C$ and $\norm{u_0}_{\cC^0}$ is of order $\eps$.  

Our aim will now be to construct, by an iterative method, a solution of the
functional equations which is close to $(u_0,v_0)$. 
 
\item	{\bf Newton's method:}
Assume for the moment that we are given a function $F:\R\to\R$ and that we
want to solve the equation $F(x)=0$. Assume further that we know an
approximation $x_0$ of the solution. By Taylor's formula, 
\[
F(x_0+\Delta) = F(x_0) + F'(x_0)\Delta + \Order{\Delta^2}, 
\]
which vanishes for $\Delta\cong-F'(x_0)^{-1}F(x_0)$. We thus expect that 
\[
x_1 = T(x_0) \defby x_0 - F'(x_0)^{-1}F(x_0)
\]
is a better approximation of the solution. Under appropriate conditions, the
sequence of $x_n=T^n(x_0)$ should converge to this solution: This is simply
Newton's method. Applying Taylor's formula to second order shows that
\[
\begin{split}
F(x_1) &= F(x_0) + F'(x_0) (x_1-x_0)
+ \frac12 F''(z)(x_1-x_0)^2 \\
&= \frac12 F''(z)\bigpar{- F'(x_0)^{-1}F(x_0)}^2,
\end{split}
\]
where $z$ lies between $x_0$ and $x_1$. Thus if $\abs{F'}\geqs C$ and
$\abs{F''}<K$, we get 
\[
\abs{F(x_1)} \leqs \frac K{2C^2} \abs{F(x_0)}^2. 
\]
This behaviour is called \defwd{quadratic convergence}. Indeed, if $M_0$ is
a small enough positive number, then the iterates of the map 
\[
M_{n+1} = \frac K{2C^2} M_n^2
\]
converge to zero very fast. In fact, we can write $M_n=(2C^2/K)\lambda_n$,
where 
\[
\lambda_{n+1} = \lambda_n^2 
\qquad
\Rightarrow
\qquad
\lambda_n = \lambda_0^{2^n},
\]
which converges to zero faster than exponentially, provided $\lambda_0<1$
(that is, provided $\abs{F(x_0)}< 2C^2/K$). Moreover, 
\[
\abs{x_n} \leqs \abs{x_0} + \sum_{j=0}^{n-1} \frac{M_j}C 
= \abs{x_0} + \frac{2C}K \sum_{j=0}^{n-1} \lambda_0^{2^n},
\]
which converges faster than a geometric series. This fast convergence can
overcome the problem of small divisors. 

\item	{\bf Solving the functional equation:}
Let us return to the functional equation. From now on, we will proceed very
heuristically, so that what follows should not be taken too literally. We
will consider a simplified equation
\[
(\cF_\w v)(\psi,\w) \defby (\cL_\w v)(\psi,\w) - g(\psi,\w+v(\psi,\w)) = 0.
\]
As a first approximation, we use $v_0(\psi,\w)= \cL_\w^{-1} g(\psi,\w)$.
Then $\norm{v_0}_{\cC^0}$ is bounded by $(M/C)\norm{g}_{\cC^r}$, where
$M=\z(r-\nu)/2$ and $C$ is the constant of the Diophantine condition.  This
shows that 
\[
(\cF_\w v_0)(\psi,\w) = g(\psi,\w) - g(\psi,\w+v_0(\psi,\w))
\]
is of order $\norm{g}_{\cC^1}\norm{v_0}_{\cC^0} \leqs (M/C)
\norm{g}_{\cC^r}^2$.  Since $\cF_\w$ is a small perturbation of $\cL_\w$,
$(\cF_\w')^{-1}$ should also be of the order $M/C$ (here $\cF_\w'$ denotes
the Fr\'echet derivative). If $K$ denotes a bound for $\cF_\w''$, the
condition for the convergence of Newton's method thus becomes 
$\norm{g}_{\cC^r}^2 < (2/K) (C/M)^3$. 

We now construct a sequence of functions $v_n$ by iterations of Newton's
method. Our estimate on $\abs{x_n}$ indicates that 
$\norm{v_n-v_0}_{\cC^0}$ can be bounded, independently of $n$, by a
constant times $(M/C)^2  \norm{g}_{\cC^r}^2$, which can be made as small as
we like by choosing $\norm{g}_{\cC^r}$ sufficiently small. Moreover, by an
argument of uniform convergence, the $v_n$ should have a continuous limit.
So it seems plausible that for any prescribed $\delta$, we can find a
continuous solution $v$ of the functional equation such that
$\norm{v}_{\cC^0}\leqs\delta$, provided $\norm{g}_{\cC^r}$ is smaller than
some function of $\delta$ and $C$. 

There are, however, several difficulties that we brushed under the carpet.
A first problem is to invert $\cF_\w'$ in a sufficiently large domain.
Since this is rather difficult, one replaces this inverse by an
approximation of the inverse. More serious is the problem of loss of
derivatives. Indeed, we only estimated the $\cC^0$-norm of $\cL_\w^{-1}(h)$
as a function of the $\cC^r$-norm of $h$, so that we are in trouble when
iterating Newton's method. To counteract this problem, one does not apply
$\cL_\w^{-1}$ directly to a $\cC^r$ function, but to an approximation of
this function which has more continuous derivatives. Thus the iterative
procedure will not converge as quickly as Newton's method, but it is still
possible to make it converge fast enough to counteract the effect of small
denominators. Finally, one has to check that the approximations $(u_n,v_n)$
{\em and} their derivatives converge uniformly. 
\qed
\end{enum}
\renewcommand{\qed}{}
\end{proof}


\subsection{Invariant Curves and Periodic Orbits}
\label{ssec_kamo}

Let us now examine some consequences of Moser's theorem. A first important
consequence it that we have obtained sufficient conditions for the stability
of elliptic fixed points. 

\begin{cor}
\label{cor_kamo}
Let $x^\star$ be an elliptic fixed point of a $\cC^4$ area-preserving map
$\Pi$ of the plane. Assume that the Jacobian matrix $\dpar\Pi x(x^\star)$
has eigenvalues $a=\e^{2\pi\icx\th}$ and $\cc a$ such that 
\begin{equation}
\label{kamo1}
\abs{a}=1 
\qquad\qquad \text{and}
\qquad\qquad
a^q \neq 1
\qquad \text{for $q=1,2,3,4$.}
\end{equation}
Assume further that the coefficient $C_1$ in the Birkhoff normal form
$\hat\z = a\z + C_1\abs{\z}^2\z + \order{\abs{\z}^4}$  (c.f.\
\eqref{nfs14}) is different from zero. Then $x^\star$ is stable, in fact
there exists a \nbh\ of $x^\star$ which is invariant under $\Pi$. 
\end{cor}
\begin{proof}
We already showed that the dynamics near $x^\star$ can be described by a map
of the form
\[
\begin{split}
\hat\ph &= \ph + \Omega(I) + f(\ph,I) \\
\hat I &= I + g(\ph,I),
\end{split}
\]
with 
\[
\Omega(I) = 2\pi\th + \frac{\e^{-2\pi\icx\th}C_1}{\icx} I^2, 
\qquad
f=\order{I^3}, \qquad g=\order{I^4}. 
\]
We should first note that due to the singularity of polar coordinates at the
origin, $f$ may only be $\cC^3$ at $I=0$. However, it is $\cC^4$
everywhere else in a \nbh\ of $I=0$ because the map $\z\mapsto\hat\z$ is
$\cC^4$ and the transformation to polar coordinates is smooth except at the
origin. For instance, we might have $f(I,\ph)=I^{7/2}\sin\ph$. 

We consider now the map in a small strip of the form $\eps\leqs I\leqs
2\eps$ (which is a small annulus in original variables). Since we exclude
the point $I=0$, $f$ and $g$ are $\cC^4$ in this domain. Replacing $I$ by
$\eps J$, we obtain a new map (in the domain $1\leqs J\leqs 2$) 
\begin{align*}
\hat\ph &= \ph + \Omega(\eps J) + \tilde f(\ph,J) &
\tilde f(\ph,J) &= f(\ph,\eps J) \\
\hat J &= J + \tilde g(\ph,J) &
\tilde g(\ph,J) &= \frac1\eps g(\ph,\eps J). 
\end{align*} 
It is easy to check that $\norm{\tilde f}_{\cC^4}$ and $\norm{\tilde
g}_{\cC^4}$ are $\order{\eps^3}$. Moreover,
\[
\Bigabs{\dtot{}{J} \Omega(\eps J)} = 
2\eps^2 \abs{C_1}J \geqs 2\eps^3\abs{C_1}.
\]
Although we did not discuss that point in detail, Moser's theorem remains
true if $W$ is a small parameter. Thus the theorem yields the existence of
invariant curves in the strip, provided we take $\eps$ small enough (i.e.,
we magnify a sufficiently small \nbh\ of the elliptic point). Since the map
leaves invariant $x^\star$ and a curve surrounding $x^\star$, continuity
implies that the interior of that curve is invariant. 
\end{proof}

Note that Condition~\eqref{kamo1} only excludes the strongest resonances
$a=\pm1$, $a=\pm\icx$ and $a=\e^{\pm2\pi\icx/3}$. The result remains true if
the argument of $a$ is another rational multiple of $2\pi$ because the
condition on $C_1$ implies the twist condition on $\Omega$. The twist
condition in turn ensures that even if the angle of rotation is rational at
$x^\star$, is will change, and thus go through Diophantine values, as one
moves away from $x^\star$. 

Let us now illustrate some further properties of maps of the form
\eqref{mos1}, by considering the well-studied \defwd{standard map}
\begin{equation}
\label{kamo2}
\begin{split}
\hat\ph &= \ph + I + \eps\sin\ph \\
\hat I &= I + \eps\sin\ph. 
\end{split}
\end{equation}
Here $\ph$ is again an angle (modulo $2\pi$), so that the phase space has
the topology of a cylinder. However, for some purposes it is useful to
consider $\ph$ as real. Then one would define, for instance, an orbit of
period $q$ as an orbit such that the point $(\ph,I)$ is mapped to 
$(\ph+2\pi p,I)$, $p\in\Z$, after $q$ iterations.  Note that the standard
map \eqref{kamo2} is area-preserving and satisfies the \defwd{twist
condition}
\begin{equation}
\label{kamo3}
\dpar{\hat\ph}I > 0.
\end{equation}
General maps satisfying this condition are called \defwd{twist maps}. 

\goodbreak

If $\eps=0$, $I$ is constant and $\ph$ increases at each iteration by $I$.
Since $\ph$ is an angle, we can distinguish between two cases:
\begin{itemiz}
\item	if $I/2\pi = p/q$ is rational, then any orbit starting on a point of
the form $(\ph,I)$ is periodic of period $q$, and makes $p$ complete
rotations after these $q$ iterations;
\item	if $I/2\pi$ is irrational, then the orbit of $(\ph,I)$ never returns
to the same point, but fills the curve of constant $I$ in a dense way; such
orbits are called \defwd{quasiperiodic}. 
\end{itemiz}
If $\eps>0$, Moser's theorem gives us sufficient conditions for the
existence of an invariant curve. Since, in the case of the standard map,
$\norm{f}_{\cC^r}=\norm{g}_{\cC^r}=1$ for all $r$, we obtain the existence
of an invariant curve with parametric equation 
\begin{equation}
\label{kamo4}
\begin{array}{rcl}
\ph &=& \psi + u(\psi,\w) \\
I &=& 2\pi\w + v(\psi,\w)
\end{array}
\qquad\qquad
0\leqs\psi\leqs2\pi
\end{equation}
for every Diophantine $\w$ of type $(C,\nu)$, provided $\eps\leqs\eps_0
C^2$ (where $\eps_0$ is an absolute constant that we called $\eps$ in the
theorem). The $\cC_1$-norms of $u$ and $v$ are bounded by some function
$\delta(\eps)$, in fact they are of order $\eps$. Invariant curves of the
form \eqref{kamo4}, which encircle the cylinder, are called
\defwd{rotational invariant curves}. 

\begin{figure}
 \centerline{\psfig{figure=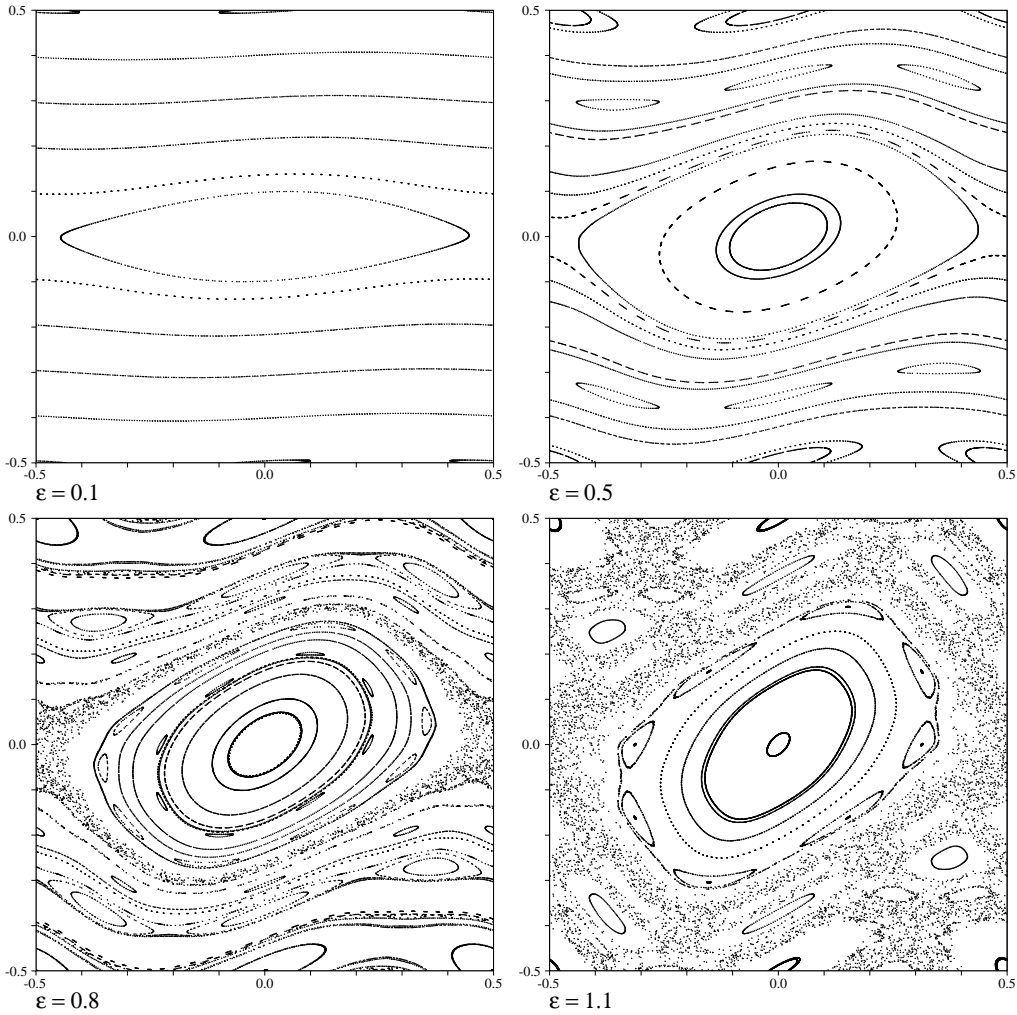,width=135mm,clip=t}}
 \vspace{3mm}
 \caption[]
 {Phase portraits of the standard map, obtained by representing several
 orbits with different initial conditions, for increasing values of the
 perturbation $\eps$.}
\label{fig_sm}
\end{figure}

\figref{fig_sm} shows that there are indeed many invariant curves for small
$\eps$. As $\eps$ grows, their number decreases, and the last one is found
to disappear for $\eps=0.97\dots$. This curve is given by \eqref{kamo4}
with $\w=(\sqrt5-1)/2$ equal to the \defwd{golden mean}, a quadratic
irrational which is particularly hard to approximate by rationals -- there
is a deep connection between this fact and the properties of continued
fractions. For $\eps$ increasing beyond $1$, the orbits become more and
more disordered, and although elliptic orbits still exist, their domains of
stability tend to shrink. 

We know by Moser's theorem that the dynamics on the invariant curve
$\eqref{kamo4}$ is given by $\psi\mapsto\psi+2\pi\w$. Thus the \nth{j}
iterate of $(\ph_0,I_0)$ has coordinates 
\begin{equation}
\label{kamo7}
\begin{split}
\ph_j &= \psi_0 + 2\pi j\w + u(\psi_0 + 2\pi j\w,\w) \\
I_j &= 2\pi\w + v(\psi_0 + 2\pi j\w,\w).
\end{split}
\end{equation}
One can associate with any orbit
$\set{(\ph_n,I_n)=\Pi^n(\ph_0,I_0)}_{n\in\Z}$ a \defwd{rotation number}
\begin{equation}
\label{kamo5}
\Omega = \frac1{2\pi} \lim_{n\to\infty} \frac1n \sum_{j=1}^n
(\ph_{j+1}-\ph_j). 
\end{equation}
There is a slight ambiguity because of the angular nature of $\ph$, which
implies that $\Omega$ is defined up to addition of an integer. For the
standard map, we can also define 
\begin{equation}
\label{kamo6}
\Omega = \frac1{2\pi} \lim_{n\to\infty} \frac1n \sum_{j=1}^n
(I_j+\eps\sin\ph_j). 
\end{equation}
For $\eps=0$, we have $\Omega=I/2\pi$. To compute the rotation number of
the invariant curve \eqref{kamo4},  we thus have to average $v(\psi_0 +
2\pi j\w,\w)$ over $j$. Since $v$ is $2\pi$-periodic in its first argument,
we use the fact that for each Fourier component $v_k(\w)\e^{\icx k\psi}$ of
$v$, 
\begin{equation}
\label{kamo8}
\sum_{j=1}^n v_k(\w)\e^{\icx k(\psi_0+2\pi j\w)} = v_k(\w)\e^{\icx k\psi_0}
\frac{1-\e^{2\pi\icx(n-1)k\w}}{1-\e^{2\pi\icx k\w}},
\end{equation} 
which is bounded uniformly in $n$. We conclude that the average of
$v(\psi_0 + 2\pi j\w,\w)$ vanishes in the limit $n\to\infty$. Hence the
rotation number $\Omega$ of any orbit on the invariant curve \eqref{kamo4}
is $\w$. The orbits of a twist map can be labelled by their rotation
number, and Moser's theorem asserts that orbits with Diophantine rotation
number belong to invariant curves for sufficiently small perturbation size
$\eps$.

Let us briefly describe other kinds of orbits appearing on \figref{fig_sm},
some of which are well understood. In particular, we have the following
result on existence of periodic orbits.

\begin{figure}
 \centerline{\psfig{figure=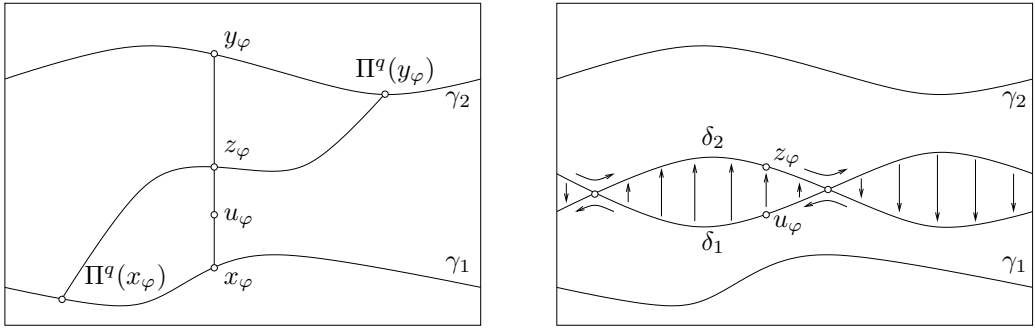,width=140mm,clip=t}}
 \figtext{
 	\writefig	3.45	1.2	$x_\ph$
 	\writefig	3.45	2.05	$u_\ph$
 	\writefig	3.45	2.95	$z_\ph$
 	\writefig	3.45	4.4	$y_\ph$
 	\writefig	1.6	1.2	$\Pi^q(x_\ph)$
 	\writefig	5.2	3.95	$\Pi^q(y_\ph)$
 	\writefig	6.4	3.6	$\gamma_2$
 	\writefig	6.4	1.4	$\gamma_1$
 	\writefig	13.8	3.6	$\gamma_2$
 	\writefig	13.8	1.4	$\gamma_1$
 	\writefig	9.8	3.05	$\delta_2$
 	\writefig	9.8	1.65	$\delta_1$
 	\writefig	10.75	1.95	$u_\ph$
 	\writefig	10.75	2.85	$z_\ph$
 }
 \captionspace
 \caption[]
 {Construction of periodic orbits with rational rotation number $p/q$
 between two invariant curves $\gamma_1$ and $\gamma_2$. $\Pi^q$ maps
 points $u_\ph$ on $\delta_1$ vertically to points $z_\ph$ on $\delta_2$.
 Their intersections are thus fixed points of $\Pi^q$. The arrows give a
 hint that one point is elliptic and the other one hyperbolic.}
\label{fig_birkhoff}
\end{figure}

\begin{theorem}[Poincar\'e--Birkhoff]
\label{thm_poincarebirkhoff}
Let $\Pi$ be an area-preserving twist map admitting two rotational invariant
curves $\gamma_1$ and $\gamma_2$, with respective rotation numbers
$\w_1<\w_2$. For every rational number $p/q\in(\w_1,\w_2)$, there exist at
least two periodic orbits of period $q$, with rotation number $p/q$, and
contained in the domain between $\gamma_1$ and $\gamma_2$. 
\end{theorem} 

\begin{proof}
Consider two points $x_\ph=(\ph,I_1)\in\gamma_1$ and
$y_\ph=(\ph,I_2)\in\gamma_2$ (\figref{fig_birkhoff}). Since
$q\w_1<p<q\w_2$, the iterate $\Pi^q(x_\ph)$ has turned by an angle smaller
than $2\pi p$, while  $\Pi^q(y_\ph)$ has turned by an angle larger than
$2\pi p$. Thus the vertical line through $x_\ph$ and $y_\ph$ must cross its
image under $\Pi^q$ at least once, say in $z_\ph$. The point $z_\ph$ is not
necessarily a fixed point of $\Pi^q$, but by construction, it is the image
of a point $u_\ph$ with the same $\ph$-coordinate. 

Consider now the curves $\delta_1=\setsuch{u_\ph}{0\leqs\ph\leqs2\pi}$ and
$\delta_2=\setsuch{z_\ph}{0\leqs\ph\leqs2\pi}$. Then
$\Pi^q(\delta_1)=\delta_2$, and the $\ph$-coordinates of all points of
$\delta_1$ are rotated by an angle $2\pi p$ (i.e., they move vertically).
Since the map is area-preserving, the area between $\gamma_1$ and
$\delta_1$ must be the same as the area between $\gamma_1$ and $\delta_2$.
Thus $\delta_1$ and $\delta_2$ must intersect at least twice. These
intersections are  fixed points of $\Pi^q$, and thus belong to two
orbits of period $q$ of $\Pi$. Since each point makes $p$ rotations during
$q$ iterations, the rotation number is $p/q$. 
\end{proof}

The assumptions of the theorem can be weakened. There exist more powerful
techniques to find periodic orbits, for instance variational ones, where one
looks for stationary points of an action functional. It can be shown that if
the periodic orbits described in the theorem are non-degenerate, then one of
them is elliptic and another one hyperbolic. Chains of elliptic and
hyperbolic points with rotation numbers $0$, $1/2$, $1/3$ and $1/4$ are
easily located on \figref{fig_sm} for $\eps=0.8$. 

As we have shown, elliptic orbits are generically stable. The dynamics in
their \nbh\ can again be described, in polar coordinates, by a twist map,
so that the same layered structure of invariant curves and elliptic and
hyperbolic orbits exists around each elliptic orbit. This self-similar
structure is most apparent for the standard map when $\eps$ is close to $1$.

Hyperbolic orbits, on the other hand, are often surrounded by chaotic
orbits. This is related to the fact that the stable and unstable manifolds
of hyperbolic orbits of given rotation number tend to intersect, producing a
so-called Smale horseshoe with complicated dynamics. 

The phase portraits on \figref{fig_sm} are quite typical for Poincar\'e
sections of two-degrees-of-freedom Hamiltonians. They also illustrate the
limits of the method of averaging. Assuming that such a phase portrait 
represents the Poincar\'e map of a Hamiltonian system with one fast
variable, the averaged system would be integrable, and thus its Poincar\'e
section would consist only of invariant curves and periodic orbits. The
chaotic orbits of the original system are smeared out by the averaging
transformation. 

Many extensions of Moser's theorem exist, for instance to multidimensional
symplectic maps. There has also been major progress in \lq\lq converse
KAM\rq\rq\ theory, which describes when and how invariant curves or tori
break up. 


\goodbreak

\subsection{Perturbed Integrable Hamiltonian Systems}
\label{ssec_kamt}

The other class of KAM theorems applies to Hamiltonian systems of the
form
\begin{equation}
\label{kamt1}
H(I,\ph) = H_0(I) + \eps H_1(I,\ph,\eps)
\end{equation}
with $n$ degrees of freedom. Let us first explain where the small
denominators come from. Assume that we want to eliminate the
$\ph$-dependence of the first-order term $H_1(I,\ph,0)$ by the Lie--Deprit
method. As we have seen, the new first order term after a canonical
transformation with generator $W_0$ will be 
\begin{align}
\label{kamt2}
K_1(J,\psi) &= H_1(J,\psi) + \poisson{H_0}{W_0}(J,\psi) \\
&= H_1(J,\psi) - \sum_{j=1}^n \Omega_j(J) \dpar{W_0}{\psi_j}(J,\psi), 
\qquad\qquad \Omega_j(J) \defby \dpar{H_0}{I_j}(J). 
\nonumber
\end{align}
Let us write $H_1$ and $W_0$ in Fourier series
\begin{equation}
\label{kamt3}
H_1(J,\psi) = \sum_{k\in\Z^n} H_{1k}(J) \e^{\icx k\cdot\psi}, 
\qquad\qquad
W_0(J,\psi) = \sum_{k\in\Z^n} W_{0k}(J) \e^{\icx k\cdot\psi}.
\end{equation}
Then we have to solve, for each $k\in\Z^n\setminus\set{0,\dots,0}$, an
equation of the form 
\begin{equation}
\label{kamt4}
\sum_{j=1}^n \Omega_j(J) W_{0k}(J) \icx k_j\e^{\icx k\cdot\psi} 
= H_{1k}(J) \e^{\icx k\cdot\psi}. 
\end{equation}
We thus take 
\begin{equation}
\label{kamt5}
W_{0k}(J) = \frac{H_{1k}(J)}{\icx k\cdot \Omega(J)} 
\qquad\qquad
\forall k\in\Z^n\setminus\set{0,\dots,0}.
\end{equation}
If $k\cdot \Omega(J) = \sum_{j=1}^n k_j \Omega_j(J) = 0$ for some $k$, then
the corresponding term of the Fourier series is \defwd{resonant} and cannot
be eliminated. But even if the $\Omega_j(J)$, which are the frequencies of
the unperturbed system, are incommensurate, the term $k\cdot \Omega(J)$ may
become arbitrarily small, unless we impose a Diophantine condition on the
$\Omega_j(J)$. The following theorem by Arnol'd gives sufficient conditions
for the existence of a change of variables, defined on some Cantor set with 
Diophantine frequencies, which eliminates the dependence of $H$ on the
angles. It thus implies the existence of invariant tori on which the motion
is quasiperiodic with Diophantine frequencies. 

\begin{theorem}[Arnol'd]
\label{thm_arnold}
Assume the Hamiltonian \eqref{kamt1} is analytic, and satisfies the
non-degeneracy condition
\begin{equation}
\label{kamt6}
\biggabs{\det\dpar{^2 H_0}{I^2}(I)} \geqs W > 0
\end{equation}
in a \nbh\ of the torus $I=I_0$. Let $\w=\Omega(I_0)\in\R^n$ satisfy the
Diophantine condition
\begin{equation}
\label{kamt7}
\bigabs{\w\cdot k} \geqs \frac C{\abs{k}^\nu} 
\qquad\qquad 
\forall k\in\Z^n\setminus\set{0,\dots,0}.
\end{equation}
where $\abs{k} = \abs{k_1}+\dots+\abs{k_n}$.
Then, if $\eps$ is sufficiently small, the Hamiltonian system \eqref{kamt1}
admits a quasiperiodic solution with frequency $\w$, i.e., this solution can
be written in the form 
\begin{equation}
\label{kamt8}
(q(t),p(t)) = F(\w_1 t,\dots,\w_n t),
\end{equation}
where $F$ is $2\pi$-periodic in all its arguments. This solution lies on an
analytic torus, which it fills densely. The distance between this torus and
the unperturbed torus $I=I_0, \ph\in\T^n$ goes to zero as $\eps\to0$. 
\end{theorem}

In this case, the invariant tori of the perturbed system can be labelled
by the frequency vector $\w$, which plays an analogous r\^ole as the
rotation number for twist maps (since the Diophantine condition only implies
ratios of components of $\w$, it does not matter whether we divide $\w$ by
$2\pi$ or not). 

Let us now examine the consequences of this result on the stability of
perturbed integrable systems. If the number of degrees of freedom is $n=2$,
the manifold of constant $H$ is three-dimensional, and contains
two-dimensional tori for $\eps$ small enough (the intersections of these
tori with a two-dimensional Poincar\'e section are invariant circles).
Since these tori act as barriers, and there are many of them, all orbits,
even those which do not start on an invariant torus, will be confined to a
relatively small portion of phase space, close to an invariant torus of the
unperturbed system. Thus the system is stable with respect to small
perturbations. 

The situation is different for $n>2$. For instance, if $n=3$, the level
sets of $H$ are five-dimensional and the invariant tori are
three-dimensional, hence they have codimension $2$. This means that tori no
longer act as barriers, and trajectories not starting on an invariant torus
may well travel a long distance in the space of actions. This phenomenon is
known as \defwd{Arnol'd diffusion}. In fact, Arnol'd has formulated the
following conjecture: for \lq\lq generic\rq\rq\ perturbations
$H_1(I,\ph,\eps)$ (in a sense to be defined), there exist a trajectory and
times $t_0$, $t_1$ such that 
\begin{equation}
\label{kamt9}
\abs{I(t_1)-I(t_0)} \geqs 1.
\end{equation} 
However, this conjecture has not been proved yet. Nekhoroshev has shown that
if $H$ is analytic, then the diffusion time must be very long: If a
trajectory satisfying \eqref{kamt9} exists, then 
\begin{equation}
\label{kamt10}
t_1 - t_0 \geqs \e^{1/\eps^\eta}, 
\qquad\qquad \eta>0.
\end{equation}
For certain trajectories, diffusion times admit even larger lower bounds,
such as $\e^{\e^{1/\eps}}$. Thus, although systems with three or more
degrees of freedom are not necessarily stable for infinite times, in
practice they are often stable on a very long time scale if the
perturbation is small. 


\section*{Bibliographical Comments}

Hamiltonian systems are discussed in detail in \cite{Arnold89}. The
averaging theorem is found in the classical textbooks \cite{GH} and
\cite{Wiggins}, as well as \cite{Verhulst}. 
The method known as Lie--Deprit series was introduced in \cite{Deprit}, and
has been extended to non-Hamiltonian systems \cite{Henrard}. 

Kolmogorov's outline of proof appeared in \cite{Kolmogorov}, and Arnol'd's
general result in \cite{Arnold63}. The first version of Moser's theorem
appeared in \cite{Moser62} and was improved by R\"ussmann \cite{Russmann70}
and Herman \cite{Herman}. Some comments on Newton's method are found in
\cite{Russmann72}. Some applications are discussed in \cite{Moser73}. 

A good review of the properties of twist maps is \cite{Meiss}. There have
been numerous developments in KAM theory. See for instance \cite{dlL01} for
a review and additional references. 


\chapter{Singular Perturbation Theory}
\label{ch_spt}

In this last chapter, we will consider so-called \defwd{slow--fast systems}
of the form 
\begin{equation}
\label{spt1}
\begin{split}
\eps\dot x &= f(x,y),\\
\dot y &= g(x,y),
\end{split}
\end{equation}
where $x\in\R^n$, $y\in\R^m$, $f$ and $g$ are of class $\cC^r$, $r\geqs2$,
in some domain $\cD\subset\R^n\times\R^m$, and $\eps$ is a small parameter.
The variable $x$ is called \defwd{fast variable} and $y$ is called
\defwd{slow variable}, because $\dot x/\dot y$ can be of order $1/\eps$. 

Equation \eqref{spt1} is an example of a \defwd{singularly perturbed}
system, because in the limit $\eps\to0$, if does not reduce to a
differential equation of the same type, but to an algebraic--differential
system 
\begin{equation}
\label{spt2}
\begin{split}
0 &= f(x,y)\\
\dot y &= g(x,y).
\end{split}
\end{equation}
Since this limiting system is not an ordinary differential equation, it is
not clear how solutions of \eqref{spt1} might be expanded into powers of
$\eps$, for instance. The set of points $(x,y)\in\cD$ for which $f(x,y)=0$
is called the \defwd{slow manifold}. If we assume that this manifold can be
represented by an equation $x=x^\star(y)$, then the dynamics of
\eqref{spt2} on the slow manifold is governed by the equation
\begin{equation}
\label{spt3}
\dot y = g(x^\star(y),y) \bydef G(y),
\end{equation}
which is easier to solve than the original equation \eqref{spt1}. One of the
questions one tries to answer in singular perturbation theory is thus: What
is the relation between solutions of the singularly perturbed system
\eqref{spt1} and those of the reduced system \eqref{spt3}? 

There is another way to study the singular limit $\eps\to0$. With respect to
the \defwd{fast time} $s=t/\eps$, \eqref{spt1} can be written as 
\begin{equation}
\label{spt4}
\begin{split}
\dtot xs &= f(x,y) \\
\dtot ys &= \eps g(x,y). 
\end{split}
\end{equation}
Taking the limit $\eps\to0$, we obtain the so-called \defwd{associated
system}
\begin{equation}
\label{spt5}
\dtot xs = f(x,y), 
\qquad\qquad y = \text{{\it constant}},
\end{equation}
in which $y$ plays the r\^ole of a parameter. 
The perturbed system in the form \eqref{spt4} can thus be considered as a
modification of the associated system \eqref{spt5} in which the \lq\lq
parameter\rq\rq\ $y$ changes slowly in time. Observe that the slow manifold
$x=x^\star(y)$ defines equilibrium points of the associated system. 

Results from regular perturbation theory (e.g.\ Corollary~\ref{cor_rpb})
show that orbits of \eqref{spt4} will remain close to those of \eqref{spt5} 
for $s$ of order $1$, i.e., for $t$ of order $\eps$, but not necessarily for
larger $t$. Our aim is to describe the dynamics for times $t$ of order $1$.

We will start by giving sufficient conditions under which \eqref{spt3} and
\eqref{spt5} indeed provide good approximations to different phases of the
motion. Then we will discuss some cases in which these conditions are
violated, and new phenomena occur. 


\section{Slow Manifolds}
\label{sec_ssm}

In order to give an idea of what to expect in the general case, let us first
discuss a solvable example. 

\begin{example}
\label{ex_ssm1}
Consider the equation
\begin{equation}
\label{ssm1}
\begin{split}
\eps\dot x &= -x + \sin(y) \\
\dot y &= 1.
\end{split}
\end{equation}
Then the slow manifold is given by the curve 
\begin{equation}
\label{ssm2}
x = x^\star(y) = \sin(y),
\end{equation}
and the dynamics reduced to the slow manifold is described by 
\begin{equation}
\label{ssm8}
x(t) = \sin(y(t)), 
\qquad\qquad y(t) = y_0 + t.
\end{equation}
The associated system for $s=t/\eps$ is 
\begin{equation}
\label{ssm3}
\dtot xs = -x + \sin(y), 
\qquad\qquad y = \text{{\it constant}},
\end{equation}
and admits the solution
\begin{equation}
\label{ssm4}
x(s) = (x_0-\sin(y)) \e^{-s} + \sin(y).
\end{equation}
Thus all solutions converge to the slow manifold $x=\sin(y)$ as
$s\to\infty$. 

\begin{figure}
 \centerline{\psfig{figure=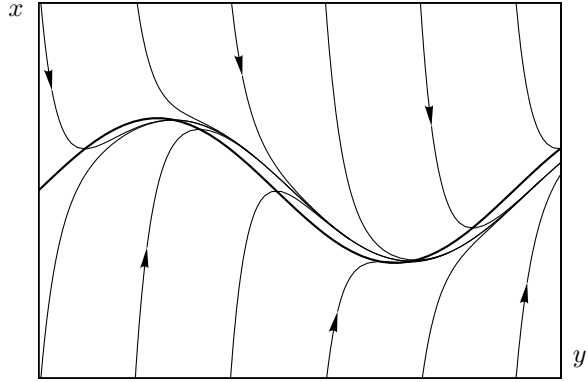,height=50mm,clip=t}}
 \figtext{
 	\writefig	11.0	0.7	$y$
 	\writefig	3.5	5.3	$x$
 }
 \captionspace
 \caption[]
 {The heavy curve represents the slow manifold \eqref{ssm2}, while the light
 curves represent solutions of the slow--fast equation \eqref{ssm1} for
 various initial conditions, and $\eps=0.07$. The solutions are quickly
 attracted by a small \nbh\ of order $\eps$ of the slow manifold.}
\label{fig_slowfast}
\end{figure}

Let us now compare the approximate solutions to the solution of the original
system \eqref{ssm1} which can, in this case, be computed by the method of
variation of the constant. It can be written in the form 
\begin{equation}
\label{ssm5}
\begin{split}
x(t) &= \bigbrak{x_0-\bar x(0,\eps)}\e^{-t/\eps}  + \bar x(t,\eps), 
\qquad
\bar x(t,\eps) = \frac{\sin(y_0+t) - \eps\cos(y_0+t)}{1+\eps^2}, \\
y(t) &= y_0 + t.
\end{split}
\end{equation}
For small $t$, we have $\bar x(t,\eps) = \sin(y_0) + \Order{\eps} +
\Order{t}$, so that 
\begin{equation}
\label{ssm6}
x(t) = \bigbrak{x_0-\sin(y_0)}\e^{-t/\eps} + \sin(y_0) + \Order{\eps} +
\Order{t},
\end{equation}
which is well approximated by the solution \eqref{ssm4} of the associated
system. As time $t$ grows, however, the drift of $y$ makes itself felt and
the associated system no longer provides a good approximation. 

On the other hand, the factor $\e^{-t/\eps}$ decreases very fast: For
$t=k\eps\abs{\log\eps}$, for instance, it is equal to $\eps^k$, and for
larger times, it goes to zero faster than any power of $\eps$. The solution
\eqref{ssm5} thus becomes very close to 
\begin{equation}
\label{ssm7}
x(t) = \frac{\sin(y(t)) - \eps\cos(y(t))}{1+\eps^2}, 
\qquad y(t) = y_0+t,
\end{equation}
which lies at a distance of order $\eps$ from the slow manifold
$x=\sin(y)$. Thus the reduced equation \eqref{ssm8} on the slow manifold
provides a good approximation of the system \eqref{ssm1} for
$t\gg\eps\abs{\log\eps}$. 
\end{example}

Tihonov's theorem gives sufficient conditions under which this behaviour
holds for general systems: Solutions are attracted, in a time of order
$\eps\abs{\log\eps}$, to a small \nbh\ of order $\eps$ of a slow manifold. 


\subsection{Tihonov's Theorem}
\label{ssec_stt}

We consider the slow--fast system
\begin{equation}
\label{stt1}
\begin{split}
\eps\dot x &= f(x,y)\\
\dot y &= g(x,y)
\end{split}
\end{equation}
with $f\in\cC^r(\cD,\R^n)$, $g\in\cC^r(\cD,\R^m)$, $r\geqs2$. Let us first
give conditions under which the slow manifold, defined implicitly by
$f(x,y)=0$, can be represented in the explicit form $x=x^\star(y)$. We
introduce the $n\times n$ matrix 
\begin{equation}
\label{stt2}
A_x(y) = \dpar fx(x,y).
\end{equation}
Assume that $f(x_0,y_0)=0$. The implicit function theorem states that if
$\det A_{x_0}(y_0)\neq 0$, then there exists a \nbh\ $\cN$ of $(x_0,y_0)$
such that all solutions of $f(x,y)=0$ in $\cN$ can be written as
$x=x^\star(y)$. Here $x^\star$ is a function of class $\cC^r$, defined in a
\nbh\ $\cN_0$ of $y_0$ in $\R^m$, and it satisfies 
\begin{equation}
\label{stt3}
\dpar{x^\star}y(y) = - A_{x^\star(y)}(y)^{-1} \dpar fy(x^\star(y),y). 
\end{equation}
We will abbreviate $A_{x^\star(y)}(y)$ by $A(y)$. From now on we will always
assume that $(x,y)$ belongs to the open set $\cN$, and that 
\begin{enum}
\item	$\det A(y)$ is bounded away from $0$ whenever $(x^\star(y),y)\in\cN$,
\item	the norms of $f$, $g$ and all their mixed partial derivatives up to
order $2$ are bounded uniformly in $\cN$.
\end{enum}

We now state a result describing how well the slow--fast system \eqref{stt1}
is approximated by its limit when $\eps\to0$. The first versions of this
result were proved independently by Tihonov and Grad\v ste\u\i n. 

\begin{theorem}
\label{thm_Tihonov}
Assume that the eigenvalues $a_j(y)$ of $A(y)$ satisfy 
\begin{equation}
\label{stt4}
\re a_j(y) \leqs -a_0 
\qquad\forall y\in\cN_0,
\quad j=1,\dots,n
\end{equation}
for some constant $a_0>0$. Let $f$ and $g$, as well as their derivatives up
to order $2$, be uniformly bounded in $\cN$. Then there exist constants 
$\eps_0$, $c_0$ to $c_4$, $K_0>0$ such that the following properties hold
for $0<\eps<\eps_0$.
\begin{enum}
\item	Any solution of \eqref{stt1} with initial condition
$(x_0,y_0)\in\cN$ such that $\norm{x_0-x^\star(y_0)}< c_0$ satisfies 
\begin{equation}
\label{stt5}
\norm{x(t) - x^\star(y(t))} \leqs c_1 \eps +
c_2\norm{x_0-x^\star(y_0)}\e^{-K_0t/\eps}
\end{equation}
for all $t\geqs0$ such that $(x(s),y(s))\in\cN$ for $0\leqs s\leqs t$. 

\item	Let $y^0(t)$ be the solution of the reduced system 
\begin{equation}
\label{stt6}
\dot y^0 = G(y^0) \defby g(x^\star(y^0),y^0)
\end{equation}
with initial condition $y^0(0)=y_0$. Let $K_1$ be a Lipschitz
constant for $G$ in $\cN_0$. Then 
\begin{equation}
\label{stt7}
\norm{y(t)-y^0(t)} \leqs c_3 \eps \e^{K_1t} +
c_4\norm{x_0-x^\star(y_0)}\e^{-K_0t/\eps} 
\end{equation}
for all $t\geqs0$ such that $(x(s),y(s))\in\cN$ for $0\leqs s\leqs t$.
\end{enum}
\end{theorem}

Roughly speaking, Relation~\eqref{stt5} implies that 
\begin{equation}
\label{stt8}
x(t) = x^\star(y(t)) + \Order{\eps}
\qquad\qquad
\text{for $t\geqs\frac1{K_0}\eps\abs{\log\eps}$,}
\end{equation}
and thus $x(t)$ reaches a \nbh\ of order $\eps$ of the slow manifold in a
time of order $\eps\abs{\log\eps}$, and stays there as long as this manifold
is attracting. If we want to determine $x(t)$ more precisely, we need to
approximate $y(t)$ as well. This is done by Relation~\eqref{stt7}, which
implies that we also have 
\begin{equation}
\label{stt9}
\begin{array}{l}
y(t) = y^0(t) + \Order{\eps} \\
x(t) = x^\star(y^0(t)) + \Order{\eps}
\end{array}
\qquad\qquad
\text{for $\frac1{K_0}\eps\abs{\log\eps}\leqs t\leqs\frac1{K_1}$.}
\end{equation}
Thus if we can solve the simpler reduced equation \eqref{stt6}, then we can
compute $x$ and $y$ up to errors of order $\eps$ and for times of order $1$.

\begin{proof}[{\sc Proof of Theorem~\ref{thm_Tihonov}}]
We will give a proof in the case $n=1$, and comment on larger $n$ at the
end of the proof. For $n=1$, $A(y)=a(y)$ is a scalar and $a(y)\leqs -a_0$
$\forall y$. 

\begin{enum}
\item	{\bf Change of variables:} 
The deviation of general solutions of \eqref{stt1} from the slow manifold
$x=x^\star(y)$ is characterized by the variable $z=x-x^\star(y)$, which
satisfies the equation
\[
\eps\dot z = \eps\dtot{}t \bigpar{x-x^\star(y)} = 
f(x^\star(y)+z,y) - \eps\dpar{x^\star}y(y) \, g(x^\star(y)+z,y). 
\]
The variables $(y,z)$ belong to
$\cN'=\setsuch{(y,z)}{(x^\star(y)+z,y)\in\cN}$. 
We will expand $f$ in Taylor series to second order in $z$, and write the
system in the form 
\[
\begin{split}
\eps\dot z &= a(y)z + b(y,z) + \eps w(y,z) \\
\dot y &= g(x^\star(y)+z,y). 
\end{split}
\]
The nonlinear term 
\[
b(y,z) = f(x^\star(y)+z,y) - a(y)z 
\]
satisfies, by Taylor's formula and our assumption on the derivatives, 
\[
\norm{b(y,z)} \leqs M\norm{z}^2 
\qquad\forall (y,z)\in\cN'
\]
for some constant $M>0$. The \defwd{drift term} is given by (see
\eqref{stt3})
\[
w(y,z) = -\dpar{x^\star}y(y) \, g(x^\star(y)+z,y) = 
a(y)^{-1} \dpar fy(x^\star(y),y) \, g(x^\star(y)+z,y).
\]
Our assumptions on $f$, $g$ and $a$ imply the
existence of a constant $W>0$ such that 
\[
\norm{w(y,z)} \leqs W 
\qquad\forall (y,z)\in\cN'.
\]
If $w$ were absent, the slow manifold $z=0$ would be invariant.  The drift
term acts like an inertial force, and may push solutions away from the slow
manifold. 

\item	{\bf Proof of \eqref{stt5}:}
The above bounds on $a$, $b$ and $w$ imply that 
\begin{align*}
\eps\dot z &\leqs -a_0z + Mz^2 + \eps W &
&\text{if $z\geqs0$} \\
\eps\dot z &\geqs -a_0z - Mz^2 - \eps W &
&\text{if $z\leqs0$} 
\end{align*}
whenever $(y,z)\in\cN$. Thus if $\ph(t)=\abs{z(t)}$, we have 
\[
\eps\dot\ph(t) \leqs -a_0\ph(t) + M\ph(t)^2 + \eps W.
\]
The function $\ph(t)$ is not differentiable at those $t$ for which $\ph=0$,
but it is left and right differentiable at these points, so that this small
inconvenience will not cause any problems. Now we take 
\[
c_0 = \frac{a_0}{2M}, \qquad
c_1 = \frac{2W}{a_0}, \qquad
\eps_0 = \frac{a_0}{2Mc_1}, 
\]
so that $\ph(0)=\abs{z(0)}<\frac{a_0}{2M}$ and $c_1\eps<\frac{a_0}{2M}$
under the hypotheses of the theorem. We define 
\[
\tau = \inf\Bigsetsuch{t>0}{\text{$\ph(t)=\frac{a_0}{2M}$ or
$(y(t),z(t))\notin\cN'$}}
\]
with the convention that $\tau=\infty$ if the set on the right-hand side is
empty. By continuity of the solutions, this definition implies that 
\[
\tau<\infty \qquad \Rightarrow \qquad
\ph(\tau)=\dfrac{a_0}{2M} \quad \text{or} \quad
(y(\tau),z(\tau))\in\partial\cN'.
\]
For $0\leqs t\leqs\tau$, we have 
\[
\eps\dot\ph(t) \leqs -\frac{a_0}2 \ph(t) + \eps W.
\]
We cannot apply our version of Gronwall's Lemma directly, because the
coefficient of $\ph(t)$ is negative. We can, however, bound the difference
between $\ph(t)$ and the value we would obtain if this relation were an
equality. We define a function $\psi(t)$ by 
\[
\ph(t) = \frac{2\eps}{a_0}W + \Bigpar{\ph(0) - \frac{2\eps}{a_0}W +
\psi(t)}\e^{-a_0t/2\eps}. 
\]
Differentiating with respect to time and using the inequality for $\dot\ph$,
we obtain that $\dot\psi(t)\leqs 0$. Since $\psi(0)=0$, we thus have 
\[
\ph(t) \leqs \frac{2\eps}{a_0}W + 
\Bigpar{\ph(0) - \frac{2\eps}{a_0}W}\e^{-a_0t/2\eps}
< \frac{a_0}{2M} \quad \forall t\in[0,\tau]. 
\]
Our definition of $c_0$, $c_1$, $\eps_0$ implies that
$\ph(t)<\frac{a_0}{2M}$ for $0\leqs t\leqs\tau$. If $\tau<\infty$, we  have
$\ph(\tau)<\frac{a_0}{2M}$ and thus necessarily
$(y(\tau),z(\tau))\in\partial\cN'$. In other words, the above upper bound
for $\ph(t)$ holds for all $t$ such that $(y(s),z(s))\in\cN'$ for
$s\in[0,t]$. Since $\ph(t)=\norm{x(t)-x^\star(y(t))}$, we have proved
\eqref{stt5}. 

\item	{\bf Proof of \eqref{stt7}:}
This part of the proof is similar to proofs in the chapter on regular
perturbation theory. We want to compare solutions of the equations 
\[
\begin{split}
\dot y &= g(x^\star(y)+z,y) \bydef G(y) + zR(y,z) \\
\dot y^0 &= g(x^\star(y^0),y^0) \bydef G(y^0).
\end{split}
\]
Taylor's formula shows that $\norm{R(y,z)}\leqs M'$ for some constant
$M'$ uniformly in $\cN'$. If $\y(t)=y(t)-y^0(t)$, we find 
\[
\dot\y = G(y^0+\y) - G(y^0) + z R(y^0+\y,z).
\]
Integrating from time $0$ to $t$, using $K_1$ as a Lipschitz constant for
$G$ and \eqref{stt5}, we obtain 
\[
\abs{\y(t)} \leqs K_1 \int_0^t \abs{\y(s)}\6s + M'(c_1 \eps +
c_2\norm{x_0-x^\star(y_0)} \e^{-K_0t/\eps}).
\]
Now the bound \eqref{stt7} follows easily from Gronwall's inequality and a
short computation. 

\item	{\bf The case $n>1$:}
In order to reduce the problem to a one-dimensional one, it is possible to
use a so-called \defwd{Liapunov function} $V:\R^n\to\R$, which decreases
along orbits of the associated system. A possibility is to use a quadratic
form 
\[
V(y,z) = z \cdot Q(y)z,
\]
for some suitable symmetric positive definite matrix $Q(y)$, depending on
$A(y)$. Then the function $\ph(t)=V(y(t),z(t))^{1/2}$ obeys a similar
equation as in the case $n=1$. 
\qed
\end{enum}
\renewcommand{\qed}{}
\end{proof}


\subsection{Iterations and Asymptotic Series}
\label{ssec_sit}

In the proof of Theorem~\ref{thm_Tihonov}, we have used the fact that the
original system
\begin{equation}
\label{sit1}
\begin{split}
\eps\dot x &= f(x,y)\\
\dot y &= g(x,y)
\end{split}
\end{equation}
is transformed, by the translation $x=x^\star(y)+z$, into the system 
\begin{equation}
\label{sit2}
\begin{split}
\eps\dot z &= f(x^\star(y)+z,y) - \eps \dpar{x^\star}y(y) g(x^\star(y)+z,y)\\
\dot y &= g(x^\star(y)+z,y).
\end{split}
\end{equation}
We can rewrite this system in the following way:
\begin{equation}
\label{sit2b}
\begin{split}
\eps\dot z &= A(y)z + b_0(y,z,\eps) + \eps w_0(y)\\
\dot y &= g_0(y,z),
\end{split}
\end{equation}
where $A(y) = \dpar fx(x^\star(y),y)$, and  $w_0(y)= -
\smash{\dpar{}y}x^\star(y)\mskip1.5mu g(x^\star(y),y)$ contains the terms
which do not vanish at $z=0$. The term $b_0$ contains all remainders of the
Taylor expansion of $f$ (at order $2$) and $g$ (at order $1$), and is
bounded in norm by a constant times $\norm{z}^2 + \eps\norm{z}$. 

We know that $z(t)$ is asymptotically of order $\eps$. The reason why $z$
does not go to zero is that the drift term $\eps w_0(y)$ may push $z$ away
from zero. If we manage to increase the order of the drift term, then $z$
may become smaller and we will have obtained a better approximation of the
solution of \eqref{sit1}.  

We could use the fact that \eqref{sit2b} is again a slow--fast system, where
$\dot z$ vanishes on a slow manifold $z=z^\star(y)$. Since $z^\star(y)$ is
in general difficult to compute, we will only use an approximation of it.
Consider the change of variables 
\begin{equation}
\label{sit3}
z = z_1 + \eps u_1(y), 
\qquad\qquad
u_1(y) = -A(y)^{-1} w_0(y).
\end{equation}
Substituting into \eqref{sit2b}, we get 
\begin{equation}
\label{sit4}
\eps\dot z_1 = A(y) z_1 + b_0(y,z_1+\eps u_1(y)) 
 - \eps^2 \dpar{u_1}y g_0(y,z_1+\eps u_1(y)).
\end{equation}
We can again extract the terms which vanish for $z=0$, which are of order
$\eps^2$. The system can thus be written as 
\begin{equation}
\label{sit5}
\begin{split}
\eps\dot z_1 &= A(y)z_1 + b_1(y,z_1,\eps) + \eps^2 w_1(y,\eps)\\
\dot y &= g_1(y,z_1,\eps),
\end{split}
\end{equation}
where $w_1$ denotes the right-hand side of \eqref{sit4} for $z_1=0$ and 
Taylor's formula shows again that $b_1(y,z,\eps)$ is of order $\norm{z}^2 +
\eps\norm{z}$. As long as the system is sufficiently differentiable, we can
repeat this procedure, increasing the order in $\eps$ of the drift term by
successive changes of variables. Thus there is a transformation 
\begin{equation}
\label{sit6}
x = x^\star(y) + \eps u_1(y) + \eps^2 u_2(y) + \dots + \eps^r u_r(y) + z_r
\end{equation}
which yields a system of the form 
\begin{equation}
\label{sit7}
\begin{split}
\eps\dot z_r &= A(y)z_r + b_r(y,z_r,\eps) + \eps^{r+1} w_r(y,\eps)\\
\dot y &= g_r(y,z_r,\eps).
\end{split}
\end{equation}
Proceeding as in the proof of Tihonov's theorem, one can show that
$z_r(t)=\Order{\eps^{r+1}}$ after a time of order $\eps\abs{\log\eps}$. This
implies that after this time, solutions of \eqref{sit1} satisfy
\begin{equation}
\label{sit8}
\begin{split}
x(t) &= x^\star(y(t)) + \eps u_1(y(t)) + \dots + \eps^r u_r(y(t)) +
\Order{\eps^{r+1}} \\
\dot y &= g(x^\star(y) + \eps u_1(y) + \dots + \eps^r u_r(y),y) +
\Order{\eps^{r+1}}
\end{split}
\end{equation}
for times of order $1$. 

One might wonder whether these iterations may be pushed all the way to
$r=\infty$ if the system is analytic, in such a way that the remainder
disappears. The answer in negative in general, because the amplitude of the
drift term $w_r$ tends to grow with $r$. We will illustrate this phenomenon
by a very simple example. 

\begin{example}
\label{ex_sit}
Consider the equation
\begin{equation}
\label{sit9}
\begin{split}
\eps\dot x &= -x + h(y) \\
\dot y &= 1,
\end{split}
\end{equation}
where we will assume that $y$ is analytic in the strip $\abs{\im y}\leqs
R$. This system can be solved exactly:
\begin{equation}
\label{sit10}
x(t) = x_0 \e^{-t/\eps} + \frac1\eps \int_0^t \e^{-(t-s)/\eps} h(y_0+s)\6s, 
\qquad y(t) = y_0 + t.
\end{equation}
Thus if $\abs{h(y)}$ is uniformly bounded by $M$, we have 
\begin{equation}
\label{sit11}
\abs{x(t)} \leqs \abs{x_0} \e^{-t/\eps} + \frac M\eps \int_0^t
\e^{-(t-s)/\eps}\6s = \abs{x_0} \e^{-t/\eps} + M (1-\e^{-t/\eps}).
\end{equation}
Thus there is no question that $x(t)$ exists and is uniformly bounded for
all $t\geqs0$, by a constant independent of $\eps$. 

Now let us pretend that we do not know the exact solution \eqref{sit10}. The
most naive way to solve \eqref{sit9} is to look for a series of the form 
\begin{equation}
\label{sit12}
x = x^0(y) + \eps x^1(y) + \eps^2 x^2(y) + \dots 
\end{equation} 
It is unlikely that all solutions admit such a series representation, but
if we obtain a particular solution of \eqref{sit9} which has this form, we
can then determine the general solution. Substituting the series
\eqref{sit12} into \eqref{sit9}, and solving order by order, we obtain 
\begin{align}
\nonumber
0 &= -x^0(y) + h(y) &
&\Rightarrow&
x^0(y) &= h(y) \\
\label{sit12b}
(x^0)'(y) &= -x^1(y) &
&\Rightarrow&
x^1(y) &= -h'(y) \\
\nonumber
(x^1)'(y) &= -x^2(y) &
&\Rightarrow&
x^2(y) &= h''(y)
\end{align}
and so on. Thus $x$ admits a formal series representation of the form 
\begin{equation}
\label{sit13}
x(t) = h(y(t)) - \eps h'(y(t)) + \dots + (-\eps)^k h^{(k)}(y(t)) + \dotsb
\end{equation}
Can this series converge? In certain cases yes, as shows the particular
case $h(y)=\sin y$ that we investigated in Example~\ref{ex_ssm1}. In that
case, the function $\bar x(t,\eps)$ admits an expansion in $\eps$ which
converges for $\abs{\eps}<1$. In general, however, the derivatives
$h^{(k)}(y)$ may grow with $k$. Cauchy's formula tells us that if $h$ is
analytic and bounded by $M$ in a disc of radius $R$ in the complex plane
around $y$, then 
\begin{equation}
\label{sit14}
h^{(k)}(y) = \frac{k!}{2\pi\icx} \int_{\abs{z-y}=R}
\frac{h(z)}{(z-y)^{k+1}} \6 z 
\qquad\Rightarrow\qquad
\abs{h^{(k)}(y)} \leqs M R^{-k} k!
\end{equation}
Thus in principle, $\abs{h^{(k)}(y)}$ may grow like $k!$, which is too fast
for the formal series \eqref{sit13} to converge. 

A more subtle method consists in trying to simplify \eqref{sit9} by
successive changes of variables. The translation $x=h(y)+z$ yields the
equation 
\begin{equation}
\label{sit15}
\eps\dot z = -z - \eps h'(y).
\end{equation}
The drift term $- \eps h'(y)$ can be further decreased by the change of
variables $z=- \eps h'(y)+z_1$, and so on. In fact, the transformation
\begin{equation}
\label{sit16}
x = \sum_{j=0}^k (-\eps)^j h^{(j)}(y) + z_k
\end{equation}
results in the equation
\begin{equation}
\label{sit17}
\eps\dot z_k = -z_k + (-\eps)^{k+1} h^{(k+1)}(y). 
\end{equation}
The solution is given by 
\begin{equation}
\label{sit18}
z_k(t) = z_k(0) \e^{-t/\eps} {-} (-\eps)^k \int_0^t \e^{-(t-s)/\eps}
h^{(k+1)}(y_0+s)\6s.
\end{equation}
Of course, the relations \eqref{sit16} and \eqref{sit18} can be obtained
directly from the exact solution \eqref{sit10} by $k$ successive
integrations by parts. The iterative method, however, also works in cases in
which we do not know an exact solution. In that case, the remainder $z_k(t)$
can be bounded as in the proof of Tihonov's theorem. 

We are looking for a particular solution admitting an expansion of the form
\eqref{sit13}. If we choose $x(0)$ in such a way that $z_k(0)=0$, we obtain
for each $k\in\Z$ a particular solution satisfying 
\begin{align}
\label{sit19}
x(t,\eps) &= \sum_{j=0}^k (-\eps)^j h^{(j)}(y(t)) + z_k(t,\eps)
\\
\label{sit20}
\abs{z_k(t,\eps)} &\leqs \eps^{k+1} M R^{-(k+1)} (k+1)!
\end{align}
where $M, R$ are independent of $k$. The first $k$ terms of this expansion
agree with the formal solution \eqref{sit13}. In general, an expression of
the form 
\begin{equation}
\label{sit21}
x(t,\eps) = \sum_{j=0}^k \eps^j u_j(t) + \eps^{k+1} r_k(t,\eps)
\qquad \forall k\in\Z,
\end{equation}
where $r_k$ is bounded for all $k$, but the series does not necessarily
converge, is called an \defwd{asymptotic series}. If the remainder grows
like in \eqref{sit20}, this asymptotic series is called of type
\defwd{Gevrey-1}. This situation is common in singularly perturbed systems:
bounded solutions exist, but they do not admit convergent series in $\eps$
(this differs from the situation we encountered in Chapter~3, where the
existence of a bounded solution is not always guaranteed). 
\end{example}

\begin{remark}
\label{rem_sit1}
We can determine at which order $k$ to stop the asymptotic expansion
\eqref{sit19} in such a way that the remainder $\eps^{k+1}\abs{z_k}$ is a
small as possible. This order will be approximately the one for which
$k!\eps^k$ is minimal. Recall Stirling's formula
\begin{equation}
\label{sit22}
k! \simeq k^k \e^{-k}.
\end{equation} 
It is convenient to optimize the logarithm of $k!\eps^k$. We have 
\begin{equation}
\label{sit23}
\dtot{}{k} \log\bigpar{k^k\e^{-k}\eps^k} = 
\dtot{}{k} \bigpar{k\log k - k + k\log\eps} = 
\log k + \log\eps,
\end{equation}
which vanishes for $\log k=-\log\eps$, i.e., $k=1/\eps$ (of course we need
to take an integer $k$ but this hardly makes a difference for small $\eps$). 
For this $k$, we have
\begin{equation}
\label{sit24}
k! \eps^k \simeq \Bigpar{\frac1\eps}^{1/\eps} \e^{-1/\eps} \eps^{1/\eps} 
= \e^{-1/\eps}.
\end{equation}
We conclude that the optimal truncation of the series is for
$k=\Order{1/\eps}$, which yields a remainder of order $\e^{-1/\eps}$. This
function goes to zero faster than any power of $\eps$ when $\eps\to0$, and
thus even though we cannot compute the solution exactly in general, we can
compute it to a relatively high degree of precision.

For general slow--fast systems one can show that if $f$ and $g$ are analytic
in some open complex domain, then the iterative scheme leading to
\eqref{sit7} can also be truncated in such a way that the remainder is
exponentially small in $\eps$. We will not prove this result here, but only
mention that it heavily relies on Cauchy's formula. 
\end{remark}

\begin{remark}
\label{rem_sit2}
There is another way which may help to understand why the asymptotic series
\eqref{sit13} does not converge for general $h$. Assume that $h$ is periodic
in $y$, and write its Fourier series in the form 
\begin{equation}
\label{sit25}
h(y) = \sum_{k\in\Z} h_k \e^{\icx k y}, 
\qquad\qquad
h_k = \frac1{2\pi} \int_0^{2\pi} h(y) \e^{-\icx k y} \6y.
\end{equation}
The formal series \eqref{sit13} will also be periodic in $y$, so that we may
look for a solution of \eqref{sit9} of the form 
\begin{equation}
\label{sit26}
x(y) = \sum_{k\in\Z} x_k \e^{\icx k y}.
\end{equation}
Substituting in the equation for $x$, we obtain 
\begin{equation}
\label{sit27}
\eps \icx k x_k = -x_k + h_k
\qquad\forall k\in\Z,
\end{equation}
and thus the periodic solution can be written as
\begin{equation}
\label{sit28}
x(y) = \sum_{k\in\Z} \frac{h_k}{1+\eps\icx k} \e^{\icx k y}.
\end{equation}
This Fourier series {\em does} converge in a \nbh\ of the real axis because
the $h_k$ decrease exponentially fast in $\abs{k}$ (to see this, shift
the integration path of $h_k$ in \eqref{sit25} by an imaginary distance).
Thus the equation admits indeed a bounded periodic solution. However, in the
plane of complex $\eps$, the function \eqref{sit28} has poles at every
$\eps=-\icx/k$ for which $h_k\neq0$. In the case $h(y)=\sin y$, there are
only two poles at $\eps=\pm\icx$, so that the series converges for
$\abs{\eps}<1$. If $h(y)$ has nonvanishing Fourier components with
arbitrarily large $k$, however, the poles accumulate at $\eps=0$, and the
radius of convergence of the series in $\eps$ is equal to zero.  
\end{remark}


\section{Dynamic Bifurcations}
\label{sec_sdb}

We examine now what happens to solutions of the slow--fast system 
\begin{equation}
\label{sdb1}
\begin{split}
\eps\dot x &= f(x,y)\\
\dot y &= g(x,y)
\end{split}
\end{equation}
in cases where the assumptions of Tihonov's theorem are violated. If the
slow manifold is represented by $x=x^\star(y)$, the main assumption of
Tihonov's theorem is that all eigenvalues of the matrix
\begin{equation}
\label{sdb2}
A(y) = \dpar fx(x^\star(y),y)
\end{equation}
have a strictly negative real part. This is equivalent to saying that
$x^\star(y)$ is an asymptotically stable equilibrium point of the associated
system 
\begin{equation}
\label{sdb3}
\dtot xs = f(x,y), 
\qquad\qquad
y = \text{{\it constant}}.
\end{equation}
A \defwd{dynamic bifurcation} occurs if the slow motion of $y$ causes some
eigenvalues of $A(y)$ to cross the imaginary axis. Then there may be a
bifurcation in the associated system \eqref{sdb3}, and the main assumption
of Tihonov's theorem is no longer satisfied. 

The two most generic cases (codimension $1$ bifurcations of equilibria) are 
\begin{itemiz}
\item	Hopf bifurcation: a pair of complex conjugate eigenvalues of $A$
cross the imaginary axis; the slow manifold $x^\star(y)$ continues to exist,
but changes its stability.
\item	Saddle--node bifurcation: an eigenvalue of $A(y)$ vanishes, and the
slow manifold $x^\star(y)$ ceases to exist (in fact, it has a fold).
\end{itemiz}

We will end this overview by discussing a few relatively simple examples of
these two dynamic bifurcations. 


\subsection{Hopf Bifurcation}
\label{ssec_shb}

\begin{example}
\label{ex_shb1}
Consider the slow--fast system
\begin{equation}
\label{shb1}
\begin{split}
\eps\dot x_1 &= y x_1 - x_2 - x_1 (x_1^2+x_2^2) \\
\eps\dot x_2 &= x_1 + y x_2 - x_2 (x_1^2+x_2^2) \\
\dot y &= 1.
\end{split}
\end{equation}
The associated system is 
\begin{equation}
\label{shb2}
\begin{split}
\dtot {x_1}s &= y x_1 - x_2 - x_1 (x_1^2+x_2^2) \\
\dtot {x_2}s &= x_1 + y x_2 - x_2 (x_1^2+x_2^2) \\
y &= \text{{\it constant}}. 
\end{split}
\end{equation}
The slow manifold is given by $x^\star(y)\equiv(0,0)$, and the matrix 
\begin{equation}
\label{shb3}
A(y) = 
\begin{pmatrix}
y & -1 \\ 1 & y
\end{pmatrix}
\end{equation}
has eigenvalues $a_{\pm}(y) = y \pm\icx$. Thus there is a Hopf bifurcation
at $y=0$. In fact the complex variable $z=x_1+\icx x_2$ satisfies 
\begin{equation}
\label{shb4}
\dtot zs = (y+\icx) z - \abs{z}^2 z,
\end{equation}
which becomes, in polar coordinates $z=r\e^{\icx\ph}$, 
\begin{equation}
\label{shb5}
\begin{split}
\dtot rs &= yr - r^3 \\
\dtot \ph s &= 1.
\end{split}
\end{equation}
If $y\leqs 0$, all solutions spiral to the origin $r=0$, while for $y>0$,
the origin is unstable, and all solutions not starting in $r=0$ are
attracted by the periodic orbit $r=\sqrt y$. 

Let us now return to the slow--fast system \eqref{shb1}. We know, by
Tihonov's theorem, that solutions starting not too far from $x=0$ at some
negative $y$ will reach a \nbh\ of order $\eps$ of the origin after a time
of order $\eps\abs{\log\eps}$. What happens as $y$ approaches zero? Since
the origin becomes unstable, intuition would probably tell us that the
solution will leave the \nbh\ of the origin as soon as $y$ becomes positive,
and approach the periodic orbit. 

This assumption can be checked by computing the solutions explicitly. In the
complex variable $z$, the slow--fast system \eqref{shb1} becomes 
\begin{equation}
\label{shb6}
\eps\dot z = (y+\icx) z - \abs{z}^2 z, 
\qquad\qquad 
\dot y = 1,
\end{equation}
and in polar coordinates,
\begin{equation}
\label{shb7}
\begin{split}
\eps\dot r &= yr - r^3 \\
\eps\dot\ph &= 1 \\
\dot y &= 1.
\end{split}
\end{equation}
The solution with initial condition $(r,\ph,y)(0)=(r_0,\ph_0,y_0)$ can be
written as 
\begin{equation}
\label{shb8}
\begin{split}
r(t) &= r_0 c(t) \e^{\alpha(t)/\eps} \\
\ph(t) &= \ph_0 + t/\eps \\
y(t) &= y_0 + t,
\end{split}
\end{equation}
where 
\begin{equation}
\label{shb9}
\begin{split}
\alpha(t) &= \int_0^t y(s) \6s = y_0 t + \frac12 t^2 \\
c(t) &= \Bigpar{1+\frac{2r_0^2}\eps \int_0^t \e^{\alpha(s)/\eps}
\6s}^{-1/2}.
\end{split}
\end{equation}
Assume that we start with $y_0<0$. Then the function $\alpha(t)$ is negative
for $0<t<-2y_0$. For these $t$, the term $\e^{\alpha(t)/\eps}$ is small,
while $c(t)\leqs1$. This implies that $r(t)$ is very small up to times $t$
close to $-2y_0$. However, the bifurcation already happens at time $t=-y_0$.

\begin{figure}
 \centerline{\psfig{figure=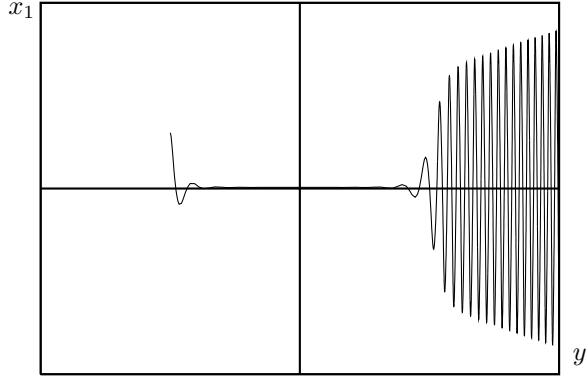,height=50mm,clip=t}}
 \figtext{
 	\writefig	11.0	0.7	$y$
 	\writefig	3.5	5.3	$x_1$
 }
 \captionspace
 \caption[]
 {A solution of Equation~\ref{shb1} starting with a negative $y_0$. The
 Hopf bifurcation occurs at $y=0$, but the solution approaches the limit
 cycle at $r=\sqrt{y}$ only after $y=-y_0$. This is the phenomenon of
 \defwd{bifurcation delay}.}
\label{fig_delay1}
\end{figure}

For $t>-2y_0$, the term $\e^{\alpha(t)/\eps}$ becomes large, but is
counteracted by $c(t)$. In fact, it is not hard to show that $r(t)$
approaches $\sqrt{y(t)}$ as $t\to\infty$, which means that the solution
\eqref{shb8} approaches the periodic orbit of the associated system. 

In summary, the solution starting at a distance $r_0$ (independent of
$\eps$) from the origin with $y_0<0$ is quickly attracted by the origin (in
accordance with Tihonov's theorem), but stays there until $y=-y_0$,
although the actual bifurcation takes place at $y=0$. Only after $y=-y_0$
will the solution approach the stable periodic orbit of the associated
system (\figref{fig_delay1}). This phenomenon is called \defwd{bifurcation
delay}. 
\end{example}

\begin{example}
\label{ex_shb2}
A particularity of the previous example is that the slow manifold does not
depend on $y$, so that there is no drift term pushing solutions away from
it. One might expect that the bifurcation delay is destroyed as soon as
$x^\star(y)$ depends on $y$. Let us consider the following modification of
\eqref{shb1}:
\begin{equation}
\label{shb10}
\begin{split}
\eps\dot x_1 &= y (x_1+y) - x_2 - (x_1+y) ((x_1+y)^2+x_2^2) \\
\eps\dot x_2 &= (x_1+y) + y x_2 - x_2 ((x_1+y)^2+x_2^2) \\
\dot y &= 1.
\end{split}
\end{equation}
The slow manifold is given by $x^\star(y)=(-y,0)$. The associated system is
the same as in Example \ref{ex_shb1}, up to a translation of $x^\star(y)$.
There is still a Hopf bifurcation at $y=0$. 

The variable $z=(x_1+y) + \icx x_2$ now obeys the equation 
\begin{equation}
\label{shb11}
\eps\dot z = (y+\icx)z - \abs{z}^2 z + \eps, 
\qquad\qquad \dot y = 1.
\end{equation}
The additional drift term stems from the $y$-dependence of the slow
manifold. 

In naive perturbation theory, one would look for a particular solution of
the form 
\begin{equation}
\label{shb12}
z = z^0(y) + \eps z^1(y) + \eps^2 z^2(y) + \dotsb
\end{equation}
Inserting into \eqref{shb11} and solving order by order, we get 
\begin{align}
\nonumber
0 &= (y+\icx) z^0 - \abs{z^0}^2 z^0 
&&\Rightarrow&
z^0(y) &= 0 \\
\label{shb13}
\dot z^0 &= (y+\icx) z^1 - \bigbrak{(z^0)^2\bar z^1 + 2\abs{z^0}^2 z^1} + 1
&&\Rightarrow&
z^1(y) &= -\frac 1{y+\icx} \\
& \dots
\nonumber
\end{align}
Thus Equation \eqref{shb11} seems to admit a particular solution 
\begin{equation}
\label{shb14}
z_{\math{p}}(t) = -\frac\eps{y(t)+\icx} + \Order{\eps^2}.
\end{equation}
To obtain the general solution, we set $z=z_{\math{p}}+\z$, substitute in
\eqref{shb11}, and use the fact that $z_{\math{p}}$ is a particular
solution. The result is an equation for $\z$ of the form 
\begin{equation}
\label{shb15}
\eps\dot\z = (y+\icx)\z + \Order{\eps^2\abs{\z}} + \Order{\eps\abs{\z}^2} +
\Order{\abs{\z}^3}.
\end{equation}
The drift term has disappeared in the transformation, so that the
right-hand side of \eqref{shb15} vanishes for $\z=0$. Although this
equation cannot be solved exactly, it should behave as in Example
\ref{ex_shb1}. To show this, let us assume that the system is started at
time $t=t_0=y_0$ in such a way that $y(0)=0$. We use the Ansatz
\begin{equation}
\label{shb16}
\z(t) = \z_0 \e^{\alpha(t,t_0)/\eps} c(t), 
\qquad
\alpha(t,t_0) = \icx (t-t_0) + \frac12 (t^2-t_0^2).
\end{equation}
If the last three terms in \eqref{shb15} were absent, the solution would be
given by \eqref{shb16} with $c(t)\equiv1$. Inserting \eqref{shb16} into
\eqref{shb15}, we obtain a differential equation for $c$, which can be used
to show that $c$ remains of order $1$ as long as $\re\alpha(t,t_0)<0$. Since
$\re\alpha(t,t_0)<0$ for $t_0<t<-t_0$, although the bifurcation happens at
$t=0$, we conclude that the bifurcation is delayed as in the previous
example. 

\begin{figure}
 \centerline{\psfig{figure=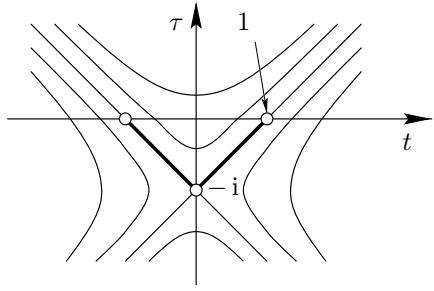,height=40mm,clip=t}}
 \figtext{
 	\writefig	9.8	2.4	$t$
 	\writefig	6.7	4.0	$\tau$
 	\writefig	7.6	4.0	$1$
 	\writefig	7.2	1.8	$-\icx$
 }
 \captionspace
 \caption[]
 {Level lines of the function $\re\alpha(t+\icx\tau)$. The function is
 small for large $\abs{\tau}$ and large for large $\abs{t}$. The integral
 $\Psi(t,t_0)$ is small only if $t_0$ and $t$ can be connected by a path on
 which $\re\alpha$ is non-increasing.}
\label{fig_delay2}
\end{figure}

There is, however, a major problem with this procedure. We already know
that the asymptotic series \eqref{shb12} is unlikely to converge, but this
does not necessarily imply that \eqref{shb11} does not admit a solution of
order $\eps$. In this case, however, \eqref{shb11} {\em does not admit a
solution of order $\eps$ for all times}. To see this, let us write the
general solution implicitly as 
\begin{equation}
\label{shb17}
z(t) = z_0 \e^{\alpha(t,t_0)/\eps} + \frac1\eps \int_{t_0}^t
\e^{\alpha(t,s)/\eps} \bigbrak{\eps - \abs{z(s)}^2 z(s)} \6s.
\end{equation}
We can try to solve this equation by iterations. Assume that
$z_0=\Order{\eps}$. As long as $z$ remains of order $\eps$, we have 
\begin{equation}
\label{shb18}
z(t) = z_0 \e^{\alpha(t,t_0)/\eps} + \int_{t_0}^t
\e^{\alpha(t,s)/\eps} \bigbrak{1 + \Order{\eps^2}} \6s.
\end{equation}
The behaviour of the solution is thus essentially contained in the integral 
\begin{equation}
\label{shb19}
\Psi(t,t_0) = \int_{t_0}^t \e^{\alpha(t,s)/\eps} \6s.
\end{equation}
As long as this integral remains of order $\eps$, $z(t)$ is of order $\eps$.
If it becomes large, however, $z$ becomes large as well (though this has to
be checked independently). $\Psi$ can be evaluated by deformation of the
integration path into the complex plane. The function
$\alpha(t)\defby\alpha(t,0)$ can be continued to the complex plane, and
satisfies 
\begin{equation}
\label{shb20}
\re\alpha(t+\icx\tau) = -\tau + \frac12 \bigbrak{t^2-\tau^2} 
= \frac12 \bigbrak{t^2 - (\tau+1)^2 + 1}.
\end{equation}
The level lines of $\re\alpha$ are hyperbolas centred at
$t+\icx\tau=-\icx$. The integral \eqref{shb19} is small if we manage to
find a path of integration connecting $t_0$ to $t$ on which $\re\alpha(t,s)
= \re\alpha(t)-\re\alpha(s)\leqs0$ for all $s$. This is only possible if
$\re\alpha(t)\leqs\re\alpha(t_0)$, but this condition is not sufficient. In
fact (see \figref{fig_delay2}),
\begin{itemiz}
\item	if $t_0\leqs-1$, such a path exists for all $t\leqs1$;
\item	if $-1\leqs t\leqs0$, such a path exists for all $t\leqs-t_0$;
\item	if $t_0\geqs0$, no such path exists for $t\geqs t_0$.
\end{itemiz}
We conclude that if $t_0\leqs0$, then $\Psi(t,t_0)$ remains small for all
times 
\begin{equation}
\label{shb21}
t \leqs \min\set{-t_0,1}.
\end{equation}
As a consequence, $z(t)$ also remains small up to this time. The term $-t_0$
is natural, it is the same as in Example \ref{ex_shb1}. However, the term
$1$, called \defwd{maximal delay} or \defwd{buffer point}, is new. It is an
effect of the drift term, and cannot be found by naive perturbation theory,
because it is, so to speak, hidden in exponentially small terms. This
example shows that one should be extremely cautious with asymptotic
expansions near dynamic bifurcations. 
\end{example}


\subsection{Saddle--Node Bifurcation}
\label{ssec_ssn}

\begin{example}
\label{ex_ssn}
Let us finally return to Van der Pol's equation, written in the form 
\begin{equation}
\label{ssn1}
\begin{split}
\eps\dot x &= y+x-\frac13x^3 \\
\dot y &= -x.
\end{split}
\end{equation}
The slow manifold is given implicitly by the equation 
\begin{equation}
\label{ssn2}
y = \frac13 x^3 - x.
\end{equation}
This equation has three solutions for $x$ if $\abs{y}<2/3$ and one if
$\abs{y}>2/3$. The associated system 
\begin{equation}
\label{ssn3}
\dtot xs = x - \frac13 x^3 + y, 
\qquad\qquad
y=\text{{\it constant}}
\end{equation} 
has three equilibria (two stable, one unstable) for $\abs{y}<2/3$ and one
(stable) equilibrium for $\abs{y}>2/3$. The points $(1,-2/3)$ and $(-1,2/3)$
are saddle--node bifurcation points. 

Let us now study the dynamics near $(1,-2/3)$. The change of variables 
\begin{equation}
\label{ssn4}
\begin{split}
x &= 1 + \x \\
y &= -\frac23 - \y
\end{split}
\end{equation}
transforms \eqref{ssn1} into the system 
\begin{equation}
\label{ssn5}
\begin{split}
\eps\dot\x &= -\y -\x^2 - \frac13 \x^3 \\
\dot\y &= 1 + \x.
\end{split}
\end{equation}
Assume that we start near the piece of slow manifold given by 
\begin{equation}
\label{ssn6}
\x = \x^\star(\y) = \sqrt{-\y} + \Order{\y}, 
\qquad \y\leqs 0. 
\end{equation}
Tihonov's theorem states that for negative $\y$ bounded away from zero,
$\x$ quickly reaches a \nbh\ of order $\eps$ of $\x^\star(\y)$. The
transformation $\x=\x^\star(\y)+\x_1$ yields the new system
\begin{equation}
\label{ssn7}
\begin{split}
\eps\dot\x_1 &= -\bigbrak{2\x^\star(\y)+(\x^\star(\y))^2}\x_1 -
\bigbrak{1+\x^\star(\y)}\x_1^2 - \frac13\x_1^3 -
\eps\dtot{\x^\star(\y)}{\y} (1+\x^\star(\y)+\x_1) \\ 
\dot\y &= 1 + \x^\star(\y) + \x_1.
\end{split}
\end{equation}
The order of the drift term can be further decreased by a translation 
\begin{equation}
\label{ssn8}
\x_1 = \x_2 + \eps u_1(\y), 
\qquad
u_1(\y) = - \dtot{\x^\star(\y)}{\y}
\frac{1+\x^\star(\y)}{2\x^\star(\y)+(\x^\star(\y))^2}.
\end{equation}
Continuing in this way, we construct an asymptotic series 
\begin{equation}
\label{ssn9}
\x = \x^\star(\y) + \eps u_1(\y) + \eps^2 u_2(\y) + \dotsb
\end{equation}
Note, however, that as $\y$ approaches $0$, we have 
\begin{equation}
\label{ssn10}
\x^\star(\y) \sim \sqrt{-\y}, 
\qquad\qquad
u_1(\y) \sim \frac1\y
\end{equation}
and further computations will show that $u_2(\y)\sim\abs{\y}^{-5/2}$. Thus
the successive terms of the asymptotic series decrease in amplitude only if 
\begin{equation}
\label{ssn11}
\frac\eps{\abs{\y}} < \sqrt{\abs\y} 
\qquad\Rightarrow\qquad
\abs{\y} > \eps^{2/3}.
\end{equation}
One says that the asymptotic series becomes \defwd{disordered} at
$\y=-\eps^{2/3}$. Following the proof of Tihonov's theorem, one can indeed
show that $\x$ behaves like the first terms of the asymptotic series for
times smaller than $-\eps^{2/3}$. For larger $\y$, the system has to be
described by other means. Since $\x$ is of order $\eps^{1/3}$ when
$\y=-\eps^{2/3}$, we introduce the scaling 
\begin{equation}
\label{ssn12}
\begin{split}
\x &= \eps^{1/3} u \\
\y &= \eps^{2/3} v,
\end{split}
\end{equation}
which transforms \eqref{ssn5} into 
\begin{equation}
\label{ssn13}
\begin{split}
\eps^{4/3} \dot u &= -\eps^{2/3} v - \eps^{2/3} u^2 - \frac13 \eps u \\
\eps^{2/3} \dot v &= 1 + \eps^{1/3} u,
\end{split}
\end{equation}
so that 
\begin{equation}
\label{ssn14}
\dtot uv = -u^2 - v + \Order{\eps^{1/3}}. 
\end{equation}
This is a perturbation of order $\eps^{1/3}$ of a solvable Riccati equation,
that can be used to describe the motion near the bifurcation point. We will
not discuss this equation in detail, but only mention that $u$ can be shown
to reach $-1$ after a \lq\lq time\rq\rq\ $v$ of order $1$. This implies that
$\x$ reaches $-\eps^{1/3}$ for $\y$ of order $\eps^{2/3}$
(\figref{fig_sndyn}a). 

\begin{figure}
 \centerline{\psfig{figure=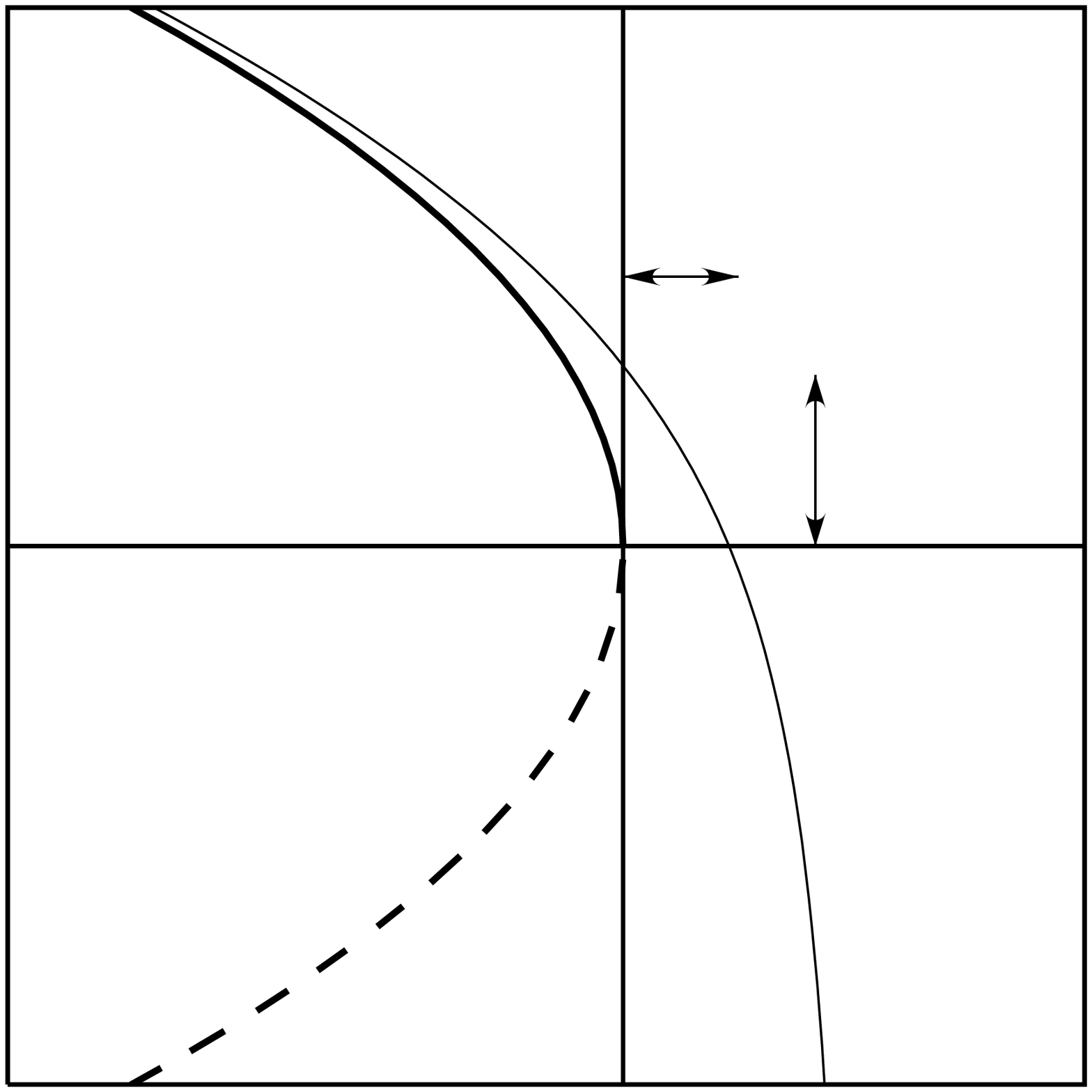,height=50mm,clip=t}
 \hspace{12mm}
 \psfig{figure=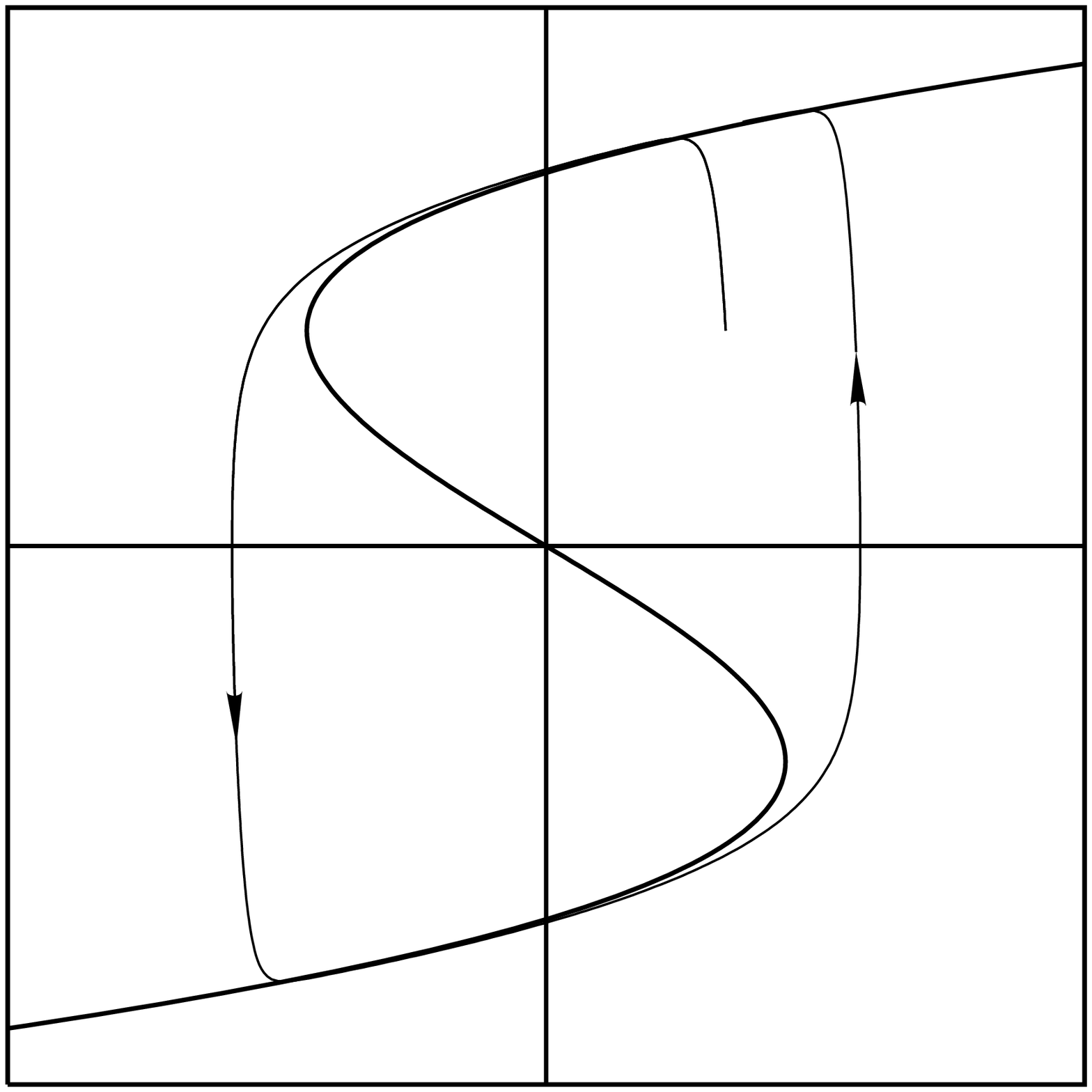,height=50mm,clip=t}}
 \figtext{
 	\writefig	1.1	5.1	(a)
 	\writefig	6.3	3.15	$\y$
 	\writefig	4.6	5.2	$\x$
 	\writefig	4.6	4.4	$\eps^{2/3}$
 	\writefig	5.5	3.3	$\eps^{1/3}$
 	\writefig	7.5	5.1	(b)
 	\writefig	12.7	3.15	$y$
 	\writefig	10.7	5.2	$x$
 }
 \captionspace
 \caption[]
 {(a) Behaviour of a solution of \eqref{ssn5} near the saddle--node
 bifurcation point. The trajectory jumps after a delay of order
 $\eps^{2/3}$. (b) This jump results in relaxation oscillations for the
 global behaviour of van der Pol's equation.}
\label{fig_sndyn}
\end{figure}

Once the bifurcation point has been passed, the solution reaches another
branch of the slow manifold, which it follows until the next bifurcation
point. The result is an asymptotically periodic motion, with alternating
slow and fast phases. In the limit $\eps\to0$, the motion follows branches
of the slow manifold and lines $y=\text{{\it constant}}$. When $\eps>0$,
this cycle is enlarged by an amount of order $\eps^{1/3}$ in the
$x$-direction and of order $\eps^{2/3}$ in the $y$-direction
(\figref{fig_sndyn}b). 
\end{example}


\section*{Bibliographical Comments}

Various approaches to singular perturbation theory have been developed,
based on asymptotic expansions \cite{MKKR,VBK}, invariant manifolds
\cite{Fe,Jones} and nonstandard analysis \cite{Benoit}. More information on
asymptotic series can be found in \cite{Wasow}. 

Tihonov's theorem was originally proved in \cite{Tihonov} and \cite{Grad},
and has been generalized to periodic orbits \cite{PR}. 

The phenomenon of delayed Hopf bifurcation has been analysed in detail by
Neishtadt \cite{Ne1,Ne2}. The saddle--node bifurcation was first considered
in \cite{Po}, and later in \cite{Hab,MKKR}. See  also \cite{LS} and
\cite{Benoit} for other dynamic bifurcations.  



\begin{thebibliography}{MMMM99}

\bibitem[Ar63]{Arnold63}
			\bibarticle{V.I.\ Arnol'd}
			{Proof of a theorem of {A}. {N}. {K}olmogorov on the
              		preservation of conditionally periodic motions under a small
              		perturbation of the {H}amiltonian}
			{Uspehi Mat. Nauk}
			{18}{1}{40}{1963}
     
\bibitem[Ar83]{Arnold}
			\bibbook{V.I.\ Arnold}
			{Geometrical Methods in the Theory of Ordinary
			Differential Equations}
			{Springer--Verlag}{New York, 1983}

\bibitem[Ar89]{Arnold89}
			\bibbook{V.I.\ Arnold}
			{Mathematical methods of classical mechanics}
			{Springer--Verlag}{New York, 1989}
      	
\bibitem[Be91]{Benoit}	
			\bibbook{E.\ Beno\^{\i}t (Ed.)}
			{Dynamic Bifurcations, Proceeding, Luminy 1990}
			{Sprin\-ger--Verlag, Lecture Notes in Mathematics
			1493}
			{Berlin, 1991}
			
\bibitem[Bo75]{Bogdanov}
			\bibarticle{R.I.\ Bogdanov}
			{Versal deformations of singular points of vector
	     		fields on the plane}
			{Functional Anal.\ Appl.}
			{9}{144}{145}{1975}
	      
\bibitem[dlL01]{dlL01}
			\bibpreprint{R.\ de\ la\ Llave}
			{A tutorial on KAM theory}
			{mp-arc 01-29}
			{2001}
			{\tt
			http://mpej.unige.ch/mp\_arc-bin/mpa?yn=01-29}

\bibitem[D69]{Deprit}	
			\bibarticle{A.\ Deprit}
			{Canonical transformations depending on a small
			parameter}
			{Celestial Mech.}{1}{12}{30}{1969}

\bibitem[F79]{Fe}	
			\bibarticle{N.\ Fenichel}
			{Geometric singular perturbation theory for ordinary
			differential equations}
			{J.\ Diff.\ Eq.}{31}{53}{98}{1979}
			
\bibitem[Gr53]{Grad}      \bibarticle{I.\,S. Grad\v ste\u\i n}
                        {Applications of A.\,M.\ Lyapunov's theory of
                        stability to the theory of differential equations
                        with small coefficients in the derivatives}
                        {Mat.\ Sbornik N.S.}{32}{263}{286}{1953}
                        
\bibitem[GH83]{GH}
			\bibbook{J. Guckenheimer, P. Holmes}
			{Nonlinear Oscillations, Dynamical Systems, 
			and Bifurcations of Vector Fields}
			{Springer-Verlag}{New York, 1983}

\bibitem[Ha79]{Hab}		
			\bibarticle{R.\ Haberman}
			{Slowly varying jump and transition phenomena 
			associated with algebraic bifurcation problems}
			{\SIAM}{37}{69}{106}{1979}
			
\bibitem[HK91]{HK}		
			\bibbook{J.\ Hale, H.\ Ko\c cak}
			{Dynamics and Bifurcations}
			{Springer--Verlag}		
			{New York, 1991}
			
\bibitem[H70]{Henrard}
			\bibarticle{J.\ Henrard}
			{On a perturbation theory using Lie transforms}
			{Celestial Mech.}{3}{107}{120}{1970}

\bibitem[H83]{Herman}
			\bibbook{M.-R.\ Herman}
			{Sur les courbes invariantes par les
			diff\'eomorphismes de l'anneau. {V}ol. 1}
			{Soci\'et\'e Math\'ematique de France}	
			{Paris, 1983}
      
\bibitem[J95]{Jones}		
			\bibtitle{C.K.R.T.\ Jones}
			{Geometric Singular Perturbation Theory}, in
			\bibbook{R.\ Johnson (Ed.)}
			{Dynamical Systems, Proc.\ (Montecatini Terme, 
			1994)}
			{Lecture Notes in Math.\ 1609}
			{Springer, Berlin, 1995} 

\bibitem[KKR98]{KKR}	
			\bibarticle{A.I.\ Khibnik, B.\ Krauskopf, C.\
			Rousseau}
			{Global study of a family of cubic Li\'enard
			equations}
			{\NL}{11}{1505}{1519}{1998}
			
\bibitem[K54]{Kolmogorov}
			\bibarticle{A.N.\ Kolmogorov}
			{On conservation of conditionally periodic motions
			for a small change in Hamilton's function}
			{Dokl. Akad. Nauk SSSR (N.S.)}
			{98}{527}{530}{1954}
			
\bibitem[LS77]{LS}		
			\bibarticle{N.R.\ Lebovitz, R.J.\ Schaar}
			{Exchange of Stabilities in Autonomous Systems I, II}
			{\SAM}{54}{229}{260}{1975}
			\bibref{\SAM}{56}{1}{50}{1977}.

\bibitem[Me92]{Meiss}	
			\bibarticle{J.D. Meiss} 
			{Symplectic maps, variational principles, and
			transport}
			{\RMP}{64}{795}{848}{1992}
			
\bibitem[MKKR94]{MKKR}	
			\bibbook{E.F.\ Mishchenko, Yu.S.\ Kolesov, 
			A.Yu.\ Kolesov, N.Kh.\ Rozov}
			{Asymptotic Methods in Singularly Perturbed Systems}
			{Consultants Bureau}{New York, 1994}
			
\bibitem[Mo62]{Moser62}
			\bibarticle{J.\ Moser}
			{On invariant curves of area-preserving mappings of
			an annulus}
			{Nachr. Akad. Wiss. G\"ottingen Math.-Phys. Kl. II}
			{1962}{1}{20}{1962}
     
\bibitem[Mo73]{Moser73}		
			\bibbook{J.\ Moser}
			{Stable and Random Motions in Dynamical Systems}
			{Princeton University Press}
			{Princeton, New Jersey, 1973, 2000}

\bibitem[Ne87]{Ne1}	
			\bibarticle{A.I.\ Neishtadt}
			{Persistence of stability loss for dynamical 
			bifurcations I}
			{\DE}{23}{1385}{1391}{1987}
			Transl. from \bibref{\DU}{23}{2060}{2067}{1987}.
			
\bibitem[Ne88]{Ne2}	
			\bibarticle{A.I.\ Neishtadt}
			{Persistence of stability loss for dynamical 
			bifurcations II}
			{\DE}{24}{171}{176}{1988}
			Transl. from \bibref{\DU}{24}{226}{233}{1988}.
			
\bibitem[PM82]{Palis}
			\bibbook{J.\ Palis, W.\ de Melo}
			{Geometric theory of dynamical systems}
			{Springer--Verlag}{New York, 1982}
      
\bibitem[PP59]{PP}
			\bibarticle{M.C.\ Peixoto, M.M.\ Peixoto}
			{Structural stability in the plane with enlarged
			boundary conditions}
			{An.\ Acad.\ Brasil.\ Ci.}
			{31}{135}{160}{1959}
			
\bibitem[PP68]{PP68}
			\bibarticle{M.M.\ Peixoto, C.C.\ Pugh}
			{Structurally stable systems on open manifolds are
			never dense}
			{Ann.\ of Math.\ (2)}{87}{423}{430}{1968}
     
\bibitem[Pe62]{Pe62}
			\bibarticle{M.M.\ Peixoto}
			{Structural stability on two-dimensional manifolds}
			{Topology}{1}{101}{120}{1962}
     
\bibitem[Pe73]{Pe73}
			\bibtitle{M.M.\ Peixoto}
			{On the classification of flows on 2-manifolds}, 
			in \bibbook{M.M.\ Peixoto (Ed.)}
 			{Dynamical systems}
			{Academic Press}{New York-London, 1973}
      
\bibitem[Po57]{Po}		
			\bibarticle{L.S.\ Pontryagin}
			{Asymptotic behavior of solutions of systems of  
			differential equations with a small parameter 
			in the derivatives of highest order}
			{Izv.\ Akad.\ Nauk SSSR.\ Ser.\ Mat.}
			{21}{605}{626}{1957}

\bibitem[PR60]{PR}		
			\bibarticle{L.S.\ Pontryagin, L.V.\ Rodygin}
			{Approximate solution of a system of ordinary 
			differential equations involving a small parameter 
			in the derivatives}
			{\Dokl}{131}{237}{240}{1960}

\bibitem[R{\"u}70]{Russmann70}
			\bibarticle{H.\ R{\"u}ssmann}
			{Kleine {N}enner. {I}. \"{U}ber invariante {K}urven
              		differenzierbarer {A}bbildungen eines {K}reisringes}
			{Nachr. Akad. Wiss. G\"ottingen Math.-Phys. Kl. II}
			{1970}{67}{105}{1970}

\bibitem[R{\"u}72]{Russmann72}
			\bibarticle{H.\ R{\"u}ssmann}
			{Kleine {N}enner. {I}{I}. {B}emerkungen zur {N}ewtonschen
              		{M}ethode}
			{Nachr. Akad. Wiss. G\"ottingen Math.-Phys. Kl. II}
			{1972}{1}{10}{1972}

\bibitem[Si65]{Silnikov}
			\bibarticle{L.P.\ \v Silnikov}
			{A case of the existence of a denumerable set of
			periodic motions}
			{Sov.\ Math.\ Dokl.}{6}{163}{166}{1965}

\bibitem[Sm66]{Smale}
			\bibarticle{S.\ Smale}
			{Structurally stable systems are not dense}
			{Amer.\ J.\ Math.}{88}{491}{496}{1966}
     
\bibitem[Ta74]{Takens}	
			\bibarticle{F.\ Takens}
			{Forced oscillations and bifurcations}
			{Comm.\ Math.\ Inst., Rijksuniversiteit Utrecht}
			{3}{1}{59}{1974}
			
\bibitem[Ti52]{Tihonov}   \bibarticle{A.\,N. Tihonov}
                        {Systems of differential equations containing small
                        parameters in the derivatives}
                        {Mat.\ Sbornik N.S.}
                        {31}{575}{586}{1952}
                        
\bibitem[VBK95]{VBK}	
			\bibbook{A.B.\ Vasil'eva, V.F.\ Butusov, L.V.\
			Kalachev}
			{The Boundary Function Method for Singular
			Perturbation Problems}
			{SIAM}
			{Philadelphia, 1995}
			
\bibitem[V96]{Verhulst}	
			\bibbook{F.\ Verhulst}
			{Nonlinear Differential Equations and Dynamical
			Systems}
			{Sprin\-ger-Verlag}
			{Berlin, 1996}

\bibitem[Wa65]{Wasow}	
			\bibbook{W.\ Wasow}
			{Asymptotic expansions for ordinary differential 
			equations}
			{Re\-print of the 1965 edition, Krieger Publishing}		
			{Huntington, NY, 1976}
					
\bibitem[Wi90]{Wiggins}
			\bibbook{S.\ Wiggins}
			{Introduction to Applied Nonlinear Dynamical Systems
			and Chaos}
			{Springer--Verlag}              
			{New York, 1990}

			
\end{thebibliography}
\end{document}